\documentclass{article}
\PassOptionsToPackage{table}{xcolor}
\usepackage{iclr2026_conference,times}

\iclrfinalcopy
\usepackage{graphicx,epstopdf}
\usepackage{enumitem,amsfonts,amssymb,mathtools,amsthm,bm,pifont,listings}
\usepackage{subcaption,hyperref,comment,amsmath}
\usepackage{array,threeparttable,booktabs,makecell,algorithm}
\usepackage{xcolor}
\def\cite{\citep}
\definecolor{customgray}{gray}{0.9}
\usepackage{arydshln}
 
\newtheorem{theorem}{Theorem}[section]

\newtheorem{lemma}[theorem]{Lemma}

\newtheorem{assumption}[theorem]{Assumption}
\newtheorem{remark}[theorem]{Remark}

\def\a{\mathbf{a}}\def\b{\mathbf{b}}\def\d{\mathbf{d}}\def\e{\mathbf{e}}\def\g{\mathbf{g}}\def\h{\mathbf{h}}\def\m{\mathbf{m}}\def\q{\mathbf{q}}\def\r{\mathbf{r}}\def\s{\mathbf{s}}\def\u{\mathbf{u}}\def\v{\mathbf{v}}\def\w{\mathbf{w}}\def\x{\mathbf{x}}\def\y{\mathbf{y}}\def\z{\mathbf{z}}\def\A{\mathbf{A}}\def\C{\mathbf{C}}\def\E{\mathbf{E}}\def\I{\mathbf{I}}\def\V{\mathbf{V}}\def\W{\mathbf{W}}\def\v{\mathbf{v}}

\def\OO{\mathcal{O}}\def\VV{\mathcal{V}}\def\XX{\mathcal{X}}\def\ZZ{\mathcal{Z}}

\def\EEE{\mathbb{E}}\def\LLL{\mathbb{L}}\def\SSS{\mathbb{S}}\def\UUU{\mathbb{U}}

\def\trans{^\mathsf{T}}\def\Rn{\mathbb{R}} \def\one{\mathbf{1}} \def\zero{\mathbf{0}}  
 \def\trans {^\mathsf{T}} \def\fro {\mathsf{F}}

\DeclareMathOperator{\st}{s.t.}
\DeclareMathOperator{\mat}{mat}
\DeclareMathOperator{\prox}{Prox}
\DeclareMathOperator{\vecc}{vec}

\DeclareMathOperator{\dom}{dom}

\DeclareMathOperator{\Prox}{\operatorname{Prox}}
\DeclareMathOperator{\argmin}{arg\,min}

\newcommand{\bfit}[1]{\textit{\textbf{#1}}}
\newcommand{\beq}{\begin{eqnarray}}
\newcommand{\eeq}{\end{eqnarray}}
\newcommand{\beqq}{\begin{equation}}
\newcommand{\eeqq}{\end{equation}}

\newcommand{\bel}{\begin{align}}
\newcommand{\eel}{\end{align}}
\newcommand{\la}{\langle}
\newcommand{\ra}{\rangle}
\newcommand{\noi}{\noindent}
\newcommand{\nn}{\nonumber}
\def\noi{\noindent}
\def\nn{\nonumber}
\def\la{\langle}
\def\ra{\rangle}

\def\Diag{{\rm{Diag}}}
\def\dist{{\rm{dist}}}
\newcommand{\step}[1]{\text{\ding{\numexpr#1+171\relax}}}

\def\ts{\textstyle}

\def\etas{\boldsymbol{\eta}}

\providecommand{\epsilons}{{\boldsymbol{\epsilon}}}
\providecommand{\deltas}{{\boldsymbol{\delta}}}

\usepackage{tikz}

\makeatother

\definecolor{lightgray}{rgb}{0.9, 0.9, 0.9}

\usepackage{url}

\definecolor{darkred}{HTML}{8B0000}

\providecommand{\normalem}{}
\renewcommand{\E}{{\mathbb E}}
\hypersetup{hidelinks}
\newcommand{\xxxx}[2]{\if\relax\detokenize\expandafter{\romannumeral-`\q#1}\relax#2\else \overset{\step{#1}}{#2}\fi}
\def\myvspace{{\vspace{-8pt}}}

\providecommand{\cmark}{\ding{52}}
\providecommand{\xmark}{\ding{56}}

\title{Adaptive Lipschitz-Free Proximal Gradient Methods for Stochastic Composite Nonconvex Optimization}

\author{Ganzhao Yuan \\
Faculty of Computer Science and Artificial Intelligence, \\
Shenzhen University of Advanced Technology (SUAT), China\\
\texttt{yuanganzhao@foxmail.com} \\
}

\begin{document}

\maketitle

\begin{abstract}
We propose {\sf ALFPG}, an Adaptive Lipschitz-Free Proximal Gradient
framework for stochastic composite nonconvex minimization.
The framework adapts its stepsizes using information collected
along the optimization trajectory, without requiring prior
knowledge of smoothness constants or performing line search.
For finite-sum problems, {\sf ALFPG-FS} combines coordinate-wise
stepsizes, Nesterov's extrapolation, and the SPIDER estimator,
with {\sf ALFPG-D} as its full-batch specialization. Under component-wise smoothness, {bounded sublevel sets, and separate lower bounds on $f$ and $h$}, {\sf ALFPG-FS} achieves expected $\epsilon$-stationarity using
$\mathcal O(N+\sqrt N \epsilon^{-2})$ component-gradient evaluations. For expectation problems, we develop {\sf ALFPG-MVR1} and {\sf ALFPG-MVR2}, where {\sf MVR} stands for Momentum-based Variance Reduction. Both use one fresh sample per iteration: {\sf ALFPG-MVR1} uses one stochastic-gradient evaluation, whereas {\sf ALFPG-MVR2} uses at most two evaluations at different points with the same sample. {\sf ALFPG-MVR1} requires smoothness of the expected loss; {\sf ALFPG-MVR2} additionally requires uniformly smooth sample losses. When the regularizer is convex or Lipschitz continuous on its
effective domain, and stochastic gradients are almost surely
bounded, the two methods achieve $\epsilon$-stationarity in
expectation using $\widetilde{\mathcal O}(\sigma^2\epsilon^{-4}+\epsilon^{-2})$
and $\widetilde{\mathcal O}(\sigma\epsilon^{-3}+\epsilon^{-2})$
stochastic-gradient evaluations, respectively.
Here, $N$ is the number of component functions,
$\sigma^2$ bounds the gradient noise variance, and
$\widetilde{\mathcal O}$ hides polylogarithmic factors.
Neither expectation variant requires prior knowledge of
$\sigma$, and both recover
$\widetilde{\mathcal O}(\epsilon^{-2})$ complexity
in the zero-noise regime. {Numerical experiments on
sparse phase retrieval, linear eigenvalue problems, and neural network
regression show that the four variants of {\sf ALFPG} generally outperform state-of-the-art proximal gradient baselines.}
\end{abstract}


\myvspace
\section{Introduction}
\myvspace

We consider composite nonconvex optimization of the form
\begin{equation}\label{eq:main}
    \min_{\x\in\mathbb{R}^n}
    F(\x):=f(\x)+h(\x),
    \qquad
    f(\x)=
    \begin{cases}
        \displaystyle \frac{1}{N}\sum_{i=1}^{N}f_i(\x),
        & \text{finite-sum setting},\\[2mm]
        \displaystyle \mathbb{E}_{\xi\sim\mathcal{D}}
        [f(\x;\xi)],
        & \text{expectation setting}.
    \end{cases}
\end{equation}
The smooth loss $f(\cdot)$ may be nonconvex, and the regularizer
$h(\cdot)$ is proper and lower semicontinuous. We assume that the required proximal subproblems admit exactly and efficiently computable global minimizers (Assumption \ref{ass:prox}). For finite-sum settings, $h(\cdot)$ can be nonconvex and non-Lipschitz under the stated {bounded-sublevel-set condition}, while for expectation settings, $h(\cdot)$ is required to be convex or Lipschitz continuous on its effective domain (Assumption~\ref{ass:regularizer}).

\paragraph{Composite Optimization and Structural Regularization.}
The composite formulation separates data fitting from
structural regularization.
The term $h(\cdot)$ can promote sparsity or encode constraints
through an extended-valued indicator function, allowing
the optimization model to incorporate prior information
about the desired solution.
Applications include sparse phase retrieval
\cite{Cai24,shechtman2014gespar},
optimization with orthogonality constraints
\cite{wen2013feasible},
and quantized model training \cite{BaiWL19}. These examples motivate methods that accommodate nonconvex losses and structural regularizers while keeping gradient evaluations and parameter tuning manageable.

\paragraph{From Adaptivity to Lipschitz-Free Updates.}
Proximal gradient methods exploit this structure by
combining gradient evaluations of $f(\cdot)$ with proximal
updates involving $h(\cdot)$.
Their stepsizes are commonly chosen using a smoothness
constant or a backtracking procedure.
Adaptive gradient and proximal gradient methods use information
collected during optimization to adjust their updates
\cite{duchi2011adaptive,ward2020adagrad,yun2021adaptive}.
For some adaptive methods, the convergence guarantee
still requires choosing the base stepsize or another algorithmic
parameter using a smoothness constant. We call an update \emph{Lipschitz-free} if it runs without line search
and its parameters can be chosen to satisfy the convergence conditions
without knowing a smoothness constant. The analysis may still assume
smoothness.

\paragraph{Challenges in Adaptive Proximal Updates.}
Adaptation for composite objectives should account for
the interaction between the loss and the regularizer.
At a composite stationary point $\x_\star$, we have
$0\in\nabla f(\x_\star)+\partial h(\x_\star)$, and
$\nabla f(\x_\star)$ can be nonzero. Consequently, the magnitude of the smooth gradient alone
need not indicate proximity to composite stationarity.
Proximal displacements incorporate the effect of $h$,
but relating them to stationarity requires controlling
the evolving metric. In the finite-sum setting, the challenge is to control how the
SPIDER estimation error accumulates between full-gradient refreshes
when extrapolation is coupled with an adaptive diagonal metric.
In the expectation setting, the challenge is to combine noise-sensitive
weights with stepsize scaling based on proximal residuals, without
knowing the smoothness constant or the noise level. The analysis must
account for the dependence between the adaptive weights and the
gradient-tracking errors.

\paragraph{Summary of Main Contributions.}
We propose Adaptive Lipschitz-Free Proximal Gradient
({\sf ALFPG}), a framework that uses proximal displacements
as a common principle for stepsize adaptation.
We measure stationarity by $\mathcal{R}(\x):=\operatorname{dist}(0,\partial F(\x))= \operatorname{dist}\!\bigl(0,\nabla f(\x)+\partial h(\x)\bigr)$, where $\partial$ denotes the limiting subdifferential.
Our main contributions are as follows.

\begin{enumerate}[leftmargin=17pt,itemsep=0.15ex,topsep=0.2ex]

\item \textbf{Adaptive Updates Using Proximal Displacements.} For finite-sum problems, {\sf ALFPG-FS} updates a diagonal
    preconditioner using metric-scaled iterate displacements and combines this metric with Nesterov extrapolation
    and the SPIDER estimator.
    Its extrapolation coefficient is coupled to the
    evolving metric to support the descent analysis.
    The full-batch specialization is denoted by {\sf ALFPG-D}.
    For expectation problems, {\sf ALFPG-MVR1} and {\sf ALFPG-MVR2} (MVR: Momentum-based Variance Reduction) combine normalized proximal displacements with
    gradient-discrepancy statistics to adapt their
    scalar stepsizes and momentum weights.
    Each iteration uses one fresh sample, with one
    stochastic-gradient evaluation for MVR1 and at most
    two evaluations for MVR2.

\item \textbf{Finite-Sum Guarantees with General Regularizers.}
Under the stated finite-sum assumptions, including {bounded
sublevel sets and separate lower bounds on $f$ and $h$}, {\sf ALFPG-FS} achieves expected $\epsilon$-stationarity using $\mathcal O(\epsilon^{-2})$
proximal iterations and $\mathcal O(N+\sqrt{N}\epsilon^{-2})$
component-gradient evaluations (Theorem~\ref{the:it:stoc}). Its deterministic variant, {\sf ALFPG-D}, achieves $\epsilon$-stationarity
using $\mathcal O(N\epsilon^{-2})$ component-gradient evaluations. Both methods combine Lipschitz-free diagonal adaptation with
extrapolation and accommodate nonconvex regularizers.

\item \textbf{Noise-Adaptive Guarantees for Expectation Problems.}
    {\sf ALFPG-MVR1} requires smoothness of the expected loss,
    whereas {\sf ALFPG-MVR2} additionally requires smoothness
    of the sample losses to control the same-sample correction. Under the respective assumptions of
Theorems~\ref{theorem:MVR1:main:paper} and~\ref{theorem:MVR2:main:paper}, including unbiased stochastic gradients with bounded variance and
almost-surely bounded norms and the regularizer condition in
Assumption~\ref{ass:regularizer}, {\sf ALFPG-MVR1} and {\sf ALFPG-MVR2} require
    $\widetilde{\mathcal{O}}
    (\epsilon^{-2}+\sigma^2\epsilon^{-4})$ and
    $\widetilde{\mathcal{O}}
    (\epsilon^{-2}+\sigma\epsilon^{-3})$
    stochastic-gradient evaluations, respectively,
    to obtain a uniformly sampled proximal iterate
    $\widehat{\x}$ satisfying
    $\mathbb{E}[\mathcal{R}(\widehat{\x})]\leq\epsilon$.
    Here, $\sigma^2$ bounds the gradient noise variance,
    and $\widetilde{\mathcal{O}}$ hides polylogarithmic factors. Both methods run without the variance bound as an input, and both
complexities reduce to $\widetilde{\mathcal{O}}(\epsilon^{-2})$
when $\sigma=0$.

\end{enumerate}

 \begin{table}[!t]
\centering
\begin{threeparttable}
\caption{Comparison among existing methods for composite nonconvex
finite-sum minimization.}
\label{tab:finite_sum}

\footnotesize
\setlength{\tabcolsep}{0.5pt}
\renewcommand{\arraystretch}{1.10}

\begin{tabular}{lccccc}
\toprule[1.2pt]
& \makecell{Lipschitz-free\\step size}
& \makecell{Regularizer\\class\tnote{c}}
& \makecell{Nesterov\\extrapolation}
& \makecell{Adaptive diagonal\\scaling}
& \makecell{Gradient-oracle\\complexity}
\\
\midrule
{\sf APG} \cite{li2015accelerated}
& \xmark & NC & \cmark & \xmark
& $\OO(N\epsilon^{-2})$ \\

{\sf SVRG-APG} \cite{li2017convergence}
& \xmark & NC & \cmark & \xmark
& Unknown\tnote{a} \\

{\sf ProxSVRG} \cite{Reddi2016}
& \xmark & CVX & \xmark & \xmark
& $\OO(N+N^{2/3}\epsilon^{-2})$ \\

{\sf SPIDER} \cite{fang2018spider}
& \xmark & $h\equiv0$ & \xmark & \xmark
& $\OO(N+\sqrt{N}\epsilon^{-2})$ \\

{\sf SpiderBoost-M} \cite{wang2019spiderboost}
& \xmark & CVX & \cmark & \xmark
& $\OO(N+\sqrt{N}\epsilon^{-2})$ \\

{\sf ProxSARAH} \cite{pham2020proxsarah}
& \xmark\tnote{b} & CVX & \xmark & \xmark
& $\OO(N+\sqrt{N}\epsilon^{-2})$ \\

{\sf SPGR} \cite{xu2019non}
& \xmark & NC & \xmark & \xmark
& $\OO(N+\sqrt{N}\epsilon^{-2})$ \\

\addlinespace[2pt]

{\sf AdaGrad-Norm} \cite{ward2020adagrad}
& \cmark & $h\equiv0$ & \xmark & \xmark
& $\OO(N\epsilon^{-2})$ \\

{\sf AGD} \cite{KavisLC22}
& \cmark & $h\equiv0$ & \cmark & \xmark
& $\OO(N\epsilon^{-2})$ \\

{\sf ADA-SPIDER} \cite{kavis2022adaptive}
& \cmark & $h\equiv0$ & \xmark & \xmark
& $\widetilde{\OO}(N+\sqrt{N}\epsilon^{-2})$ \\

{\sf AC-PG} \cite{alfpglan2024projected}
& \cmark & IND & \xmark & \xmark
& $\OO(N\epsilon^{-2})$\tnote{d} \\

{\sf AdaPGNC} \cite{alfpgye2025adapgnc}
& \cmark & CVX & \xmark & \xmark
& $\OO(N\epsilon^{-2})$ \\

\midrule
\rowcolor{gray!18}[\tabcolsep][\tabcolsep]
\textbf{{\sf ALFPG-D}} [ours]
& \cmark & NC & \cmark & \cmark
& $\OO(N\epsilon^{-2})$ \\

\rowcolor{gray!18}[\tabcolsep][\tabcolsep]
\textbf{{\sf ALFPG-FS}} [ours]
& \cmark & NC & \cmark & \cmark
& $\OO(N+\sqrt{N}\epsilon^{-2})$ \\
\bottomrule[1.2pt]
\end{tabular}

\begin{tablenotes}[flushleft]
\footnotesize
\setlength{\itemsep}{1pt}
\item[a] This work established linear convergence under a KL condition,
but no general stationarity oracle complexity is provided.
\item[b] This algorithm does not qualify as using adaptive step sizes
in our strict sense, as it relies on the Lipschitz constant.
\item[c] NC: possibly nonconvex; CVX: convex; $h\equiv0$: no regularizer;
IND: $h(\x)$ is the indicator function of a convex set.
\item[d] AC-PG uses objective values and full gradients
($N$ component-gradient evaluations per gradient) on a compact convex set.
\end{tablenotes}
\end{threeparttable}
\end{table}

\begin{table}[!t]
\centering
\begin{threeparttable}
\caption{Comparison among existing methods for composite nonconvex
expectation minimization.}
\label{tab:expectation}

\footnotesize
\setlength{\tabcolsep}{1.4pt}
\renewcommand{\arraystretch}{1.10}

\begin{tabular}{lccccc}
\toprule[1.2pt]
& \makecell{Lipschitz-free\\step size}
& \makecell{Regularizer\\class\tnote{a}}
& \makecell{Function\\$f(\cdot)$\tnote{b}}
& \makecell{Maximum\\batch size}
& \makecell{Gradient-oracle\\complexity}
\\
\midrule
{\sf SPIDER} \cite{fang2018spider}
& \xmark & $h\equiv0$ & $\LLL_+$
& $\OO(1+\sigma^2\epsilon^{-2})$ & $\OO(\epsilon^{-3})$ \\

{\sf SpiderBoost} \cite{wang2019spiderboost}
& \xmark & CVX & $\LLL_+$
& $\OO(1+\sigma^2\epsilon^{-2})$ & $\OO(\epsilon^{-3})$ \\

{\sf RSPG} \cite{ghadimi2016mini}
& \xmark & CVX & $\LLL$
& $\OO(1+\sigma^2\epsilon^{-2})$ & $\OO(\epsilon^{-4})$ \\

{\sf ProxSARAH} \cite{pham2020proxsarah}
& \xmark & CVX & $\mathbb{MS}$
& $\OO(1+\sigma^2\epsilon^{-2})$ & $\OO(\epsilon^{-3})$ \\

{\sf SMM} \cite{davis2019stochastic}
& \xmark & CVX & $\LLL_-$
& $\OO(1)$ & $\OO(\epsilon^{-4})$ \\

{\sf ProxHSGD} \cite{tran2022hybrid}
& \xmark & CVX & $\mathbb{MS}$
& $\OO(1+\sigma^3\epsilon^{-1})$ & $\OO(\epsilon^{-3})$ \\

{\sf MB-SPG} \cite{xu2019non}
& \xmark & NC & $\LLL$
& $\OO(1+\sigma^2\epsilon^{-2})$ & $\OO(\epsilon^{-4})$ \\

{\sf SPGR} \cite{xu2019non}
& \xmark & NC & $\LLL_+$
& $\OO(1+\sigma^2\epsilon^{-2})$ & $\OO(\epsilon^{-3})$ \\

{\sf STORM} \cite{cutkosky2019momentum}
& \xmark & $h\equiv0$ & $\LLL_+$
& $\OO(1)$ & $\OO(\epsilon^{-3})$ \\

{\sf PStorm} \cite{xu2023momentum}
& \xmark & CVX & $\mathbb{MS}$
& $\OO(1)$ & $\OO(\epsilon^{-3})$ \\

{\sf ProxGen} \cite{yun2021adaptive}
& \xmark & NC & $\LLL$
& $\OO(1+\sigma^2\epsilon^{-2})$ & $\OO(\epsilon^{-4})$ \\

{\sf SuperADAM-}$\LLL$ \cite{huang2021super}
& \xmark & IND & $\LLL$
& $\OO(1)$ & $\widetilde{\OO}(\epsilon^{-4})$ \\

{\sf SuperADAM-}$\LLL_+$ \cite{huang2021super}
& \xmark & IND & $\LLL_+$
& $\OO(1)$ & $\widetilde{\OO}(\epsilon^{-3})$ \\

\addlinespace[2pt]

{\sf AdaGrad-Norm} \cite{ward2020adagrad}
& \cmark & $h\equiv0$ & $\LLL$
& $\OO(1)$ & $\widetilde{\OO}(\epsilon^{-2}+\sigma^2\epsilon^{-4})$ \\

{\sf STORM+} \cite{levy2021stormplus}
& \cmark & $h\equiv0$ & $\LLL_+$
& $\OO(1)$ & $\widetilde{\OO}(\epsilon^{-2}+\sigma\epsilon^{-3})$ \\

{\sf OptEMA} \cite{alfpgyuan2026optema}
& \cmark & $h\equiv0$ & $\LLL$
& $\OO(1)$ & $\widetilde{\OO}(\epsilon^{-2}+\sigma^2\epsilon^{-4})$ \\

{\sf AC-SPG} \cite{alfpglan2024projected}
& \cmark & IND & $\LLL_+$
& $\OO(1+\sigma^2\epsilon^{-2})$
& $\OO(\epsilon^{-2}+\sigma^2\epsilon^{-4})$\tnote{c} \\

{\sf AC-VR-SPG} \cite{alfpglan2024projected}
& \cmark & IND & $\LLL_+$
& $\OO(1+\sigma^2\epsilon^{-2})$
& $\OO(\epsilon^{-2}+\sigma\epsilon^{-3})$\tnote{c} \\

\midrule
\rowcolor{gray!18}[\tabcolsep][\tabcolsep]
\textbf{{\sf ALFPG-MVR1}} [ours]
& \cmark & CVX/LC & $\LLL$
& $\OO(1)$
& $\widetilde{\OO}(\epsilon^{-2}+\sigma^2\epsilon^{-4})$ \\

\rowcolor{gray!18}[\tabcolsep][\tabcolsep]
\textbf{{\sf ALFPG-MVR2}} [ours]
& \cmark & CVX/LC & $\LLL_+$
& $\OO(1)$
& $\widetilde{\OO}(\epsilon^{-2}+\sigma\epsilon^{-3})$ \\
\bottomrule[1.2pt]
\end{tabular}

\begin{tablenotes}[flushleft]
\footnotesize
\setlength{\itemsep}{1pt}
\item[a] NC: $h(\x)$ is possibly nonconvex; CVX: $h(\x)$ is convex;
IND: $h(\x)$ is the indicator function of a convex set;
CVX/LC: $h(\x)$ is convex or Lipschitz continuous on the effective domain.
\item[b] $\LLL$: $f(\x)$ is $L$-smooth;
$\LLL_+$: $f(\x;\xi)$ is component-wise $L$-smooth;
$\LLL_-$: $f(\x)$ is weakly convex;
$\mathbb{MS}$: the mean-squared smoothness condition
displayed in the related-work discussion.
\item[c] Both methods use noise-dependent batches and query stochastic
function values and gradients on compact convex sets.
\end{tablenotes}
\end{threeparttable}
\end{table}

\myvspace
 \section{Related Work}
 \myvspace

We review proximal gradient methods for nonconvex composite optimization
in finite-sum and expectation settings, focusing on extrapolation,
adaptive stepsizes, variance reduction, and their oracle complexity
guarantees.

\paragraph{Proximal Gradient Methods for Nonconvex Composite Optimization.}
Proximal gradient methods combine gradient evaluations of the smooth
loss with proximal updates for the regularizer. For a nonconvex smooth
loss and a convex regularizer, {\sf RSPG} gives stochastic-gradient
oracle bounds in terms of a generalized projected-gradient mapping
\cite{ghadimi2016mini}. \citet{xu2019non} establish subdifferential-residual
guarantees for possibly nonconvex and nonsmooth regularizers.
In the expectation setting with a fixed positive noise level, their fixed-batch {\sf MB-SPG} method uses
batches of order $\epsilon^{-2}$ and achieves
$\mathcal{O}(\epsilon^{-4})$ stochastic-gradient complexity, and their
recursive-estimator method {\sf SPGR} achieves
$\mathcal{O}(\epsilon^{-3})$ under additional sample-wise smoothness.
{\sf ALFPG} studies Lipschitz-free stepsize adaptation within this
class of composite problems.

\paragraph{Extrapolation in Nonconvex Proximal Optimization.}
Accelerated proximal gradient methods achieve improved
objective-value convergence rates in convex optimization
\cite{beck2009fast}.
Methods such as {\sf iPiano} \cite{ochs2014ipiano}, {\sf APG} \cite{li2015accelerated}, and {\sf iPALM} \cite{pock2016inertial} incorporate inertial or extrapolation steps into nonconvex proximal
optimization. \citet{li2017convergence} analyze {\sf APGnc} and {\sf SVRG-APGnc}.
For the latter, they establish a linear convergence rate for the
expected objective gap at the outer iterates under a global
Kurdyka--{\L}ojasiewicz condition with exponent $1/2$.

\paragraph{Variance Reduction for Finite Sums.}
{\sf SVRG} \cite{NIPS2013TongZ} and {\sf SAGA} \cite{defazio2014saga} use finite-sum structure to reduce
gradient-estimation variance. Their proximal extensions, {\sf ProxSVRG} and {\sf ProxSAGA}, have
complexity guarantees for nonconvex smooth losses with nonsmooth
convex regularizers \cite{Reddi2016}.
{\sf SARAH} \cite{nguyen2017sarah} and {\sf SPIDER} \cite{fang2018spider} use recursive gradient estimators,
and {\sf PAGE} \cite{li2021page} combines recursive updates with probabilistic refreshes.
For convex regularizers, {\sf Prox-SpiderBoost}, {\sf Prox-SpiderBoost-M} \cite{wang2019spiderboost},
and {\sf ProxSARAH} \cite{pham2020proxsarah} achieve $\mathcal{O}(N+\sqrt{N}\epsilon^{-2})$ component-gradient complexity for a proximal-gradient mapping criterion. {\sf SPGR} \cite{xu2019non} attains the same dependence on $N$ and
$\epsilon$ while permitting nonconvex regularizers and directly
controlling a subdifferential residual. The dynamic stepsize schedule of {\sf ProxSARAH} depends on the
smoothness constant $L$.

\paragraph{Adaptive Gradients and Proximal Preconditioning.}
{\sf AdaGrad} \cite{duchi2011adaptive} and {\sf Adam} \cite{KingmaB14} adapt their updates using historical
gradient information. \citet{reddi2018convergence} identify convergence failures of {\sf Adam} and develop {\sf AMSGrad}. Nonconvex convergence guarantees are established for {\sf AdaGrad-Norm} \cite{ward2020adagrad} and classes of {\sf Adam}-type methods \cite{chen2019convergence}. {\sf AdaSPIDER} \cite{kavis2022adaptive} combines Lipschitz-free stepsizes with recursive
variance reduction for smooth nonconvex finite sums.
{\sf ProxGen} \cite{yun2021adaptive} combines adaptive preconditioning with momentum and allows
nonconvex regularizers. {\sf SUPER-ADAM} \cite{huang2021super} uses
adaptive matrices for constrained nonconvex optimization. The convergence guarantees for both methods
impose smoothness-dependent parameter restrictions. \citet{alfpglan2024projected} develop {\sf AC-PG}, {\sf AC-SPG},
and {\sf AC-VR-SPG}, which use function values and gradients to estimate
local curvature and select scalar stepsizes over compact convex sets.
\citet{alfpgye2025adapgnc} develop {\sf AdaPGNC}, which uses local upper
and lower curvature estimates to select scalar stepsizes for deterministic
composite problems with convex $h$.
{\sf ALFPG} adapts its stepsizes using proximal displacements. Its
expectation variants also use gradient discrepancies and draw one
fresh sample per iteration.

\paragraph{Recursive Momentum for Expectation Problems.}
{\sf STORM} \cite{cutkosky2019momentum} uses a recursive correction based on the same sample evaluated at two consecutive iterates.
{\sf STORM+} \cite{levy2021stormplus} gives a parameter-free choice of
the stepsize and momentum weight. {\sf Ada-STORM} \cite{jiang2024adaptive} achieves an
$\mathcal{O}(T^{-1/3})$ rate for smooth nonconvex objectives without
uniform bounds on gradients or function values, using an initialization
batch of order $T^{1/3}$. Among proximal variants, {\sf ProxSARAH} \cite{pham2020proxsarah} uses
periodic snapshot batches. At a fixed noise level, the $\mathcal{O}(\epsilon^{-3})$ guarantee for {\sf ProxHSGD}
\cite{tran2022hybrid} uses an accuracy-dependent initialization batch,
and {\sf PStorm} \cite{xu2023momentum} achieves the same order with
constant batches throughout. {\sf ALFPG-MVR1} uses an exponential moving average, and
{\sf ALFPG-MVR2} adds a same-sample correction. Both draw one fresh
sample per iteration, including initialization. Their convergence
guarantees assume almost-surely bounded stochastic gradients.

\paragraph{Smoothness Assumptions and Noise Adaptivity.} Uniform sample-wise smoothness implies mean-squared smoothness,
which in turn implies smoothness of the expected loss when stochastic
gradients are unbiased. The mean-squared condition is
$\mathbb{E}_{\xi}\big[
\|\nabla f(\x;\xi)-\nabla f(\y;\xi)\|^2\big]
\leq L^2\|\x-\y\|^2$ for all $\x,\y$.
It is used in the analyses of {\sf ProxSARAH} \cite{pham2020proxsarah}, {\sf ProxHSGD} \cite{tran2022hybrid},
and {\sf PStorm} \cite{xu2023momentum}.
Our {\sf ALFPG-MVR1} result assumes smoothness of the expected loss.
The {\sf ALFPG-MVR2} result additionally assumes uniform sample-wise
smoothness. Variance-sensitive guarantees and recovery of deterministic rates
have also been studied in smooth adaptive optimization
\cite{ward2020adagrad,cutkosky2019momentum}. {\sf OptEMA} \cite{alfpgyuan2026optema} uses gradient statistics
to adapt a first- or second-moment averaging coefficient. When $h\equiv0$,
both {\sf OptEMA} and {\sf ALFPG-MVR1} can use one fresh sample per
iteration and, under their respective assumptions, require
$\widetilde{\OO}(\epsilon^{-2}+\sigma^2\epsilon^{-4})$
stochastic-gradient evaluations to obtain an iterate $\widehat{\x}$
satisfying $\mathbb{E}\|\nabla f(\widehat{\x})\|\leq\epsilon$.


\paragraph{Stationarity Measures and Oracle Complexity.}
We measure stationarity by
$\mathcal{R}(\x):=\operatorname{dist}(0,\partial F(\x))$.
Related guarantees use proximal-gradient mappings or, for weakly convex
objectives, Moreau-envelope gradients \cite{davis2019stochastic}.
A proximal-gradient mapping bound can imply a residual bound at the
corresponding proximal point, with suitable control of stochastic
estimation errors.
Gradient oracle complexity counts all individual gradient evaluations,
including initialization, refreshes, and evaluations at different points
using the same sample. We count component-gradient evaluations for
finite sums and stochastic-gradient evaluations for expectation problems.

Tables~\ref{tab:finite_sum} and~\ref{tab:expectation}
compare representative methods for the finite-sum and expectation
settings, respectively, in terms of their regularizer assumptions,
stepsize requirements, and oracle complexity.

\textbf{Notation.}
For a nonempty set $S$, write
$\|S\|:=\dist(\zero,S)=\inf_{u\in S}\|u\|$.
The symbols $\odot$ and $\div$ denote element-wise multiplication
and division, respectively, and $\min(\x)$ denotes the smallest
entry of $\x$. We use $[\x_t]_i$ to denote the $i$th entry of $\x_t$, with $t$ indexing iterations. For $\v\geq\zero$, define the weighted seminorm
$\|\x\|_{\v}:=\sqrt{\la\v,\x\odot\x\ra}$.
The notation $\widetilde{\OO}(\cdot)$ hides polylogarithmic factors.
Further notation, technical preliminaries, and relevant lemmas
appear in Appendix~\ref{app:sect:notations}.

 \myvspace
\section{The Proposed Algorithms}
\label{sect:alg}
\myvspace
 
This section presents the Adaptive Lipschitz-Free Proximal Gradient ({\sf ALFPG}) framework for Problem~\eqref{eq:main}. {\sf ALFPG-FS} addresses finite-sum problems, with {\sf ALFPG-D} as its full-batch specialization. {\sf ALFPG-MVR1} and {\sf ALFPG-MVR2} address expectation problems.

The algorithms require the following exact-proximal-computation assumption.

\begin{assumption}\label{ass:prox}
(Exact proximal computations) For $\gamma>0$, $\a\in\mathbb{R}^n$, and
$\v\in\mathbb{R}_{++}^n$, define $\Prox_{\gamma h}(\a;\v)
:=\arg\min_{\x\in\mathbb{R}^n}
\left\{h(\x)+\frac{1}{2\gamma}\|\x-\a\|_{\v}^{2}\right\}$. For every $\a\in\mathbb{R}^n$ and $\gamma>0$, we assume that a global minimizer of this subproblem can be computed exactly and efficiently: (i) for every $\v\in\mathbb{R}_{++}^n$ in {\sf ALFPG-FS} and its full-batch specialization {\sf ALFPG-D}; (ii) for $\v=\mathbf{1}$ in {\sf ALFPG-MVR1} and {\sf ALFPG-MVR2}. A minimizer is selected whenever the solution is nonunique. When $\v=\one$, we denote $\Prox_{\gamma h}(\a;\v)$ as $\Prox_{\gamma h}(\a)$.
\end{assumption}

\begin{remark}
\ding{182} Assumption \ref{ass:prox} is commonly employed in nonconvex proximal gradient methods. \ding{183} When $\v=\one$, the diagonal preconditioner reduces to the identity preconditioner. Assumption \ref{ass:prox} holds for certain regularizers $h(\x)$. Common examples include capped-$\ell_1$ penalty \cite{zhang2010analysis}, log-sum penalty \cite{candes2008enhancing}, minimax concave penalty \cite{zhang2010nearly}, Geman penalty \cite{geman1995nonlinear}, $\ell_p$ regularization with $p\in\{0,\tfrac{1}{2},\tfrac{2}{3},1\}$, and indicator functions for cardinality constraints, orthogonality constraints in matrices, and rank constraints in matrices. \ding{184} When $\v$ is a general positive vector, the variable-metric operator can still be evaluated for certain coordinate-wise separable regularizers $h(\x)$. Examples include $\ell_p$ regularization with $p\in\{0,\tfrac{1}{2},\tfrac{2}{3},1\}$ (with or without bound constraints) \cite{yun2021adaptive}, and the W-shaped regularizer \cite{BaiWL19}.

\end{remark}

\begin{algorithm}[!t]

\caption{Proposed {\sf ALFPG-FS} for Finite-Sum Settings.}

\small

\textbf{Input:} $\x_0\in\dom(F)$. Let $\underline{\rm{v}},\alpha>0$, $\beta\geq 0$, and $\theta \in[0,1)$. Choose integers $q,b\geq 1$ (take $q=b=\lceil\sqrt{N}\rceil$).

\textbf{Initialize:} Let $\x_{-1}=\y_{0}=\x_0$. Set $\v_{-1} =\v_{0} = \underline{\rm{v}}\one$, $\sigma_{-1}=\theta$.

For $t=0,\ldots,T$:

\begin{enumerate}[leftmargin=17pt,itemsep=0.15ex,topsep=0.2ex]

\item If $\mathrm{mod}(t,q)\ne0$, draw a mini-batch
$\mathcal I_t$ of $b$ indices uniformly with replacement from
$\{1,\ldots,N\}$. Set 
\beq \label{eq:SPIDER:v}
\g_t = \left\{
  \begin{array}{ll}
     \nabla f(\y_t) , & \hbox{$\textrm{mod}(t,q)=0$;} \\
   \g_{t-1} + \nabla f(\y_t;\mathcal{I}_t) - \nabla f(\y_{t-1};\mathcal{I}_{t}), & \hbox{else.}
  \end{array}
\right.
\eeq
When $q=1$, every iteration uses the full gradient,
giving {\sf ALFPG-D}.


\item Let $\x_{t+1} \in\Prox_h(\y_t - \g_t \div {\v_t};\v_t)$, $\d_t := \x_{t+1}-\x_{t}$.

\item Let $\s_t:= \alpha \|\r_t\|_2^2\cdot\one + \beta \r_t \odot \r_t$, where $\r_t := \v_t \odot \d_t $.

\item Set $\v_{t+1} = \ts \sqrt{ \v_{0} \odot \v_{0} + \ts \sum_{i=0}^{t} \s_i}$.

\item Let $\sigma_t := \theta (1-   \sigma_{t-1}) \cdot \min( \v_t \div \v_{t+1})$.

\item Set $\y_{t+1} = \x_{t+1} + \sigma_t \d_t$.

\end{enumerate}

\label{alg:main}

\end{algorithm}

\subsection{{\sf ALFPG-FS} for Finite-Sum Settings}

{\sf ALFPG-FS} combines a displacement-based diagonal metric, metric-coupled Nesterov-type extrapolation, and the SPIDER
estimator for finite-sum composite optimization. When $q=1$, it specializes to the full-batch method
{\sf ALFPG-D}. Both methods are specified in
Algorithm~\ref{alg:main}.

\paragraph{Design and Motivation of Algorithm~\ref{alg:main}.} {At the extrapolated point $\y_t$, the algorithm computes the full
gradient when $\mathrm{mod}(t,q)=0$. Otherwise, it applies the SPIDER
recursion~\eqref{eq:SPIDER:v} \cite{fang2018spider} using the same
mini-batch $\mathcal I_t=(i_{t,1},\ldots,i_{t,b})$ at $\y_t$ and
$\y_{t-1}$. The indices are sampled uniformly with replacement, and $\nabla f(\y;\mathcal I_t):=b^{-1}\sum_{j=1}^{b}\nabla f_{i_{t,j}}(\y)$. After the proximal step, the metric-scaled displacement $\r_t=\v_t\odot(\x_{t+1}-\x_t)$ determines the next
metric through $[\v_{t+1}]_j^2=\underline{\rm v}^{\,2}+\alpha\sum_{i=0}^{t}\|\r_i\|_2^2
+\beta\sum_{i=0}^{t} [\r_i]_j^2$. The $\alpha$-term accumulates overall motion, and the $\beta$-term
accumulates coordinate-wise motion. Together, they determine the
coordinate-wise inverse stepsizes. The extrapolation coefficient $\sigma_t$ is coupled to the change
in the metric to support the descent analysis.}

\paragraph{Relationship to AdaGrad-Norm and Main Novelty of Algorithm~\ref{alg:main}.} Consider $h\equiv0$, $\beta=\theta=0$, and $q=1$. Then
$\y_t=\x_t$, $\g_t=\nabla f(\x_t)$, and $\v_t=c_t\one$ for some $c_t>0$. The proximal step implies $\r_t=-\g_t$, so the updates reduce to $c_{t+1}^{\,2}=c_t^{\,2}+\alpha\|\g_t\|_2^2$, $\x_{t+1}=\x_t-\frac{1}{c_t}\g_t$. After rescaling the accumulator and base stepsize, {\sf AdaGrad-Norm}
follows the same accumulation rule but uses $c_{t+1}$ in the current
update \cite{ward2020adagrad}. In this special case, the two methods
differ in whether the current gradient is included in the denominator
of the current step. For general $h$, $\r_t$ reflects the proximal step and need not equal $-\g_t$. The proposed framework combines this displacement-based metric with metric-coupled extrapolation and, for $q>1$, SPIDER estimation.

\paragraph{Additional Remarks on Algorithm \ref{alg:main}.} \ding{182} The metric update is equivalently written
$\v_{t+1}=\sqrt{\v_t\odot\v_t+\s_t}$. The choice $\beta=0$ gives a scalar metric, while $\beta>0$ permits different
coordinate-wise stepsizes. \ding{183} The generalized proximal step permits
nonsmooth and nonconvex $h$ under Assumption~\ref{ass:prox}, which
requires an exactly computable global minimizer of the relevant
diagonal-metric subproblem. \ding{184} The parameter choice $\alpha>0$ ensures the uniform
metric-ratio bound
$\max(\v_t)/\min(\v_t)\leq
\kappa:=1+\sqrt{\beta/\alpha}$
(Lemma~\ref{lemma:prop:w}(b)).

\begin{algorithm}[!t]
\caption{{Proposed {\sf ALFPG-MVR1} and {\sf ALFPG-MVR2} for Expectation Settings.}}
\label{alg:mvr}
\small

\textbf{Input:}
$\theta>0$, $\x_1\in\dom(h)$.

\textbf{Initialize:}
$\x_0=\x_1$, $\m_0=\zero$, $\rho_{-1}=\rho_0=1$, and $(p,\lambda)=
\begin{cases}
(1/2,1/2), & \text{MVR1},\\[1mm]
(2/3,1/4), & \text{MVR2}.
\end{cases}$.

For $t=1,\ldots,T$:
\begin{enumerate}[leftmargin=17pt,itemsep=0.15ex,topsep=0.2ex]

\item Draw a fresh sample $\xi_t\sim\mathcal D$. At $t=1$, set $\xi_0:=\xi_1$.

\item Compute $\g_t:=\nabla f(\x_t;\xi_t)$. For MVR2, also compute $\g_t^-:=\nabla f(\x_{t-1};\xi_t)$. \\
    At $t=1$, set $\g_0:=\g_1$; MVR2 reuses $\g_1^-=\g_1$. 

\item Set $\rho_t :=
\min_{1\le k\le t}
\left(
\frac{
1+\max_{1\le i\le k}\chi_i^2
}{
1+\sum_{i=1}^k\chi_i^2
}
\right)^{p}$, where $\chi_t:=
\begin{cases}
\|\g_{t-1}-\g_t\|, & \text{MVR1}, \\[1mm]
\|\g_{t-1}-\g_t^-\|, & \text{MVR2}.
\end{cases}$.

\item Let $\alpha_t:=\rho_{t-1}$, and $\m_t :=
\begin{cases}
\alpha_t\g_t + (1-\alpha_t)\m_{t-1},
& \text{MVR1},
\\[1mm]
\alpha_t\g_t + (1-\alpha_t) (\m_{t-1} +  \g_t-\g_t^- ),
& \text{MVR2}.
\end{cases}$.

\item Set $\beta_t:= \left(1+\sum_{i=1}^{t-1}\|\r_i\|^2 \right)^{-1/2}$, and $\gamma_t =\theta \beta_t \alpha_t^{\lambda}$.

\item Let $\x_{t+1} \in \prox_{\gamma_t h}
\bigl(\x_t-\gamma_t\m_t\bigr)$, and $\r_t:= \frac{\x_t-\x_{t+1}}{\gamma_t}$.

\end{enumerate}


\end{algorithm}

\subsection{{\sf ALFPG-MVR1} and {\sf ALFPG-MVR2} for Expectation Settings}

{\sf ALFPG-MVR1} uses an Exponential Moving Average (EMA) estimator
with one stochastic-gradient evaluation per iteration \cite{KingmaB14},
whereas {\sf ALFPG-MVR2} uses a STORM-type recursive estimator
with at most two evaluations on the same fresh sample
\cite{cutkosky2019momentum}.
The analysis requires expected-loss smoothness for the former
and additionally uniform sample-wise smoothness for the latter.
Both methods are specified in Algorithm~\ref{alg:mvr}.

\paragraph{Design and Motivation of Algorithm~\ref{alg:mvr}.} {Algorithm~\ref{alg:mvr} combines a stochastic gradient estimator,
a proximal step, and two adaptive factors:
$\alpha_t$ controls gradient mixing, while $\beta_t$ scales
the proximal step.
For {\sf ALFPG-MVR2}, $\chi_t$ compares gradients from two
samples at the same point $\x_{t-1}$ and thus reflects sampling noise.
To retain one stochastic-gradient evaluation per iteration, {\sf ALFPG-MVR1} instead uses
$\chi_t=\|\g_{t-1}-\g_t\|$ as a proxy, reusing the stored gradient.
This proxy reflects both sampling noise and iterate motion. We use a running maximum to set the discrepancy scale
and a running minimum to keep the weight nonincreasing.
Specifically,
$\rho_t:=
\min_{1\leq k\leq t}
\bigl(
\frac{1+\max_{1\leq i\leq k}\chi_i^2}
     {1+\sum_{i=1}^{k}\chi_i^2}
\bigr)^p$. For {\sf ALFPG-MVR2}, when the squared differences remain
on a comparable scale $\sigma^2$, this ratio suggests the
heuristic behavior
$\rho_t\approx\bigl(\frac{1+\sigma^2}{1+t\sigma^2}\bigr)^p$.
This interpretation is not an assumption used in the analysis.
Lemma~\ref{lemma:weight} controls the variation of $\rho_t$
and $1/\rho_t$ and provides a contraction inequality for
gradient tracking.
We use the lagged weight $\alpha_t=\rho_{t-1}$ and compute
$\beta_t$ from past proximal residuals, so both the averaging
coefficient and the stepsize
$\gamma_t=\theta\beta_t\alpha_t^\lambda$
are determined before sampling $\xi_t$.}

\paragraph{Relationship to Existing Methods and Main Novelty of Algorithm~\ref{alg:mvr}.} {\sf AdaGrad-Norm} \cite{ward2020adagrad} uses a scalar stepsize
denominator based on $\sum_{i=1}^{t}\|\g_i\|^2$. Our noise-sensitive weight $\rho_t$ uses accumulated squared
gradient differences: $\sum_{i=1}^{t}\|\g_i-\g_{i-1}\|^2$ for {\sf ALFPG-MVR1}
and $\sum_{i=1}^{t}\|\g_{i-1}-\g_i^-\|^2$
for {\sf ALFPG-MVR2}. For each variant, the accumulated sum and the running maximum
$\max_{1\leq i\leq t}\chi_i^2$ enter the definition of $\rho_t$
through the running-minimum rule. The separate factor
$\beta_t=(1+\sum_{i=1}^{t-1}\|\r_i\|^2)^{-1/2}$
uses proximal residuals.
These quantities serve different purposes:
gradient differences reflect sampling noise and, for
{\sf ALFPG-MVR1}, iterate motion, whereas proximal residuals
incorporate the regularizer into stepsize adaptation.
The {\sf ALFPG-MVR1} estimator is an exponential moving average
using one stochastic-gradient evaluation per iteration.
{\sf ALFPG-MVR2} uses a recursive correction evaluated at
$\x_t$ and $\x_{t-1}$ with the same fresh sample, as in
momentum-based variance reduction methods
\cite{cutkosky2019momentum,tran2022hybrid}.
Our design couples these estimators with the averaging
coefficient $\alpha_t=\rho_{t-1}$ based on gradient discrepancies
and the proximal stepsize
$\gamma_t=\theta\beta_t\alpha_t^\lambda$,
which also incorporates the accumulated proximal residuals.

\paragraph{Additional Remarks on Algorithm~\ref{alg:mvr}.} \ding{182} We use $(p,\lambda)=(1/2,1/2)$ for
{\sf ALFPG-MVR1} and $(2/3,1/4)$ for {\sf ALFPG-MVR2}.
Since $\alpha_t=\rho_{t-1}$ and $\beta_t$ depends only on past
proximal residuals, $\gamma_t=\theta\beta_t\alpha_t^\lambda$
is determined before sampling $\xi_t$.
\ding{183} {\sf ALFPG-MVR1} uses one stochastic-gradient
evaluation per iteration, whereas {\sf ALFPG-MVR2} uses
at most two evaluations with the same fresh sample.
Neither variant requires $L$ or the noise level $\sigma$
as an algorithmic input.
\ding{184} Both variants construct $\chi_t$ without additional
gradient evaluations.
{\sf ALFPG-MVR1} uses $\chi_t=\|\g_{t-1}-\g_t\|$,
which reflects both sampling noise and iterate motion.
{\sf ALFPG-MVR2} reuses
$\g_t^-=\nabla f(\x_{t-1};\xi_t)$ from its recursive correction
to form $\chi_t=\|\g_{t-1}-\g_t^-\|$, comparing the stored
and current-sample gradients at the same point $\x_{t-1}$. \ding{185} When $\sigma=0$, {\sf ALFPG-MVR2} has $\chi_t=0$ and $\rho_t=1$.
For {\sf ALFPG-MVR1}, the discrepancy reduces to
$\chi_t=\|\nabla f(\x_t)-\nabla f(\x_{t-1})\|
\leq L\|\x_t-\x_{t-1}\|$.
It therefore reflects only gradient variation along the iterates
and tends to zero as the iterate displacements vanish,
consistent with its interpretation as an inexpensive proxy
for the same-point discrepancy. Both variants attain the
$\widetilde{\mathcal O}(T^{-1/2})$ expected stationarity rate
in this noise-free case under their respective assumptions.

 \myvspace
\section{Oracle Complexity for {\sf ALFPG-D} and {\sf ALFPG-FS}}
\label{sect:iter}
 
\myvspace

This section analyzes the oracle complexity of {\sf ALFPG-D} and {\sf ALFPG-FS}. In {\sf ALFPG-FS}, $\x_t$ depends on the mini-batches sampled before iteration $t$.
We use $\EEE[\cdot]$ to denote expectation over the algorithmic randomness and make the following assumptions.

\begin{assumption}\label{ass:compact}
{ For every finite $a\in\mathbb R$, the sublevel set $\{\x\in\mathbb R^n:F(\x)\leq a\}$ is bounded.} Furthermore, $f_\star:=\inf_{\x}f(\x)>-\infty$ and
$h_\star:=\inf_{\x}h(\x)>-\infty$.
\end{assumption}

\begin{assumption}\label{ass:f}
Each $f_{j}(\cdot)$ is $L$-smooth, meaning that $\|\nabla f_{j}(\x) - \nabla f_{j}(\tilde{\x})\| \leq L \|\x-\tilde{\x}\|$ for all $\x,\tilde{\x}\in\mathbb R^n$ and $j\in[N]$. Consequently, $f(\cdot)$ is also $L$-smooth.

\end{assumption}

\begin{remark}
\ding{182} The bounded-sublevel-set condition in
Assumption~\ref{ass:compact} holds when the effective domain
of $F$ is bounded, and also allows an unbounded effective domain.
For example, $f_j\equiv0$ and $h(\x)=\|\x\|_1$ satisfy
Assumption~\ref{ass:compact} on $\mathbb R^n$.
\ding{183} Lemma~\ref{lemma:level-localization} bounds the
iterate norms uniformly over the iteration index.
This bound is used in the adaptive-metric estimates below. \ding{184} Assumption \ref{ass:f} is a standard requirement in the convergence analysis of nonconvex algorithms.
\end{remark}

\myvspace
\subsection{Convergence Analysis of {\sf {ALFPG-D}}}
\label{ALFPG:FS}
\myvspace

{A proof overview is provided in Appendix~\ref{subsect:D:overview}.} The following theorem establishes the iteration complexity.


\begin{theorem} \label{the:it}
(Proof in Section \ref{app:the:it}, \rm{Iteration Complexity}). {Under Assumptions~\ref{ass:prox}, \ref{ass:compact}, and~\ref{ass:f},
let $\{\x_t\}_{t=0}^{T+1}$ be generated by {\sf ALFPG-D}. Then $\ts \frac1{T+1}\sum_{t=0}^{T}
\operatorname{dist}(0,\partial F(\x_{t+1})) \leq\frac{C}{\sqrt{T+1}}$, where $C>0$ is a constant independent of $T$. Consequently, there exists $\bar t\in[T+1]$ such that
$\|\nabla f(\x_{\bar t})+\partial h(\x_{\bar t})\|\leq\epsilon$, provided $T+1\geq C^2\epsilon^{-2}$.}
\end{theorem}

\begin{remark}

Theorem~\ref{the:it} gives $\mathcal O(\epsilon^{-2})$ proximal iterations and $\mathcal O(N\epsilon^{-2})$ component-gradient evaluations, without requiring the
smoothness constant $L$ as an algorithmic input.

\end{remark}

\subsection{Convergence Analysis of {\sf ALFPG-FS}}

{A proof overview is provided in
Appendix~\ref{subsect:FS:overview}.}
The following theorem establishes expected stationarity
and component-gradient complexity bounds.


The following theorem gives the expected stationarity and gradient oracle complexity guarantees.

\begin{theorem}\label{the:it:stoc}
(Proof in Appendix~\ref{app:the:it:stoc}, \rm{Iteration and Oracle Complexity})
Under Assumptions~\ref{ass:prox}, \ref{ass:compact}, and~\ref{ass:f},
let $\{\x_{t}\}_{t=0}^{T+1}$ be generated by {\sf ALFPG-FS}. Then


\begin{enumerate}[label=\textbf{(\alph*)}, leftmargin=18pt,
itemsep=2pt, topsep=1pt, parsep=0pt, partopsep=0pt]


\item We have $\frac1{T+1}\sum_{t=0}^{T}\EEE[ \operatorname{dist}(0,\partial F(\x_{t+1}))]
\leq\frac{{C_{q,b}}}{\sqrt{T+1}}$, where ${C_{q,b}}$ is defined explicitly in (\ref{eq:third:AEPG:stoc}), which is independent of $T$ and may depend on $q,b$. Consequently, there exists $\bar t\in\{1,\ldots,T+1\}$ such that $\EEE[\operatorname{dist}(0,\partial F(\x_{\bar t}))]\leq\epsilon$, provided $T+1\geq {C_{q,b}^2}\epsilon^{-2}$. The same expected bound holds for $\x_J$, where $J$ is sampled uniformly
from $\{1,\ldots,T+1\}$ independently of the algorithmic randomness.

\item {If $(q-1)/b\leq1$, the bound in part~(a) holds with $\overline{C}$ in place of $C_{q,b}$, where $\overline{C}$ is given explicitly in (\ref{eq:CCCC}) and is independent of $T,q,b$.}

\item Take $b=q=\lceil\sqrt N\rceil$. With the remaining algorithmic parameters fixed, at most $N+3\sqrt N\max\{1,\lceil\overline{C}^{\,2}\epsilon^{-2}\rceil\}=\mathcal O(N+\sqrt N\epsilon^{-2})$ component-gradient evaluations suffice to obtain $\EEE[\operatorname{dist}(0,\partial F(\x_J))]\leq\epsilon$.

\end{enumerate}

\end{theorem}

\begin{remark}
\ding{182} {\sf AdaSPIDER}~\cite{kavis2022adaptive} attains
$\widetilde{\mathcal O}(N+\sqrt N\epsilon^{-2})$
component-gradient complexity for smooth finite-sum problems.
With $q=b=\lceil\sqrt N\rceil$,
Theorem~\ref{the:it:stoc}(c) gives
$\mathcal O(N+\sqrt N\epsilon^{-2})$ complexity for
{\sf ALFPG-FS}, without logarithmic factors in $N$ or
$\epsilon^{-1}$, while allowing nonconvex regularizers
under the stated assumptions.
\ding{183} The deterministic proof uses the recursive bound
in~\eqref{eq:resursive:bounding:1}.
The stochastic proof instead controls $\EEE[\VV_t]$
through~\eqref{eq:resursive:bounding:2},
accounting for the randomness of the adaptive metric.
\end{remark}

\section{Oracle Complexity for {\sf ALFPG-MVR1} and {\sf ALFPG-MVR2}}
\label{sect:iter2}

We analyze {\sf ALFPG-MVR1} and {\sf ALFPG-MVR2} under the following assumptions. We fix a deterministic initial point
$\x_1\in\dom(h)$.

\begin{assumption}\label{ass:oracle}
(Unbiased stochastic gradients with bounded variance) Let $\{{\mathcal F}_t\}_{t\geq0}$ be the natural filtration generated by the i.i.d. samples $\xi_1,\xi_2,\ldots$. At iteration $t$, $\x_t$, $\x_{t-1}$, $\alpha_t$, and all quantities constructed before sampling are ${\mathcal F}_{t-1}$-measurable, while $\xi_t$ is independent of ${\mathcal F}_{t-1}$ and distributed according to $\mathcal D$. For every deterministic $\x$,
$\E[\nabla f(\x;\xi)]=\nabla f(\x)$, $\E\|\nabla f(\x;\xi)-\nabla f(\x)\|^2\leq\sigma^2.$

\end{assumption}

\begin{assumption}\label{ass:bounded-oracle}
(Bounded stochastic gradients) There exists ${{\widehat G}}>0$ such that $\|\nabla f(\x;\xi)\|\leq{{\widehat G}}$ almost surely for every $\x\in\mathbb R^n$. Consequently, Jensen’s inequality gives $\|\nabla f(\x)\|\leq {{\widehat G}}$ for all $\x$.
\end{assumption}

\begin{assumption}\label{ass:smooth}
(Smoothness of the expected loss and lower boundedness) The expected loss $f$ is differentiable and $L$-smooth:
$\|\nabla f(\x)-\nabla f(\y)\| \leq L\|\x-\y\|, ~\forall\x,\y\in\mathbb R^n.$
Moreover, $F_\star:=\inf_{\x\in\mathbb R^n}F(\x)>-\infty,~\Delta:=F(\x_1)-F_\star<\infty.$
\end{assumption}

\begin{assumption}\label{ass:regularizer}
(Regularizer: convexity or Lipschitz continuity) The function $h:\mathbb R^n\to(-\infty,+\infty]$ is proper, lower semicontinuous, and proximable. At least one of the following conditions holds: (i) $h$ is convex. (ii) There exists $\hat H\geq0$ such that $|h(\x)-h(\y)|\leq\hat H\|\x-\y\|,~\forall\x,\y\in\dom(h)$.
\end{assumption}

\begin{remark}
Convex regularizers may be extended-valued and need not be Lipschitz. The second alternative permits nonconvex regularizers, including indicator functions of nonempty closed sets. Lipschitz continuity on $\dom(h)$ need not bound ambient subgradients.
\end{remark}

\paragraph{Notation and constants.} Define $\s_t:=\m_t-\nabla f(\x_t)$ and
$\r_t:=(\x_t-\x_{t+1})/\gamma_t$. Set $\gamma_0=0$,
$\r_0=\zero$, $\alpha_0=1$, and $\s_0=-\nabla f(\x_1)$, and write
$\E_t[\cdot]:=\E[\cdot\mid\mathcal F_{t-1}]$.
We measure composite stationarity by $\mathcal{R}(\x):=\dist\!\left(0,\nabla f(\x)+\partial h(\x)\right)$. Constants may depend on $L$, $\theta$, ${{\widehat G}}$, $K_h$, and $\Delta$,
but are independent of $T$ and $\sigma$. Here $K_h$ is specified in Lemmas~\ref{lemma:MVR1:basic:bound} and~\ref{lemma:MVR2:basic:bound}.

\subsection{Convergence Analysis of {\sf ALFPG-MVR1}} \label{subsect:analysis:A}

For {\sf ALFPG-MVR1}, the tracking estimator is an exponential
moving average, and the analysis uses smoothness of the expected
loss. 

{An overview of the proof strategy is provided
in Section~\ref{subsect:MVR1:overview}.} The following theorem establishes a noise-adaptive convergence rate for {\sf ALFPG-MVR1}.

\begin{theorem}\label{theorem:MVR1:main:paper}
(Noise-adaptive composite stationarity of {\sf ALFPG-MVR1}, proof in Theorem~\ref{theorem:MVR1}) Suppose Assumptions \ref{ass:prox}, \ref{ass:oracle}, \ref{ass:bounded-oracle},
\ref{ass:smooth}, and~\ref{ass:regularizer} hold, and run {\sf ALFPG-MVR1}. Write $\ell_T:=\ln(e+T)$. Then there exists a constant $C_R>0$, independent of
$T$ and $\sigma$, such that, for every $T\geq 1$,
\beq \label{eq:main-rate:A}
\ts \frac1T
\sum_{t=1}^T\mathbb E[\mathcal{R}(\x_{t+1})] \leq
C_R\left(\frac{\ell_T^{5/4}}{T^{1/2}}+\frac{ \sigma^{1/2} \ell_T}{T^{1/4}}\right).
\eeq
Consequently, if $R\sim{\rm Unif}\{1,\ldots,T\}$ is independent of the algorithmic randomness, then
$\mathbb E[\mathcal{R}(\x_{R+1})] = \widetilde{\mathcal O}
(T^{-1/2}+\sigma^{1/2}T^{-1/4})$. Since {\sf ALFPG-MVR1} uses one stochastic-gradient
evaluation per iteration, its stochastic-gradient oracle complexity for obtaining
$\mathbb E[\mathcal{R}(\x_{R+1})]\leq\epsilon$ is $\widetilde{\mathcal O}(\epsilon^{-2}+\sigma^2\epsilon^{-4})$.

\end{theorem}

\subsection{Convergence Analysis of {\sf ALFPG-MVR2}}
\label{subsect:analysis:I}

For {\sf ALFPG-MVR2}, we additionally assume uniform
smoothness of the sample losses.

\begin{assumption}\label{ass:individual-smooth}
(Individual smoothness for {\sf ALFPG-MVR2}) For almost every $\xi$, the sample loss $f(\cdot;\xi)$ is differentiable and $L$-smooth:
$\|\nabla f(\x;\xi)-\nabla f(\y;\xi)\| \leq L\|\x-\y\|, ~\forall\x,\y\in\mathbb R^n.$
\end{assumption}

{An overview of the proof strategy is provided
in Section~\ref{subsect:MVR2:overview}.} The following theorem establishes a noise-adaptive convergence rate for {\sf ALFPG-MVR2}.

\begin{theorem}\label{theorem:MVR2:main:paper}
(Noise-adaptive composite stationarity of {\sf ALFPG-MVR2}, proof in Theorem~\ref{app:theorem:MVR2})
Suppose Assumptions \ref{ass:prox}, \ref{ass:oracle}, \ref{ass:bounded-oracle},
\ref{ass:smooth}, \ref{ass:regularizer}, and~\ref{ass:individual-smooth}
hold, and run {\sf ALFPG-MVR2}. Write $\ell_T:=\ln(e+T)$. Then there exists a constant $C_R>0$, independent of $T$ and $\sigma$, such that, for every $T\geq1$,
\beq\label{eq:main-rate}
\ts \frac1T
\sum_{t=1}^T\mathbb E[\mathcal{R}(\x_{t+1})] \leq C_R  \left(\frac{\ell_T}{T^{1/2}}+\frac{\sigma^{1/3}\ell_T}{T^{1/3}}\right).
\eeq
Consequently, if $R\sim{\rm Unif}\{1,\ldots,T\}$
is independent of the algorithmic randomness, then
$\mathbb E[\mathcal{R}(\x_{R+1})]
=\widetilde{\mathcal O}
(T^{-1/2}+\sigma^{1/3}T^{-1/3})$. Since {\sf ALFPG-MVR2} uses at most two stochastic-gradient
evaluations per iteration, its stochastic-gradient oracle complexity for obtaining
$\mathbb E[\mathcal{R}(\x_{R+1})]\leq\epsilon$ is $\widetilde{\mathcal O}(\epsilon^{-2}+\sigma\epsilon^{-3})$.

\end{theorem}

\begin{remark}

\ding{182} Theorems~\ref{theorem:MVR1:main:paper} and~\ref{theorem:MVR2:main:paper} give the gradient oracle complexities $\widetilde{\mathcal O}(\sigma^2\epsilon^{-4}+\epsilon^{-2})$ and $\widetilde{\mathcal O}(\sigma\epsilon^{-3}+\epsilon^{-2})$, respectively. Each iteration uses one fresh sample; {\sf ALFPG-MVR1} uses one gradient evaluation and {\sf ALFPG-MVR2} uses at most two. \ding{183} The logarithmic factors arise from the adaptive-weight and motion bounds. Both complexities reduce to $\widetilde{\mathcal O}(\epsilon^{-2})$ for $\sigma=0$. These results use a bounded stochastic-gradient oracle and either alternative in Assumption~\ref{ass:regularizer}. \end{remark}

\myvspace
\section{Experiments}
\myvspace
We evaluate the ALFPG methods on sparse phase retrieval, linear eigenvalue problems, and neural network regression with box-constrained $\ell_1$ regularization, using eight benchmark instances for each task. 

$\blacktriangleright$ \textbf{Compared Methods}. \textbf{(i)} In the finite-sum setting, we compare {\sf ALFPG-FS} with {\sf ProxSARAH} \cite{pham2020proxsarah}, {\sf SpiderBoost}, {\sf SpiderBoost-M} \cite{wang2019spiderboost}, and {\sf SGP-SPIDER} (a projected subgradient/SPIDER baseline related to \cite{yang2020proxsgd}) for sparse phase retrieval and neural network regression. \textbf{(ii)} In the deterministic setting, we compare {\sf ALFPG-D} with {\sf OptM} \cite{wen2013feasible}, {\sf APG} \cite{li2015accelerated}, and {\sf FOForth} \cite{Gau2018} on linear eigenvalue problems. For {\sf FOForth}, different retraction strategies are employed to handle the orthogonality constraint, resulting in several variants: {\sf FOForth-GR}, {\sf FOForth-P}, and {\sf FOForth-QR}. Similarly, for {\sf OptM}, both QR and Cayley retraction strategies are utilized, giving rise to two variants: {\sf OptM-QR} and {\sf OptM-Cayley}. It is worth noting that both {\sf FOForth} and {\sf OptM} incorporate the Barzilai-Borwein non-monotonic line search in their implementations. \textbf{(iii)} In the expectation setting, we compare {\sf ALFPG-MVR1} and {\sf ALFPG-MVR2} with
{\sf SpiderBoost} \cite{wang2019spiderboost},
{\sf ProxSARAH} \cite{pham2020proxsarah},
{\sf ProxHSGD} \cite{tran2022hybrid},
{\sf MB-SPG} \cite{xu2019non},
{\sf ProxGen} \cite{yun2021adaptive},
{\sf PStorm} \cite{xu2023momentum}, and {\sf Prox-AdaGrad-Norm} on sparse phase retrieval, linear eigenvalue problems, and neural network regression. {\sf Prox-AdaGrad-Norm} combines the scalar stepsize of
{\sf AdaGrad-Norm} \cite{ward2020adagrad} with a proximal
update for $h$: $\eta_t= \eta \big(b_0^2+\sum_{i=0}^{t}\|\nabla f(\x_i)\|^2\big)^{-1/2}$, $\x_{t+1}\in
\operatorname{prox}_{\eta_t h}\big(\x_t-\eta_t\nabla f(\x_t)\big)$, where $\eta>0$ and $b_0>0$.

$\blacktriangleright$ \textbf{Main Findings}. \textbf{(i)} In the finite-sum and deterministic comparisons, extrapolated {\sf ALFPG-FS} and {\sf ALFPG-D} variants often reduce the objective more rapidly than their non-extrapolated counterparts. Furthermore, {\sf ALFPG-FS} and {\sf ALFPG-D} generally outperform the other methods in attaining low objective values. \textbf{(ii)} In the expectation setting, {\sf ALFPG-MVR1} and {\sf ALFPG-MVR2} often achieve rapid objective reduction and competitive final objective values relative to the tested stochastic baselines. 

$\blacktriangleright$ \textbf{Further Details}. Appendix~\ref{app:sect:exp} provides the problem formulations, experimental settings, complete results, and implementation details needed to interpret these comparisons. 

$\blacktriangleright$ \textbf{MATLAB Code}. The MATLAB code used to generate the reported results is included in the supplemental material.

\myvspace

\section{Conclusions}

\myvspace
This paper develops {\sf ALFPG}, a family of Lipschitz-free proximal gradient methods for composite nonconvex optimization. {\sf ALFPG-FS} combines a displacement-based diagonal metric, metric-coupled extrapolation, and the SPIDER estimator. The analysis establishes stationarity guarantees for this combination under exact proximal computation and {bounded sublevel sets}. With a full gradient at every iteration, {\sf ALFPG-FS} reduces to {\sf ALFPG-D}. Under the stated finite-sum assumptions, {\sf ALFPG-D} uses $\mathcal{O}(N/\epsilon^2)$ component-gradient evaluations to reach $\epsilon$-stationarity, and, when $f$ and $h$ are also bounded below, {\sf ALFPG-FS} uses $\mathcal{O}(N+\sqrt N \epsilon^{-2})$ evaluations to reach the same stationarity criterion in expectation. For expectation problems, {\sf ALFPG-MVR1} and {\sf ALFPG-MVR2} attain $\widetilde{\mathcal O}(\sigma^2\epsilon^{-4}+\epsilon^{-2})$ and $\widetilde{\mathcal O}(\sigma\epsilon^{-3}+\epsilon^{-2})$ gradient oracle complexities under their respective assumptions. Their bounds adapt to the noise level and recover $\widetilde{\mathcal O}(\epsilon^{-2})$ when $\sigma=0$. Across
experiments on sparse phase retrieval, linear eigenvalue problems, and neural network
regression, the four {\sf ALFPG} variants generally outperform the state-of-the-art proximal gradient
baselines.


\section*{LLM Usage}
A large language model (LLM) was used to assist with language editing, improving the presentation of proofs, and reviewing mathematical derivations, with the authors retaining full responsibility for the content and correctness of the paper.

\normalem

\bibliographystyle{iclr2026_conference}
\bibliography{mybib}
\onecolumn
\clearpage
\appendix

{\huge Appendix}

The organization of the appendix is as follows:

Appendix \ref{app:sect:notations} provides notations, technical preliminaries, and relevant lemmas.

Appendix~\ref{app:sect:iter} offers the proofs of the results in Section~\ref{sect:iter}.

Appendix~\ref{app:sect:iter2} contains the proofs of the results in Section~\ref{sect:iter2}.



{Appendix~\ref{app:sect:exp} provides the experiment results.}

\section{Notations, Technical Preliminaries, and Relevant Lemmas}
\label{app:sect:notations}

\subsection{Notations}

{Bold lowercase and uppercase letters denote real-valued
vectors and matrices, respectively. Square brackets with a subscript denote vector components. The following notation is used throughout the paper.}

\begin{itemize}[leftmargin=12pt,itemsep=0.2ex]

\item $[n]$: The set $\{1,2,\ldots,n\}$.

\item {$[\x]_i$, $[\x_t]_i$: The $i$th components of $\x$ and $\x_t$,
respectively, where $t$ denotes the iteration index.}

\item {$A_{ij}$: The entry in row $i$ and column $j$
of the matrix $\A$.}

\item {$\zero$, $\one$: The all-zero and all-one vectors
of the appropriate dimension.}

\item $\|\mathbf{x}\|$: Euclidean norm, defined as
$\|\mathbf{x}\|=\|\mathbf{x}\|_2
=\sqrt{\la\mathbf{x},\mathbf{x}\ra}$.

\item $\la\a,\b\ra$: The Euclidean inner product of $\a$ and $\b$.

\item $\x\pm\y$: Element-wise addition or subtraction of $\x$ and $\y$.

\item $\x\odot\y$: Element-wise multiplication of $\x$ and $\y$.

\item $\x\div\y$: Element-wise division of $\x$ by $\y$,
provided that all components of $\y$ are nonzero.

\item {$\sqrt{\x}$, $\x^p$: Element-wise square roots
and powers, whenever these operations are well defined.}

\item $\la\a,\b\ra_{\v}$: The weighted Euclidean inner product,
defined by $\la\a,\b\ra_{\v}=\la\a\odot\b,\v\ra$, where $\v>\zero$.

\item $\|\x\|_{\v}$: The weighted Euclidean seminorm, defined by
$\|\x\|_{\v}=\sqrt{\la\x\odot\x,\v\ra}$, where $\v\geq\zero$.

\item $\a\leq\b$: Component-wise inequality, meaning
${[\a]_i\leq[\b]_i}$ for every $i$.

\item {$\min(\x)$: The smallest entry of $\x$, defined by
$\min(\x):=\min_{i\in[n]}[\x]_i$.}

\item {$\max(\x)$: The largest entry of $\x$, defined by
$\max(\x):=\max_{i\in[n]}[\x]_i$.}

\item $\dom(h)$: The effective domain
$\{\x:h(\x)<+\infty\}$ of the proper function $h$.

\item $\iota_{\Omega}(\x)$: The indicator function of $\Omega$, equal to $0$ for $\x\in\Omega$ and $+\infty$ otherwise.


\item $\dist(\x,S)$: The distance from $\x$ to a set $S$,
defined by $\dist(\x,S):=\inf_{\u\in S}\|\x-\u\|$.

\item {$\|S\|$: Shorthand for
$\dist(\zero,S)=\inf_{\u\in S}\|\u\|$.
In particular,
$\|\partial h(\x)\|=\dist(\zero,\partial h(\x))$.}

\item $\mathcal{R}(\x)$: The composite stationarity
residual with $\mathcal{R}(\x)
:=\dist\!\bigl(\zero,\nabla f(\x)+\partial h(\x)\bigr)$.

\item {$\Prox_{\gamma h}(\a;\v)$: The variable-metric
proximal mapping, defined for $\gamma>0$ and $\v>\zero$ by $\Prox_{\gamma h}(\a;\v)
:=\argmin_{\x\in\Rn^n}
\left\{h(\x)+\frac{1}{2\gamma}\|\x-\a\|_{\v}^{2}\right\}$. The mapping can be set-valued when $h$ is nonconvex.
For $\v=\one$, we write $\Prox_{\gamma h}(\a)$ or
$\prox_{\gamma h}(\a)$.}

\item $\EEE[v]$: Expected value of the random variable $v$.

\item {$\mathcal F_t$: The filtration generated by
the algorithmic randomness through iteration $t$.
In the expectation setting, we write
$\E_t[v]:=\EEE[v\mid\mathcal F_{t-1}]$.}

\item $\{A_i\}_{i=0}^{\infty}$, $\{B_i\}_{i=0}^{\infty}$:
Sequences indexed by the integers $i=0,1,2,\ldots$.

\item $\A\trans$: The transpose of the matrix $\A$.



\item $b$: The mini-batch size parameter of {\sf ALFPG-FS}.

\item $q$: The full-gradient refresh period of
{\sf ALFPG-FS}. 


\end{itemize}

\subsection{Technical Preliminaries}

We introduce key concepts from nonsmooth analysis, focusing on the Fr\'{e}chet subdifferential and the limiting (Fr\'{e}chet) subdifferential \cite{Mordukhovich2006,Rockafellar2009,bertsekas2015convex}. Let $F: \Rn^n \rightarrow (-\infty,+\infty]$ be an extended real-valued, not necessarily convex function. The domain of $F(\cdot)$ is defined as $\text{dom}(F):= \{\x\in\Rn^n: |F(\x)|<+\infty\}$. The Fr\'{e}chet subdifferential of $F$ at $\x\in\text{dom}(F)$, denoted as $\hat{\partial}F(\x)$, is given by $$\hat{\partial}{F}(\x) := \{\mathbf{v}\in\Rn^n: \lim_{\mathbf{z} \rightarrow \x} \inf_{\mathbf{z}\neq \x} \frac{ {F}(\mathbf{z}) - {F}(\x) - \la \mathbf{v},\mathbf{z}-\x \ra  }{\|\mathbf{z}-\x\|}\geq 0\}.$$

\noi The limiting subdifferential of $F(\cdot)$ at $\x\in\text{dom}({F})$, denoted $\partial{F}(\x)$, is defined as: $$\partial{F}(\x):= \{\mathbf{v}\in \Rn^n: \exists \x_k \rightarrow \x, {F}(\x_k)  \rightarrow {F}(\x), \mathbf{v}^k \in\hat{\partial}{F}(\x_k) \rightarrow \mathbf{v},\forall k\}.$$
 \noi It is important to note that $\hat{\partial}{F}(\x) \subseteq \partial{F}(\x)$. If $F(\cdot)$ is continuously differentiable at $\x$, then $\hat{\partial}{F}(\x) = \partial{F}(\x) = \{\nabla F(\x)\}$, where $\nabla F(\x)$ represents the gradient of $F(\cdot)$ at $\x$. For convex function $F(\cdot)$, both $\hat{\partial}{F}(\x)$ and $\partial{F}(\x)$ reduce to the classical subdifferential for convex functions: $\hat{\partial}{F}(\x) = \partial{F}(\x) = \{\mathbf{v}\in\Rn^n: F(\mathbf{z})-F(\x)-\la\mathbf{v},\mathbf{z}-\x \ra\geq 0,\forall \mathbf{z}\in\Rn^n\}$.

\subsection{Relevant Lemmas}

We collect elementary inequalities used in the subsequent convergence proofs.

\begin{lemma} \label{lemma:pyth}
(\rm{Pythagoras Relation}) For any vectors $\a,\v,\x,\x_+ \in \Rn^n$ with $\v\geq \zero$, we have: $$\tfrac{1}{2}\|\x - \a\|_{\v}^2 - \tfrac{1}{2}\|\x_+ - \a \|_{\v}^2 = \tfrac{1}{2}\| \x  -  \x_+ \|_{\v}^2 +   \la \a - \x_{+}, \x_{+} - \x\ra_{\v}.$$
\end{lemma}

\begin{lemma} \label{lemma:frac:4}
For all $a,b\geq0$ and $c,d>0$, we have: $\frac{a+b}{c+d} \leq \max(\tfrac{a}{c},\tfrac{b}{d})$.

\begin{proof}
 We consider two cases: (\bfit{i}) $\frac{a}{c}\geq \frac{b}{d}$. We derive: $b\leq \frac{ad}{c}$, leading to $\frac{a+b}{c+d}\leq \frac{a+\frac{ad}{c}}{c+d} = \frac{a}{c}\cdot \frac{c+d}{c+d}=\frac{a}{c}$. (\bfit{ii}) $\frac{a}{c}<\frac{b}{d}$. We have: $a \leq \frac{bc}{d}$, resulting in $\frac{a+b}{c+d}\leq \frac{\frac{bc}{d}+b}{c+d}=\frac{b}{d}\cdot \frac{c+d}{c+d}=\frac{b}{d}$.

\end{proof}
\end{lemma}

\begin{lemma} \label{lemma:quad:abc}
Assume $a x^2\leq bx + c$, where $b,c,x\geq 0$ and $a>0$. Then, we have: $x\leq \sqrt{{c}/{a}} + {b}/{a}$.

\begin{proof}

Given the quadratic inequality $a x^2\leq bx + c$, we have $\tfrac{b-\sqrt{b^2+4a c}}{2a}\leq x\leq \tfrac{b+\sqrt{b^2+4a c}}{2a }$. Since $x\geq 0$, we have $0 \leq x\leq \tfrac{b+\sqrt{b^2+4 a c}}{2a} \leq \tfrac{b+ b+2\sqrt{a c}}{2a} = {b}/{a}+\sqrt{c/a}$, where the last inequality uses $\sqrt{a+d}\leq \sqrt{a}+\sqrt{d}$ for all $a,d\geq0$.

\end{proof}
\end{lemma}

\begin{lemma} \label{lemma:remove:log}

Assume that $\{A_i\}_{i=0}^n$ and $\{B_i\}_{i=0}^{n+1}$ are two non-negative sequences with $A_0\leq A_1\leq \ldots \leq A_n$. Then, we have: $$\ts \sum_{t=0}^n A_t (B_t - B_{t+1})\leq [\max_{i=0}^n A_i]\cdot[\max_{j=0}^n B_j].$$

\begin{proof}

We have: $\ts \sum_{t=0}^n A_t (B_t - B_{t+1} ) = \ts  [\sum_{t=1}^n (A_t - A_{t-1}) B_t ] +  A_0 B_0 - A_n B_{n+1}  \overset{\step{1}}{\leq} \ts    [\sum_{t=1}^n ( A_t - A_{t-1}) B_t] + A_0 B_0 \overset{\step{2}}{\leq} \ts [\sum_{t=1}^n ( A_t - A_{t-1}) ]\cdot [\max_{j=0}^n B_j]  +  A_0 [\max_{j=0}^n B_j] = \ts A_n [\max_{j=0}^n B_j]$, where step \step{1} uses $A_n,B_{n+1}\geq 0$; step \step{2} uses $\{A_t\}_{t=0}^n$ is non-decreasing.

\end{proof}

\end{lemma}

\begin{lemma} \label{lemma:online}

Let $\{A_t\}_{t=0}^{\infty}$ be a sequence of nonnegative real numbers, and let $c>0$. Then, for every $t\geq 0$,
\beq
\ts \sum_{i=0}^t \tfrac{A_i}{\sqrt{c +  \sum_{j=0}^i A_j}} \leq 2 \sqrt{ c + \sum_{i=0}^t A_i  }
\eeq

\begin{proof}

This lemma extends the result of Lemma 5 in \cite{McMahanS10}.

We define $S_t:= c+\sum_{i=0}^t A_i$.

Initially, we define $h(x):=  \tfrac{x}{\sqrt{y}} + 2 \sqrt{y-x} - 2 \sqrt{y}$, where $y> x\geq 0$. We have $h'(x) = y^{-1/2}-(y-x)^{-1/2}  \leq 0$. Therefore, $h(x)$ is non-increasing for all $x\geq 0$. Given $h(0)=0$, it holds that
\beq \label{eq:xx:bound}
h(x):=  \tfrac{x}{\sqrt{y}} + 2 \sqrt{y-x} - 2 \sqrt{y} \leq 0.
\eeq

We complete the proof of the lemma using mathematical induction.

\paragraph{Part (a).} The lemma holds for $t = 0$, since $\tfrac{A_0}{\sqrt{c+A_0}} \leq 2 \sqrt{c+A_0}$.

\paragraph{Part (b).} Now, fix some $t$ and assume that the lemma holds for $t-1$. We proceed as follows:
\beq
  \ts \sum_{i=0}^t {A_i}/{\sqrt{S_i}} &=&   \ts {A_t}/{\sqrt{S_t}}+  \sum_{i=0}^{t-1} {A_i}/{\sqrt{S_i}}  ~~\overset{\step{1}}{\leq}~~ {A_t}/{\sqrt{S_t}}+ 2 \sqrt{S_{t-1}} \nn\\
  &\overset{\step{2}}{=}&\ts  {A_t}/{\sqrt{S_t}}+ 2 \sqrt{S_t - A_t } ~~\overset{\step{3}}{\leq}~~ 2 \sqrt{S_t}, \nn
\eeq
\noi where step \step{1} uses the inductive hypothesis that the conclusion of this lemma holds for $t-1$; step \step{2} uses $S_t := c+\sum_{i=0}^t A_i$; step \step{3} uses Inequality (\ref{eq:xx:bound}) with $x=A_t\geq 0$ and $y=S_t>0$.

\end{proof}
\end{lemma}

\begin{lemma}
\label{lemma:online:pp}
Let $\{A_t\}_{t=0}^{\infty}$ be a sequence of nonnegative real numbers, and let $c>0$ and $p\in(0,1]$. Then, for every $t\geq 0$,
\beq
\ts \sum_{i=0}^t \tfrac{A_i}{c + \sum_{j=0}^i A_j} \leq -1+\tfrac{1}{p} \left(1+ \tfrac{1}{c}\sum_{i=0}^t A_i\right)^p.\nn
\eeq

\begin{proof}

We define $S_t:= \sum_{i=0}^t A_i$.

We define $g(x) = \tfrac{x}{1+x}-\tfrac{1}{p}(1+x)^p+1$ and $f(x) = x^{p}$, where $x\geq 0$ and $p\in(0,1]$.

First, we have $g(0)=1-1/p \leq 0$ and $g'(x)= (1+x)^{-2} -(1+x)^{p-1}  = \tfrac{1-(1+x)^{p+1}}{(1+x)^2}\leq 0$. We derive:
\beq \label{eq:gxgx:pp}
g(x) \leq 0.
\eeq

Second, since $f(x)$ is concave, we have $\forall x,y > 0,\,f(x) \leq f(y) + \la x-y, f'(y) \ra$, leading to
\beq
\ts \forall x,y>0,\,x^{p}  \leq  y^{p} - p (y-x) \cdot y^{p-1}.  \label{eq:xy:sqrt4}
\eeq

We finish the proof of the lemma using mathematical induction.

\paragraph{Part (a).} We first consider $t=0$. We have
\beq
\ts A_0 / (c + S_0) - \tfrac{1}{p} (1+S_0/c)^{p} +1 & \overset{\step{1}}{=}& A_0 / (c + A_0) - \tfrac{1}{p} (1+A_0/c)^{p}+1 \nn\\
&\overset{}{=}& \tfrac{A_0/c}{1 + A_0/c}  - \tfrac{1}{p} (1+A_0/c)^{p} +1 ~~\overset{\step{2}}{\leq}~~ 0,\nn
\eeq
\noi where step \step{1} uses $A_0=S_0$; step \step{2} uses Inequality (\ref{eq:gxgx:pp}) with $x=A_0/c$. Thus, the claim holds for $t=0$.

\paragraph{Part (b).} Now, fix some $t$ and assume that the lemma holds for $t-1$. We derive:
\beq
\ts \sum_{i=0}^t {A_i}/ (c+S_i) &= & \ts {A_t}/ (c+S_t) + \sum_{i=0}^{t-1} {A_i}/ (c+S_i) \nn\\
&\overset{\step{1}}{\leq}& \ts -1+ {A_t}/ (c+S_t) + \tfrac{1}{p} c^{-p} (c+ S_{t-1})^{p}\nn\\
&\overset{\step{2}}{=}& \ts -1+ {A_t}/ (c+S_t) + \tfrac{1}{p} c^{-p} (c+ S_{t} -A_t)^{p}\nn\\
&\overset{\step{3}}{\leq}& \ts-1+  {A_t}/ (c+S_t) + \tfrac{1}{p} c^{-p} \{ (c+ S_t)^{p} -  p A_t  (c + S_t)^{p-1} \} \nn\\
&\overset{}{=}&\ts -1+ \tfrac{1}{p} c^{-p} (c+ S_t)^{p} +  A_t (c + S_t)^{p-1}\cdot \{  (c + S_t)^{-p}  -  c^{-p} \} \nn\\
&\overset{\step{4}}{\leq}&\ts -1+ \tfrac{1}{p} c^{-p} (c+ S_t)^{p} + 0, \nn
\eeq
\noi where step \step{1} uses the inductive hypothesis that the conclusion of this lemma holds for $t-1$; step \step{2} uses $S_t := \sum_{i=0}^t A_i$; step \step{3} uses Inequality (\ref{eq:xy:sqrt4}) with $x=c+S_t-A_t>0$ and $y=c+S_t> 0$; step \step{4} uses $(c+S_t)^{-p}\leq c^{-p}$ for all $p\in(0,1]$.

\end{proof}
\end{lemma}

\begin{lemma}  \label{lemma:bound:zzz}

Let $\{Z_t\}_{t=0}^{\infty}$ be a non-negative sequence satisfying $\ts Z_{t+1} \leq a + b \sqrt{ \max_{i=0}^tZ_{i}  }$ for all $t\geq 0$, where $a,b\geq 0$. It follows that: $Z_t \leq \overline{\rm{Z}} := \max\left(Z_0, \left(\tfrac{1}{2} b + \tfrac{1}{2} \sqrt{b^2 + 4a}\right)^2\right)$ for all $t\geq 0$. Furthermore, an alternative valid upper bound for $Z_t$ is given by $\overline{\rm{Z}}_+ := \max(Z_0,2 b^2 +2 a)$.

\begin{proof}

We define $M_t := \max_{0\leq i\leq t} Z_i$ for all $t\geq 0$.

\paragraph{Part (a).} For all $t\geq 0$, we have:
\beq
M_{t+1} \overset{\step{1}}{=} \max (M_t,Z_{t+1}) \overset{\step{2}}{\leq}    \max (M_t,a+b\sqrt{M_t}),\nn
\eeq
\noi where step \step{1} uses the definition of $M_t$; step \step{2} uses the assumed recursion $Z_{t+1}\leq a+b \sqrt{M_t}$.

\paragraph{Part (b).} We establish a fixed-point upper bound that the sequence $M_t$ cannot exceed such that $M_t \leq \overline{\rm{M}}$. For $\overline{\rm{M}}$ to be a valid upper bound, it should satisfy the recurrence relation: $\overline{\rm{M}} = a + b \sqrt{\overline{\rm{M}}}$. This is because if $M_t \leq \overline{\rm{M}}$, then we have: $M_{t+1}\leq\max\{M_t,a+b\sqrt{M_t}\}\leq\overline{\rm M}$, which would imply by induction that $M_t \leq \overline{\rm{M}}$ for all $t\geq 0$. Solving the quadratic equation $\overline{\rm{M}} = a + b \sqrt{\overline{\rm{M}}}$ yields a positive root $(\tfrac{1}{2} (b + \sqrt{b^2 + 4a}))^2 := c$. Taking into account the case $M_0$, for all $t\geq 0$, we have $M_t\leq \max(M_0,c) = \max(Z_0,c):= \overline{\rm{Z}}$.

\paragraph{Part (c).} Let $\overline{\rm Z}=\max(Z_0,c)$. Since $c=a+b\sqrt c$, we have $a+b\sqrt{\overline{\rm Z}}\leq\overline{\rm Z}$. If $M_t\leq\overline{\rm Z}$, then
\[
M_{t+1}\leq\max\{M_t,a+b\sqrt{M_t}\}\leq\overline{\rm Z}.
\]
As $M_0=Z_0\leq\overline{\rm Z}$, induction proves the bound for every $t$.

\paragraph{Part (d).} Finally, we have: $c := \left(\tfrac{1}{2} b +\tfrac{1}{2} \sqrt{b^2 + 4a} \right)^2 \leq (\tfrac{1}{2} (b + b + 2 \sqrt{a}))^2 = (b+\sqrt{a})^2\leq 2b^2 + 2a$. Hence, $\overline{\rm{Z}}_+ := \max(Z_0, 2 b^2 +2 a)$ is also a valid upper bound for $Z_t$.
 \end{proof}

 \end{lemma}

\begin{lemma} \label{eq:log:diff}
Let $b_1,\ldots,b_T\geq0$. Then:

(a) For every $t=1,\ldots,T$, $\frac{b_t}{1+\sum_{i=1}^t b_i}
\leq
\ln\left(1+\sum_{i=1}^t b_i\right) - \ln\left(1+\sum_{i=1}^{t-1}b_i\right)$.

(b) We have $\sum_{t=1}^T \frac{b_t}{1+\sum_{i=1}^t b_i}\leq \ln\left(1+\sum_{i=1}^T b_i\right)$.

\begin{proof}

Define $B_0:=1$ and $B_t:=1+\sum_{i=1}^t b_i$.
Then $B_t-B_{t-1}=b_t\geq0$.

\paragraph{Part (a).} Since $u\mapsto 1/u$ is decreasing on $(0,\infty)$, we have: $\ln B_t-\ln B_{t-1} = \int_{B_{t-1}}^{B_t}\frac{1}{u}\,du
\geq \frac{B_t-B_{t-1}}{B_t} =
\frac{b_t}{B_t}$.

\paragraph{Part (b).} Summing (a) over $t=1,\ldots,T$ gives $\sum_{t=1}^T\frac{b_t}{B_t} \leq \sum_{t=1}^T(\ln B_t-\ln B_{t-1}) = \ln B_T-\ln B_0 = \ln B_T$.

\end{proof}

\end{lemma}

\section{Proofs for Section~\ref{sect:iter}}
\label{app:sect:iter}

The analysis proceeds in three stages.
We first establish bounds for the adaptive metric and
extrapolation coefficients, together with a Lyapunov descent
inequality shared by both variants.
We then derive the deterministic stationarity guarantee
for {\sf ALFPG-D}.
Finally, we control the SPIDER estimation error and take
expectations to obtain the expected stationarity guarantee
for {\sf ALFPG-FS}.

\subsection{Overview of the proof strategy for {\sf ALFPG-D}}
\label{subsect:D:overview}

Proximal optimality, smoothness, and the metric-dependent
extrapolation rule yield the core inequality
(Lemma~\ref{prop:H}):
\beq \label{eq:D:overview:descent}
\ts \ZZ_{t+1}-\ZZ_t
\leq
-c_1\underbrace{
\tfrac{\|\r_t\|_2^2}{\min(\v_t)}
}_{\SSS_t}
+c_2\underbrace{
\tfrac{\|\r_t\|_2^2}{\min(\v_t)^2}
}_{\UUU_t},
\eeq
\noi where $\r_t=\v_t\odot(\x_{t+1}-\x_t)$, $\ZZ_t\geq 0$, and $c_1,c_2>0$. The proof proceeds in four steps.

\begin{enumerate}[label=(\roman*), leftmargin=18pt,
itemsep=0.3ex, topsep=0.3ex, parsep=0pt]

\item The metric update gives the cumulative bounds $\sum_{t=0}^{T}\SSS_t\leq s_1\VV_{T+1}$, $\sum_{t=0}^{T}\UUU_t\leq s_2\sqrt{\VV_{T+1}}$, where $\VV_{T+1}$ is the terminal accumulator (Lemma~\ref{lemma:AB}).

\item Multiplying~\eqref{eq:D:overview:descent}
by $\min(\v_t)$ and applying summation by parts uses
the negative term to control accumulator growth: $\VV_{T+1}
\leq w_1+w_2\max_{0\leq t\leq T}\ZZ_t$, where $w_1,w_2>0$ (see~\eqref{eq:upper:bound:VVT}).

\item Summing~\eqref{eq:D:overview:descent} and dropping
the negative term gives $\ZZ_{T+1}
\leq \ZZ_0+c_2s_2\sqrt{\VV_{T+1}}
\leq \ZZ_0+c_2s_2
\sqrt{w_1+w_2\max_{0\leq t\leq T}\ZZ_t}$. This bound yields uniform bounds on $\ZZ_t$
and $\VV_t$ (Lemma~\ref{lemma:smmable}).

\item Proximal optimality leads to $\frac{1}{T+1}\sum_{t=0}^{T}
\operatorname{dist}(0,\partial F(\x_{t+1}))
\leq \frac{C}{\sqrt{T+1}}$, where $C>0$.

\end{enumerate}

\subsection{Overview of the proof strategy for {\sf ALFPG-FS}}
\label{subsect:FS:overview}

Lemma~\ref{prop:H} gives the core pathwise inequality
\beq\label{eq:FS:overview:descent}
\ts \widehat{\ZZ}_{t+1}-\widehat{\ZZ}_t
\leq
-\langle\d_t,\e_t\rangle-c_1 \underbrace{
\tfrac{\|\r_t\|_2^2}{\min(\v_t)}
}_{\SSS_t} + c_2 \underbrace{
\tfrac{\|\r_t\|_2^2}{\min(\v_t)^2}
}_{\UUU_t},
\eeq
where $\e_t:=\g_t-\nabla f(\y_t)$, $\d_t:=\x_{t+1}-\x_t$, and $\ZZ_t=\mathbb E[\widehat{\ZZ}_t] \geq 0$. The proof proceeds in four steps.

\begin{enumerate}[label=(\roman*), leftmargin=18pt,
itemsep=0.3ex, topsep=0.3ex, parsep=0pt]

\item Summing the SPIDER variance recursion over refresh
epochs and applying the motion bounds gives
$\mathbb E[\sum_{t=0}^{T}\|\e_t\|_2^2]
\leq \frac{L^2(q-1)u_2}{b}
\mathbb E[\sqrt{\VV_{T+1}}]$,
where $u_2$ is defined in Lemma~\ref{lemma:AB2}
(see~\eqref{eq:fs:summed:error}).

\item The cumulative metric bounds, Young's inequality,
and Jensen's inequality give
$\frac{\alpha+\beta}{2c_1}\ZZ_{T+1}
+\mathbb E[\VV_{T+1}]
\leq A+B\sqrt{\mathbb E[\VV_{T+1}]}$,
with $A,B>0$. Hence $\ZZ_t$ and $\mathbb E[\VV_t]$ are uniformly bounded
(Lemma~\ref{lemma:fi:preresults}).

\item The accumulator identity and Cauchy--Schwarz
control the cumulative residuals through
$\mathbb E[\VV_{T+1}]$. Proximal optimality gives the averaged expected stationarity bound
$C_{q,b}/\sqrt{T+1}$.

\item For $(q-1)/b\leq1$, a metric-threshold argument gives the bound
$\overline C/\sqrt{T+1}$, with $\overline C>0$ (see~\eqref{eq:CCCC}).
Setting $q=b=\lceil\sqrt N\rceil$ yields
$\mathcal O(N+\sqrt N\epsilon^{-2})$
component-gradient complexity.

\end{enumerate}

\subsection{Initial Theoretical Analysis}

{We first establish the metric bounds, the extrapolation upper bound,
and the Lyapunov descent inequality directly from the updates.
These yield localization under bounded sublevel sets.
The resulting iterate bound is then used to prove metric growth
and the extrapolation lower bound, followed by the accumulator estimates.}

\begin{lemma}\label{lemma:prop:w}
(\rm{Properties of $\v_t$}) We define $M_t :=  \sum_{i=0}^{t} \|\r_i\|_2^2 \in \Rn_+$. For all $t\geq 0$, we have:

\begin{enumerate}[label=\textbf{(\alph*)}, leftmargin=18pt, itemsep=2pt, topsep=1pt, parsep=0pt, partopsep=0pt]

\item $\sqrt{ \underline{\rm{v}}^2  +  \alpha M_{t}} \leq \v_{t+1} \leq \sqrt{ \underline{\rm{v}}^2 + (\alpha+\beta)M_{t}} := \VV_{t+1}$.

\item $\tfrac{\max(\v_t)}{\min(\v_t)} \leq \kappa:= 1 + \sqrt{ \beta/\alpha}$.

\end{enumerate}

\begin{proof}

We define $M_t:= \sum_{i=0}^{t} \|\r_i\|_2^2 \in \Rn_+$, where $\r_t := \v_t\odot \d_t$, and $\d_t:= \x_{t+1}-\x_{t}$.

We define $\s_t:= \alpha \|\r_t\|_2^2 \cdot \mathbf{1} + \beta \r_t\odot \r_t\in \Rn^n$.
The update gives $\v_{t+1}\geq\v_t\geq\underline{\rm v}\one$ coordinatewise.

\paragraph{Part (a).} Given $\v_{t+1}= \ts \sqrt{ \v_{0} \odot \v_{0} + \ts \sum_{i=0}^{t} \s^i}$, we derive:
\beq
\ts \v_{t+1} \odot \v_{t+1} = \ts  {(\v_{0})\odot (\v_{0}) + \alpha  \underbrace{\ts \sum_{i=0}^{t} \|\r_i\|_2^2}_{ \ts =    M_t } \cdot \one  + \beta \underbrace{\ts \sum_{i=0}^{t} \r_i \odot \r_i }_{ \ts \leq  M_t \cdot \one}}. \nn
\eeq
\noi This results in the following lower and upper bounds for $\v_{t+1}$ for all $t\geq 0$:
\beq
\ts \sqrt{ \underline{\rm{v}}^2 + \alpha M_t } \leq \v_{t+1} \leq \sqrt{\underline{\rm{v}}^2  +  (\alpha+\beta)M_t}.\nn
\eeq

\paragraph{Part (b).} For all $t\geq 1$, we derive:
\beq \label{eq:kappa:1}
\ts \tfrac{\max(\v_{t})}{\min(\v_{t})}  ~~\overset{\step{1}}{\leq} ~~ \ts  \sqrt{\tfrac{ \underline{\rm{v}}^2 + (\alpha+\beta) M_{t-1}  }{  \underline{\rm{v}}^2 + \alpha M_{t-1} } } ~~{\leq}~~ \sqrt{  \tfrac{\alpha+\beta}{\alpha}} ~~\leq~~ 1+  \sqrt{\beta/\alpha},
\eeq
\noi where step \step{1} uses \textbf{Part (a)}. Note that Inequality (\ref{eq:kappa:1}) also holds for $t=0$ since $\v_0=\underline{\rm{v}}\one$.

\end{proof}

\end{lemma}

With the metric bounds in place, we establish the extrapolation upper bound and the coupling needed for descent.

\begin{lemma}\label{lemma:prop:theta}
(\rm{Properties of $\sigma_t$}) For all $t\geq 0$, we have the following results.

\begin{enumerate}[label=\textbf{(\alph*)}, leftmargin=18pt, itemsep=2pt, topsep=1pt, parsep=0pt, partopsep=0pt]

\item $0\leq \sigma_t \leq \theta$.

\item $(\sigma_{t-1}-1)\v_t + \sigma_t \v_{t+1} \leq -(1-\theta)^2 \v_t$.

\end{enumerate}

\begin{proof}

We define $\sigma_t := \ts \theta (1-   \sigma_{t-1}) \cdot \min( \v_t \div \v_{t+1})$, where $t\geq 0$.

\paragraph{Part (a).} We now prove that $\sigma_t \in [0,\theta]$. We complete the proof using mathematical induction. First, we consider $t=0$, we have:
\beq
\sigma_0 ~~=~~ \ts \min(\v_0\div \v_{1}) \theta (1-\sigma_{-1}) ~~\overset{\step{1}}{=}~~\ts \min(\v_0\div \v_{1}) \theta (1-\theta) ~~\overset{\step{2}}{\leq}~~ \theta (1-\theta) ~~\overset{\step{3}}{\leq} ~~\theta,\nn
\eeq
\noi where step \step{1} uses $\sigma_{-1}=\theta$; step \step{2} uses $\min(\v_t\div \v_{t+1})\in(0,1]$; step \step{3} uses $\theta\in[0,1)$. Second, we fix some $t$ and assume that $\sigma_{t-1}\in [0,\theta]$. We analyze the following term for all $t\geq 1$:
\beq
\sigma_t := \ts \theta (1-   \sigma_{t-1}) \cdot \min( \v_t \div \v_{t+1}).\nn
\eeq
\noi Given $\sigma_{t-1}\in[0,\theta]$, $\min( \v_t \div \v_{t+1})\in(0,1]$, and $\theta\in[0,1)$, we conclude that $\sigma_t\in [0,\theta]$.

\paragraph{Part (b).} For all $t\geq 0$, we derive the following results:
\beq
 (\sigma_{t-1}-1)\v_t + \sigma_t \v_{t+1} &\overset{\step{1}}{=}  & (\sigma_{t-1}-1)\v_t + \theta (1-   \sigma_{t-1}) \cdot \min(      \v_t \div \v_{t+1}   )  \v_{t+1} \nn\\
 &\overset{\step{2}}{\leq}& (\sigma_{t-1}-1)\v_t +  \theta (1-   \sigma_{t-1}) \v_t   ~~\overset{}{=} ~~   - (1-\theta) (1-\sigma_{t-1}) \v_t \nn\\
 &\overset{\step{3}}{ \leq} & \ts -(1-\theta)  (1-   \theta)\v_t  ~~ \overset{}{=} ~~ \ts - (1-\theta)^2 \v_t,  \nn
    \eeq
    \noi where step \step{1} uses the choice for $\sigma_t$ for all $t\geq 0$; step \step{2} uses $\min( \a\div \v) \v\leq \a$ for all $\a,\v\in\Rn^n$ with $\v>\zero$; step \step{3} uses $\sigma_{t-1}\leq \theta < 1$ for all $t\geq 0$.

\end{proof}

\end{lemma}

We let $\bar{\x} \in \arg \min_{\x}F(\x)$, where $F(\x):= f(\x)+h(\x)$. Such a minimizer exists by lower semicontinuity and bounded sublevel sets. We define the following sequence: $$\ZZ_t:=\EEE[\widehat{\ZZ}_t],~\widehat{\ZZ}_t := F(\x_{t}) - F(\bar{\x})+  \tfrac{1}{2} \|\x_{t} - \x_{t-1}\|_{\sigma_{t-1} (\v_t + L)}^2.$$

We now combine proximal optimality, smoothness, and the extrapolation inequality to obtain a pathwise Lyapunov estimate. Keeping the gradient-estimation error explicit allows the same estimate to serve both the deterministic and stochastic analyses.

\begin{lemma}\label{prop:H}
(\rm{Properties of $\ZZ_t$}) We define
$c_1 := \tfrac{1}{2}(\tfrac{1-\theta}{\kappa})^2$, $c_2:= \tfrac{3L}{2}$. We have almost surely:
\beq
\widehat{\ZZ}_{t+1} - \widehat{\ZZ}_t       \leq    \la \d_{t}, \nabla f(\y_t) - \g_t \ra + c_2\UUU_t   - c_1 \SSS_t ,
\eeq
\noi where $\SSS_t := \tfrac{\|\r_t\|_2^2}{ \min(\v_t)  }$, $\UUU_t := \tfrac{\|\r_t\|_2^2}{ \min(\v_t)^2}$.

Before rescaling the quadratic terms, we also have
\begin{equation}\label{eq:local:raw}
\widehat{\ZZ}_{t+1}-\widehat{\ZZ}_t
\leq-\langle\e_t,\d_t\rangle
+\tfrac12(1+2\theta)L\|\d_t\|^2
-\tfrac12(1-\theta)^2\|\d_t\|_{\v_t}^2,
\qquad \e_t:=\g_t-\nabla f(\y_t).
\end{equation}

\begin{proof}

Let $\bar{\x} \in \arg \min_{\x}F(\x)$, where $F(\x):= f(\x)+h(\x)$.

Define $\a_t := \y_t - \g_t \div \v_t$, and $Q_t:=\la \d_{t}, \nabla f(\y_t) - \g_t\ra$, where $\d_t := \x_{t+1}-\x_{t}$.

Define $\widehat{\ZZ}_t := F(\x_{t}) - F(\bar{\x})+  \tfrac{1}{2} \|\x_{t} - \x_{t-1}\|_{\sigma_{t-1} (\v_t + L)}^2$, and $\XX^{t}:=   \tfrac{1}{2} \|\x_{t} - \x_{t-1}\|_{\sigma_{t-1} (\v_t + L)}^2$.

Using the optimality of $\x_{t+1} \in \arg \min_{\x} h(\x) + \tfrac{1}{2} \|\x - \a_t\|_{\v_t}^2$, we have the following inequality:
\begin{align}
\ts h(\x_{t+1}) + \tfrac{1}{2}\|\x_{t+1} - \a_t\|_{\v_t}^2\leq h(\x_{t}) + \tfrac{1}{2}\|\x_{t} - \a_t\|_{\v_t}^2. \label{eq:opt:2}
\end{align}

Given $f(\x)$ is $L$-smooth, we have:
\beq \label{eq:L:smooth}
f(\x_{t+1}) \leq f(\x_{t}) + \la\x_{t+1}-\x_{t},\nabla f(\x_{t})\ra + \tfrac{L}{2}\|\x_{t+1}-\x_{t}\|_2^2.
\eeq
\noi Adding Inequalities (\ref{eq:opt:2}) and (\ref{eq:L:smooth}) together yields:
\begin{align} \label{eq:FFXX}
&~  F(\x_{t+1}) - F(\x_{t}) - \tfrac{L}{2}\|\x_{t+1}-\x_{t}\|_2^2\nn \\
\overset{}{\leq} &~ \ts \la \x_{t+1}-\x_{t}, \nabla f(\x_{t})\ra + \tfrac{1}{2}\|\x_{t} - \a_t \|_{\v_t}^2 - \tfrac{1}{2}\|\x_{t+1} - \a_t \|_{\v_t}^2 \nn\\
\overset{\step{1}}{=} &~ \ts  \la \x_{t+1}-\x_{t}, \nabla f(\x_{t})\ra + \tfrac{1}{2}\| \x_{t}  -  \x_{t+1} \|_{\v_t}^2 +   \la \a_t - \x_{t+1}, \x_{t+1} - \x_{t} \ra_{\v_t} \nn\\
\overset{\step{2}}{=} &~ \ts \la \x_{t+1}-\x_{t}, \nabla f(\x_{t}) -  {\g_t} \ra+ \tfrac{1}{2} \|\x_{t}  -  \x_{t+1}\|_{\v_t}^2  + \la \y_t - \x_{t+1}, \x_{t+1} - \x_{t}\ra_{\v_t}  \nn\\
\overset{\step{3}}{=} &~ \ts Q_t + \la \x_{t+1}-\x_{t}, \nabla f(\x_{t}) - \nabla f(\y_t)\ra +  \tfrac{1}{2}\|\x_{t}  -  \y_t\|_{\v_t}^2  - \tfrac{1}{2}\|\y_t - \x_{t+1}\|_{\v_t}^2\nn\\
\overset{\step{4}}{\leq} &~ \ts Q_t + L \|\x_{t+1}-\x_{t}\| \|\x_{t} - \y_t\| + \tfrac{1}{2}\|\x_{t}  -  \y_t\|_{\v_t}^2  - \tfrac{1}{2}\|\y_t - \x_{t+1}\|_{\v_t}^2 \nn\\
\overset{\step{5}}{=} &~ Q_t +  \ts  \sigma_{t-1} L\|\x_{t+1}-\x_{t}\| \|\x_{t} - \x_{t-1}\|  \nn\\
&~ \ts + \tfrac{ (\sigma_{t-1})^2}{2}\|\x_{t} - \x_{t-1}\|_{\v_t}^2  - \tfrac{1}{2}\|\x_{t+1}-\x_{t} - \sigma_{t-1} (\x_{t} - \x_{t-1})\|_{\v_t}^2 \nn\\
\overset{}{=} &~  Q_t +  \sigma_{t-1} L\|\x_{t+1}-\x_{t}\| \|\x_{t} - \x_{t-1}\|  - \tfrac{1}{2}\|\x_{t+1} - \x_{t}\|_{\v_t}^2 + \sigma_{t-1} \la \x_{t+1}-\x_{t}, \x_{t} - \x_{t-1}\ra_{\v_t} \nn\\
\overset{\step{6}}{\leq}&~ \ts Q_t + \tfrac{\sigma_{t-1} L}{2} \|\x_{t+1}-\x_{t}\|_2^2 + \tfrac{\sigma_{t-1} L}{2} \|\x_{t} - \x_{t-1}\|_2^2 \nn\\
&~ \ts - \tfrac{1}{2}\|\x_{t+1} - \x_{t}\|_{\v_t}^2 + \tfrac{\sigma_{t-1}}{2}  \|\x_{t+1}-\x_{t}\|_{\v_t}^2  + \tfrac{\sigma_{t-1}}{2} \|\x_{t} - \x_{t-1}\|_{\v_t}^2  \nn \\
\overset{}{=}&~  Q_t + \ts - \tfrac{1}{2}\|\x_{t+1} - \x_{t}\|_{\v_t}^2 + \underbrace{\tfrac{1}{2} \|\x_{t} - \x_{t-1}\|_{ \sigma_{t-1} (\v_t + L)}^2}_{:= \XX^{t}} + \tfrac{1}{2} \|\x_{t+1} - \x_{t}\|_{\sigma_{t-1} (\v_t + L)}^2,
\end{align}
\noi where step \step{1} uses the Pythagoras Relation as in Lemma \ref{lemma:pyth} that $\x_+=\x_{t+1}$, $\x=\x_{t}$, $\a = \a_t$, and $\v=\v_t$; step \step{2} uses $\a_t := \y_t - {\g_t}\div \v_t$; step \step{3} uses the definition of $Q_t$, and the Pythagoras Relation as in Lemma \ref{lemma:pyth} that $\x_+=\x_{t+1}$, $\x=\x_{t}$, $\a = \y_t$, and $\v=\v_t$; step \step{4} uses $L$-smoothness of $f(\cdot)$; step \step{5} uses $\y_{t} - \x_{t} = \sigma_{t-1} (\x_{t} - \x_{t-1})$; step \step{6} uses $a b \leq \frac{a^2}{2} + \frac{b^2}{2}$ for all $a,b\in\Rn$, and $\la \a,\b \ra_{\v} \leq \frac{1}{2}\|\a\|_{\v}^2+\frac{1}{2}\|\b\|_{\v}^2$ for all $\a,\b,\v\in\Rn^n$ with $\v\geq \zero$.

\noi Define $\widehat{\ZZ}_t := F(\x_{t}) - F(\bar{\x})+  \tfrac{1}{2} \|\x_{t} - \x_{t-1}\|_{\sigma_{t-1} (\v_t + L)}^2$. Given Inequality (\ref{eq:FFXX}), we have the following inequalities for all $t\geq 0$ almost surely:
\begin{align}
~& \ts\widehat{\ZZ}_{t+1} - \widehat{\ZZ}_t   - Q_t  \nn\\
\leq~ & \ts \tfrac{L}{2} \|\x_{t+1}-\x_{t}\|_2^2- \tfrac{1}{2} \|\x_{t+1} - \x_{t}\|_{\v_t}^2 + \tfrac{1}{2} \|\x_{t+1}-\x_{t}\|_{\sigma_{t-1} (\v_t + L)}^2  + \XX^{t+1} \nn\\
\overset{\step{1}}{=}~ & \ts  \tfrac{L}{2} \|\x_{t+1}-\x_{t}\|_2^2- \tfrac{1}{2} \|\x_{t+1} - \x_{t}\|_{\v_t}^2   + \tfrac{1}{2} \| \x_{t+1} - \x_{t}\|_{ [\sigma_{t-1} \v_{t} + \sigma_{t-1} L + \sigma_t \v_{t+1} + \sigma_t L ] }^2  \nn\\
=~ & \ts \tfrac{L + \sigma_{t-1} L + \sigma_t L }{2} \|\x_{t+1}-\x_{t}\|_2^2  + \tfrac12\langle\sigma_{t-1}\v_t+\sigma_t\v_{t+1}-\v_t,\d_t\odot\d_t\rangle  \nn\\
\overset{\step{2}}{\leq} ~&  \tfrac{3L}{2}\|\x_{t+1}-\x_{t}\|_2^2  - \tfrac{1}{2}(1-\theta)^2 \|\x_{t+1}-\x_{t}\|_{\v_t}^2 \nn\\
\overset{\step{3}}{\leq} ~& \ts \underbrace{\ts\tfrac{3L}{2}}_{:= c_2} \cdot \underbrace{\tfrac{1}{\min(\v_t)^2} \|\v_t\odot (\x_{t+1}-\x_{t})\|_{2}^2}_{:= \UUU_t}  - \underbrace{ \tfrac{1}{2}(\tfrac{1-\theta}{\kappa})^2}_{:= c_1}   \cdot \underbrace{\tfrac{1}{  \min(\v_t) } \| \v_t\odot (\x_{t+1}-\x_{t})\|_{2}^2}_{ := \SSS_t}, \nn
\end{align}
\noi where step \step{1} uses the definition of $\XX^{t}:=   \tfrac{1}{2} \|\x_{t} - \x_{t-1}\|_{\sigma_{t-1} (\v_t + L)}^2$; step \step{2} uses $\sigma_t\leq 1$ for all $t\geq 0$, and $\sigma_t \v_{t+1} + (\sigma_{t-1}-1)\v_t \leq -(1-\theta)^2 \v_t$ for all $t\geq 0$ as shown in Lemma \ref{lemma:prop:theta}\textbf{(b)}; step \step{3} uses the following two inequalities for all $\d\in\Rn^n$ with $\d=\x_{t+1}-\x_{t}$:
\begin{align}
\|\d\|_{2}^2 ~~\leq~~& \tfrac{1}{\min(\v_t)^2} \|\v_t\odot \d\|_{2}^2,\nn\\
\|\d\|_{\v_t}^2\kappa^2 ~~\geq~~& {\|\d\|_{\v_t}^2\frac{\max(\v_t)^2}{\min(\v_t)^2}
~~\geq~~\|\d\|_2^2\frac{\max(\v_t)^2}{\min(\v_t)}}~~\geq~~{\frac{\|\v_t\odot\d\|_2^2}{\min(\v_t)}}.\nn
\end{align}

The equality preceding step \step{2}, together with
$0\leq\sigma_{t-1},\sigma_t\leq\theta$ and
Lemma~\ref{lemma:prop:theta}\textbf{(b)}, gives
\eqref{eq:local:raw} before replacing $(1+2\theta)L/2$ by $c_2$.

\end{proof}

\end{lemma}

{For each fixed $q,b$, the following lemma bounds the objective values and iterate norms uniformly over the iteration index.}

{Define $\ddot{\theta}:=(1-\theta)^2/2$, $\ddot{L}:=(1+2\theta)L\bigl(\tfrac12+2(q-1)\bigr)$, $\breve{m}:=\max\{\underline{\rm v},\ddot{L}/\ddot{\theta}\}$, and $P=\frac{\breve{m}^2-\underline{\rm v}^{\,2}}{\alpha\underline{\rm v}^{\,2}}$.}

{\noi Define $R:=\max_{F(\x)\leq F(\x_0)+{\ddot L}P}\|\x-\x_0\|$.
The sublevel set contains $\x_0$ and is compact by
Assumption~\ref{ass:compact} and lower semicontinuity of $F$.
Hence the maximum is attained and finite.} 

{\noi Define $\ddot{G}:=(\kappa \breve{m}+L)(R+\theta\sqrt P)
+\|\nabla f(\x_0)\|
+2L(1+2\theta)\sqrt{(q-1)P}$.}

{Let $\underline{H}:=\min_{\z}\{h(\z)+\tfrac{\underline{\rm v}}4
\|\z-\x_0\|^2\}$, which is finite by Assumption~\ref{ass:prox}. }

\begin{lemma}[Boundedness of Objective Values and Iterates]
\label{lemma:level-localization}
Suppose Assumptions~\ref{ass:prox}, \ref{ass:compact}, and
\ref{ass:f} hold. Let Algorithm~\ref{alg:main} start from a fixed
$\x_0\in\operatorname{dom}(F)$ with refresh period $q$.
For all $t\geq 0$, we have the following results almost surely:

\begin{enumerate}[label=\textbf{(\alph*)}, leftmargin=18pt,
itemsep=2pt, topsep=1pt, partopsep=0pt]

\item $F(\x_t)\leq\overline{\rm F}$, where $\overline{\rm F}:=F(\x_0)+\ddot{L} P+ \ddot{L} \left(R+\frac{4 \ddot{G}}{\underline{\rm v}}
+2\sqrt{\frac{h(\x_0)-\underline{H}}{\underline{\rm v}}}\right)^2$.

\item $\|\x_t\|\leq\overline{\rm x}$, where $\overline{\rm x}:=\max\left\{1,\ \max_{F(\x)\leq\overline{\rm F}}\|\x\|\right\}$.

\end{enumerate}

\begin{proof} Let $\e_t:=\g_t-\nabla f(\y_t)$, $\d_t:=\x_{t+1}-\x_t$, and
$m_t:=\min(\v_t)$. {Set $\d_{-1}=\zero$.}

\paragraph{Bounding the SPIDER error.} Let $k_t:=q\lfloor t/q\rfloor$ be the last refresh index at or before $t$. Since $\e_{k_t}=\zero$, we have
\beq \label{eq:et2}
\ts \|\e_t\|^2 &=& \ts \left\|\sum_{s=k_t+1}^{t}(\e_s-\e_{s-1})\right\|^2
\overset{\step{1}}{\leq}
4L^2\left(\sum_{s=k_t+1}^{t}\|\y_s-\y_{s-1}\|\right)^2
\notag\\
&\overset{\step{2}}{\leq}& \ts 4L^2(q-1)\sum_{s=k_t+1}^{t}\|\y_s-\y_{s-1}\|^2
~~\overset{\step{3}}{\leq}~~ 4L^2(q-1)(1+2\theta)^2\sum_{s=0}^{t-1}\|\d_s\|^2,
\eeq
\noi where step \step{1} uses the triangle inequality, the SPIDER
update~\eqref{eq:SPIDER:v}, and component-wise $L$-smoothness; step \step{2} uses Cauchy--Schwarz and $t-k_t\leq q-1$; step \step{3} uses the extrapolation identity
$\y_s-\y_{s-1}=(1+\sigma_{s-1})\d_{s-1}-\sigma_{s-2}\d_{s-2}$,
the bound $0\leq\sigma_j\leq\theta$, and the triangle inequality
in the product Euclidean norm after extending the sum to
$s=1,\ldots,t$.

\paragraph{Summed descent inequality.} For every $T\geq 0$, we have
\beq
\ts \sum_{t=0}^{T-1}\|\e_t\|^2 \overset{\step{1}}{\leq}
4L^2(q-1)^2\sum_{s=1}^{T-1}\|\y_s-\y_{s-1}\|^2 \overset{\step{2}}{\leq}
4L^2(q-1)^2(1+2\theta)^2\sum_{t=0}^{T-1}\|\d_t\|^2,\nn
\eeq
\noi where step \step{1} uses the refresh-interval estimate in
step \step{2} of~\eqref{eq:et2}, and the fact that each squared displacement
is counted at most $q-1$ times within its refresh interval; step \step{2} uses the extrapolation bound from step \step{3} of~\eqref{eq:et2}.

Using the telescoping sum of Inequality (\ref{eq:local:raw}), $\widehat{\ZZ}_0=F(\x_0)-F(\bar{\x})$, and
$\widehat{\ZZ}_T\geq F(\x_T)-F(\bar{\x})$ yields:
\beq \label{eq:local:summed}
\ts F(\x_T)-F(\x_0) &{\leq}& \ts
\sum_{t=0}^{T-1}\left(
-\langle\e_t,\d_t\rangle
+\tfrac12(1+2\theta)L\|\d_t\|^2
-\ddot{\theta}\|\d_t\|_{\v_t}^2\right)
\notag\\
&\overset{\step{1}}{\leq}& \ts \sum_{t=0}^{T-1}(\ddot{L}-\ddot{\theta} m_t)\|\d_t\|^2,
\eeq
\noi where step \step{1} uses Cauchy--Schwarz, the bound on the summed
SPIDER error, and $\|\d_t\|_{\v_t}^2\geq m_t\|\d_t\|^2$.

\paragraph{Part (a). Boundedness of the objective values.}
By the definition of $\breve{m}$, the terms in~\eqref{eq:local:summed}
with $m_t\geq \breve{m}$ are nonpositive. We therefore estimate the
displacements at indices with $m_t<\breve{m}$.

First, for every $t\geq 1$, Lemma~\ref{lemma:prop:w}\textbf{(a)}
gives $m_t^2=\min_j [\v_t]_j^2
\geq \underline{\rm v}^{\,2}
+\alpha\sum_{s=0}^{t-1}\|\r_s\|^2$. Since $\r_s=\v_s\odot\d_s$ and
${[\v_s]_j}\geq\underline{\rm v}$ for every $j$, we have $\|\r_s\|^2
=\sum_j [\v_s]_j^2\cdot [\d_s]_j^2
\geq\underline{\rm v}^{\,2}\|\d_s\|^2$. Thus, whenever $m_t<\breve{m}$,
\beq
\ts \sum_{s=0}^{t-1}\|\d_s\|^2
~~\overset{\step{1}}{\leq}~~\frac{m_t^2-\underline{\rm v}^{\,2}}
     {\alpha\underline{\rm v}^{\,2}}~~\overset{\step{2}}{\leq}~~\frac{\breve{m}^2-\underline{\rm v}^{\,2}}{\alpha\underline{\rm v}^{\,2}} =P, \nn
\eeq
\noi where step \step{1} uses the two preceding bounds and also holds
at $t=0$, since the sum is empty and $m_0=\underline{\rm v}$; step \step{2} uses $m_t<\breve{m}$.

For every $t$ with $m_t<\breve{m}$, applying~\eqref{eq:local:summed}
with $T=t$ and using $\ddot{L}-\ddot{\theta} m_s\leq \ddot{L}$ gives
\beq \label{eq:Fxt:bound}
\ts F(\x_t)-F(\x_0)
\leq \sum_{s=0}^{t-1}(\ddot{L}-{\ddot{\theta}m_s})\|\d_s\|^2
\leq \ddot{L}\sum_{s=0}^{t-1}\|\d_s\|^2
\leq \ddot{L} P.
\eeq
In particular, whenever $m_t<\breve{m}$, the current iterate $\x_t$
belongs to the sublevel set defining $R$. 

We next bound $\|\d_t\|$ for such indices. This also controls
the step at which the metric first reaches the threshold,
when $m_t<\breve{m}$ and $m_{t+1}\geq \breve{m}$. We further derive:
\begin{align*}
\|\v_t\odot(\y_t-\x_0)-\g_t\|
&\overset{\step{1}}{\leq}
\|\v_t\odot(\y_t-\x_0)\|
+\|\nabla f(\y_t)-\nabla f(\x_0)\|
+\|\nabla f(\x_0)\|+\|\e_t\|\\
&\overset{\step{2}}{\leq}
(\kappa m_t+L)\|\y_t-\x_0\|
+\|\nabla f(\x_0)\|+\|\e_t\|\\
&\overset{\step{3}}{\leq}
(\kappa \breve{m} + L)\|\y_t-\x_0\|
+\|\nabla f(\x_0)\|+\|\e_t\|\\
&\overset{\step{4}}{\leq}
(\kappa \breve{m}+L)(R+\theta\sqrt P)
+\|\nabla f(\x_0)\|
+2L(1+2\theta)\sqrt{(q-1)P}
=:\ddot{G},
\end{align*}
\noi where step \step{1} uses $\g_t=\nabla f(\y_t)+\e_t$, and the triangle
inequality; step \step{2} uses
$\|\v_t\odot(\y_t-\x_0)\|\leq\max(\v_t)\|\y_t-\x_0\|$,
$\max(\v_t)\leq\kappa m_t$ from Lemma~\ref{lemma:prop:w}\textbf{(b)}, and
$\|\nabla f(\y_t)-\nabla f(\x_0)\|\leq L\|\y_t-\x_0\|$; step \step{3} uses $m_t<\breve{m}$; step \step{4} uses $\|\y_t-\x_0\|\leq R+\theta\sqrt P$,
which follows from the extrapolation update,
$\|\x_t-\x_0\|\leq R$, and $\|\d_{t-1}\|\leq\sqrt P$,
and $\|\e_t\|\leq 2L(1+2\theta)\sqrt{(q-1)P}$
from Inequality (\ref{eq:et2}).

Moreover, we have:
\begin{align*}
\underline{H}+\tfrac{\underline{\rm v}}4\|\x_{t+1}-\x_0\|^2
&\overset{\step{1}}{\leq}
h(\x_{t+1})+\tfrac12\|\x_{t+1}-\x_0\|_{\v_t}^2\\
&\overset{\step{2}}{\leq}
h(\x_0)+\langle\v_t\odot(\y_t-\x_0)-\g_t,\x_{t+1}-\x_0\rangle\\
&\overset{\step{3}}{\leq}
h(\x_0)+\ddot{G}\|\x_{t+1}-\x_0\|,
\end{align*}
\noi where step \step{1} uses the definition of $\underline{H}$ and
$\v_t\geq\underline{\rm v}\one$; step \step{2} uses proximal optimality by comparing
$\x_{t+1}$ with $\x_0$ and expanding the quadratic; step \step{3} uses Cauchy--Schwarz and the definition of $\ddot{G}$.

By Inequality \eqref{eq:local:summed}, we have
\beq
F(\x_T) &\overset{\step{1}}{\leq}& \ts F(\x_0)+\ddot{L}\sum_{\substack{0\leq t<T,~m_t<\breve{m}}}\|\d_t\|^2 \nn \\
&\overset{\step{2}}{\leq}& \ts F(\x_0)+\ddot{L}\left[
P+\left(R+\frac{4\ddot{G} }{\underline{\rm v}}
+2\sqrt{\frac{h(\x_0)-\underline{H}}{\underline{\rm v}}}\right)^2
\right]
~~=:~~\overline{\rm F}, \nn
\eeq
\noi where step \step{1} uses $\ddot{L}-\ddot{\theta} m_t\leq 0$ for $m_t\geq \breve{m}$ and
$\ddot{L}-\ddot{\theta} m_t\leq \ddot{L}$ for the remaining indices; step \step{2} uses the monotonicity of $m_t$. If the sum is
nonempty, let $j$ be its last index. Then
\beq
\ts \sum_{\substack{0\leq t<T\\m_t<\breve{m}}}\|\d_t\|^2
=\sum_{t=0}^{j-1}\|\d_t\|^2+\|\d_j\|^2
\leq P+\|\d_j\|^2, \nn
\eeq
where the last inequality uses $m_j<\breve{m}$. Solving the preceding
quadratic inequality at $t=j$ and using $\|\x_j-\x_0\|\leq R$
gives
\beq \ts \|\d_j\|
\leq \|\x_{j+1}-\x_0\|+\|\x_j-\x_0\|
\leq R+\frac{4\ddot{G} }{\underline{\rm v}}
+2\sqrt{\frac{h(\x_0)-\underline{H}}{\underline{\rm v}}}. \nn
\eeq
If the sum is empty, step \step{2} holds directly.

\paragraph{Part (b). Boundedness of the iterates.}
The sublevel set $\{\x:F(\x)\leq\overline{\rm F}\}$ contains
$\x_0$ and is compact by Assumption~\ref{ass:compact} and lower
semicontinuity. Since every iterate belongs to this set,
$\|\x_T\|\leq\overline{\rm x}$ for all $T\geq 0$, where
\beq
\ts \overline{\rm x}
~~:=~~\max\left\{1,\ \max_{F(\x)\leq\overline{\rm F}}\|\x\|\right\}.\nn
\eeq
\end{proof}
\end{lemma}

\begin{lemma}[Metric growth and extrapolation lower bound]\label{lemma:metric-growth}
Define $\overline{\rm x}$ in Lemma~\ref{lemma:level-localization}. For all $t\geq0$, we have:
\begin{enumerate}[label=\textbf{(\alph*)}, leftmargin=18pt, itemsep=2pt, topsep=1pt, parsep=0pt, partopsep=0pt]
\item $\tfrac{\min(\v_{t+1})}{\min(\v_{t})} \leq \grave{\kappa} := 1+  2 \kappa \overline{\rm{x}} \sqrt{ \alpha +\beta}$.

\item $\theta(1-\theta)/(\kappa\grave{\kappa})\leq\sigma_t$.
\end{enumerate}
\begin{proof}
The iterate bound in Lemma~\ref{lemma:level-localization} is now available for both claims.

\paragraph{Part (a).} For all $t\geq 0$, we derive the following results:
\beq
\ts \tfrac{ \min(\v_{t+1})  }{ \min(\v_t)  } &\overset{\step{1}}{=}& \ts  \sqrt{\tfrac{ \min( \v_t\odot \v_t + \alpha \| \r_t\|_2^2 \cdot \one  + \beta \r_t \odot \r_t )  }{ \min(\v_t\odot \v_t)  }} ~~\overset{\step{2}}{\leq} ~~ \ts  \sqrt{ \tfrac{ \min(\v_t\odot \v_t) + (\alpha+\beta) \|\r_t \|_2^2    }{ \min(\v_t\odot \v_t)  }}   \nn\\
&\overset{\step{3}}{\leq} & \ts  \sqrt{ \tfrac{ \min(\v_t\odot \v_t) + (\alpha+\beta) \max(\v_t)^2\|\x_{t+1}-\x_{t}\|_2^2    }{ \min(\v_t\odot \v_t)  }}  ~~\overset{\step{4}}{\leq} ~~ \ts  \sqrt{ \tfrac{ \min(\v_t\odot \v_t) + 4 (\alpha+\beta) \max(\v_t)^2 \overline{\rm{x}}^2   }{ \min(\v_t\odot \v_t)  }}   \nn\\
&\overset{\step{5}}{\leq} & \ts  \sqrt{ 1 + 4 (\alpha +\beta) \kappa^2 \overline{\rm{x}}^2 }  ~~\overset{}{\leq} ~~ \ts  1+  2 \kappa \overline{\rm{x}} \sqrt{ \alpha +\beta} := \grave{\kappa},  \nn
\eeq
\noi where step \step{1} uses the update rule for $\v_{t+1}$ that $\v_{t+1} = \sqrt{\v_t \odot\v_t + \s_t}$; step \step{2} uses the fact that $\r_t\odot \r_t \leq \|\r_t\|_2^2$; step \step{3} uses $\r_t:= \v_t\odot \d_t$ with $\d_t:= \x_{t+1}-\x_{t}$; step \step{4} uses $\|\x_{t+1}-\x_{t}\|\leq \|\x_{t+1}\|+\|\x_{t}\|\leq 2 \overline{\rm{x}}$ {from Lemma~\ref{lemma:level-localization}}; step \step{5} uses $\max(\v_t)/\min(\v_t)\leq \kappa$ for all $t$.

\paragraph{Part (b).} We now establish the lower bound for $\sigma_t$. For all $t\geq 0$, we have:
\beq
\sigma_t &:=& \theta (1-   \sigma_{t-1}) \cdot \min( \v_t \div \v_{t+1}   ) ~~\overset{\step{1}}{\geq} ~~ \theta (1-   \theta) \cdot \min( \v_t \div \v_{t+1}   ) \nn\\
&\overset{\step{2}}{\geq}& \theta (1-   \theta) \cdot \tfrac{\min(\v_t)}{\max(\v_{t+1})  }  ~~\overset{}{=}~~ \theta (1-   \theta) \cdot  \tfrac{\min(\v_{t+1})}{\max(\v_{t+1})} \cdot \tfrac{\min(\v_t)}{\min(\v_{t+1})}  ~~\overset{\step{3}}{\geq}~~ \theta (1-   \theta) \cdot  \tfrac{1}{\kappa} \cdot \tfrac{1}{\grave{\kappa}}, \nn
\eeq
\noi where step \step{1} uses $\sigma_{t-1}\leq \theta$; step \step{2} uses $\min(\a\div \b)\geq \tfrac{\min(\a)}{\max(\b)}$ for all $\a\geq \zero$ and $\b>\zero$; step \step{3} uses $\tfrac{\min(\v_{t+1})}{ \max(\v_{t+1})}\geq \tfrac{1}{\kappa}$ and $\tfrac{\min(\v_t)}{\min(\v_{t+1})} \geq \tfrac{1}{\grave{\kappa}}$ for all $t$, as shown in Lemma~\ref{lemma:prop:w}\textbf{(b)} and Part~\textbf{(a)} of this lemma, respectively.
\end{proof}
\end{lemma}

To sum the Lyapunov inequality, we need cumulative bounds for its two adaptive terms. The next lemma expresses both sums in terms of the terminal accumulator, preparing the boundedness arguments that follow.

\begin{lemma} \label{lemma:AB}
We define $\VV_t$ as in Lemma \ref{lemma:prop:w}. We have the following results.

\begin{enumerate}[label=\textbf{(\alph*)}, leftmargin=18pt, itemsep=2pt, topsep=1pt, parsep=0pt, partopsep=0pt]

\item $\sum_{t=0}^T \SSS_t \leq s_1 \VV_{T+1}$, where $s_1 := \tfrac{2\grave{\kappa}}{\alpha}$.

\item $\sum_{t=0}^T \UUU_t \leq s_2 \sqrt{\VV_{T+1}}$, where $s_2 := \tfrac{4 \grave{\kappa}^2}{\alpha}\underline{\rm{v}}^{-1/2}$.

\end{enumerate}

\begin{proof}

We define $\VV_{t+1}:=\sqrt{ \underline{\rm{v}}^2 + (\alpha+\beta)M_{t}}$, where $M_t :=  \sum_{i=0}^{t} \|\r_i\|_2^2$.

\paragraph{Part (a).} We derive the following results:
\beq
\ts \sum_{t=0}^T \tfrac{ \|\r_t\|_2^2}{\min(\v_t)} &\overset{}{=}& \ts  \sum_{t=0}^T \tfrac{ \|\r_t\|_2^2 }{   { \min(\v_{t+1}) } } \cdot { \tfrac{{\min(\v_{t+1})}}{{\min(\v_{t})}} }   ~~\overset{\step{1}}{\leq}~~ \ts \ts \grave{\kappa} \cdot \sum_{t=0}^T \tfrac{ \|\r_t\|_2^2 }{   \sqrt{ \underline{\rm{v}}^2 + \alpha  \sum_{j=0}^t \|\r_j\|_2^2 } } \nn \\
&\overset{\step{2}}{\leq}&  \ts {\tfrac{2\grave{\kappa}}{\alpha}} \cdot  \sqrt{ \underline{\rm{v}}^2 + \alpha \sum_{t=0}^T \|\r_t\|_2^2 }~~\overset{}{\leq}~~ \ts \underbrace{{\tfrac{2\grave{\kappa}}{\alpha}}}_{:= s_1} \cdot \underbrace{\ts \sqrt{\underline{\rm{v}}^2 + (\alpha+\beta)  \sum_{t=0}^T \|\r_t\|_2^2}}_{:= \VV_{T+1}},
\eeq

\noi where step \step{1} uses {Lemma~\ref{lemma:metric-growth}\textbf{(a)}} that $\min(\v_{t+1})/\min(\v_t)\leq \grave{\kappa}$, and Lemma \ref{lemma:prop:w}\textbf{(a)}; step \step{2} uses Lemma \ref{lemma:online}.

\paragraph{Part (b).} We have the following results:
\[\begin{array}{rcl}
\ts  \sum_{t=0}^T \tfrac{ \|\r_t\|_2^2}{\min(\v_t)^2} & \overset{}{=} & \ts \sum_{t=0}^T \tfrac{\|\r_t\|_2^2}{\min(\v_{t+1})^2 } \cdot \left(\tfrac{\min(\v_{t+1})}{ \min(\v_t)}\right)^2~~ \overset{\step{1}}{\leq} ~~ \ts  \tfrac{\grave{\kappa}^2}{\alpha} \cdot \sum_{t=0}^T \tfrac{\|\r_t\|_2^2}{\underline{\rm{v}}^2/\alpha + \sum_{j=0}^t \|\r_j\|_2^2 } \\
& \overset{\step{2}}{\leq} & \ts \tfrac{\grave{\kappa}^2}{\alpha} \cdot \tfrac{1}{1/4} \cdot \left( 1 + \tfrac{1}{\underline{\rm{v}}^2/\alpha} \sum_{t=0}^T \|\r_t\|_2^2 \right)^{1/4}  ~~ \overset{}{=} ~~ \ts \tfrac{4 \grave{\kappa}^2}{\alpha \cdot \underline{\rm{v}}^{1/2} } \cdot\left(\underline{\rm{v}}^2 + \alpha \sum_{t=0}^T \|\r_t\|_2^2 \right)^{1/4} \\
& \overset{\step{3}}{\leq} & \ts \underbrace{\ts \tfrac{4 \grave{\kappa}^2}{\alpha \cdot \underline{\rm{v}}^{1/2} }}_{:= s_2} \cdot\underbrace{\ts \left(\underline{\rm{v}}^2 + (\alpha+\beta) \sum_{t=0}^T \|\r_t\|_2^2 \right)^{1/4}}_{:= \sqrt{\VV_{T+1}} } ,
\end{array}\]
\noi where step \step{1} uses Lemma \ref{lemma:prop:w}\textbf{(a)} and {Lemma~\ref{lemma:metric-growth}\textbf{(a)}}; step \step{2} uses Lemma \ref{lemma:online:pp} with $p=1/4$, $A_t = \|\r_t\|_2^2$, and $c =\underline{\rm{v}}^2/\alpha$; step \step{3} uses $\beta\geq0$.

\end{proof}

\end{lemma}

\subsection{Analysis for {\sf {ALFPG-D}}}

This subsection provides the convergence analysis of {\sf {ALFPG-D}}.

For the full-batch method, the gradient-estimation error vanishes. Combining the preceding descent
and summation estimates yields coupled bounds for the Lyapunov sequence and the metric accumulator; the next lemma closes these bounds uniformly in the iteration horizon.

\begin{lemma}\label{lemma:smmable}
(\rm{Boundedness of $\ZZ_t$ and $\VV_t$}) Define $w_1:=\underline{\rm{v}} +   \tfrac{\alpha +\beta}{c_1} \cdot c_2 s_1$, and $w_2:=\tfrac{\alpha +\beta}{c_1}$. Define $c_1$ and $c_2$ in Lemma \ref{prop:H}. Define $s_1$ and $s_2$ in Lemma \ref{lemma:AB}. For all $t\geq0$, we have:

\begin{enumerate}[label=\textbf{(\alph*)}, leftmargin=18pt, itemsep=2pt, topsep=1pt, parsep=0pt, partopsep=0pt]
\item {$\ZZ_t\leq\overline{\ZZ}$, where $\overline{\ZZ}:=2c_2^2s_2^2w_2+2\ZZ_0+2c_2s_2\sqrt{w_1}$.}

\item {$\VV_t\leq\overline{\rm v}$, where $\overline{\rm v}:=w_1+w_2\overline{\ZZ}$. }


\end{enumerate}

\end{lemma}

\begin{proof}

We define $\VV_{t+1}:=\sqrt{ \underline{\rm{v}}^2 + (\alpha+\beta)M_{t}}$, where $M_t :=  \sum_{i=0}^{t} \|\r_i\|_2^2$.

We define $\r_t:= \v_t\odot\d_t$, where $\d_t:= \x_{t+1}-\x_{t}$.

Initially, for the full-batch, deterministic setting where $\nabla f(\y_t)=\g_t$, we obtain from Lemma \ref{prop:H} that
\beq\label{eq:first:ZZ:D}
0\leq \ZZ_t - \ZZ_{t+1}+   \tfrac{c_2 \|\r_t\|_2^2}{\min(\v_t)^2}  -   \tfrac{ c_1 \|\r_t\|_2^2}{\min(\v_t)}
\eeq
\noi Multiplying both sides of Inequality (\ref{eq:first:ZZ:D}) by $\min(\v_t)$ yields:
\[\begin{array}{rcl}
0 \leq - c_1 \|\r_t\|_2^2 +  \min(\v_t)(\ZZ_t - \ZZ_{t+1})  +  \tfrac{c_2 \|\r_t\|_2^2 }{\min(\v_t)}.
\end{array}\]

\noi Summing this inequality over $t$ from $t=0$ to $T$, we obtain:
\beq \label{eq:sum:delta::D}
\ts  0 &\leq& \ts - c_1 \sum_{t=0}^T \|\r_t\|_2^2 + \sum_{t=0}^T \min(\v_t)(\ZZ_t - \ZZ_{t+1} ) + c_2 \sum_{t=0}^T \tfrac{\|\r_t\|_2^2}{\min(\v_t)}          \nn\\
& \overset{\step{1}}{\leq} & \ts  - c_1 \sum_{t=0}^T \|\r_t\|_2^2 + (\max_{t=0}^T \ZZ_t) \cdot \min(\v_T) +  c_2 s_1 \VV_{T+1}   \nn\\
& \overset{\step{2}}{\leq} & \ts  - \tfrac{c_1}{\alpha +\beta} \left( \VV_{T+1}^2 -\underline{\rm{v}}^2 \right)  + (\max_{t=0}^T \ZZ_t) \cdot \VV_{T+1} + c_2 s_1 \VV_{T+1},
\eeq
\noi where step \step{1} uses Lemma \ref{lemma:remove:log} with $A_i = \min(\v_i)$ for all $i\in[T]$ with $A_1\leq A_2 \leq \ldots \leq A_T$, and $B_j = \ZZ_j$ for all $j\geq 0$, and Lemma \ref{lemma:AB} that $\sum_{t=0}^T \SSS_t \leq s_1 \VV_{T+1}$; step \step{2} uses the definition of $\VV_T$, along with the facts that $\v_t\leq \v_{t+1} \leq \VV_{t+1}$ and $\ZZ_i\geq 0$.

Notably, the inequality $\sum_{t=0}^T \min(\v_t) (\ZZ_t-\ZZ_{t+1})\leq \min(\v^T) \sum_{t=0}^T (\ZZ_t - \ZZ_{t+1})$ does not hold in general, as the sequence $\{\ZZ_t\}_{t=0}^T$ is not necessarily monotonic. To address this, we instead employ the alternative inequality which can be implied by Lemma \ref{lemma:remove:log}: $\sum_{t=0}^T \min(\v_t) (\ZZ_t-\ZZ_{t+1})\leq \min(\v_T) \max_{t=0}^T \ZZ_t$.

In view of Inequality (\ref{eq:sum:delta::D}), we have the following quadratic inequality for all $T\geq 0$:
\beq
\ts \tfrac{c_1}{\alpha +\beta} (\VV_{T+1})^2 \leq \left(c_2 s_1 + \max_{t=0}^T \ZZ_t \right) \cdot \VV_{T+1} + \tfrac{c_1}{\alpha +\beta}\underline{\rm{v}}^2.\nn
\eeq
Applying Lemma \ref{lemma:quad:abc} with $a=\tfrac{c_1}{\alpha +\beta}$, $b=c_2 s_1 + \max_{t=0}^T \ZZ_t$, $c=\tfrac{c_1}{\alpha +\beta}\underline{\rm{v}}^2$, and $x=\VV_{T+1}$ yields:
\beq \label{eq:upper:bound:VVT}
\VV_{T+1} &\leq& \ts  \sqrt{c/a} + {b}/{a} ~~=~~ \ts \underline{\rm{v}} + \tfrac{\alpha +\beta}{c_1} \cdot \left(c_2 s_1 + \max_{t=0}^T \ZZ_t\right)\nn\\
&=& \ts \underbrace{\ts \underline{\rm{v}} +   \tfrac{\alpha +\beta}{c_1} \cdot c_2 s_1}_{:= w_1} + \underbrace{ \tfrac{\alpha +\beta}{c_1} }_{:= w_2}\cdot \max_{t=0}^T \ZZ_t,
\eeq
The upper bound for $\VV_{T+1}$ is established in Inequality (\ref{eq:upper:bound:VVT}), but it depends on the unknown variable $(\max_{t=0}^T \ZZ_t)$.

\paragraph{Part (a).} We now show that $(\max_{t=0}^T \ZZ_t)$ is always bounded above by a finite constant $\overline{\ZZ}$. Dropping the negative term $-\tfrac{ c_1\|\r_t\|_2^2}{\min(\v_t)}$ on the right-hand side of Inequality (\ref{eq:first:ZZ:D}) and summing over $t$ from $t=0$ to $T$ yields:
\beq  \label{eq:upper:bound:HHT}
\ts \ZZ_{T+1}    &\leq& \ts \ZZ_0 + c_2 \sum_{t=0}^T \tfrac{ \|\r_t\|_2^2}{\min(\v_t)^2} \nn\\
& \overset{\step{1}}{\leq} & \ts \ZZ_0 + c_2 s_2 \sqrt{\VV_{T+1}} \nn\\
& \overset{\step{2}}{\leq} & \ts\ZZ_0 + c_2 s_2 \sqrt{ w_1 + w_2 \max_{t=0}^T \ZZ_t   }  \nn\\
& \overset{\step{3}}{\leq} & \ts \underbrace{ \ts \ZZ_0 + c_2 s_2 \sqrt{ w_1 }}_{:= \dot{a}}    + \underbrace{ c_2 s_2 \sqrt{  w_2}}_{:= \dot{b}}\cdot\sqrt{ \max_{t=0}^T \ZZ_t   }   \label{eq:resursive:bounding:1} \\
& \overset{\step{4}}{\leq} & \ts \max(\ZZ_{0}, 2 \dot{b}^2 + 2 \dot{a} )  \nn\\
& \overset{\step{5}}{=} & \ts 2 \dot{b}^2 + 2 \dot{a}  ,
\eeq
\noi where step \step{1} uses Lemma \ref{lemma:AB}\textbf{(b)}; step \step{2} uses Inequality (\ref{eq:upper:bound:VVT}); step \step{3} uses $\sqrt{a+b}\leq \sqrt{a}+\sqrt{b}$ for all $a,b\geq 0$; step \step{4} uses Lemma \ref{lemma:bound:zzz} with $a=\dot{a}$ and $b=\dot{b}$; step \step{5} uses the fact that $\dot{a}\geq \ZZ_0$. Thus one admissible explicit choice is
\begin{equation}\label{eq:det:lyapunov:bound}
\overline{\ZZ}:=2\dot b^{2}+2\dot a
=2c_2^2s_2^2w_2+2\ZZ_0+2c_2s_2\sqrt{w_1}.
\end{equation}

\paragraph{Part (b).} Define
\begin{equation}\label{eq:det:metric:bound}
\begin{aligned}
\overline{\rm v}
&:=w_1+w_2\overline{\ZZ}\\
&=w_1+w_2\bigl(2c_2^2s_2^2w_2+2\ZZ_0
                 +2c_2s_2\sqrt{w_1}\bigr),
\end{aligned}
\end{equation}
where $w_2=(\alpha+\beta)/c_1$ and
$w_1=\underline{\rm v}+w_2c_2s_1$.
Equation~\eqref{eq:upper:bound:VVT} and Part~\textbf{(a)} give
\[
\VV_{T+1}\leq w_1+w_2\max_{0\leq t\leq T}\ZZ_t
\leq w_1+w_2\overline{\ZZ}=\overline{\rm v}.
\]
The same bound holds for $\VV_0=\underline{\rm v}$.

\end{proof}

\subsubsection{Proof of Theorem \ref{the:it}}
\label{app:the:it}

The metric boundedness established above makes the squared proximal displacements summable. We first record this bound and then use proximal optimality to convert it into the stationarity estimate of Theorem \ref{the:it}.

\begin{proof}

Define $\d_t:=\x_{t+1}-\x_{t}$ and $\r_t:=\v_t \odot \d_t$.

Using $\overline{\rm v}$ from \eqref{eq:det:metric:bound}, set $\ts \breve{\rm{X}}:=\frac1\alpha\big[
\big( \overline{\rm v}/\underline{\rm v}\big)^2-1\big]$, and $C:=2(L+\overline{\rm v}) \breve{\rm{X}}^{1/2}$.

\paragraph{Part (a).} We have the following inequalities:
\beq \label{eq:bnd:xxt}
\ts \sum_{t=0}^T \|\x_{t+1}-\x_{t}\|_2^2& \overset{\step{1}}{\leq} & \ts \tfrac{1}{\underline{\rm{v}}^2}{\sum_{t=0}^T\min(\v_t)^2\|\x_{t+1}-\x_{t}\|_2^2} \nn\\
& \overset{\step{2}}{\leq} & \ts \tfrac{1}{\underline{\rm{v}}^2}\sum_{t=0}^T \|\v_t\odot (\x_{t+1}-\x_{t})\|_{2}^2 \nn\\
& \overset{\step{3}}{=} & \ts \tfrac{1}{\underline{\rm{v}}^2} \sum_{t=0}^T \|\r_t\|_2^2 = \tfrac{1}{\underline{\rm{v}}^2} M_T\nn\\
& \overset{\step{4}}{\leq} & \ts\tfrac{1}{\underline{\rm{v}}^2}\tfrac{1}{\alpha} (\overline{\rm{v}}^2 - \underline{\rm{v}}^2) := \breve{\rm{X}},
\eeq
\noi where step \step{1} uses $\v_t\geq \underline{\rm{v}}$; step \step{2} uses ${\min(\v)\|\d\|\leq \|\v\odot\d\|}$ for all $\v,\d\in\Rn^n$ with $\v\geq \zero$; step \step{3} uses the definition of $\r_t:= \v_t\odot \d_t$; step \step{4} uses Lemma \ref{lemma:prop:w}\textbf{(a)} that $\underline{\rm{v}}^2  +  \alpha M_{t} \leq (\v_{t+1})^2\leq \overline{\rm{v}}^2$ for all $t$.

\paragraph{Part (b).} First, by the first-order necessary condition of $\x_{t+1}$ that $\x_{t+1}\in\Prox_h(\y_t-\g_t\div\v_t;\v_t)$, we have:
\beq \label{eq:first:order:cond}
\zero \in \partial h(\x_{t+1}) +  \g_t +  \v_t \odot (\x_{t+1}-\y_t).
\eeq
\noi Second, we obtain:
\beq \label{eq:subg:opt}
\ts \|\nabla f(\x_{t+1})  + \partial h(\x_{t+1})   \| &\overset{\step{1}}{{\leq}}& \| \ts \nabla f(\x_{t+1})  - \nabla f(\y_t) -  \v_t \odot (\x_{t+1}-\y_t) \| \nn\\
& \overset{\step{2}}{\leq} & \ts L \| \y_t - \x_{t+1} \| + \max(\v_t) \|\y_t -\x_{t+1}\| \nn\\
& \overset{\step{3}}{=} & \ts (L+\max(\v_t)) \cdot \| \x_{t} + \sigma_{t-1} (\x_{t} - \x_{t-1}) -  \x_{t+1} \| \nn\\
& \overset{\step{4}}{\leq} & \ts (L+ {\overline{\rm{v}}}) \cdot ( \| \x_{t}-  \x_{t+1} \|  + \|\x_{t} - \x_{t-1}\|),
\eeq
\noi where step \step{1} uses the inclusion in (\ref{eq:first:order:cond}) with $\g_t=\nabla f(\y_t)$, as in {\sf {ALFPG-D}}; step \step{2} uses the triangle inequality, the $L$-smoothness of $f(\cdot)$, and $\|\v_t\odot \a\|\leq \max(\v_t)\|\a\|$ for all $\a \in \Rn^n$; step \step{3} uses $\y_t=\x_{t}+\sigma_{t-1}(\x_{t}-\x_{t-1})$; step \step{4} uses $\v_t\leq \overline{\rm{v}}$, the triangle inequality, and $\sigma_{t-1}\leq 1$ for all $t$.

\noi Third, we obtain the following results:
\beq \label{eq:third:AEPG}
&&\ts \sum_{t=0}^{T} \|\partial h(\x_{t+1}) +  \nabla f(\x_{t+1}) \|_2^2 \nn\\
& \overset{\step{1}}{\leq} & \ts 2(L+ {\overline{\rm{v}}})^2 \cdot \sum_{t=0}^T   ( \|  \x_{t+1}- \x_{t}  \|_2^2  + \|\x_{t} - \x_{t-1}\|_2^2 ) \nn\\
& \overset{}{=} & \ts 2  (L+ {\overline{\rm{v}}})^2 \cdot \{\sum_{t=0}^T    \|  \x_{t+1}- \x_{t}  \|_2^2  +   \sum_{t=-1}^{T-1} \|\x_{t+1} - \x_{t}\|_2^2 \} \nn\\
& \overset{}{=} & \ts 2 (L+ {\overline{\rm{v}}})^2 \cdot \{ \| \x_{-1} - \x_0\|_2^2 - \| \x_{T+1} - \x_{T} \|_2^2 + 2 \sum_{t=0}^T    \| \x_{t+1}- \x_{t}  \|_2^2   \} \nn\\
& \overset{\step{2}}{\leq} & \ts {4} (L+ {\overline{\rm{v}}})^2 \cdot \sum_{t=0}^T   \| \x_{t+1}- \x_{t}  \|_2^2   \nn\\
& \overset{\step{3}}{\leq} & {4 (L+ {\overline{\rm{v}}})^2} \cdot {\breve{\rm{X}}} ~~{=~~C^2},
\eeq
\noi where step \step{1} uses Inequality (\ref{eq:subg:opt}); step \step{2} uses the choice $\x_{-1}=\x_0$ as shown in Algorithm \ref{alg:main}, and $- \|\x_{T+1} - \x_{T}\|^2 \leq 0$; step \step{3} uses Inequality (\ref{eq:bnd:xxt}).

Finally, using the inequality $\|\a\|_2^2\geq \tfrac{1}{T+1}(\|\a\|_1)^2$ for all $\a\in\Rn^{T+1}$, we deduce from Inequality (\ref{eq:third:AEPG}) that
\beq
\ts \tfrac{1}{T+1} \sum_{t=0}^{T} \|\partial h(\x_{t+1}) +  \nabla f(\x_{t+1}) \| {~~\leq~~\tfrac{C}{\sqrt{T+1}}}.\nn
\eeq

\end{proof}

\subsection{Analysis for {\sf ALFPG-FS}}

This subsection provides the convergence analysis of {\sf ALFPG-FS}.

We fix $q\geq 1$. For all $t\geq 0$, we denote $r_{t}:= \lfloor\tfrac{t}{q}\rfloor+1$ \footnote{For example, if $q=3$ and $t \in \{0,1,2,3,4,5,6,7,8,9\}$, then the corresponding values of $r_t$ are $\{1,1,1,2,2,2,3,3,3,4\}$.}, leading to $ (r_{t} - 1)q \leq t \leq r_{t} q  - 1$.

We first recall the SPIDER variance recursion from \cite{fang2018spider}. It relates the estimation error
accumulated between refreshes to the motion of the extrapolated iterates.

\begin{lemma} \label{eq:spider:estimator}
(\rm{Lemma 1 in} \cite{fang2018spider}) The SPIDER estimator produces a stochastic gradient $\g_t$ that, for all $t$ with ${(r_t-1)q+1} \leq t \leq r_t q - 1$, we have:  $\EEE[\|\g_t- \nabla f(\y_t)\|_2^2] - \EEE[\|\g_{t-1} - \nabla f(\y_{t-1})\|_2^2]
\leq \tfrac{L^2}{b} Y_{t-1}$, where $Y_i:=  \EEE [  \|\y_{i+1} - \y_{i}\|_2^2 ]$.

At a refresh iteration, $\g_t=\nabla f(\y_t)$ and the estimation error is zero.
\end{lemma}


Telescoping the SPIDER recursion within one refresh interval bounds the current estimation error. Combining that bound with the common Lyapunov inequality gives the stochastic descent estimate below.

\begin{lemma} \label{lemma:ZZ:ZZ}

For any positive constant $\phi$, we define $c_1 := \tfrac{1}{2}(\tfrac{1-\theta}{\kappa})^2$, $c'_2:= \tfrac{(3+\phi)L}{2}$, $c_3:= \tfrac{L}{2 \phi } \tfrac{q}{b}$. We define $Y_i:=  \EEE [  \|\y_{i+1} - \y_{i}\|_2^2 ]$. For all $t$ with $(r_t-1)q \leq t\leq r_t q-1$, we have:

 \begin{enumerate}[label=\textbf{(\alph*)}, leftmargin=18pt, itemsep=2pt, topsep=1pt, parsep=0pt, partopsep=0pt]

\item $\EEE[\|  \g_t - \nabla f(\y_t)    \|_2^2] \leq \tfrac{L^2}{b} \sum_{i= (r_{t}-1)q }^{t-1}Y_i$.

\item $\ZZ_{t+1} - \ZZ_t +  \EEE[c_1 \SSS_t - c'_2 \UUU_t ]   \leq   \ts  \tfrac{c_3}{q}\sum_{i=(r_{t}-1)q}^{t-1} Y_i$.
\end{enumerate}
\begin{proof}

We define $c_1 := \tfrac{1}{2}(\tfrac{1-\theta}{\kappa})^2$, $c'_2:= \tfrac{(3+\phi)L}{2}$, and $c_3:= \tfrac{L}{2 \phi }\tfrac{q}{b}$, where $\phi>0$ can be any constant.

We define $\ZZ_t := \EEE [F(\x_{t}) - F(\bar{\x})+\tfrac12\|\x_{t}-\x_{t-1}\|_{\sigma_{t-1}(\v_t+L)}^2]$.

We define $Y_i:=  \EEE [  \|\y_{i+1} - \y_{i}\|_2^2 ]$.

We define $\SSS_t := \tfrac{\|\r_t\|_2^2}{ \min(\v_t)  }$, and $\UUU_t := \tfrac{\|\r_t\|_2^2}{ \min(\v_t)^2}$.

\paragraph{Part (a).} Telescoping the inequality $\EEE [\|\g_t- \nabla f(\y_t)\|_2^2] - {\EEE[\|\g_{t-1} - \nabla f(\y_{t-1})\|_2^2]} \leq \ts  \tfrac{L^2}{b}\EEE [ \|\y_t - \y_{t-1}\|_2^2 ] $ (as stated in Lemma \ref{eq:spider:estimator}) over $t$ from $(r_{t}-1)q+1$ to $t$, where $t\leq r_{t} q -1$, we obtain:
\beq \label{eq:approximation:vf}
&& \ts \EEE[\|  \g_t - \nabla f(\y_t)    \|_2^2] \nn\\
&\overset{}{\leq} &  \ts \EEE[\|\g^{ (r_t-1)q} - \nabla f (\y^{  (r_t-1)q })\|_2^2] + \tfrac{L^2}{b} \sum_{i= (r_{t}-1)q+1 }^{t}\EEE[ \|\y_{i} - \y_{i-1}\|_2^2 ]     \nn\\
&\overset{\step{1}}{=} &  \ts 0 + \tfrac{L^2}{b} \sum_{i= (r_{t}-1)q }^{t-1} \underbrace{\ts \EEE[ \|\y_{i+1} - \y_{i}\|_2^2 ]}_{:= Y_i},
\eeq
\noi where step \ding{172} uses $\g_{j} = \nabla f (\y_{j})$ when $j$ is a multiple of $q$. Notably, Inequality (\ref{eq:approximation:vf}) holds for every $t$ of the form $t = (r_{t}-1)q$, since at these points we have $\g_t=\nabla f(\y_t)$.

\paragraph{Part (b).} For every $t$, Lemma~\ref{prop:H} and Young's inequality give
\beq
\ZZ_{t+1}-\ZZ_t+c_1\EEE[\SSS_t]
&\leq& \ts \frac{3L}{2}\EEE[\UUU_t]
+\frac{\phi L}{2}\EEE\|\d_t\|^2
+\frac1{2\phi L}\EEE\|\nabla f(\y_t)-\g_t\|^2\notag\\
&\leq& \ts c'_2\EEE[\UUU_t]
+\frac{c_3}{q}\sum_{i=(r_t-1)q}^{t-1}Y_i,\label{eq:to:be:tel}
\eeq
where $\phi>0$, $c'_2=(3+\phi)L/2$, and $c_3=Lq/(2\phi b)$. The second inequality uses $\|\d_t\|^2\leq\UUU_t$ and Part (a). We may take $L>0$ without loss of generality by enlarging a zero smoothness bound.

\end{proof}

\end{lemma}

The stochastic descent estimate involves motion at the extrapolated points rather than at the proximal
iterates. The next lemma converts this motion into the adaptive displacement sums already bounded
in Lemma \ref{lemma:AB}.

\begin{lemma}\label{lemma:AB2}

We define $V_t:=\min(\v_t)$, and $Y_i:=  \EEE [  \|\y_{i+1} - \y_{i}\|_2^2 ]$. We have:

\begin{enumerate}[label=\textbf{(\alph*)}, leftmargin=18pt, itemsep=2pt, topsep=1pt, parsep=0pt, partopsep=0pt]

\item ${\EEE\!\left[\sum_{t=0}^{T}V_t\|\y_{t+1}-\y_t\|_2^2\right]} \leq  u_1 \EEE[\VV_{T+1}]$, where $u_1:= (8+2\grave{\kappa}) s_1$.

\item $\sum_{t=0}^{T} Y_t \leq  u_2 \EEE[\sqrt{\VV_{T+1}}]$, where $u_2 := 10s_2$.

\end{enumerate}

\begin{proof}

Define $V_t:=\min(\v_t)$, and $Y_i:=  \EEE [  \|\y_{i+1} - \y_{i}\|_2^2 ]$.

First, for all $\tau>0$, we have the following results:
\beq \label{eq:v:y:y:1}
\| \y_{t+1} - \y_t\|_2^2 &\overset{\step{1}}{=} & \| (\x_{t+1} + \sigma_t \d_t) - (\x_{t} + \sigma_{t-1} \d_{t-1})\|_2^2 \nn\\
&=& \| (1+\sigma_t) \d_t - \sigma_{t-1} \d_{t-1}  \|_2^2 \nn\\
&\overset{\step{2}}{\leq} &  [(1+\tau) \| (1+\sigma_t) \d_t \|_2^2 + (1+1/\tau) \|\sigma_{t-1} \d_{t-1}  \|_2^2 ] \nn\\
&\overset{\step{3}}{\leq} &  4 (1+\tau) \| \d_t \|_2^2 + (1+1/\tau) \|\d_{t-1}  \|_2^2  , \eeq
\noi where step \step{1} uses $\y_{t+1} = \x_{t+1} + \sigma_t \d_t$; step \step{2} uses $\|\a+\b\|_2^2\leq (1+\tau)\|\a\|_2^2+(1+1/\tau)\|\b\|_2^2$ for all $\tau>0$; step \step{3} uses $\sigma_t\leq \theta<1$.

Second, we obtain the following inequalities:
\beq \label{eq:v:y:y:2}
\min(\v_t) \| \y_{t+1} - \y_t\|_2^2 &\overset{\step{1}}{{\leq}} & \min(\v_t) \cdot [ 4 (1+\tau) \|\d_t \|_2^2 + (1+1/\tau) \|\d_{t-1}  \|_2^2 ] \nn\\
&\overset{\step{2}}{\leq} & 4 (1+\tau) \min(\v_t) \|\d_t \|_2^2 + (1+1/\tau) \grave{\kappa} \min(\v_{t-1}) \|\d_{t-1}  \|_2^2 \nn\\
&\overset{\step{3}}{\leq} &  (4+\grave{\kappa}) \left(\min(\v_t) \|\d_t \|_2^2 +  \min(\v_{t-1}) \|\d_{t-1}  \|_2^2\right),
\eeq
\noi where step \step{1} uses Inequality (\ref{eq:v:y:y:1}); step \step{2} uses  $\min(\v_t)\leq \min(\v_{t-1})\grave{\kappa}$, as shown in {Lemma~\ref{lemma:metric-growth}\textbf{(a)}}; step \step{3} uses the choice $\tau = \tfrac{\grave{\kappa}}{4}$. Note that Inequality \eqref{eq:v:y:y:2} also holds at $t=0$ with $\v_{-1}=\v_0$ and $\d_{-1}=\zero$.

\paragraph{Part (a).} We have the following inequalities:
\beq
\ts\sum_{t=0}^{T} \min(\v_t) \| \y_{t+1} - \y_t\|_2^2 &\overset{\step{1}}{\leq} & \ts  (4+\grave{\kappa})  \left( \sum_{t=0}^{T} \min(\v_t) \|\d_t \|_2^2 + \sum_{t=0}^{T}\min(\v_{t-1}) \|\d_{t-1}  \|_2^2\right) \nn\\
&\overset{}{=} & \ts  (4+\grave{\kappa})  \left( \sum_{t=0}^{T} \min(\v_t) \|\d_t \|_2^2 + \sum_{t=-1}^{T-1}\min(\v_{t}) \|\d_{t}  \|_2^2\right) \nn\\
&\overset{\step{2}}{=} & \ts  (4+\grave{\kappa})  \left( \sum_{t=0}^{T} \min(\v_t) \|\d_t \|_2^2 + \sum_{t=0}^{T-1}\min(\v_{t}) \|\d_{t}  \|_2^2\right) \nn\\
&\overset{}{\leq} & \ts (8+2\grave{\kappa})  \left( \sum_{t=0}^{T} \min(\v_t) \|\d_t \|_2^2 \right) \nn\\
&\overset{\step{3}}{\leq} & \ts (8+2\grave{\kappa}) \sum_{t=0}^{T} \min(\v_{t}) \cdot \tfrac{1}{\min(\v_t)^2}\| \v_t\odot \d_t \|_2^2 \nn\\
&\overset{\step{4}}{\leq} & \ts \underbrace{\ts(8+2\grave{\kappa}) \cdot s_1}_{:= u_1}\cdot \VV_{T+1}, \eeq
\noi where step \step{1} uses Inequality (\ref{eq:v:y:y:2}); step \step{2} uses $\d_{-1} = \x_{0}-\x_{-1}=\zero$; step \step{3} uses $\min(\v_t)\|\d_t\|\leq \|\v_t\odot \d_t\|$; step \step{4} uses $\sum_{t=0}^T \SSS_t \leq s_1 \VV_{T+1}$ with $\SSS_t := \tfrac{\|\r_t\|_2^2}{ \min(\v_t)  }$, as shown in Lemma \ref{lemma:AB}\textbf{(a)}.

\paragraph{Part (b).} We obtain the following inequalities:
\beq
 \ts\sum_{t=0}^{T} \| \y_{t+1} - \y_t\|_2^2 &\overset{\step{1}}{\leq} & \ts \sum_{t=0}^{T} \left( 8\| \d_t \|_2^2 + 2 \|\d_{t-1}\|_2^2 \right) \nn\\
&\overset{\step{2}}{=} & \ts 8 \left( \sum_{t=0}^{T}\| \d_t \|_2^2\right) + 2 \left(\sum_{t=1}^{T} \|\d_{t-1}\|_2^2 \right) \nn\\
&\overset{}{=} & \ts 8 \left( \sum_{t=0}^{T}\| \d_t \|_2^2\right) + 2 \left(\sum_{t=0}^{T-1} \|\d_{t}\|_2^2 \right) \nn\\
&\overset{}{\leq} & \ts 10 \left( \sum_{t=0}^{T}\| \d_t \|_2^2\right)\nn\\
&\overset{\step{3}}{\leq } & \ts 10 \left( \sum_{t=0}^{T} \tfrac{1}{\min(\v_t)^2} \| \v_t\odot \d_t \|_2^2\right)\nn\\
&\overset{\step{4}}{\leq} & \ts \underbrace{10 \cdot s_2}_{:= u_2} \cdot \sqrt{\VV_{T+1}} , \nn
\eeq
\noi where step \step{1} uses Inequality (\ref{eq:v:y:y:1}) with $\tau=1$; step \step{2} uses $\d_{-1} = \x_{0}-\x_{-1}=\zero$; step \step{3} uses $\min(\v_t)\|\d_t\|\leq \|\v_t\odot \d_t\|$; step \step{4} uses $\sum_{t=0}^T \UUU_t \leq s_2 \sqrt{\VV_{T+1}}$ with $\UUU_t := \tfrac{\|\r_t\|_2^2}{ \min(\v_t)^2}$, as shown in Lemma \ref{lemma:AB}\textbf{(b)}.

Taking expectations in the two pathwise bounds proves the stated results.
\end{proof}

\end{lemma}

We next count how often each motion term occurs in the sums over a refresh interval. In particular, part (c) controls the total SPIDER error without introducing the metric-growth factor appearing in the weighted bound of part (a).


\begin{lemma}\label{lemma:sum:y:two}

We define $Y_i:=  \EEE [  \|\y_{i+1} - \y_{i}\|_2^2 ]$, $V_i:= \min(\v_i)$, and $q':= q \grave{\kappa}^{q-1}$. For all $t$ with $(r_t-1)q \leq t \leq r_t q-1$, we have:

\begin{enumerate}[label=\textbf{(\alph*)}, leftmargin=13pt, itemsep=2pt, topsep=1pt, parsep=0pt, partopsep=0pt]

\item $\sum_{j={(r_t-1)q}}^t  [ V_j \sum_{i=(r_{j}-1)q}^{j-1} Y_i ] \leq q' \sum_{i=(r_{t}-1)q}^{t-1} V_i  Y_i$.

\item $\sum_{j=(r_{t}-1)q}^{t}  [\sum_{i=(r_{j}-1)q}^{j-1} Y_i] \leq (q-1) \sum_{i=(r_{t}-1)q}^{t-1} Y_i$.

\item $\sum_{t=0}^{T}  [\sum_{i=(r_{t}-1)q}^{t-1} Y_i] \leq (q-1)\sum_{i=0}^{T}Y_i$.

\end{enumerate}

\begin{proof}

We define $V_t:=\min(\v_t)$, where $\{V_j\}_{j=0}^{\infty}$ is non-decreasing. We define $q' := q \grave{\kappa}^{q-1}$.

\paragraph{Bounding the Term $\left(t -  (r_t -1) q\right)$.} For any integer $t\geq 0$, we derive the following inequalities:
\beq \label{eq:length:tt}
t -  (r_t -1) q  ~~ \overset{\step{1}}{=}  ~~  t -  ( \lfloor\tfrac{t}{q}\rfloor+1  -1) q ~~=~~  t -  \lfloor\tfrac{t}{q}\rfloor  q ~~\overset{\step{2}}{\leq} ~~ q - 1,
\eeq
\noi where step \step{1} uses $r_t := \lfloor\tfrac{t}{q}\rfloor+1$; step \step{2} uses the fact that $t -  \lfloor\tfrac{t}{q}\rfloor  q\leq q-1$ for all integer $t\geq 0$ and $q\geq 1$.

\paragraph{Part (a).} For any $t$ with $t\geq (r_t-1)q$, we have the following results:
\beq \label{eq:vv:rt1:q:rt1q}
\tfrac{ \min(\v_t)}{ \min (\v_{ (r_t-1)q })} ~~= ~~\tfrac{ \min(\v_{(r_t-1)q+1})}{ \min (\v_{ (r_t-1)q })} \cdot \tfrac{ \min(\v_{(r_t-1)q+2})}{ \min (\v_{ (r_t-1)q+1 })}   \ldots \cdot \tfrac{ \min(\v_t)}{ \min (\v_{ t-1 })}~~\overset{\step{1}}{\leq} ~~ \grave{\kappa}^{q-1},
\eeq
\noi where step \step{1} uses the fact that the product length is at most $\left( [t] - [(r_{t}-1)q+1] + 1  \right)$ and Inequality (\ref{eq:length:tt}).

For all $t$ with $(r_t-1)q \leq t \leq r_t q-1$, we have:
\[{\begin{array}{rcl}
\ts \sum_{j={(r_t-1)q}}^t \left( V_j \cdot  \sum_{i=(r_{j}-1)q}^{j-1} Y_i  \right) &\overset{\step{1}}{\leq} & \ts \sum_{j={(r_t-1)q}}^t  \left( V_j \cdot  \sum_{i=(r_{j}-1)q}^{t-1}  Y_i \right) \\
&\overset{\step{2}}{=} & \ts \sum_{j={(r_t-1)q}}^t \left( V_j \cdot  \sum_{i=(r_{t}-1)q}^{t-1}  Y_i \right) \\
&\overset{\step{3}}{=} & \ts \sum_{i=(r_{t}-1)q}^{t-1}  \left(  Y_i \cdot   \sum_{j={(r_t-1)q}}^t V_j  \right) \\
&\overset{\step{4}}{\leq} & \ts \sum_{i=(r_{t}-1)q}^{t-1}  \left(  Y_i \cdot   \sum_{j={(r_t-1)q}}^t V_t  \right) \\
&\overset{\step{5}}{\leq} & \ts q V_t \sum_{i=(r_{t}-1)q}^{t-1}   Y_i \\
&\overset{\step{6}}{\leq} & \ts q (\grave{\kappa}^{q-1} V_{ (r_t-1)q}) \sum_{i=(r_{t}-1)q}^{t-1}   Y_i \\
&\overset{\step{7}}{\leq} & \ts \underbrace{\ts q \grave{\kappa}^{q-1}}_{\ts := q'}\cdot \sum_{i=(r_{t}-1)q}^{t-1} V_i  Y_i,
\end{array}}\]
\noi where step \step{1} uses $j\leq t$ for all $j\in [ (r_t-1)q,t]$; step \step{2} uses $r_j=r_t$ for all $j\in [ (r_t-1)q,t]$ with $t\in [(r_t-1)q , r_t q-1]$; step \step{3} uses the fact that $\sum_{j=\underline{j}}^{\overline{j}} ( \a_j \sum_{i=\underline{i}}^{ \overline{i}} \b_i) = \sum_{i=\underline{i}}^{ \overline{i}} ( \b_i \sum_{j=\underline{j}}^{  \overline{j}} \a_j )$ for all $\underline{i}\leq\overline{i}$ and $\underline{j}\leq \overline{j}$; step \step{4} uses $V_j\leq V_t$ as $j\leq t$; step \step{5} uses $t- (r_t-1)q +1 \leq q$; step \step{6} uses Inequality (\ref{eq:vv:rt1:q:rt1q}); step \step{7} uses $i\geq (r_t-1)q$.

\paragraph{Part (b).} For all $t$ with $(r_t-1)q \leq t \leq r_t q-1$, we have:
\[{\begin{array}{rcl}
\ts \sum_{j=(r_{t}-1)q}^{t}  \sum_{i=(r_{j}-1)q}^{j-1} Y_i &\overset{\step{1}}{=} & \ts \sum_{i=(r_{t}-1)q}^{t-1} (t-i) Y_i \\
&\overset{\step{2}}{\leq} & \ts \left( [t-1]- [(r_t-1)q] + 1\right) \cdot \sum_{i=(r_{t}-1)q}^{t-1}  Y_i  \\
&\overset{\step{3}}{\leq} & \ts (q-1)\sum_{i=(r_{t}-1)q}^{t-1}Y_i, \end{array}}\]
\noi where step \step{1} uses basic reduction; step \step{2} uses $i\geq (r_{t}-1)q$; step \step{3} uses Inequality (\ref{eq:length:tt}).

\paragraph{Part (c).} Each $Y_i$ can appear only at later indices in its own epoch, at most $q-1$ times. Thus, including a possibly incomplete final epoch, $\sum_{t=0}^{T}\sum_{i=(r_t-1)q}^{t-1}Y_i
\leq(q-1)\sum_{i=0}^{T}Y_i$.

\end{proof}
\end{lemma}

The preceding estimates control the accumulated SPIDER error through the extrapolated motion.
We now combine this control with pathwise descent and apply Jensen's inequality only after taking
expectations, obtaining a closed first-moment bound for the random metric accumulator.

\begin{lemma}\label{lemma:fi:preresults}
(\rm{Boundedness of $\ZZ_t$ and ${\VV_t}$}) Under Assumptions~\ref{ass:prox}, \ref{ass:compact}, and~\ref{ass:f},
let the sequence be generated by {\sf ALFPG-FS} from a fixed
$\x_0\in\operatorname{dom}(F)$. There exist finite constants $\overline{\ZZ}>0$ and $\overline{\rm v}>0$, independent of $t$, such that, for all $t\geq0$,

\textbf{(a)} $\ZZ_t :=\EEE[\widehat{\ZZ}_t ] \leq\overline{\ZZ}$.

\textbf{(b)} $\EEE[\VV_t]\leq\overline{\rm v}$.

\end{lemma}

\begin{proof}

Let $\e_t=\g_t-\nabla f(\y_t)$.

Define $\VV_{t+1}:=\sqrt{ \underline{\rm{v}}^2 + (\alpha+\beta)M_{t}}$, where $M_t :=  \sum_{i=0}^{t} \|\r_i\|_2^2$.

Define $\SSS_t := \tfrac{\|\r_t\|_2^2}{ \min(\v_t)  }$, $\UUU_t := \tfrac{\|\r_t\|_2^2}{ \min(\v_t)^2}$. Let $c_1:=\frac12\left(\frac{1-\theta}{\kappa}\right)^2$. Define $Y_i:=  \EEE [  \|\y_{i+1} - \y_{i}\|_2^2]$.

Define $u_2=10s_2$, where $s_2$ and $u_2$ are defined in
Lemmas~\ref{lemma:AB}(b) and~\ref{lemma:AB2}(b), respectively.

Define $\ZZ_t=\EEE[\widehat{\ZZ}_t]$, where $\widehat{\ZZ}_t
:=F(\x_{t})-F(\bar{\x}) +\frac12\|\x_{t}-\x_{t-1}\|_{\sigma_{t-1}(\v_t+L)}^2$. It holds that $\widehat{\ZZ}_0=\ZZ_0$.

\paragraph{Bounding the term $\EEE[\sum_{t=0}^{T}\|\e_t\|^2]$.} We derive:
\beq\label{eq:fs:summed:error}
\ts \EEE[\sum_{t=0}^{T}\|\e_t\|^2] &\xxxx{1}{\leq}& \ts \frac{L^2}{b}
\sum_{t=0}^{T}\sum_{i=(r_t-1)q}^{t-1}Y_i\nn\\
&\xxxx{2}{\leq}& \ts \frac{L^2(q-1)}b\sum_{i=0}^{T}Y_i~~\xxxx{3}{\leq}~~ \ts \frac{L^2(q-1)u_2}{b} \EEE[\sqrt{\VV_{T+1}}],
\eeq
\noi where step \step{1} uses Lemma~\ref{lemma:ZZ:ZZ}(a); step \step{2} uses Lemma~\ref{lemma:sum:y:two}(c); step \step{3} uses Lemma~\ref{lemma:AB2}(b).

\paragraph{Pathwise bounds for the adaptive sums.} We derive:
\beq \label{eq:fs:lower:S1}
\ts \SSS_t~~:=~~ \frac{\|\r_t\|^2}{\min(\v_t)}
~~\geq~~\frac{\|\r_t\|^2}{\VV_t}~~\xxxx{1}{=}~~ \frac{(\VV_{t+1}-\VV_t)(\VV_{t+1}+\VV_t)}{(\alpha+\beta)\VV_t}~~\geq~~\frac{2}{\alpha+\beta}(\VV_{t+1}-\VV_t),
\eeq
\noi where step \step{1} uses $\VV_{t+1}^2-\VV_t^2= (\alpha+\beta)(M_t-M_{t-1})=(\alpha+\beta)\|\r_t\|^2$. Summing this inequality and invoking Lemma~\ref{lemma:AB}(b) give
the following bounds on every sample path:
\begin{equation}\label{eq:fs:pathwise:sums}
\ts \sum_{t=0}^{T}\SSS_t\geq\frac{2}{\alpha+\beta}(\VV_{T+1}-\underline{\rm v}),
\qquad
\sum_{t=0}^T \UUU_t\leq s_2\sqrt{\VV_{T+1}}.
\end{equation}

\paragraph{Pathwise descent inequality.} We have from Lemma \ref{prop:H}:
\beq
\widehat{\ZZ}_{t+1}-\widehat{\ZZ}_t &\leq& \ts -\langle\d_t,\e_t\rangle
+\frac{3L}{2}\UUU_t-c_1\SSS_t \nn\\
&\overset{\step{1}}{\leq}& \ts \tfrac{L}{2}\|\d_t\|^2 + \tfrac{1}{2L}\|\e_t\|^2 + \frac{3L}{2}\UUU_t-c_1\SSS_t \nn\\
&\overset{\step{2}}{\leq}& \ts \tfrac{L}{2} \UUU_t + \tfrac{1}{2L}\|\e_t\|^2 + \frac{3L}{2}\UUU_t-c_1\SSS_t, \nn
\eeq
\noi where step \step{1} uses Young's inequality; step \step{2} uses $\|\d_t\|^2\leq\UUU_t$. Summing these pathwise inequalities yields:
\beq \label{eq:fs:pathwise:descent}
\ts \widehat{\ZZ}_{T+1}+c_1\sum_{t=0}^{T}\SSS_t
&\leq& \ts \ZZ_0+2L\sum_{t=0}^{T}\UUU_t
+\frac1{2L}\sum_{t=0}^{T}\|\e_t\|^2\nn\\
&\xxxx{1}{\leq}& \ts \ZZ_0 + 2Ls_2 \sqrt{\VV_{T+1}}
+\frac1{2L}\sum_{t=0}^{T}\|\e_t\|^2, \nn
\eeq
\noi where step \step{1} uses the second inequality in (\ref{eq:fs:pathwise:sums}). Taking expectations yields:
\beq\label{eq:fs:expected:descent}
\ts \ZZ_{T+1}
&\leq&\ts  -c_1\EEE \left[\sum_{t=0}^{T}\SSS_t\right] +\ZZ_0+ \EEE[2Ls_2 \sqrt{\VV_{T+1}} + \frac1{2L}\sum_{t=0}^{T}\|\e_t\|^2] \nn\\
&\xxxx{1}{\leq}&\ts - \tfrac{2c_1 }{\alpha+\beta}(\EEE[\VV_{T+1}] -\underline{\rm v}) +\ZZ_0+ \EEE[2Ls_2 \sqrt{\VV_{T+1}}] + \frac1{2L} \frac{L^2(q-1)u_2}{b} \EEE[\sqrt{\VV_{T+1}} ], \nn
\eeq
where step \step{1} uses Inequality (\ref{eq:fs:pathwise:sums}).

\paragraph{Bounding the term $\EEE[\VV_{T+1}]$.} Set $x:=\EEE[\VV_{T+1}]$ and $y:=\ZZ_{T+1}$. Using Jensen's inequality yields:
\beq\label{eq:resursive:bounding:2}
\tfrac{\alpha+\beta}{2c_1}\cdot y + x~~\leq~~ A+B\sqrt{x},
\eeq
\noi where $A:=\underline{\rm v} +\frac{ (\alpha+\beta) \ZZ_0}{2c_1}$, and $B:=\frac{(\alpha+\beta)}{2c_1}\cdot (2Ls_2 + \frac{L^2(q-1)u_2}{2L b})$.

Dropping the non-negative term $y$ in Inequality (\ref{eq:resursive:bounding:2}) and using $B\sqrt{x}\leq(x+B^2)/2$, we obtain:
\beq
\ts ~~x:=~~\EEE[\VV_{T+1}] ~~\leq ~~\overline{\rm v}~~:=~~2A+B^2.\nn
\eeq

Dropping the non-negative term $x$ and using $x\leq \overline{\rm v}$, we conclude that
\beq
\ts y~~:=~~\ZZ_{T+1} ~~\leq~~ \overline{\ZZ}~~:=~~ \tfrac{2c_1}{\alpha+\beta} \left(A+B \sqrt{\overline{\rm v}}\right).\nn
\eeq

\end{proof}

  \subsubsection{Proof of Theorem \ref{the:it:stoc}}
\label{app:the:it:stoc}

We now convert the first-moment metric bound into expected stationarity. The proof controls the displacement, estimator-error, and metric-scaled residual terms separately, and then counts all component-gradient evaluations.

\begin{proof}
Write
$\d_t:=\x_{t+1}-\x_t$, $\r_t:=\v_t\odot\d_t$, and
$\e_t:=\g_t-\nabla f(\y_t)$, with
$\d_{-1}=\r_{-1}=\zero$.

Let $s_2=\frac{4\grave{\kappa}^{2}}{\alpha\sqrt{\underline{\rm v}}}$, which is defined in Lemma~\ref{lemma:AB}(b).

Let $u_2=10s_2$, which is defined in Lemma~\ref{lemma:AB2}(b).

Let $\mathcal{R}(\x):=\operatorname{dist}(0,\partial F(\x))$.

\paragraph{Part (a).} We derive:
\beq \label{eq:fs:thm:motion:revised}
\ts \EEE\sum_{t=0}^{T}\|\d_t\|_2^2~~\xxxx{1}{\leq}~~
s_2\EEE[\sqrt{\VV_{T+1}}]
~~\xxxx{2}{\leq}~~ s_2\sqrt{\EEE[\VV_{T+1}]} ~~\leq ~~s_2\sqrt{\overline{\rm v}}~~:=~~D.
\eeq
where step \step{1} uses Lemma~\ref{lemma:AB}(b) that
$\sum_{t=0}^{T}\UUU_t\leq s_2\sqrt{\VV_{T+1}}$; step \step{2} uses Jensen's inequality.

Furthermore, we derive:
\beq\label{eq:fs:thm:error:revised}
\ts \EEE[\sum_{t=0}^{T}\|\e_t\|_2^2]~~\xxxx{1}{\leq}~~ \frac{L^2(q-1)u_2}{b} \EEE[\sqrt{\VV_{T+1}}]~~{\leq}~~ \frac{L^2(q-1)u_2}{b} \sqrt{\overline{\rm v}}~~:=E.
\eeq
\noi where step \step{1} uses Inequality (\ref{eq:fs:summed:error}).

The optimality condition for the proximal subproblem implies
\[
\zero\in\partial h(\x_{t+1})+\g_t
+\v_t\odot(\x_{t+1}-\y_t).
\]
Since $\y_t=\x_t+\sigma_{t-1}\d_{t-1}$,
$0\leq\sigma_{t-1}\leq\theta$, and
$\sigma_{t-1}\v_t\leq\theta\v_{t-1}$ coordinatewise
for $t\geq1$, smoothness of $f$ yields
\begin{align}\label{eq:fs:thm:residual:revised}
\mathcal{R}(\x_{t+1})
&\leq\bigl\|\nabla f(\x_{t+1})-\g_t
-\v_t\odot(\x_{t+1}-\y_t)\bigr\| \notag\\
&\leq L\bigl(\|\d_t\| + \theta\|\d_{t-1}\|\bigr)+\|\e_t\|+\|\r_t\|+\theta\|\r_{t-1}\|.
\end{align}
For $t=0$, the terms with index $-1$ vanish.

By Cauchy--Schwarz and \eqref{eq:fs:thm:motion:revised}--\eqref{eq:fs:thm:error:revised},
\beq
\ts \EEE[\sum_{t=0}^{T}\|\d_t\|] \leq\sqrt{(T+1)D}, \qquad \EEE[\sum_{t=0}^{T}\|\e_t\|] \leq\sqrt{ (T+1) E}.\nn
\eeq
The accumulator identity
$\VV_{T+1}^2-\underline{\rm v}^{2}
=(\alpha+\beta)\sum_{t=0}^{T}\|\r_t\|_2^2$
also gives
\beq
\ts \EEE[\sum_{t=0}^{T}\|\r_t\|_2]~~\leq~~\sqrt{\frac{T+1}{\alpha+\beta}}\,
\EEE\bigl[\sqrt{\VV_{T+1}^2-\underline{\rm v}^{2}}\bigr]~~\leq~~\sqrt{\frac{T+1}{\alpha+\beta}}\,\EEE[\VV_{T+1}]
\leq\sqrt{\frac{T+1}{\alpha+\beta}}\,\overline{\rm v}.\nn
\eeq
Summing~\eqref{eq:fs:thm:residual:revised} over
$t=0,\ldots,T$, and noting that the shifted sums do not
exceed their unshifted counterparts, we obtain
\beq \label{eq:third:AEPG:stoc}
\ts \frac{1}{T+1}\sum_{t=0}^{T}\EEE[\mathcal{R}(\x_{t+1})]
\leq\frac{{C_{q,b}}}{\sqrt{T+1}},\qquad
{C_{q,b}}:=L(1+\theta)\sqrt D+\sqrt{E}
+\frac{(1+\theta)\overline{\rm v}}{\sqrt{\alpha+\beta}}.
\eeq
This constant is finite and independent of $T$. At least one $\bar t\in\{1,\ldots,T+1\}$ satisfies $\EEE[\mathcal{R}(\x_{\bar t})]\leq {C_{q,b}}/\sqrt{T+1}$. If $J$ is uniform on $\{1,\ldots,T+1\}$ and independent of the
algorithm, then $\EEE[\mathcal{R}(\x_J)]=(T+1)^{-1}\sum_{t=0}^{T}\EEE[\mathcal{R}(\x_{t+1})]$. Hence both claims follow when $T+1\geq {C_{q,b}}^2\epsilon^{-2}$.

\paragraph{Part (b).} Assume $f_*=\inf f>-\infty$, $h_*=\inf h>-\infty$,
and $(q-1)/b\leq1$. Set
\beq
\ts m_t=\min(\v_t),~\Delta_0=F(\x_0)-f_*-h_*,~M=\underline{\rm v}+\frac{20L}{(1-\theta)^2},~P=\frac{M^2-\underline{\rm v}^2}{\alpha\underline{\rm v}^2}.\nn
\eeq
Since $m_0<M$, the index $\ell=\max\{0\leq t\leq T:m_t<M\}$ is well defined.
The metric update gives
$m_\ell^2\geq\underline{\rm v}^2+
\alpha\underline{\rm v}^2\sum_{t<\ell}\|\d_t\|^2$, so
$\sum_{t<\ell}\|\d_t\|^2\leq P$.

\noindent$\quad\bullet$\enspace
\emph{SPIDER error before the metric threshold.}
Let $I_t=\mathbf 1_{\{m_t<M\}}=\mathbf 1_{\{t\leq\ell\}}$.
Since $I_t$ is known before sampling and $I_t\leq I_{t-1}$ for $t\geq1$,
the conditional variance recursion at nonrefresh steps gives
\beq
\ts \EEE[I_t\|\e_t\|^2]~~\leq~~\EEE[I_{t-1}\|\e_{t-1}\|^2]
+\frac{L^2}{b}\EEE[I_t\|\y_t-\y_{t-1}\|^2]. \nn
\eeq
The error is zero at refreshes, so summing within refresh intervals
counts each motion term at most $q-1$ times. By the extrapolation rule,
$1+2\theta\leq3$, and $(q-1)/b\leq1$, we obtain
\beq
\ts \EEE [ \sum_{t=0}^{\ell}\|\e_t\|^2] ~~\leq~~\frac{L^2(q-1)}b\EEE\sum_{t=1}^{\ell}\|\y_t-\y_{t-1}\|^2
\leq9L^2P. \nn
\eeq

\noindent$\quad\bullet$\enspace
\emph{Last step below the threshold.}
Summing \eqref{eq:local:raw} up to $\ell-1$ and applying
Cauchy--Schwarz bounds the expected sum of the error inner products by
$3LP$. Since $(1+2\theta)L/2\leq2L$, this gives
$\EEE[F(\x_\ell)-f_*-h_*]\leq\Delta_0+5LP$.
Proximal comparison with $\x_\ell$, Young's inequality, smoothness, and
$\max(\v_\ell)\leq\kappa M$ yield
\beq
\ts \|\d_\ell\|^2\leq
\frac{4(h(\x_\ell)-h_*)}{\underline{\rm v}}
+\frac{4}{\underline{\rm v}^{\,2}}
\bigl((\kappa M+L)\|\d_{\ell-1}\|+\|\nabla f(\x_\ell)\|
+\|\e_\ell\|\bigr)^2. \nn
\eeq
Using $\|\nabla f(\x_\ell)\|^2\leq2L(f(\x_\ell)-f_*)$ and
$(u+v+w)^2\leq3(u^2+v^2+w^2)$, we obtain
\beq
\ts \EEE\sum_{t=0}^{\ell}\|\d_t\|^2
\leq D_0:=
P+\frac{12\bigl[(\kappa M+L)^2+9L^2\bigr]P}{\underline{\rm v}^{\,2}}
+\left(\frac4{\underline{\rm v}}+\frac{24L}{\underline{\rm v}^{\,2}}\right)
(\Delta_0+5LP). \nn
\eeq

\noindent$\quad\bullet$\enspace
\emph{Full-interval bounds.}
Taking $I_t=1$ in the SPIDER summation gives
$\EEE\sum_{t=0}^T\|\e_t\|^2\leq9L^2\EEE\sum_{t=0}^T\|\d_t\|^2$.
Summing \eqref{eq:local:raw} and applying Cauchy--Schwarz gives
\beq
\ts
\frac{(1-\theta)^2}{2}\EEE\sum_{t=0}^{T}m_t\|\d_t\|^2
\leq\Delta_0+5L\EEE\sum_{t=0}^{T}\|\d_t\|^2. \nn
\eeq
For $m_t\geq M$, we have $5L\leq(1-\theta)^2m_t/4$.
Absorbing these terms and using $m_t\geq\underline{\rm v}$ yields
\beq
\ts
\EEE\sum_{t=0}^T\|\d_t\|^2
\leq\frac{1}{\underline{\rm v}}\EEE\sum_{t=0}^Tm_t\|\d_t\|^2
\leq D^\star:=\frac{4(\Delta_0+5LD_0)}{(1-\theta)^2\underline{\rm v}}. \nn
\eeq
Consequently $\EEE\sum_{t=0}^T\|\e_t\|^2\leq9L^2D^\star$.
The accumulator identity gives
$\VV_{t+1}-\VV_t\leq(\alpha+\beta)\kappa^2m_t\|\d_t\|^2/2$, hence
\beq
\ts \EEE[\VV_{T+1}]
\leq\underline{\rm v}\left(1+\frac{(\alpha+\beta)\kappa^2D^\star}{2}\right). \nn
\eeq
Summing \eqref{eq:fs:thm:residual:revised} as in part~(a) and using
$1+\theta\leq2$ gives
\beq \label{eq:CCCC}
\ts \frac{1}{T+1}\sum_{t=0}^{T}\EEE[\mathcal{R}(\x_{t+1})]
\leq\frac{\overline{C}}{\sqrt{T+1}},\qquad
\overline{C}:=5L\sqrt{D^\star}
+\frac{\underline{\rm v}}{\sqrt{\alpha+\beta}}
\bigl(2+(\alpha+\beta)\kappa^2D^\star\bigr). 
\eeq
Thus $\overline{C}$ depends only on $L,\Delta_0,\alpha,\beta,\theta,\underline{\rm v}$.
For the independent uniform index $J$ of part~(a),
$\EEE[\mathcal{R}(\x_J)]\leq\overline{C}/\sqrt{T+1}$.

\paragraph{Part (c).} Let $b$ denote the mini-batch size and $q$ the refresh frequency of {\sf ALFPG-FS}. Set $H:=T+1$. When $\operatorname{mod}(t,q)=0$,
the full gradient $\nabla f(\y_t)$ is computed using $N$
component-gradient evaluations. This occurs
$F:=\lceil H/q\rceil$ times. Otherwise, the SPIDER update
evaluates the mini-batch at both $\y_t$ and $\y_{t-1}$,
using at most $2b$ component-gradient evaluations.
This occurs $H-F$ times. Hence the total number of
component-gradient evaluations satisfies
\begin{align*}
NF+2b(H-F)
&\overset{\step{1}}{\leq}
N+\left(\frac{N}{q}
       +2b\left(1-\frac1q\right)\right)H\\
&\overset{\step{2}}{\leq} N+3\sqrt N\,H,
\end{align*}
where step \step{1} bounds $NF$ using $F\leq H/q+1$ and bounds $2b(H-F)$ using $F\geq H/q$. Step \step{2} uses $q=b=\lceil\sqrt N\rceil$, leading to $N/q\leq\sqrt N$
and $2b(1-1/q)=2(q-1)\leq 2\sqrt N$. By part~(b), taking $H=\max\{1,\lceil \overline{C}^{\,2}\epsilon^{-2}\rceil\}$ ensures $\EEE[\mathcal{R}(\x_J)]\leq\epsilon$. Substitution into $N+3\sqrt N\,H$ gives the claimed count.

\end{proof}

\section{Proof for Section \ref{sect:iter2}}
\label{app:sect:iter2}

We first establish adaptive-weight estimates shared by
{\sf ALFPG-MVR1} and {\sf ALFPG-MVR2}.
For each variant, we then combine tracking-error estimates,
proximal descent, and cumulative moment bounds to control
iterate motion and establish the stationarity guarantees.
Throughout this appendix, Assumptions~\ref{ass:prox},
\ref{ass:oracle}, \ref{ass:bounded-oracle}, \ref{ass:smooth},
and~\ref{ass:regularizer} hold.
The analysis of {\sf ALFPG-MVR2} additionally uses
Assumption~\ref{ass:individual-smooth}.

Define $\s_t:=\m_t-\nabla f(\x_t)$, and $\r_t:=(\x_t-\x_{t+1})/\gamma_t$.

\subsection{Overview of the proof strategy for {\sf ALFPG-MVR1}}
\label{subsect:MVR1:overview}

Combining proximal descent with the EMA tracking-error
recursion gives
(Lemma~\ref{lemma:MVR1:descent})
\begin{equation}\label{eq:MVR1:overview:descent}
Q_{t+1}-Q_t
+\frac{\gamma_t}{2}\|\r_t\|^2
+\frac{\alpha_t}{2}\|\s_t\|^2
\leq
2\dot{\mathcal E}_{t+1}+2\mathcal E_{t+1},
\end{equation}
where $Q_t:=F(\x_t)+2\|\s_t\|^2-C_f\ln(1+\sum_{i=1}^{t-1}\|\r_i\|^2)$,
$\dot{\mathcal E}_t
=\alpha_t^2\|\g_t-\nabla f(\x_t)\|^2$,
and $\mathcal E_t$ is the conditionally mean-zero cross term. Define $V_T:=\sum_{t=1}^T\alpha_t^2\chi_t^2$, $H_T:=\sum_{t=1}^T \tfrac{\gamma_t^2}{\alpha_t} \|\r_t\|^2$, $U_T:=\sum_{t=1}^{T}\alpha_t\|\s_t\|^2$. The proof proceeds in four steps.

\begin{enumerate}[label=(\roman*), leftmargin=18pt,
itemsep=0.3ex, topsep=0.3ex, parsep=0pt]

\item Self-normalized sums and summation by parts yield: $V_T=\widetilde{\mathcal O}(1)$, $H_T=\widetilde{\mathcal O}(1)$, $\mathbb E[\sum_{t=1}^{T}\dot{\mathcal E}_t] =\widetilde{\mathcal O}(1)$, $\mathbb E[U_T]=\widetilde{\mathcal O}(1)$, and $\mathbb E[U_T^2]=\widetilde{\mathcal O}(1)$ (Lemma~\ref{lemma:MVR1:logTerm}).

\item Proximal descent gives $\EEE[W_T]=\widetilde{\mathcal O}(1)$, where $W_T:=\sum_{t=1}^{T}\gamma_t\|\r_t\|^2$. This results in $\alpha_T\sum_{t=1}^{T}\|\r_t\|^2
\leq W_T/\theta+W_T^2/(4\theta^2)$,
hence
$\mathbb E[\alpha_T\sum_{t=1}^{T}\|\r_t\|^2]
=\widetilde{\mathcal O}(1)$
(Lemma~\ref{lemma:MVR1:bridge}).

\item Smoothness of the expected loss and the noise
variance bound, together with the motion estimate, yield
$\mathbb E[\alpha_T^{-1}]
=\widetilde{\mathcal O}(1+\sigma\sqrt T)$. 

\item Proximal optimality gives
$\operatorname{dist}(0,\partial F(\x_{t+1}))
\leq(1+L\theta)\|\r_t\|+\|\s_t\|$. This further leads $\tfrac1T\sum_{t=1}^T\mathbb E[\mathcal{R}(\x_{t+1})] \leq\widetilde{\mathcal O}
(T^{-1/2}+\sigma^{1/2}T^{-1/4})$.

\end{enumerate}

\subsection{Overview of the proof strategy for {\sf ALFPG-MVR2}}
\label{subsect:MVR2:overview}

Combining the rescaled tracking-error recursion with
proximal descent gives
(Lemma~\ref{lemma:MVR2:descent})
\begin{equation}\label{eq:MVR2:overview:descent}
\mathcal Q_{t+1}-\mathcal Q_t
+\frac{\gamma_t}{2}\|\r_t\|^2
+\frac{\sqrt{\alpha_t}}{2}\|\s_t\|^2
\leq
4\dot{\mathcal E}_{t+1}+2\mathcal E_{t+1},
\end{equation}
where $\mathcal Q_t:=F(\x_t)+2P_t-C_fH_{t-1}$, $P_t:=\|\s_t\|^2/\sqrt{\alpha_t}$, $H_t:=\sum_{i=1}^t\frac{\gamma_i^2}{\sqrt{\alpha_i}}\|\r_i\|^2$, $\dot{\mathcal E}_t
=\alpha_t^{3/2}\|\g_t^--\nabla f(\x_{t-1})\|^2$, and $\mathcal E_t$ is the conditionally mean-zero cross term. Define $V_T:=\sum_{t=1}^T \alpha_t^{3/2} \chi_t^2$, $Q_T:=\sum_{t=1}^T\dot{\mathcal E}_t$, and $U_T:=\sum_{t=1}^{T}\sqrt{\alpha_t}\|\s_t\|^2$.
The proof proceeds in five steps.

\begin{enumerate}[label=(\roman*), leftmargin=18pt,
itemsep=0.3ex, topsep=0.3ex, parsep=0pt]

\item Taking expectations and summing
\eqref{eq:MVR2:overview:descent}
bounds the cumulative weighted residuals and tracking errors
in terms of normalized motion and weighted sampling noise (Lemma~\ref{lemma:MVR2:s:recursion}).

\item Self-normalized sums and summation by parts yield: $V_T= \widetilde{\mathcal O}(1)$, $\mathbb E[Q_{T}]=\widetilde{\mathcal O}(1)$, $\EEE[H_T]= \widetilde{\mathcal O}(1)$, and $\mathbb E[U_T]= \widetilde{\mathcal O}(1)$ (Lemma~\ref{lemma:MVR2:logTerm}).

\item Proximal descent gives $\EEE[W_T]=\widetilde{\mathcal O}(1)$ and $\EEE[W_T^2]=\widetilde{\mathcal O}(1)$, where $W_T:=\sum_{t=1}^{T}\gamma_t\|\r_t\|^2$. This results in $\sqrt{\alpha_T}\sum_{t=1}^{T}\|\r_t\|^2
\leq W_T/\theta+W_T^2/(4\theta^2)$,
hence
$\mathbb E[\sqrt{\alpha_T}\sum_{t=1}^{T}\|\r_t\|^2]
=\widetilde{\mathcal O}(1)$
(Lemma~\ref{lemma:MVR2:bridge}).

\item Smoothness of the expected loss and the noise
variance bound, together with the motion estimate, yield $\mathbb E[\alpha_T^{-1/2}]
\leq \widetilde{\mathcal O}(1+ T^{1/3} \sigma^{2/3})$.

\item Proximal optimality gives
$\operatorname{dist}(0,\partial F(\x_{t+1}))
\leq(1+L\theta)\|\r_t\|+\|\s_t\|$. This further leads to $\tfrac1T\sum_{t=1}^T\mathbb E[\mathcal{R}(\x_{t+1})] \leq\widetilde{\mathcal O}(T^{-1/2}+\sigma^{1/3}T^{-1/3})$.

\end{enumerate}

\subsection{Shared Preliminary Estimates}

Both expectation variants use the same running-minimum construction for the averaging weight. We first bound its variation; these estimates allow shifted weights to be compared and preserve contraction in the tracking-error recursions.

\begin{lemma}\label{lemma:weight}
(Adaptive-weight regularity) For $t\geq1$, define $\rho_t=\min_{1\leq k\leq t}\left(\frac{1+\max_{1\leq i\leq k}\chi_i^2}{1+\sum_{i=1}^k\chi_i^2}\right)^{p}$, where ${p}\in(0,1)$ and $\rho_{-1}=\rho_0:=1$. Then

\begin{enumerate}[label=\textbf{(\alph*)}, leftmargin=18pt,
itemsep=2pt, topsep=1pt, parsep=0pt, partopsep=0pt]

\item $0\leq \frac1{\rho_t}-\frac1{\rho_{t-1}}\leq {p}$.

\item $\rho_{t-1}-\rho_t \leq {p}\rho_{t-1}^{1+1/{p}}$.

\item $(1-\rho_t)\sqrt{\frac{\rho_{t-1}}{\rho_t}} \leq1-(1-{p})\rho_t$.

\end{enumerate}

\end{lemma}

\begin{proof}

\noi Define $D_t:=\sum_{i=1}^t\chi_i^2$, $N_t:=\max_{1\leq i\leq t}\chi_i^2$, where $D_0=N_0:=0$.

\noi Define $u_t:=\frac{1+D_t}{1+N_t}\geq1$, $U_t:=\max_{1\leq k\leq t}u_k$, where $u_0=U_0:=1$.

\noi Since $\left(\frac{1+N_k}{1+D_k}\right)^{p}=u_k^{-{p}}$ and $f(x) =x^{-{p}}$ is decreasing on $(0,\infty)$, we have $\rho_t=U_t^{-{p}}$.

\paragraph{Inverse-weight increment.}  On the one hand, if $\chi_t^2\leq N_{t-1}$, then $N_t=N_{t-1}$ and $u_t-u_{t-1}=\frac{\chi_t^2}{1+N_{t-1}}\leq1$. On the other hand, if $\chi_t^2>N_{t-1}$, then $N_t=\chi_t^2$ and $u_t
= \frac{1+D_{t-1}+\chi_t^2}{1+\chi_t^2} = 1+\frac{D_{t-1}}{1+\chi_t^2} \leq 1+\frac{D_{t-1}}{1+N_{t-1}} \leq u_{t-1}+1$. Thus, $u_t\leq u_{t-1}+1\leq U_{t-1}+1$. Since $U_t=\max\{U_{t-1},u_t\}$, we have
\begin{equation}\label{eq:U-increment}
0\leq U_t-U_{t-1}\leq1.
\end{equation}

\paragraph{Part (a).} Because $f(x) =x^{p}$ is increasing and concave on $[1,\infty)$, we derive:
\beq
0 ~~\leq~~U_t^{p}-U_{t-1}^{p}~~\leq~~{p}U_{t-1}^{{p}-1}(U_t-U_{t-1})
~~\leq ~~ p,\nn
\eeq
\noi where we used $U_{t-1}\geq1$, ${p}-1<0$, and
\eqref{eq:U-increment}. Using $\rho_t=U_t^{-{p}}$, we have:
\beq
\ts 0~~\leq~~\tfrac1{\rho_t}-\tfrac1{\rho_{t-1}}~~\leq~~ {p}.\nn
\eeq

\paragraph{Part (b).} Using the integral representation of $x\mapsto x^{-{p}}$
and \eqref{eq:U-increment}, we obtain
\beq
\ts \rho_{t-1}-\rho_t &=& \ts U_{t-1}^{-{p}}-U_t^{-{p}}~~= ~~
{p}\int_{U_{t-1}}^{U_t}z^{-{p}-1}\,dz  \nn\\
&\leq& {p}(U_t-U_{t-1})U_{t-1}^{-{p}-1} ~~\leq~~ \ts {p}U_{t-1}^{-{p}-1}~~=~~ {p}\rho_{t-1}^{1+1/{p}}.\nn
\eeq

\noi This proves $\rho_{t-1} - \rho_t \leq {p} \rho_{t-1}^{1 + 1/{p}}$.

\paragraph{Part (c).} Using Part (a), we have the following equivalent inequalities:
\beq
\ts \frac{\rho_{t-1}-\rho_t}
{\rho_t\rho_{t-1}}
\leq {p} ~~\Leftrightarrow~~~~
\rho_{t-1}(1-{p}\rho_t)\leq\rho_t ~~\Leftrightarrow~~
\frac{\rho_{t-1}}{\rho_t}
\leq
\frac1{1-{p}\rho_t}.\nn
\eeq
Since $0<\rho_t\leq1$ and ${p}\in(0,1)$, it follows that
\beq
\ts (1-\rho_t) \sqrt{\frac{\rho_{t-1}}{\rho_t}}~~\leq~~\frac{1-\rho_t}{\sqrt{1-{p}\rho_t}}~~\leq~~\frac{1-\rho_t}{1-{p}\rho_t}~~=~~1-(1-{p})\rho_t-\frac{{p}(1-{p})\rho_t^2}{1-{p}\rho_t}~~\leq~~1-(1-{p})\rho_t.\nn
\eeq

\end{proof}


\subsection{Analysis for {\sf ALFPG-MVR1}}
\label{subsect:result:MVR1}

\paragraph{Notation.} Let $p=\lambda=1/2$. Write $\mathbb E_t[\cdot]=\mathbb E[\cdot\mid\mathcal F_{t-1}]$.

Define $\epsilons_t:=\g_t-\nabla f(\x_t)$, $\s_t:=\m_t-\nabla f(\x_t)$, $\r_t:=\frac{\x_t-\x_{t+1}}{\gamma_t}$, $A_t:=1+\sum_{i=1}^t\|\r_i\|^2$, $D_t:=\sum_{i=1}^t\chi_i^2$, $N_t:=\max_{1\leq i\leq t}\chi_i^2$, where $A_0=1$ and $D_0=N_0=0$.

Set $\g_0:=\g_1$, consistently with $\x_0=\x_1$ and
$\xi_0=\xi_1$, so $\chi_1=0$.

Set $\alpha_t=\rho_{t-1}$, $\beta_t=A_{t-1}^{-1/2}$, $\gamma_t=\theta\sqrt{\alpha_t}\beta_t$. Both $\alpha_t$ and $\gamma_t$ are $\mathcal F_{t-1}$-measurable.

Set $\alpha_0=1$, $\gamma_0=0$, $\r_0=\zero$, and $\s_0=-\nabla f(\x_1)$.

Write $U_T:=\sum_{t=1}^T\alpha_t\|\s_t\|^2$, $W_T:=\sum_{t=1}^T\gamma_t\|\r_t\|^2$, $H_T:=\sum_{t=1}^T\frac{\gamma_t^2}{\alpha_t}\|\r_t\|^2$, $M_T:=\sum_{t=1}^T\|\r_t\|^2=A_T-1$, $S_T:=\sum_{t=1}^T\|\s_t\|^2$, $\ell_T:=\ln(e+T)$.

We begin the MVR1 analysis with pathwise bounds for the estimator, the proximal residual, and successive adaptive weights. The convex and domain-Lipschitz cases require different residual arguments, but lead to the common estimates used below.

\begin{lemma}\label{lemma:MVR1:basic:bound}
(Preliminary bounds for MVR1) Let $K_h:=\begin{cases}
\tfrac{1}{\theta}\|\x_1-\widehat{\x}_1\|,& h\text{ convex},\\
\hat H,&\text{otherwise},
\end{cases}$, $\widehat{\x}_1:=\prox_{\theta h}(\x_1)$, $\widehat{C}:=2{{\widehat G}}+ 2 K_h$, $\hat{R}:=1+\widehat{C}$, and $\kappa:=1+\hat{R}$. Almost surely, for every $t\geq1$,

(a) $\|\m_t\|\leq{{\widehat G}}$, $\|\s_t\|\leq2{{\widehat G}}$.

(b) $\|\r_t\|\leq\widehat{C} t$, $\frac{\|\r_t\|^2}{A_{t-1}}\leq \hat{R}^2$.

(c) $\frac{\alpha_t}{\alpha_{t+1}}\leq\frac32$, $\frac{\beta_t}{\beta_{t+1}}\leq\kappa$.
\end{lemma}

\begin{proof}

Define $\widehat{\x}_1:=\prox_{\theta h}(\x_1)$.

\paragraph{Bounding the term $\|\r_1\|$ when $h(\cdot)$ is convex.} For the first iteration, we have $\gamma_1=\theta$ since $\alpha_1=\beta_1=1$. By the nonexpansiveness of the proximal operator, $\|\x_2-\widehat{\x}_1\| = \|\prox_{\theta h}(\x_1-\theta\m_1) -\prox_{\theta h}(\x_1)\| \leq
\theta\|\m_1\|$. Therefore,
\beq
\|\r_1\| = \frac{\|\x_1-\x_2\|}{\theta} \leq \frac{\|\x_1-\widehat{\x}_1\|}{\theta} +\frac{\|\widehat{\x}_1-\x_2\|}{\theta} \leq K_h+\|\m_1\| \leq K_h+{{\widehat G}} . \label{eq:bound:r1}
\eeq

\paragraph{Part (a-i).} The stochastic momentum estimator is updated as $\m_t=(1-\alpha_t)\m_{t-1}+\alpha_t\g_t$, where $\alpha_t\in(0,1]$. Since the stochastic oracle is bounded, $\|\g_t\|\leq {{\widehat G}}$ almost surely. Therefore, $\m_t$ is a convex combination of previous
bounded estimators and stochastic gradients. By induction,
\beq
\|\m_t\| \leq (1-\alpha_t)\|\m_{t-1}\| +\alpha_t\|\g_t\| \leq {{\widehat G}}.\nn
\eeq

\paragraph{Part (a-ii).} By the definition $\s_t=\m_t-\nabla f(\x_t)$, and Jensen's inequality together with the bounded stochastic oracle, we have $\|\nabla f(\x_t)\|
=\|\mathbb E_t[\g_t]\|\leq \mathbb E_t\|\g_t\|
\leq {{\widehat G}}$. Hence,
\beq
\ts \|\s_t\| \leq \|\m_t\|+\|\nabla f(\x_t)\|
\leq2{{\widehat G}}.\nn
\eeq

\paragraph{Part (b-i).} We next establish an upper bound for $\|\r_t\|$ when $h$ is convex. Assume $t\geq 2$. Note that $\x_{t+1}
= \prox_{\gamma_t h}(\x_t-\gamma_t\m_t)$, and $\r_t=\frac{\x_t-\x_{t+1}}{\gamma_t}$. The first-order optimality condition of the proximal subproblem gives $\r_t-\m_t\in\partial h(\x_{t+1})$. Similarly, using the proximal update at iteration $t-1$, we have $\r_{t-1}-\m_{t-1}\in\partial h(\x_t)$.

Since $h$ is convex, its subdifferential is monotone. Therefore,
\[
\left\langle
(\r_t-\m_t)-(\r_{t-1}-\m_{t-1}),
\x_{t+1}-\x_t
\right\rangle
\geq0 .
\]
Using $\x_{t+1}-\x_t=-\gamma_t\r_t$, we obtain
\[
\left\langle
(\r_t-\m_t)-(\r_{t-1}-\m_{t-1}),
\r_t
\right\rangle
\leq0 .
\]
Expanding the inner product gives
\[{\begin{array}{rcl}
\|\r_t\|^2 ~~\leq~~ \la \r_{t-1}+\m_t-\m_{t-1}, \r_t \ra ~~\xxxx{1}{\leq} ~~\|\r_t\| \|\r_{t-1}+\m_t-\m_{t-1}\| \leq \|\r_t\| \left( \|\r_{t-1}\| + 2 {{\widehat G}} \right),
\end{array}}\]
\noi where step \step{1} uses the Cauchy--Schwarz inequality. If $\r_t=\zero$, the claimed bound is immediate; otherwise, dividing by $\|\r_t\|$ gives:
\beq \label{eq:rt:1}
\|\r_t\| \leq \|\r_{t-1}\| + 2 {{\widehat G}} .
\eeq
Applying this recursion repeatedly gives
\beq
\|\r_t\| ~~\leq~~ \|\r_1\|+2(t-1){{\widehat G}} ~~\leq~~K_h + {{\widehat G}} + 2(t-1){{\widehat G}} ~~  ~~\leq~~ \widehat C t .
\eeq
Note that the inequality above also holds for $t=1$.

We next establish an upper bound for $\|\r_t\|$ when $h(\cdot)$ is possibly nonconvex and domain-Lipschitz. Global optimality of the proximal update, with $\x_t$ as a comparison point, gives
\[
h(\x_{t+1})
+\langle\m_t,\x_{t+1}-\x_t\rangle
+\frac{\|\x_{t+1}-\x_t\|^2}{2\gamma_t}
\leq h(\x_t).
\]
Since both $\x_t$ and $\x_{t+1}$ belong to $\dom(h)$,
we obtain
\[
\begin{aligned}
\frac{\|\x_{t+1}-\x_t\|^2}{2\gamma_t}
&\leq
h(\x_t)-h(\x_{t+1})
-\langle\m_t,\x_{t+1}-\x_t\rangle\\
&\leq
(K_h+\|\m_t\|)
\|\x_{t+1}-\x_t\|.
\end{aligned}
\]
If $\x_{t+1}=\x_t$, then $\r_t=\zero$.
Otherwise, dividing by $\|\x_{t+1}-\x_t\|$ gives
\beq \label{eq:rt:2}
\|\r_t\| = \frac{\|\x_{t+1}-\x_t\|}{\gamma_t} \leq 2(K_h+\|\m_t\|) \leq 2(K_h+{{\widehat G}}) =\widehat C.
\eeq
Combining both cases gives $\|\r_t\|\leq \widehat{C} t$ for all $t\geq 1$.

\paragraph{Part (b-ii).} In view of Inequalities (\ref{eq:rt:1}) and (\ref{eq:rt:2}), we have: $\|\r_t\|\leq\|\r_{t-1}\|+\widehat C$, for all $t\geq1$. Combining with $A_{t-1} \geq 1+\|\r_{t-1}\|^2$, we have:
\beq
\frac{\|\r_t\|^2}{A_{t-1}} ~~\leq~~ \frac{(\|\r_{t-1}\|+\widehat C)^2}{\|\r_{t-1}\|^2+1}~~\leq ~~ 1 + \widehat{C}^2 ~~ \leq~~\hat R^2.\nn
\eeq

\paragraph{Part (c-i).} Lemma~\ref{lemma:weight}(a) with $p=1/2$ gives $\frac1{\rho_t}-\frac1{\rho_{t-1}}\leq\frac12$. Therefore,
\beq
\ts \frac{\rho_{t-1}}{\rho_t} \leq 1+\frac{\rho_{t-1}}2 \leq \frac32. \nn
\eeq

\paragraph{Part (c-ii).} By definition, $\beta_t = A_{t-1}^{-1/2}$. Thus,
\beq
\ts \frac{\beta_t}{\beta_{t+1}} ~~=~~ \sqrt{\frac{A_t}{A_{t-1}}} ~~=~~ \sqrt{1+\frac{\|\r_t\|^2}{A_{t-1}}} ~~\leq~~ \sqrt{1+\hat R^2} ~~\leq \kappa.\nn
\eeq
\end{proof}

The preliminary bounds control the size of the quantities entering the EMA recursion. We now separate its tracking error into contraction, iterate drift, and fresh sampling noise, retaining the conditionally mean-zero cross term for the later moment estimates.

\begin{lemma}\label{lemma:MVR1:s:recursion}
(Tracking-error recursion for MVR1) Define $\mathcal E_t = 2\alpha_t
\langle\b_t,\boldsymbol\epsilon_t\rangle$, $\dot{\mathcal E}_t=\alpha_t^2 \|\boldsymbol\epsilon_t\|^2$, and $\ddot{\mathcal E}_t = \frac{2 L^2\gamma_{t-1}^2}{\alpha_{t-1}}
\|\r_{t-1}\|^2$, where $\b_t=(1-\alpha_t)
(\m_{t-1}-\nabla f(\x_t))$.

(a) Then almost surely, $\|\s_t\|^2 \leq (1-\alpha_t)\|\s_{t-1}\|^2 + \ddot{\mathcal E}_t + \dot{\mathcal E}_t + \mathcal E_t$.

(b) $\mathbb E_t[\mathcal E_t]=0$, and $\mathcal E_t^2 \leq 16{{\widehat G}}^2 \dot{\mathcal E}_t$.

\end{lemma}

\begin{proof}

Define $\epsilons_t:=\g_t-\nabla f(\x_t)$, $\s_t=\m_t-\nabla f(\x_t)$.

Let $\b_t := (1-\alpha_t) (\m_{t-1}-\nabla f(\x_t))$, and $\deltas_t:=(1-\alpha_t)(\nabla f(\x_{t-1})-\nabla f(\x_t))$.

\paragraph{Exact decomposition of the tracking error.} Recall that $\m_t=(1-\alpha_t)\m_{t-1}+\alpha_t\g_t$. We derive:
\beq \label{eq:MVR1:tracking:decomposition}
\s_t &=& (1-\alpha_t) (\m_{t-1}-\nabla f(\x_t)) +\alpha_t\boldsymbol\epsilon_t \nn\\
&=& (1-\alpha_t)\s_{t-1} + (1-\alpha_t)(\nabla f(\x_{t-1})-\nabla f(\x_t))
+\alpha_t \epsilons_t  ~~=~~ \b_t+\alpha_t\boldsymbol\epsilon_t.
\eeq

\paragraph{Bounding the deterministic drift $\|\deltas_t\|$.} By the $L$-smoothness of $f$,
\beq \label{eq:delta:drift}
\|\deltas_t\| ~~\leq~~(1-\alpha_t)L\|\x_t-\x_{t-1}\|~~\xxxx{1}{\leq}~~(1-\alpha_t)L\gamma_{t-1} \|\r_{t-1}\|~~{\leq}~~L\gamma_{t-1} \|\r_{t-1}\|.
\eeq
\noi where step \step{1} uses $\x_t-\x_{t-1} = -\gamma_{t-1}\r_{t-1}$.

\paragraph{Part (a).} Since $\b_t = (1-\alpha_t)\s_{t-1} + \alpha_t ( \tfrac{1}{\alpha_t}\deltas_t)$, due to the convexity of $\|\cdot\|_2^2$, we have:
\beq
\|\b_t\|^2 &\leq & (1-\alpha_t)\|\s_{t-1}\|^2 + \alpha_t \|( \tfrac{1}{\alpha_t} \deltas_t)\|^2~~\xxxx{1}{\leq} ~~ (1-\alpha_t)\|\s_{t-1}\|^2 + \tfrac{1}{\alpha_t} L^2 \gamma_{t-1}^2 \|\r_{t-1}\|^2\nn\\
&\xxxx{2}{\leq} &(1-\alpha_t)\|\s_{t-1}\|^2 + \tfrac{2L^2 \gamma_{t-1}^2}{\alpha_{t-1}}  \|\r_{t-1}\|^2, \label{eq:com1}
\eeq
\noi where step \step{1} uses Inequality (\ref{eq:delta:drift}); step \step{2} uses $\frac1{\alpha_t} \leq \tfrac{2}{\alpha_{t-1}}$, which can be implied by Lemma~\ref{lemma:MVR1:basic:bound} that $\frac{\alpha_{t-1}}{\alpha_t}\leq
\frac32<2$.

On the other hand, we have:
\beq
\|\s_t\|^2 = \|\b_t + \alpha_t \epsilons_t\|^2 = \|\b_t\|^2 + \dot{\mathcal{E}}_t + \mathcal{E}_t.\label{eq:com2}
\eeq

Combining Inequality (\ref{eq:com1}) and Inequality (\ref{eq:com2}) yields:
\beq
\ts \|\s_t\|^2 \leq (1-\alpha_t) \|\s_{t-1}\|^2 + \tfrac{2L^2\gamma_{t-1}^2}{\alpha_{t-1}}\|\r_{t-1}\|^2 + \dot{\mathcal E}_t+\mathcal E_t.\nn
\eeq

\paragraph{Part (b).} Since $\b_t$ is $\mathcal F_{t-1}$-measurable and $\mathbb E_t[\boldsymbol\epsilon_t]=0$, we have $\mathbb E_t[\mathcal E_t]
= 2\alpha_t \langle\b_t, \mathbb E_t[\boldsymbol\epsilon_t]\rangle =0$. Finally, $\mathcal E_t^2 = 4\alpha_t^2
\langle\b_t,\boldsymbol\epsilon_t\rangle^2$. By Cauchy--Schwarz,
\beq
\mathcal E_t^2 ~~\leq~~ 4\alpha_t^2 \|\b_t\|^2 \|\boldsymbol\epsilon_t\|^2~~\xxxx{1}{\leq}~~ 4\alpha_t^2 \cdot 4 {{\widehat G}}^2 \|\boldsymbol\epsilon_t\|^2~~=~~ 16{{\widehat G}}^2 \dot{\mathcal E}_t.\nn
\eeq
\noi where step \step{1} uses $\|\b_t\| \leq (1-\alpha_t)
(\|\m_{t-1}\|+\|\nabla f(\x_t)\|) \leq 2{{\widehat G}}$.

\end{proof}

To connect tracking accuracy with optimization progress, we combine proximal descent with the preceding error recursion. The logarithmic correction absorbs the normalized motion term and yields a joint Lyapunov inequality.

\begin{lemma} \label{lemma:MVR1:descent}
(Proximal descent and a joint Lyapunov inequality for MVR1) Let $C_f:=\theta^2\kappa^2(2 L+1)^2$ and
$Q_t:=F(\x_t)+2\|\s_t\|^2-C_f\ln A_{t-1}$.
Then, almost surely,
\begin{align}
F(\x_{t+1})-F(\x_t)+\frac{\gamma_t}{2}\|\r_t\|^2
&\leq\frac{\alpha_t}{2}\|\s_t\|^2
+\frac{1+L}{2}\frac{\gamma_t^2}{\alpha_t}\|\r_t\|^2,
\label{eq:A-rev:descent}\\
Q_{t+1}-Q_t+\frac{\gamma_t}{2}\|\r_t\|^2
+\frac{\alpha_t}{2}\|\s_t\|^2
&\leq2\dot{\mathcal E}_{t+1}+2\mathcal E_{t+1}.
\label{eq:one-step-descent}
\end{align}

\end{lemma}

\begin{proof}

Define $\s_t:=\m_t-\nabla f(\x_t)$ and
$\r_t:=(\x_t-\x_{t+1})/\gamma_t$.

Let $A_t:=1+\sum_{i=1}^t\|\r_i\|^2$, where $A_0:=1$.

Let $\x_{t+1} \in\arg\min_{\y} \left\{h(\y)+ \la \m_{t}, \y-\x_t\ra + \frac{1}{2\gamma_t}\|\y-\x_t\|^2\right\}$.

By the optimality of $\x_{t+1}$, we have:
\beq
h(\x_{t+1}) \le h(\x_t)-\la \m_{t}, \x_{t+1}-\x_t\ra-\frac{1}{2\gamma_t}\|\x_{t+1}-\x_t\|^2.\nn
\eeq
\noi By $L$-smoothness of $f$,
\beq
f(\x_{t+1})\le f(\x_t)+\la \x_{t+1}-\x_t, \nabla f(\x_t)  \ra + \frac{L}{2}\|\x_{t+1}-\x_t\|^2.\nn
\eeq

\paragraph{Part (a).} Adding the two inequalities gives
\beq
 \ts \ts F(\x_{t+1}) - F(\x_t) &\le& \ts \la \x_{t+1}-\x_t,\nabla f(\x_t) - \m_{t} \ra - \frac{1}{2\gamma_t} \|\x_{t+1}-\x_t\|^2 + \frac{L}{2}\|\x_{t+1}-\x_t\|^2 \nn \\
&=& \ts \la \gamma_t \r_t , \s_t \ra - \frac{\gamma_t}{2} \|\r_t\|^2 + \frac{L\gamma_t^2}{2}\|\r_t\|^2 \nn \\
&\xxxx{1}{\le}& \ts \frac{\alpha_t}{2}\|\s_t\|^2 + \frac{\gamma_t^2}{2\alpha_t}\|\r_t\|^2 - \frac{\gamma_t}{2} \|\r_t\|^2 + \frac{L\gamma_t^2}{2}\|\r_t\|^2  \nn\\
&\xxxx{}{\le}& \ts \frac{\alpha_t}{2}\|\s_t\|^2  - \frac{\gamma_t}{2} \|\r_t\|^2 + \frac{ (L+1)\gamma_t^2}{2\alpha_t}\|\r_t\|^2,  \nn
\eeq
\noi where step \step{1} uses weighted Young inequality $\gamma_t\langle\r_t,\s_t\rangle
\leq \frac{\alpha_t}{2}\|\s_t\|^2 + \frac{\gamma_t^2}{2\alpha_t}\|\r_t\|^2$.

\paragraph{Part (b).} We further derive:
\beq
&& \ts F(\x_{t+1}) - F(\x_t) + \tfrac{\gamma_t}{2} \|\r_t\|^2 + \tfrac{\alpha_t}{2} \|\s_t\|^2\nn\\
& \leq & \ts \alpha_t \|\s_t\|^2   + \frac{ (L+1)\gamma_t^2}{2\alpha_t}\|\r_t\|^2\nn\\
&\xxxx{1}{\leq} & \ts 2 \alpha_{t+1} \|\s_t\|^2   + \frac{ (L+1)\gamma_t^2}{2\alpha_t}\|\r_t\|^2\nn\\
&\xxxx{2}{\leq} & \ts 2 \left( \|\s_{t}\|^2 - \|\s_{t+1}\|^2 + \tfrac{2 L^2 \gamma_t^2}{\alpha_{t}} \|\r_t\|^2 + \dot{\mathcal{E}}_{t+1} + \mathcal{E}_{t+1} \right)   + \frac{ (L+1)\gamma_t^2}{2\alpha_t}\|\r_t\|^2\nn\\
&\xxxx{ }{=} & \ts 2 \left( \|\s_{t}\|^2 - \|\s_{t+1}\|^2 + \dot{\mathcal{E}}_{t+1} + \mathcal{E}_{t+1} \right)   +  \left( \frac{ (L+1)}{2} + 4 L^2 \right)\cdot \left( \tfrac{\gamma_t^2}{\alpha_t} \|\r_t\|^2 \right) ,  \label{eq:dec:222}
\eeq
\noi where step \step{1} uses $\tfrac{\alpha_t}{\alpha_{t+1}}=\tfrac{\rho_{t-1}}{\rho_t}\leq 1 + p = \tfrac{3}{2}<2$; step \step{2} uses Lemma \ref{lemma:MVR1:s:recursion}(a) that $\alpha_t \|\s_{t-1}\|^2  \leq  \|\s_{t-1}\|^2 - \|\s_t\|^2 + \ddot{\mathcal{E}}_t + \dot{\mathcal{E}}_t + \mathcal{E}_t$.

Furthermore, we have:
\beq
&& \big( \tfrac{ (L+1)}{2} + 4 L^2 \big)\cdot \big( \tfrac{\gamma_t^2}{\alpha_t} \|\r_t\|^2 \big) ~~\xxxx{1}{\leq} ~~ \big( 2L+1 \big)^2\cdot \big( \theta^2 \beta_t^2  \|\r_t\|^2 \big) \nn\\
&\xxxx{2}{\leq} & \big( 2L+1 \big)^2 \cdot \big( \theta^2 \kappa^2 \beta_{t+1}^2  \|\r_t\|^2 \big) ~~=~~\big( 2L+1 \big)^2 \cdot \big( \theta^2 \kappa^2   \big) \cdot \tfrac{\|\r_t\|^2}{  1 + \sum_{i=1}^t \|\r_i\|^2} \nn\\
&\xxxx{3}{\leq} & \big( 2L+1 \big)^2 \cdot \big( \theta^2 \kappa^2   \big) \cdot \big( \ln(A_{t}) - \ln (A_{t-1})\big) = C_f \big( \ln(A_{t}) - \ln (A_{t-1})\big), \label{eq:tail:bound}
\eeq
\noi where step \step{1} uses $\gamma_t=\theta \sqrt{\alpha_t} \beta_t$; step \step{2} uses $\tfrac{\beta_t}{\beta_{t+1}}\leq \kappa$; step \step{3} uses Lemma \ref{eq:log:diff}(a).

Define $Q_t:=F(\x_t)+2\|\s_t\|^2-C_f\ln A_{t-1}$. In view of Inequalities (\ref{eq:dec:222}) and (\ref{eq:tail:bound}), we have:
\beq
Q_{t+1} - Q_t + \tfrac{\gamma_t}{2} \|\r_t\|^2 + \tfrac{\alpha_t}{2} \|\s_t\|^2 \leq 2\dot{\mathcal{E}}_{t+1} + 2\mathcal{E}_{t+1}.\nn
\eeq

\end{proof}

The one-step estimates must now be converted into cumulative bounds. We first control the discrepancy and normalized motion sums, and then derive the first two moments of the weighted tracking error; the second-moment bound will be needed for the motion estimate.

\begin{lemma}\label{lemma:MVR1:logTerm}

(Logarithmic bounds for MVR1) For every $T\geq1$, we have:

(a) $V_T:=\sum_{t=1}^T\alpha_t^2\chi_t^2 \leq C_v\ell_T$ almost surely, where $C_v= (1+4{{\widehat G}}^2)^2 \ln (e+4{{\widehat G}}^2 )$.

(b) $H_T:=\sum_{t=1}^T \tfrac{\gamma_t^2}{\alpha_t} \|\r_t\|^2\leq C_h\ell_T$ almost surely, where $C_h=3 \kappa^2 \theta^2 \ln (e+\widehat{C}^2)$.

(c) $\mathbb E[\sum_{t=1}^{T}\dot{\mathcal E}_t]=\mathbb E[\sum_{t=1}^{T}\alpha_t^2\|\boldsymbol\epsilon_t\|^2]\leq C_e\ell_T$, where $C_e:={{\widehat G}}^2+C_v$.

(d) Let $U_T:=\sum_{t=1}^T\alpha_t\|\s_t\|^2$. We have: $\mathbb E[U_T] \leq C_u \ell_T$, where $C_u=4 L^2C_h + 4 C_v$.

(e) $\mathbb E[U_T^2]\leq C'_u\ell_T^2
\leq C_u^2\ell_T^2$, where $C'_u:=\bigl(12{{\widehat G}}^2+4L^2C_h+2C_v\bigr)C_u$. Here, the simpler choice $C'_u=C_u^2$ is also valid.

\end{lemma}

\begin{proof}

Define $D_t=\sum_{i=1}^t\chi_i^2$ and $N_t=\max_{1\leq i\leq t}\chi_i^2$, where $\chi_t=\|\g_{t-1}-\g_t\|$.

Define $V_T=\sum_{t=1}^T\alpha_t^2\chi_t^2$. Let $A_t:=1+\sum_{i=1}^t\|\r_i\|^2$. Let $c_\chi:=(1+4{{\widehat G}}^2)^2$.

Define $C_v= c_\chi \ln (e+4{{\widehat G}}^2 )$. Denote $\ell_T = \ln(e+T)$.

\paragraph{Bounding the term $\sum_{t=2}^{T+1}\alpha_t^2\chi_t^2$.} We derive:
\beq
&& \ts \sum_{t=2}^{T+1}\alpha_t^2\chi_t^2 ~~\xxxx{1}{\le}~~ \ts (1+4{{\widehat G}}^2)\sum_{t=2}^{T+1} \frac{\chi_t^2}{1+D_{t-1}}~~\xxxx{2}{\le}~~\ts (1+4{{\widehat G}}^2)^2 \sum_{t=2}^{T+1}\frac{\chi_t^2}{1+D_t} \nn \\
&\xxxx{3}{\le}&\ts c_\chi\sum_{t=2}^{T+1}\left(\ln(1+D_t)-\ln(1+D_{t-1})\right) ~~\xxxx{}{\le}~~\ts c_\chi  \ln (1 + D_{T+1}) ~~\xxxx{4}{\leq}~~\ts c_\chi\ln (1 + T 4 {{\widehat G}}^2 ), \nn
\eeq
\noi where step \step{1} uses $\alpha_t^2 \le \frac{1+N_{t-1}}{1+D_{t-1}} \le \frac{1+4{{\widehat G}}^2}{1+D_{t-1}}$; step \step{2} uses $\frac{1+D_t}{1+D_{t-1}}=1+\frac{\chi_t^2}{1+D_{t-1}} \le 1+4{{\widehat G}}^2$, which follows from $\chi_t^2\le 4{{\widehat G}}^2$; step \step{3} uses Lemma \ref{eq:log:diff}(a); step \step{4} uses $D_{T+1} = \sum_{t=2}^{T+1}\chi_t^2 \le 4{{\widehat G}}^2 T$.

\paragraph{Part (a).} We derive:
\beq
\ts \sum_{t=1}^T \alpha_t^2 \chi_t^2 &\leq & \ts \sum_{t=1}^T \tfrac{1+N_{t-1}}{1+D_{t-1}} \chi_t^2 ~~ \xxxx{1}{\leq} ~~  \sum_{t=1}^T \tfrac{1+4 {{\widehat G}}^2 }{1+D_{t}} \cdot \tfrac{1+D_{t}}{1+D_{t-1}} \cdot  \chi_t^2 \ts \nn \\
&\xxxx{2}{\leq} &\ts (1+4 {{\widehat G}}^2)^2 \sum_{t=1}^T \tfrac{\chi_t^2}{1+D_{t}}  ~~ \xxxx{3}{\leq} ~~  c_{\chi} \sum_{t=1}^T  ( \ln (1+ D_t) - \ln (1+ D_{t-1})) \ts \nn \\
& \xxxx{4}{=} & \ts c_{\chi} \ln (1+ D_T) ~~ \leq ~~  c_{\chi} \ln (1+ T 4 {{\widehat G}}^2 ) ~~ \xxxx{5}{\leq} ~~  c_{\chi} \ln (e+4{{\widehat G}}^2 ) \ln(e+T),\nn
\eeq
\noi where step \step{1} uses $N_{t-1} \leq 4 {{\widehat G}}^2$; step \step{2} uses $\tfrac{1+D_t}{1+D_{t-1}} = 1+\frac{\chi_t^2}{1+D_{t-1}} \leq 1+4{{\widehat G}}^2$; step \step{3} uses Lemma \ref{eq:log:diff}(a); step \step{4} uses $D_0=0$; step \step{5} uses $\ln(1+ab)\leq\ln(e+a)\ln(e+b)$ for all $a,b\geq 0$.

\paragraph{Part (b).} We derive:
\beq
\ts H_T &=& \ts \sum_{t=1}^T \frac{\gamma_t^2}{\alpha_t} \|\r_t\|^2~~\xxxx{1}{=}~~\theta^2\sum_{t=1}^T \beta_t^2 \|\r_t\|^2 ~~\xxxx{2}{\leq}~~\kappa^2 \theta^2\sum_{t=1}^T \beta_{t+1}^2 \|\r_t\|^2 \nn\\
& \xxxx{3}{=} & \ts \kappa^2 \theta^2\sum_{t=1}^T \tfrac{\|\r_t\|^2}{1 + \sum_{i=1}^t \|\r_i\|^2}~~\xxxx{4}{\leq}~~\ts \kappa^2 \theta^2 \sum_{t=1}^T (\ln (A_t) - \ln (A_{t-1})) \nn\\
& \xxxx{4}{ \leq } & \ts  \kappa^2 \theta^2 \ln (A_T)~~\xxxx{5}{ \leq } ~~ \kappa^2 \theta^2 \ln ( 1 + \widehat{C}^2 \sum_{t=1}^T t^2 ) ~~\xxxx{6}{ \leq } ~~ \kappa^2 \theta^2 \ln ( 1 + \widehat{C}^2 T^3 ) \nn\\
&\xxxx{}{ \leq }& ~~ 3 \kappa^2 \theta^2 \ln (e+\widehat{C}^2) \ln (e+T)~~=~~C_h\ell_T,\nn
\eeq
\noi where step \step{1} uses $\gamma_t=\theta\sqrt{\alpha_t}\beta_t$; step \step{2} uses Lemma~\ref{lemma:MVR1:basic:bound} (c) that $\frac{\beta_t}{\beta_{t+1}}\leq\kappa$; step \step{3} uses $\beta_t =A_{t-1}^{-1/2}$, where $A_{t-1} = 1+\sum_{i=1}^{t-1}\|\r_i\|^2$; step \step{4} uses Lemma \ref{eq:log:diff}(a); step \step{5} uses Lemma~\ref{lemma:MVR1:basic:bound}(b) that $\|\r_t\|\leq\widehat C t$, leading to $A_T=1+\sum_{t=1}^T\|\r_t\|^2\leq1+\widehat C^{\,2}\sum_{t=1}^Tt^2$; step \step{6} uses $\sum_{t=1}^Tt^2=\frac{T(T+1)(2T+1)}6
\leq T^3$.

\paragraph{Part (c).} The first sample needs to be handled separately, because $\g_0=\g_1$ implies $\chi_1=0$. Since $\alpha_1=1$, unbiasedness at the fixed initial point gives
\[
\mathbb E[\dot{\mathcal E}_1]~~=~~\mathbb E\|\epsilons_1\|^2~~
=~~\mathbb E\|\g_1\|^2-\|\nabla f(\x_1)\|^2~~\leq~~{{\widehat G}}^2.
\]
For $t\geq2$, $\nabla f(\x_t)-\g_{t-1}$ and $\alpha_t$ are $\mathcal F_{t-1}$-measurable. Hence
\beq
\mathbb E_t[\chi_t^2]~~=~~\|\nabla f(\x_t)-\g_{t-1}\|^2+\mathbb E_t\|\epsilons_t\|^2
~~\geq~~\mathbb E_t\|\epsilons_t\|^2.\nn
\eeq
Using the tower property before summing yields
\beq
\ts \mathbb E[\sum_{t=1}^T\dot{\mathcal E}_t] ~~\leq~~{{\widehat G}}^2+\mathbb E\sum_{t=2}^T[\alpha_t^2\chi_t^2]~~\leq~~{{\widehat G}}^2+C_v\ell_T\leq({{\widehat G}}^2+C_v)\ell_T~~=~~C_e\ell_T.\nn
\eeq

\paragraph{Part (d).} Applying Lemma~\ref{lemma:MVR1:s:recursion}(a)
at iteration $t+1$ gives
\beq
\ts \|\s_{t+1}\|^2 &\leq & \ts (1-\alpha_{t+1})\|\s_t\|^2 +\ddot{\mathcal E}_{t+1} +\dot{\mathcal E}_{t+1} +\mathcal E_{t+1} \nn\\
&\xxxx{1}{\leq} & \ts (1- \tfrac{1}{2}\alpha_{t})\|\s_t\|^2 +\ddot{\mathcal E}_{t+1} +\dot{\mathcal E}_{t+1} +\mathcal E_{t+1},\nn
\eeq
\noi where step \step{1} uses $\tfrac{\alpha_t}{\alpha_{t+1}}\leq \tfrac{3}{2}<2$. Consequently,
\beq
\alpha_t\|\s_t\|^2 &\le& 2\left(
\|\s_t\|^2-\|\s_{t+1}\|^2
+\ddot{\mathcal E}_{t+1}
+\dot{\mathcal E}_{t+1}
+\mathcal E_{t+1}
\right).\nn
\eeq
Summing over $t=1,\ldots,T$ and taking the expectation of both sides, we obtain
\beq
\ts \EEE[U_T] &\le & \ts \EEE\left[ 2\left( \|\s_1\|^2-\|\s_{T+1}\|^2
+\sum_{t=2}^{T+1}\ddot{\mathcal E}_t +\sum_{t=2}^{T+1}\dot{\mathcal E}_t
+\sum_{t=2}^{T+1}\mathcal E_t \right)\right] \nn\\
&\xxxx{1}{\le} & \ts  2\left( {{\widehat G}}^2 + \EEE[\sum_{t=2}^{T+1}\ddot{\mathcal E}_t +\sum_{t=2}^{T+1}\dot{\mathcal E}_t
 ] \right)  \nn\\
&\xxxx{2}{\le} & \ts 2\left( {{\widehat G}}^2 + 2L^2C_h\ell_T + \EEE[  \sum_{t=2}^{T+1}\alpha_t^2\chi_t^2 ] \right)  \nn\\
&\xxxx{3}{\le} & \ts  2\left( {{\widehat G}}^2 + 2L^2C_h\ell_T + c_\chi\ln(e+4 {{\widehat G}}^2 ) \ell_T \right)  \nn\\
&\xxxx{}{\le} & \ts  2\left( {{\widehat G}}^2 + 2L^2C_h + c_\chi\ln(e+4 {{\widehat G}}^2 ) \right) \ell_T \nn\\
&\xxxx{}{\le} & \ts  \left( 4 L^2C_h + 4 c_\chi\ln(e +4 {{\widehat G}}^2 ) \right) \ell_T = (4 L^2C_h + 4 C_v) \ell_T, \nn
\eeq
\noi where step \step{1} $\mathbb E\|\s_1\|^2 = \mathbb E\|\m_1-\nabla f(\x_1)\|^2= \mathbb E\|\g_1-\nabla f(\x_1)\|^2
= \mathbb E\|\g_1\|^2-\|\nabla f(\x_1)\|^2 \le {{\widehat G}}^2$ (which  follows from $\m_1=\g_1$), and $\mathbb E\big[\sum_{t=2}^{T+1}\mathcal E_t\big] = \sum_{t=2}^{T+1}\mathbb E\big[\mathbb E_t[\mathcal E_t]\big]=0$ (which follows from Lemma~\ref{lemma:MVR1:s:recursion}(b) that $\mathbb E_t[\mathcal E_t]=0$ and the tower property); step \step{2} uses $\ddot{\mathcal E}_t = \frac{2L^2\gamma_{t-1}^2}{\alpha_{t-1}} \|\r_{t-1}\|^2$, leading to $\sum_{t=2}^{T+1}\ddot{\mathcal E}_t =
2L^2 \sum_{t=1}^{T} \frac{\gamma_t^2}{\alpha_t}\|\r_t\|^2 = 2L^2H_T \le 2L^2C_h\ell_T$; step \step{3} uses $\EEE[  \sum_{t=2}^{T+1}\alpha_t^2\chi_t^2 ] \leq c_\chi\ln (1 + T 4 {{\widehat G}}^2 ) \leq c_\chi\ln(e+4 {{\widehat G}}^2 ) \ell_T$.

\paragraph{Part (e).}
Let $q_t:=\|\s_t\|^2$, $w_t:=\alpha_tq_t$, $U_0:=0$,
and
$\eta_t:=\ddot{\mathcal E}_t+\dot{\mathcal E}_t+\mathcal E_t$.
By Lemma~\ref{lemma:MVR1:basic:bound}(a),
$0\leq q_t,w_t\leq4{{\widehat G}}^2$.
For $t\geq2$, conditional unbiasedness gives
\[
\mathbb E_t[\chi_t^2]
=\|\nabla f(\x_t)-\g_{t-1}\|^2
+\mathbb E_t\|\boldsymbol\epsilon_t\|^2
\geq\mathbb E_t\|\boldsymbol\epsilon_t\|^2.
\]
Since $U_{t-1}\alpha_t^2$ is
$\mathcal F_{t-1}$-measurable, $U_0=0$, and
$\mathbb E_t[\mathcal E_t]=0$, parts (a)--(b) imply
\beq \label{eq:MVR1:U2:weighted-errors}
\ts \mathbb E\left[\sum_{t=1}^T U_{t-1}\eta_t\right]~~\leq~~
\mathbb E\left[
U_T\left(2L^2H_T+V_T\right)\right] ~~\leq~~
(2L^2C_h+C_v)\ell_T\,\mathbb E[U_T].
\eeq
Here the initial noise has zero weight, and
$\sum_{t=1}^T\ddot{\mathcal E}_t
=2L^2H_{T-1}\leq2L^2H_T$, with $H_0:=0$.

The tracking recursion yields
\beq
\ts w_t\leq(1-\alpha_t)(q_{t-1}-q_t)+\eta_t.
\eeq
Set $P_t:=U_{t-1}(1-\alpha_t)$.
Since $U_t$ is nondecreasing and $\alpha_t$ is
nonincreasing, $P_t$ is nondecreasing and $P_1=0$.
Summation by parts therefore gives
\beq
\ts \sum_{t=1}^T P_t(q_{t-1}-q_t)
=-P_Tq_T+\sum_{t=1}^{T-1}(P_{t+1}-P_t)q_t
\leq4{{\widehat G}}^2U_T.
\eeq
Multiplying the recursion by $U_{t-1}$, summing,
and using \eqref{eq:MVR1:U2:weighted-errors}, we obtain
\beq
\ts \mathbb E \left[\sum_{t=1}^T U_{t-1}w_t\right]
\leq
\bigl[4{{\widehat G}}^2+(2L^2C_h+C_v)\ell_T\bigr]
\mathbb E[U_T].
\eeq

Finally,
$U_T^2=\sum_{t=1}^T w_t^2
+2\sum_{t=1}^T U_{t-1}w_t$
and $w_t^2\leq4{{\widehat G}}^2w_t$.
Using part (d) and $\ell_T\geq1$, we conclude that
\begin{align*}
\mathbb E[U_T^2]
&\leq
\bigl[12{{\widehat G}}^2+(4L^2C_h+2C_v)\ell_T\bigr]
\mathbb E[U_T]\\
&\leq
(12{{\widehat G}}^2+4L^2C_h+2C_v)C_u\ell_T^2
=C'_u\ell_T^2.
\end{align*}
Moreover,
$2C_v-12{{\widehat G}}^2
\geq2+4{{\widehat G}}^2+32{{\widehat G}}^4\geq0$.
Thus
$12{{\widehat G}}^2+4L^2C_h+2C_v\leq C_u$,
which proves $C'_u\leq C_u^2$.

\end{proof}

We next use proximal descent to transfer the tracking-error bounds to the descent energy. A pathwise comparison with the residual accumulator then controls the weighted cumulative motion, which is the remaining ingredient for the stationarity proof.

\begin{lemma}
\label{lemma:MVR1:bridge} (Descent energy and motion bridge for MVR1) Let $\Delta:=F(\x_1)-F_\star$, $\ell_T:=\ln(e+T)$. Let $C_h$ and $C_u$ be defined in
Lemma~\ref{lemma:MVR1:logTerm}. Then, for every integer $T\geq1$,
\begin{enumerate}[
label=\textbf{(\alph*)},
leftmargin=18pt,
itemsep=2pt,
topsep=1pt]
\item
$\mathbb E[W_T]\leq C_w\ell_T$ and
$\mathbb E[W_T^2]\leq C_w^2\ell_T^2$, where $C_w:=2\Delta+C_u+(1+L)C_h$.
\item
$\mathbb E[\alpha_T M_T]\leq C_m\ell_T^2$, where $C_m:=\frac{C_w}{\theta}
+\frac{C_w^2}{4\theta^2}$.
\end{enumerate}
\end{lemma}

\begin{proof}

Let $C_w:=2\Delta+C_u+(1+L)C_h$.

\paragraph{Part (a-i).} Summing the first inequality of
Lemma~\ref{lemma:MVR1:descent} over $t=1,\ldots,T$ gives
\[
F(\x_{T+1})-F(\x_1)+\frac12W_T
\leq
\frac12U_T+\frac{1+L}{2}H_T.
\]
Since $F(\x_{T+1})\geq F_\star$, we obtain
\begin{equation}
\label{eq:MVR1:bridge:pathwise-descent}
0\leq W_T
\leq 2\Delta+U_T+(1+L)H_T.
\end{equation}
Taking expectations and applying
Lemma~\ref{lemma:MVR1:logTerm}(b),(d), we have
\beq
\ts \mathbb E[W_T]~~\leq~~ 2\Delta+C_u\ell_T+(1+L)C_h\ell_T
~~\leq~~ \bigl(2\Delta+C_u+(1+L)C_h\bigr)\ell_T
~~=~~C_w\ell_T,
\eeq
where the second inequality uses $\ell_T\geq1$.

\paragraph{Part (a-ii).} Set the deterministic quantity
\[
a_T:=2\Delta+(1+L)C_h\ell_T\geq0.
\]
Since $H_T\leq C_h\ell_T$ almost surely,
\eqref{eq:MVR1:bridge:pathwise-descent} yields
$0\leq W_T\leq a_T+U_T$.
Using Lemma~\ref{lemma:MVR1:logTerm}(d),(e),
we derive
\beq
\mathbb E[W_T^2] &\leq& \ts a_T^2+2a_T\mathbb E[U_T]+\mathbb E[U_T^2]~~\leq~~ a_T^2+2a_TC_u\ell_T+C_u^2\ell_T^2 \nn \\
&=& (a_T+C_u\ell_T)^2~~=~~\bigl[ 2\Delta+\bigl(C_u+(1+L)C_h\bigr)\ell_T\bigr]^2
~~\leq~~C_w^2\ell_T^2,\nn
\eeq
where the last inequality again uses $\ell_T\geq1$.
This proves both estimates in part (a).

\paragraph{Part (b).} Let $A_t:=1+\sum_{i=1}^t\|\r_i\|^2$, with $A_0:=1$.

The running-minimum definition of $\rho_t$ ensures
that $\alpha_t=\rho_{t-1}$ is nonincreasing.
Combining this property with
$\gamma_t=\theta\sqrt{\alpha_t}/\sqrt{A_{t-1}}$,
we obtain
\beq \label{eq:MVR1:bridge:root-motion}
\ts W_T &=& \ts \theta\sum_{t=1}^T
\sqrt{\alpha_t}\,
\frac{A_t-A_{t-1}}{\sqrt{A_{t-1}}}
~~\geq~~\theta\sqrt{\alpha_T}\sum_{t=1}^T
\frac{A_t-A_{t-1}}{\sqrt{A_{t-1}}}\nn\\
&\geq& \ts 2\theta\sqrt{\alpha_T}
\sum_{t=1}^T\bigl(\sqrt{A_t}-\sqrt{A_{t-1}}\bigr)
~~=~~ \ts 2\theta\sqrt{\alpha_T}\bigl(\sqrt{1+M_T}-1\bigr).
\eeq
Here the second inequality follows from
\beq
\ts \frac{A_t-A_{t-1}}{\sqrt{A_{t-1}}}~~=~~
\bigl(\sqrt{A_t}-\sqrt{A_{t-1}}\bigr)
\left(1+\sqrt{\frac{A_t}{A_{t-1}}}\right)
~~\geq~~2\bigl(\sqrt{A_t}-\sqrt{A_{t-1}}\bigr),
\eeq
because $A_t\geq A_{t-1}$. Rearranging~\eqref{eq:MVR1:bridge:root-motion} gives
\[
\sqrt{\alpha_T(1+M_T)}
\leq
\sqrt{\alpha_T}+\frac{W_T}{2\theta}.
\]
Both sides are nonnegative.
Squaring and subtracting $\alpha_T$ therefore yields
\begin{equation}
\label{eq:MVR1:bridge:pathwise-motion}
\alpha_T M_T
\leq
\frac{\sqrt{\alpha_T}}{\theta}W_T
+\frac{W_T^2}{4\theta^2}
\leq
\frac{W_T}{\theta}
+\frac{W_T^2}{4\theta^2},
\end{equation}
where the last inequality uses $0<\alpha_T\leq1$.
Taking expectations and applying part (a),
we conclude that
\beq
\ts \mathbb E[\alpha_T M_T] &\leq&
\ts \frac{1}{\theta}\mathbb E[W_T] +\frac{1}{4\theta^2}\mathbb E[W_T^2]\nn\\
&\leq& \ts \frac{C_w}{\theta}\ell_T + \frac{C_w^2}{4\theta^2}\ell_T^2
~~\leq~~ \left(
\frac{C_w}{\theta}
+\frac{C_w^2}{4\theta^2}
\right)\ell_T^2
~~=~~C_m\ell_T^2,\nn
\eeq
where the third inequality uses $\ell_T\geq1$.

\end{proof}

The preceding lemmas bound the weighted tracking error and cumulative motion. We finish by estimating the inverse adaptive weight in terms of the noise variance, removing the weights by
Cauchy–Schwarz, and applying proximal optimality at the uniformly sampled output.

\begin{theorem} \label{theorem:MVR1}
(Noise-adaptive composite stationarity of {\sf ALFPG-MVR1}) Suppose Assumptions \ref{ass:prox}, \ref{ass:oracle}, \ref{ass:bounded-oracle},
\ref{ass:smooth}, and~\ref{ass:regularizer} hold, and run {\sf ALFPG-MVR1}. Write $\ell_T:=\ln(e+T)$. Then there exists a constant $C_R>0$, independent of
$T$ and $\sigma$, such that, for every $T\geq 1$,
\beq
\ts \frac1T
\sum_{t=1}^T\mathbb E[\mathcal{R}(\x_{t+1})] \leq
C_R\left(\frac{\ell_T^{5/4}}{T^{1/2}}+\frac{ \sigma^{1/2} \ell_T}{T^{1/4}}\right).\nn
\eeq
Consequently, if $R\sim{\rm Unif}\{1,\ldots,T\}$ is independent of the algorithmic randomness, then
$\mathbb E[\mathcal{R}(\x_{R+1})] = \widetilde{\mathcal O}
(T^{-1/2}+\sigma^{1/2}T^{-1/4})$. Since {\sf ALFPG-MVR1} uses one stochastic-gradient
evaluation per iteration, its stochastic-gradient oracle complexity for obtaining
$\mathbb E[\mathcal{R}(\x_{R+1})]\leq\epsilon$ is $\widetilde{\mathcal O}(\epsilon^{-2}+\sigma^2\epsilon^{-4})$.

\end{theorem}

\begin{proof}

Define $D_t:=\sum_{i=1}^t\chi_i^2$, $N_t:=\max_{1\leq i\leq t}\chi_i^2$, $S_T:=\sum_{t=1}^T\|\s_t\|^2$, with $D_0=N_0=0$.

Define $B_T:=1+3(T-1)\sigma^2+2L^2C_h\ell_T$. Let $\chi_1=0$.

Define $C_u$ and $C_m$ be the constants in
Lemmas~\ref{lemma:MVR1:logTerm} and~\ref{lemma:MVR1:bridge}, and set $C_R:=\big((1+L\theta)\sqrt{C_m}+\sqrt{C_u}\big)  \big(3 + 2 L^{1/2} C_h^{1/4}\big)$.

\paragraph{Bounding the term $\mathbb E[\alpha_T^{-1}]$.} For $t\geq2$, set
$\z_t:=\nabla f(\x_t)-\g_{t-1}$.
Since $\z_t$ is $\mathcal F_{t-1}$-measurable and
$\mathbb E_t[\boldsymbol\epsilon_t]=0$, $\mathbb E_t[\chi_t^2]=\|\z_t\|^2+\mathbb E_t\|\boldsymbol\epsilon_t\|^2$. Moreover, $\z_t=\nabla f(\x_t)-\nabla f(\x_{t-1})
-\boldsymbol\epsilon_{t-1}$. By $L$-smoothness, the variance bound, and
$\x_t-\x_{t-1}=-\gamma_{t-1}\r_{t-1}$, we obtain
\beq
\mathbb E[\chi_t^2] &\leq& \ts \mathbb E[ 2\|\nabla f(\x_t)-\nabla f(\x_{t-1})\|^2 + 2 \|\boldsymbol\epsilon_{t-1}\|^2 + \|\boldsymbol\epsilon_t\|^2]\nn\\
&\leq& \ts 2L^2\mathbb{E}[\gamma_{t-1}^2\|\r_{t-1}\|^2] +3\sigma^2.\nn
\eeq
\noi Since $\chi_1=0$, summing over $t=2,\ldots,T$ yields:
\beq \label{eq:theorem:MVR1:gradient-differences}
\mathbb E[D_T] &= & \ts \mathbb E[\sum_{t=1}^T \chi_t^2 ] \leq \ts 3(T-1)\sigma^2 + 2L^2\mathbb E\left[
\sum_{t=1}^{T-1}\gamma_t^2\|\r_t\|^2 \right]  \nn\\
&\leq& \ts  \ts 3(T-1)\sigma^2+2L^2\mathbb E[H_T] ~~\leq~~ 3(T-1)\sigma^2+2L^2C_h\ell_T ~~=~~-1 + B_T.\label{eq:DT}
\eeq
We further derive for $T\geq 2$:
\beq \label{eq:theorem:MVR1:inverse-weight}
\ts \mathbb E[\alpha_T^{-1}] ~~\xxxx{1}{\leq} ~~\mathbb E\sqrt{1+D_T}~~\xxxx{2}{\leq}~~\sqrt{1+\mathbb E[D_T]} \leq B_T^{1/2},
\eeq
\noi where step \step{1} uses the running-minimum definition of $\alpha_T$ for all $T\geq 1$ that $\alpha_T^{-2} = \max_{1\leq k\leq T-1}
\frac{1+D_k}{1+N_k}
\leq 1+D_{T-1}
\leq 1+D_T$; step \step{2} uses Jensen's inequality; step \step{3} uses Inequality (\ref{eq:DT}). Inequality \eqref{eq:theorem:MVR1:inverse-weight} follows directly from $\alpha_1=1$ and $D_1=0$.

\paragraph{Bounding the composite residual.}
Global optimality of the proximal update implies $\r_t-\m_t\in\partial h(\x_{t+1})$. Consequently, $L$-smoothness and $\s_t=\m_t-\nabla f(\x_t)$ give
\beq \label{eq:theorem:MVR1:residual}
\ts \mathcal{R}(\x_{t+1}) &:=& \ts \operatorname{dist}\bigl(0,\nabla f(\x_{t+1})+\partial h(\x_{t+1})\bigr)\nn\\
&\leq&\ts \|\nabla f(\x_{t+1})+\r_t-\m_t\|~~\leq~~\|\r_t\|+\|\s_t\| + \|\nabla f(\x_{t+1})-\nabla f(\x_t)\| \nn\\
&\leq& \ts (1+L\gamma_t)\|\r_t\|+\|\s_t\|~~\leq~~
(1+L\theta)\|\r_t\|+\|\s_t\|,
\eeq
where $\gamma_t=\theta\sqrt{\alpha_t}\beta_t\leq\theta$.

\paragraph{Removing the adaptive weights.} We have:
\beq \label{eq:thethe:1}
\ts \mathbb E[\sqrt{M_T}]~~=~~ \mathbb E \big[
\sqrt{\alpha_T M_T}\sqrt{\alpha_T^{-1}} \big]~~\xxxx{1}{\leq}~~\sqrt{\mathbb E[\alpha_T M_T]}
\sqrt{\mathbb E[\alpha_T^{-1}]} ~~\xxxx{2}{\leq}~~ \sqrt{C_m}\,\ell_T \cdot B_T^{1/4},
\eeq
\noi where step \step{1} uses Cauchy--Schwarz inequality; step \step{2} uses Lemma~\ref{lemma:MVR1:bridge}(b) and Inequality \eqref{eq:theorem:MVR1:inverse-weight};

On the other hand, we derive:
\beq \label{eq:thethe:2}
\ts \mathbb E[\sqrt{S_T}]~~\xxxx{1}{\leq}~~ \mathbb E \big[
\sqrt{U_T}\sqrt{\alpha_T^{-1}}\big]~~\xxxx{2}{\leq}~~\sqrt{\mathbb E[U_T]}\sqrt{\mathbb E[\alpha_T^{-1}]} ~~\leq~~\sqrt{C_u}\,\ell_T^{1/2}B_T^{1/4},
\eeq
\noi where step \step{1} uses $\alpha_T S_T \leq \sum_{t=1}^T\alpha_t\|\s_t\|^2=U_T$ since $\alpha_t\geq\alpha_T$ for $t\leq T$; step \step{2} uses Cauchy--Schwarz inequality; step \step{3} uses Lemma~\ref{lemma:MVR1:logTerm}(d).

\paragraph{Establishing the noise-adaptive rate.}
Since $R$ is uniform and independent of the algorithm, we derive:
\beq
&&\ts \frac1T
\sum_{t=1}^T\mathbb E[\mathcal{R}(\x_{t+1})]
~~\xxxx{1}{\leq}~~ \frac{1+L\theta}{T}
\mathbb E\sum_{t=1}^T\|\r_t\|+\frac1T
\mathbb E\sum_{t=1}^T\|\s_t\| \nn \\
&\xxxx{2}{\leq}& \ts \frac{1+L\theta}{\sqrt T}\mathbb E[\sqrt{M_T}]+\frac1{\sqrt T}\mathbb E[\sqrt{S_T}]
~~\xxxx{3}{\leq}~~\frac{(1+L\theta)\sqrt{C_m}\ell_T+\sqrt{C_u}\ell_T^{1/2}}{\sqrt T} B_T^{1/4} \nn \\
&\leq& \ts \left((1+L\theta)\sqrt{C_m}+\sqrt{C_u}\right) \cdot \frac{\ell_T}{\sqrt T}B_T^{1/4} \nn\\
&\xxxx{4}{\leq}& \ts \underbrace{\left((1+L\theta)\sqrt{C_m}+\sqrt{C_u}\right)  \left(3 + 2 L^{1/2} C_h^{1/4} \right) }_{:=C_R} \cdot \frac{\ell_T}{\sqrt T} \cdot\left(\ell_T^{1/4}+T^{1/4}\sigma^{1/2} \right), \nn
\eeq
where step \step{1} uses Inequality \eqref{eq:theorem:MVR1:residual}; step \step{2} uses
Cauchy--Schwarz inequality; step \step{3} uses Inequalities (\ref{eq:thethe:1}) and (\ref{eq:thethe:2}); step \step{4} uses $B_T^{1/4} \leq (1+2L^2C_h\ell_T)^{1/4} + \bigl(3(T-1)\sigma^2\bigr)^{1/4}\leq \left(3 + 2 L^{1/2} C_h^{1/4} \right) \left(\ell_T^{1/4}+T^{1/4}\sigma^{1/2} \right)$.

\end{proof}

\subsection{Analysis for {\sf ALFPG-MVR2}}
\label{subsect:result:MVR2}

\paragraph{Notation.} Let $p=2/3$ and $\lambda=1/4$. Write $\mathbb E_t[\cdot]=\mathbb E[\cdot\mid\mathcal F_{t-1}]$.

\noi Define $\epsilons_t=\g_t-\nabla f(\x_t)$, $\epsilons_t^-=\g_t^--\nabla f(\x_{t-1})$, $\s_t=\m_t-\nabla f(\x_t)$, $\r_t=(\x_t-\x_{t+1})/\gamma_t$, $A_t=1+\sum_{i=1}^t\|\r_i\|^2$, $D_t=\sum_{i=1}^t\chi_i^2$, $N_t=\max_{1\leq i\leq t}\chi_i^2$, where $A_{-1}=A_0=1$, $D_0=N_0=0$.

Set $\alpha_0=1$, $\gamma_0=0$, $\r_0=\mathbf0$, and $\s_0=-\nabla f(\x_1)$. Set $\g_0=\g_1^-=\g_1$, leading to $\chi_1=0$.

Define $P_t=\frac{\|\s_t\|^2}{\sqrt{\alpha_t}}$, $U_T=\sum_{t=1}^T\alpha_tP_t=\sum_{t=1}^T\sqrt{\alpha_t}\|\s_t\|^2$, $W_T=\sum_{t=1}^T\gamma_t\|\r_t\|^2$, $H_T=\sum_{t=1}^T\frac{\gamma_t^2}{\sqrt{\alpha_t}}\|\r_t\|^2$, $M_T=\sum_{t=1}^T\|\r_t\|^2=A_T-1$, $S_T=\sum_{t=1}^T\|\s_t\|^2$. Set $\ell_T=\ln(e+T)$.

For MVR2, the same-sample correction changes the estimator bounds, so the MVR1 preliminary estimates cannot be reused directly. We first control the estimator increments and proximal residuals; in the convex case, the initial accumulator ratio is handled separately from subsequent iterations.

\begin{lemma}\label{lemma:MVR2:basic:bound}
(Preliminary bounds for MVR2) Let $K_h:=\begin{cases}
\theta^{-1}\|\x_1-\widehat{\x}_1\|,&h\text{ convex},\\
\hat H,&\text{otherwise},
\end{cases}$, $\widehat{\x}_1:=\prox_{\theta h}(\x_1)$. Almost surely, for every $t\geq1$ and $T\geq1$,

\begin{enumerate}[label=\textbf{(\alph*)},leftmargin=18pt,
itemsep=2pt,topsep=1pt]
\item $\|\m_t\|\leq2{{\widehat G}}/\alpha_t$ and
$\|\m_t-\m_{t-1}\|\leq4{{\widehat G}}$.
\item $\|\r_t\|\leq C_r t$, where $C_r = 2K_h+4{{\widehat G}}$.

\item If $h$ is convex, $A_t/A_{t-1}\leq\kappa $ for all $t\geq1$, where $\kappa:=\max(2+16{{\widehat G}}^2,1 + C_r^2)$.

\end{enumerate}
\end{lemma}

\begin{proof}

Let $\alpha_t=\rho_{t-1}$. Define $C_r := 2K_h+4{{\widehat G}}$.

\paragraph{Part (a).}
The stochastic momentum estimator is updated as
$\m_t=(1-\alpha_t)(\m_{t-1}-\g_t^-)+\g_t$.
Since the stochastic oracle is bounded, we derive:
\beq
\|\m_t\| &\leq& \ts (1-\alpha_t)\|\m_{t-1}\|+(2-\alpha_t){{\widehat G}}\nn\\
&\xxxx{1}{\leq}& \ts (1-\alpha_t)\frac{2{{\widehat G}}}{\alpha_t}+(2-\alpha_t){{\widehat G}} ~~=~~\frac{2{{\widehat G}}}{\alpha_t}-\alpha_t{{\widehat G}} ~~\leq~~\frac{2{{\widehat G}}}{\alpha_t},\nn
\eeq
\noi where step \step{1} uses $\alpha_t\leq\alpha_{t-1}$ and
$\|\m_{t-1}\|\leq2{{\widehat G}}/\alpha_{t-1}$; the initial case follows
from $\m_0=\zero$.
Furthermore,
\begin{eqnarray*}
\|\m_t-\m_{t-1}\|
&\leq&\alpha_t\|\m_{t-1}\|+(2-\alpha_t){{\widehat G}}
~\leq~4{{\widehat G}}.\nn
\end{eqnarray*}

\paragraph{Part (b-i).}
Suppose that $h$ is convex.
Since $\gamma_1=\theta$ and $\m_1=\g_1$,
nonexpansiveness of the proximal operator gives
\begin{eqnarray*}
\|\r_1\|
&\leq&\frac{\|\x_1-\widehat{\x}_1\|}{\theta}
 +\frac{\|\widehat{\x}_1-\x_2\|}{\theta}
~\leq~K_h+{{\widehat G}}.\nn
\end{eqnarray*}
For $t\geq2$, proximal optimality implies
$\r_t-\m_t\in\partial h(\x_{t+1})$ and
$\r_{t-1}-\m_{t-1}\in\partial h(\x_t)$.
Monotonicity of $\partial h$, together with
$\x_{t+1}-\x_t=-\gamma_t\r_t$, yields
\beq
\|\r_t\|^2~~\leq~~\langle\r_{t-1}+\m_t-\m_{t-1},\r_t\rangle~~\leq~~\bigl(\|\r_{t-1}\|+4{{\widehat G}}\bigr)\|\r_t\|.\nn
\eeq
Therefore $\|\r_t\|\leq\|\r_{t-1}\|+4{{\widehat G}}$, and hence $\|\r_t\|\leq C_r \cdot t$.

\paragraph{Part (b-ii).}
Suppose that $h$ is nonconvex and Lipschitz on $\dom(h)$.
Global optimality of the proximal update, with $\x_t$ as a comparison
point, gives
\begin{eqnarray*}
\frac{\|\x_{t+1}-\x_t\|^2}{2\gamma_t}
&\leq&h(\x_t)-h(\x_{t+1})
 -\langle\m_t,\x_{t+1}-\x_t\rangle\nn\\
&\leq&(K_h+\|\m_t\|)\|\x_{t+1}-\x_t\|.\nn
\end{eqnarray*}
If $\x_{t+1}=\x_t$, then $\r_t=\zero$.
Otherwise, dividing by $\|\x_{t+1}-\x_t\|$ gives
\beq
\|\r_t\|~~\leq~~2(K_h+\|\m_t\|)
~~\leq~~2K_h+4{{\widehat G}} t^{2/3} ~~\leq~~C_r \cdot t.\nn
\eeq

\paragraph{Part (c).} For $t\geq 2$, we derive:
\beq
\tfrac{A_t}{A_{t-1}} ~~\xxxx{1}{=}~~ \tfrac{A_{t-1} + \|\r_t\|^2}{A_{t-1}} ~~{\leq}~~ 1 + \tfrac{\|\r_t\|^2}{A_{t-1}} ~~\xxxx{2}{\leq}~~ 1 + \tfrac{(\|\r_{t-1}\| + 4{{\widehat G}} )^2}{1 + \|\r_{t-1}\|^2} ~~\xxxx{3}{\leq}~~ 1 +  1 + 16{{\widehat G}}^2,\nn
\eeq
\noi where step \step{1} uses the definition of $A_t:=1+\sum_{i=1}^t\|\r_i\|^2$; step \step{2} uses $\|\r_t\|\leq\|\r_{t-1}\|+4{{\widehat G}}$; step \step{3} uses $(a+b)^2\leq (1+a^2)(1+b^2)$ for all $a,b$.
For $t=1$, we have: $\tfrac{A_t}{A_{t-1}} = 1 + \|\r_1\|^2 \leq 1 + C_r^2$. Taking both cases into consideration, we have $\tfrac{A_t}{A_{t-1}} \leq \kappa:=\max(2+16{{\widehat G}}^2,1 + C_r^2)$.

\end{proof}

We now exploit individual smoothness to control the difference of two gradients evaluated with the same sample. Rescaling the tracking error produces a contractive recursion with an explicit
innovation term and a conditionally mean-zero cross term.

\begin{lemma}\label{lemma:MVR2:s:recursion}
(Tracking-error recursion for MVR2) Retain the notation introduced above, and define
\[
\begin{gathered}
\boldsymbol\delta_t
:=\epsilons_t-\epsilons_t^-,
\qquad
\boldsymbol\zeta_t
:=\frac{\boldsymbol\delta_t+\alpha_t\epsilons_t^-}
        {\alpha_t^{1/4}},
\qquad
\varpi_t:=(1-\alpha_t)
\left(\frac{\alpha_{t-1}}{\alpha_t}\right)^{1/4},
\\
\mathcal E_t
:=2\varpi_t\langle\z_{t-1},\boldsymbol\zeta_t\rangle,
\qquad
\mathcal V_t:=\|\boldsymbol\zeta_t\|^2,
\qquad
\dot{\mathcal E}_t
:=\alpha_t^{3/2}\|\epsilons_t^-\|^2,\quad \ddot{\mathcal E}_t
:=11 L^2
\frac{\gamma_{t-1}^2}{\sqrt{\alpha_{t-1}}}
\|\r_{t-1}\|^2.
\end{gathered}
\]
Then, almost surely, for every $t\geq1$:

\begin{enumerate}[label=\textbf{(\alph*)},leftmargin=18pt,
itemsep=2pt,topsep=1pt]
\item
$\varpi_t^2\leq1-\alpha_t$, and
$\mathcal V_t\leq\ddot{\mathcal E}_t
+2\dot{\mathcal E}_t$.
In particular,
\begin{equation}
P_t
=\varpi_t^2P_{t-1}+\mathcal V_t+\mathcal E_t
\leq(1-\alpha_t)P_{t-1}
+\ddot{\mathcal E}_t
+2\dot{\mathcal E}_t+\mathcal E_t.
\label{eq:MVR2:P-recursion}
\end{equation}

\item $\mathbb E_t[\mathcal E_t]=0$, and $\mathcal E_t^2\leq4P_{t-1}\mathcal{V}_t$.

\end{enumerate}
\end{lemma}

\begin{proof}

Define $\s_t=\m_t-\nabla f(\x_t)$, and $\z_t=\s_t/\alpha_t^{1/4}$, where $\z_0=\zero$ and $P_0=0$.

Define the random increment $\boldsymbol\eta_t:=\epsilons_t-(1-\alpha_t)\epsilons_t^-
=\boldsymbol\delta_t+\alpha_t\epsilons_t^-$. Let $\alpha_t=\rho_{t-1}$.

\paragraph{Exact decomposition of the tracking error.}

Using the MVR2 update, we obtain
\begin{align}
\s_t
&=(1-\alpha_t)(\m_{t-1}-\g_t^-)
  +\g_t-\nabla f(\x_t) \notag\\
&=(1-\alpha_t)\s_{t-1}
  +\epsilons_t-(1-\alpha_t)\epsilons_t^- \notag\\
&=(1-\alpha_t)\s_{t-1}+\boldsymbol\eta_t.
\label{eq:MVR2:tracking:decomposition}
\end{align}
Dividing by $\alpha_t^{1/4}$ gives
\[
\z_t=\varpi_t\z_{t-1}+\boldsymbol\zeta_t.
\]
At $t=1$, $\alpha_1=1$ and
$\boldsymbol\delta_1=\zero$, so this identity also
holds under the convention $\z_0=\zero$.

\paragraph{Bounding the term $\mathcal V_t:=\|\boldsymbol\zeta_t\|^2$.} By individual smoothness and smoothness of the expected loss,
\beq
\ts \|\deltas_t\|&=& \ts \|\g_t-\g_t^- -\nabla f(\x_t)+\nabla f(\x_{t-1})\|~~\leq~~\|\g_t-\g_t^-\|+\|\nabla f(\x_t)-\nabla f(\x_{t-1})\| \nn \\
&\leq& \ts 2L\|\x_t-\x_{t-1}\|~~=~~2L\gamma_{t-1}\|\r_{t-1}\|.\nn
\eeq
The last equality also holds at $t=1$ because
$\x_0=\x_1$ and $\gamma_0=0$.
Consequently,
\beq
\|\boldsymbol\eta_t\|^2~~=~~ \|\boldsymbol\delta_t+\alpha_t\epsilons_t^-\|^2~~\leq~~2\|\boldsymbol\delta_t\|^2 + 2\alpha_t^2\|\epsilons_t^-\|^2~~\leq~~8L^2\gamma_{t-1}^2\|\r_{t-1}\|^2 + 2\alpha_t^2\|\epsilons_t^-\|^2.\nn
\eeq
Dividing both sides by $\sqrt{\alpha_t}$ and using the definition of $\VV_t=\frac{\|\etas_t\|^2}{\sqrt{\alpha_t}}$ yields:
\beq \label{eq:MVR2:innovation-bound}
\mathcal V_t~~\leq ~~ \ts \frac{8L^2\gamma_{t-1}^2}{\sqrt{\alpha_t}}
\|\r_{t-1}\|^2 +2\alpha_t^{3/2}\|\epsilons_t^-\|^2
~~\xxxx{1}{\leq}~~ \ts 8L^2\sqrt{\frac53}
\frac{\gamma_{t-1}^2}{\sqrt{\alpha_{t-1}}}
\|\r_{t-1}\|^2 +2\dot{\mathcal E}_t
~~\leq ~~\ddot{\mathcal E}_t+2\dot{\mathcal E}_t,
\eeq
\noi where step \step{1} uses $\alpha_{t-1}/\alpha_t = {\rho_{t-2}}/{\rho_{t-1}} \leq1 + 2/3$ by Lemma~\ref{lemma:weight}(a).

\paragraph{Part (a).} For all $t\geq 2$, we derive:
\beq \label{eq:varpi2}
\ts \varpi_t^2~~\xxxx{1}{=}~~(1-\alpha_t)^2\sqrt{\frac{\alpha_{t-1}}{\alpha_t}}~~\xxxx{2}{\leq}~~\frac{(1-\alpha_t)^2}{\sqrt{1-\frac23\alpha_t}}~~\xxxx{3}{\leq}~~1-\alpha_t,
\eeq
\noi where step \step{1} uses the definition of $\varpi$; step \step{2} uses $\frac1{\alpha_t}-\frac1{\alpha_{t-1}}\leq\frac23$, leading to $\frac{\alpha_{t-1}}{\alpha_t}
\leq\frac1{1-\frac23\alpha_t}$; step \step{3} uses $(1-a)\leq\sqrt{1-\tfrac{2}{3}a}$ for $a\in[0,1]$. Inequality (\ref{eq:varpi2}) also holds at $t=1$ because $\varpi_1=0$.

Expanding the squared norm of the rescaled recursion yields
\beq
P_t &=&\|\varpi_t\z_{t-1}+\boldsymbol\zeta_t\|^2~~=~~\varpi_t^2P_{t-1}
  +\|\boldsymbol\zeta_t\|^2+2\varpi_t\langle\z_{t-1},\boldsymbol\zeta_t\rangle \nn \\
&=& \varpi_t^2P_{t-1}+\mathcal V_t+\mathcal E_t ~~\xxxx{1}{\leq} ~~ (1-\alpha_t)P_{t-1}+\ddot{\mathcal E}_t + 2\dot{\mathcal E}_t+\mathcal E_t,\nn
\eeq
where step \step{1} uses Inequality \eqref{eq:MVR2:innovation-bound}, $\varpi_t^2\leq 1-\alpha_t$, and $P_{t-1}\geq 0$.

\paragraph{Part (b).}
By conditional unbiasedness, $\mathbb E_t[\epsilons_t]=\zero$, and $\mathbb E_t[\epsilons_t^-]=\zero$. Since $\alpha_t$ is $\mathcal F_{t-1}$-measurable,
\[
\mathbb E_t[\boldsymbol\zeta_t]
=\alpha_t^{-1/4}
\mathbb E_t[
\epsilons_t-(1-\alpha_t)\epsilons_t^-]
=\zero.
\]
Moreover, $\varpi_t$ and $\z_{t-1}$ are
$\mathcal F_{t-1}$-measurable. Therefore,
\[
\mathbb E_t[\mathcal E_t]
=2\varpi_t
\langle\z_{t-1},
\mathbb E_t[\boldsymbol\zeta_t]\rangle
=0.
\]
Finally, Cauchy--Schwarz and $\varpi_t^2\leq1$ give
\[
\mathcal E_t^2
=4\varpi_t^2
\langle\z_{t-1},\boldsymbol\zeta_t\rangle^2
\leq4\varpi_t^2P_{t-1}\mathcal V_t
\leq4P_{t-1}\mathcal V_t.
\]
\end{proof}

The rescaled recursion is now coupled to proximal descent. The resulting Lyapunov inequality tracks both objective progress and estimator accuracy, with the normalized motion accumulated in
the correction term.

\begin{lemma}
\label{lemma:MVR2:descent}
(Proximal descent and a joint Lyapunov inequality for MVR2) Let $C_f:=(1+L)/2+22 L^2$ and
$\mathcal Q_t:=F(\x_t)+2P_t-C_fH_{t-1}$.
Then, almost surely,
\begin{align}
F(\x_{t+1})-F(\x_t)+\frac{\gamma_t}{2}\|\r_t\|^2
&\leq\frac{\sqrt{\alpha_t}}2\|\s_t\|^2
 +\frac{1+L}{2}\frac{\gamma_t^2}{\sqrt{\alpha_t}}\|\r_t\|^2,
\label{eq:MVR2:descent}\\
\mathcal Q_{t+1}-\mathcal Q_t+\frac{\gamma_t}{2}\|\r_t\|^2
 +\frac{\sqrt{\alpha_t}}2\|\s_t\|^2
&\leq4\dot{\mathcal E}_{t+1}+2\mathcal E_{t+1}.
\label{eq:MVR2:one-step-descent}
\end{align}
\end{lemma}

\begin{proof}

By global optimality of the proximal update and $L$-smoothness of $f$,
\[
\begin{aligned}
h(\x_{t+1})&\leq h(\x_t)
 -\langle\m_t,\x_{t+1}-\x_t\rangle
 -\frac{\|\x_{t+1}-\x_t\|^2}{2\gamma_t},\\
f(\x_{t+1})&\leq f(\x_t)
 +\langle\nabla f(\x_t),\x_{t+1}-\x_t\rangle
 +\frac L2\|\x_{t+1}-\x_t\|^2.
\end{aligned}
\]

\paragraph{Part (a).}
Adding the preceding two inequalities and using
$\x_{t+1}-\x_t=-\gamma_t\r_t$, we derive:
\begin{eqnarray*}
F(\x_{t+1})-F(\x_t)
&\leq&\gamma_t\langle\r_t,\s_t\rangle
 -\frac{\gamma_t}{2}\|\r_t\|^2
 +\frac{L\gamma_t^2}{2}\|\r_t\|^2\nn\\
&\xxxx{1}{\leq}&\frac{\sqrt{\alpha_t}}2\|\s_t\|^2
 +\frac{\gamma_t^2}{2\sqrt{\alpha_t}}\|\r_t\|^2
 -\frac{\gamma_t}{2}\|\r_t\|^2
 +\frac{L\gamma_t^2}{2}\|\r_t\|^2\nn\\
&\leq&\frac{\sqrt{\alpha_t}}2\|\s_t\|^2
 -\frac{\gamma_t}{2}\|\r_t\|^2
 +\frac{1+L}{2}\frac{\gamma_t^2}{\sqrt{\alpha_t}}\|\r_t\|^2,\nn
\end{eqnarray*}
\noi where step \step{1} uses the weighted Young inequality
$\gamma_t\langle\r_t,\s_t\rangle
\leq\frac{\sqrt{\alpha_t}}2\|\s_t\|^2
+\frac{\gamma_t^2}{2\sqrt{\alpha_t}}\|\r_t\|^2$.

\paragraph{Part (b).} Since $\sqrt{\alpha_t}\|\s_t\|^2=\alpha_tP_t$, part (a) gives
\beq
&& F(\x_{t+1})-F(\x_t)+\frac{\gamma_t}{2}\|\r_t\|^2 +\frac{\sqrt{\alpha_t}}2\|\s_t\|^2 ~~\leq~~ \ts \alpha_tP_t +\frac{1+L}{2}\frac{\gamma_t^2}{\sqrt{\alpha_t}}\|\r_t\|^2 \nn \\
&\xxxx{1}{\leq}& 2\alpha_{t+1}P_t
 +\frac{1+L}{2}\frac{\gamma_t^2}{\sqrt{\alpha_t}}\|\r_t\|^2~~\xxxx{2}{\leq}~~2(P_t-P_{t+1})
 +C_f\frac{\gamma_t^2}{\sqrt{\alpha_t}}\|\r_t\|^2
 +4\dot{\mathcal E}_{t+1}+2\mathcal E_{t+1},\nn
\eeq
\noi where step \step{1} uses $\alpha_t/\alpha_{t+1}\leq1+2/3<2$;
step \step{2} uses Lemma~\ref{lemma:MVR2:s:recursion}(a) and
$\ddot{\mathcal E}_{t+1}
=11 L^2 \gamma_t^2\|\r_t\|^2/\sqrt{\alpha_t}$.
Finally,
$H_t-H_{t-1}=\gamma_t^2\|\r_t\|^2/\sqrt{\alpha_t}$.
Substituting the definition of $\mathcal Q_t$ proves
\eqref{eq:MVR2:one-step-descent}.
\end{proof}

We next establish the cumulative estimates needed to sum the MVR2 descent inequality. The proof
treats the convex and domain-Lipschitz motion bounds separately and controls second moments
without assuming that the rescaled tracking error is uniformly bounded.

\begin{lemma}\label{lemma:MVR2:logTerm}

(Logarithmic bounds for MVR2) Then, for every $T\geq1$,

\begin{enumerate}[label=\textbf{(\alph*)},leftmargin=18pt,
itemsep=2pt,topsep=1pt]
\item $V_T:=\sum_{t=1}^T \alpha_t^{3/2} \chi_t^2 \leq C_v\ell_T$ almost surely, where $C_v:=2 (1+4{{\widehat G}}^2)^2 \ln(e+4{{\widehat G}}^2)$.

\item Let $Q_T:=\sum_{t=1}^T\dot{\mathcal E}_t$. We have: $\mathbb E[Q_{T}]\leq C_q \ell_T$ and
$\mathbb E[Q_{T}^2]\leq C_q^2\ell_T^2$, where $C_q:= 8{{\widehat G}}^2+2C_v$.

\item Let $H_T:=\sum_{t=1}^T\frac{\gamma_t^2}{\sqrt{\alpha_t}}\|\r_t\|^2$. We have: $\EEE[H_T]\leq C_h\ell_T$ and $\EEE[H_T^2]\leq C_h^2\ell_T^2$, where $C_h:=\begin{cases}
\dot{C}_h,&h\text{ convex}\\
\ddot{C}_h,&\text{otherwise}
\end{cases}$, $\dot{C}_h:=6\theta^2 \kappa \ln( e + C_r)$, $\ddot{C}_h:=\theta^2\bigl[C_A+B_h+80(11 L^2\theta^2+2C_q)\bigr]$, $C_A:=3+\ln(1+C_r^2)$, and $B_h:=8(K_h+{{\widehat G}})^2$.

\item Let $U_T:=\sum_{t=1}^T\sqrt{\alpha_t}\|\s_t\|^2=\sum_{t=1}^T\alpha_tP_t$. Then $\mathbb E[U_T]\leq C_u\ell_T$ and $\mathbb E[U_T^2]\leq C_u^2\ell_T^2$. Here, $C_u:=20(11 L^2 C_h+2C_q)$.

\end{enumerate}

\end{lemma}

\begin{proof}

Define $\chi_t:= \|\g_{t-1} - \g_t^-\|$. Let $c_\chi:=(1+4{{\widehat G}}^2)^2$.

Define $\alpha_t = \rho_{t-1}$, where $\rho_t :=
\min_{1\le k\le t}\left(\frac{1+\max_{1\le i\le k}\chi_i^2}{1+\sum_{i=1}^k\chi_i^2}\right)^{2/3}$.

Define $D_t:=\sum_{i=1}^t\chi_i^2$, and $N_t:=\max_{1\leq i\leq t}\chi_i^2$.

\paragraph{Part (a).} We derive:
\[{\begin{array}{rcl}
\ts V_T &\leq& \ts  V_{T+1} ~~\xxxx{1}{=}~~ \sum_{t=2}^{T+1}\alpha_t^{3/2}\chi_t^2~~\xxxx{2}{\leq}~~(1+4{{\widehat G}}^2)\sum_{t=2}^{T+1}
             \frac{\chi_t^2}{1+D_{t-1}}\\
&\xxxx{3}{\leq}& \ts (1+4{{\widehat G}}^2)^2\sum_{t=2}^{T+1}
             \frac{\chi_t^2}{1+D_t}~~\xxxx{4}{\leq}~~c_\chi \big(\ln(1+D_{T+1})-\ln(1+D_{1})\big)\\
&\xxxx{}{\leq}&   c_\chi \ln(1+4{{\widehat G}}^2 (T+1))
~~\xxxx{5}{\leq}~~\underbrace{ 2 c_\chi \ln(e+4{{\widehat G}}^2)  }_{:=C_v}\ell_T.
\end{array}}\]
\noi where step \step{1} uses $\chi_1=0$; step \step{2} uses $\alpha_t^{3/2}\leq(1+N_{t-1})/(1+D_{t-1})$ and $N_{t-1}\leq 4{{\widehat G}}^2$; step \step{3} uses $(1+D_t)/(1+D_{t-1})\leq1+4{{\widehat G}}^2$; step \step{4} uses $\chi_t^2/(1+D_t)\leq\ln(1+D_t)-\ln(1+D_{t-1})$; step \step{5} uses $\ln(1+ a (T+1)) \leq \ln(e+a)\ln(e+T+1)\leq 2\ln(e+a)\ln(e+T)$ for all $a\geq 0$.

\paragraph{Part (b-i).}
Let $q_t:=\dot{\mathcal E}_t
=\alpha_t^{3/2}\|\epsilons_t^-\|^2$ and
$Q_t:=\sum_{i=1}^t q_i$, where $Q_0=0$.
Since $0<\alpha_t\leq1$ and
$\|\epsilons_t^-\|\leq2{{\widehat G}}$, we have
\[
0\leq q_t\leq4{{\widehat G}}^2,
\qquad
q_t^2\leq4{{\widehat G}}^2q_t.
\]

For $t\geq2$, recall that
$\chi_t=\|\epsilons_{t-1}-\epsilons_t^-\|$.
Since $\epsilons_{t-1}$ is
$\mathcal F_{t-1}$-measurable and
$\mathbb E_t[\epsilons_t^-]=\mathbf0$, we obtain
\[
\mathbb E_t[\chi_t^2]
=
\|\epsilons_{t-1}\|^2
+\mathbb E_t\|\epsilons_t^-\|^2
\geq
\mathbb E_t\|\epsilons_t^-\|^2.
\]
Because $\alpha_t$ is $\mathcal F_{t-1}$-measurable, $\mathbb E_t[q_t]\leq\alpha_t^{3/2}\mathbb E_t[\chi_t^2]$ for all $t\geq2$.

We obtain:
\beq
V_{T+1}&\leq&
(1+4{{\widehat G}}^2)^2\ln(1+D_{T+1})~~\xxxx{1}{\leq}~~
(1+4{{\widehat G}}^2)^2\ln(1+4{{\widehat G}}^2T)\nn\\
&\leq& (1+4{{\widehat G}}^2)^2\ln(e+4{{\widehat G}}^2)\ell_T.
\nn
\eeq
\noi where step \step{1} uses $D_{T+1}=\sum_{t=2}^{T+1}\chi_t^2
\leq4{{\widehat G}}^2T$ since $\chi_1=0$ and $\chi_t^2\leq4{{\widehat G}}^2$. Taking expectations, summing the conditional comparison over
$t=2,\ldots,T+1$, and bounding $q_1$ separately therefore gives
\beq
\ts \mathbb E[Q_T] &\leq& \ts \mathbb E[Q_{T+1}]~~\leq~~ 4{{\widehat G}}^2+
\mathbb E\left[
\sum_{t=2}^{T+1}\alpha_t^{3/2}\chi_t^2
\right]\nn\\
&\leq& \ts 4{{\widehat G}}^2+\frac{C_v}{2}\ell_T~~\leq~~
\frac{C_v}{2}
+\frac{C_v}{2}\ell_T
\leq C_v\ell_T. \nn
\eeq

\paragraph{Part (b-ii).}
The identity
$Q_t^2-Q_{t-1}^2=2Q_{t-1}q_t+q_t^2$,
together with $Q_0=0$, gives
\beq
\ts \mathbb E[Q_{T+1}^2] = 2\mathbb E\left[
\sum_{t=1}^{T+1}Q_{t-1}q_t \right] + \mathbb E\left[
\sum_{t=1}^{T+1}q_t^2 \right].
\nn
\eeq
For $t\geq2$, the measurability of
$Q_{t-1}\alpha_t^{3/2}$ gives
\beq
\ts \mathbb E[Q_{t-1}q_t]~~=~~\mathbb E\left[
Q_{t-1}\alpha_t^{3/2} \mathbb E_t\|\epsilons_t^-\|^2
\right]~~\leq~~ \mathbb E\left[Q_{t-1}\alpha_t^{3/2}\mathbb E_t[\chi_t^2]\right]~~=~~ \mathbb E\left[Q_{t-1}\alpha_t^{3/2}\chi_t^2\right].
\nn
\eeq
Consequently,
\beq
\label{eq:MVR2:noise-second-moment}
\mathbb E[Q_{T+1}^2]&\xxxx{1}{\leq}& \ts 4{{\widehat G}}^2\mathbb E[Q_{T+1}]+2\mathbb E\left[ \sum_{t=2}^{T+1} Q_{t-1}\alpha_t^{3/2}\chi_t^2\right] ~~\xxxx{2}{\leq}~~
4{{\widehat G}}^2\mathbb E[Q_{T+1}] +2\mathbb E[Q_{T+1}V_{T+1}]\nn \\
&\xxxx{3}{\leq}&\ts (4{{\widehat G}}^2+2C_v\ell_T)
\mathbb E[Q_{T+1}]~~\xxxx{4}{\leq}~~ (4{{\widehat G}}^2+2C_v\ell_T)(4{{\widehat G}}^2+C_v\ell_T) \nn \\
&\leq& \ts 2(4{{\widehat G}}^2+C_v\ell_T)^2
\leq C_q^2\ell_T^2,
\eeq
\noi where step \step{1} uses $q_t^2\leq4{{\widehat G}}^2q_t$ and $Q_0=0$; step \step{2} uses $Q_{t-1}\leq Q_{T+1}$; step \step{3} uses $V_{T+1}\leq C_v\ell_T$; step \step{4} uses $\mathbb E[Q_{T+1}] \leq4{{\widehat G}}^2+C_v\ell_T$.

\paragraph{Part (c-i).} We consider the case when $h(\cdot)$ is convex. We have almost surely:
\[{\begin{array}{rcl}
\ts H_T &=& \ts \sum_{t=1}^T
\frac{\gamma_t^2}{\sqrt{\alpha_t}}\|\r_t\|^2~~\xxxx{1}{=} ~~\theta^2 \sum_{t=1}^T \tfrac{\|\r_t\|^2}{A_{t-1}} ~~\xxxx{2}{\leq}~~\theta^2 \kappa \sum_{t=1}^T \tfrac{\|\r_t\|^2}{A_{t}} ~~\xxxx{3}{\leq} ~~ \theta^2 \kappa \ln(A_T)\\
&\xxxx{4}{\leq}& \ts \theta^2 \kappa \ln( 1 + T C_r^2 T^2 )~~\leq~~\theta^2 \kappa \ln( e + T^3 ) \ln( e + C_r^2 )~~\leq~~6\theta^2 \kappa \ln( e + C_r)\ell_T,
\end{array}}\]
\noi where step \step{1} uses $\gamma_t=\theta\alpha_t^{1/4}\beta_t= \theta\alpha_t^{1/4} A_{t-1}^{-1/2}$; step \step{2} uses $A_t/A_{t-1}\leq\kappa$; step \step{3} uses $\frac{\|\r_t\|^2}{A_t}\leq \ln A_t-\ln A_{t-1}$; step \step{4} uses $A_T = 1 + \sum_{t=1}^T \|\r_t\|^2$ and $\|\r_t\|\leq C_r t$. Letting $C_h=\dot{C}_h=6\theta^2 \kappa \ln( e + C_r)$, we directly conclude that
\begin{equation}
\mathbb E[H_T^2]\leq C_h^2\ell_T^2,
\qquad
\mathbb E[H_T]\leq C_h\ell_T.\nn
\end{equation}

\paragraph{Part (c-ii).} We consider the case when $h(\cdot)$ is Lipschitz on $\dom(h)$ with constant $K_h$ (convexity is not required). Define $A_t=1+\sum_{i=1}^t\|\r_i\|^2$, and $P_t=\frac{\|\s_t\|^2}{\sqrt{\alpha_t}}$, $C_r:=2K_h+4{{\widehat G}}$, $C_A:=3+\ln(1+C_r^2)$, $B_h:=8(K_h+{{\widehat G}})^2$, $C_v:=2(1+4{{\widehat G}}^2)^2\ln(e+4{{\widehat G}}^2)$, $C_q:=8{{\widehat G}}^2+2C_v$.

Define $v_t:=\frac{\mathcal V_t}{A_{t-1}}$, $e_t:=\frac{\mathcal E_t}{A_{t-1}}$, $\overline B_t:=\sum_{i=1}^t v_i$, $\overline Z_t:=\sum_{i=1}^t e_i$, with $\overline B_0=\overline Z_0=0$.

\noindent$\quad\bullet$\enspace
\emph{A telescoping potential.} Define $\Psi_t:=\frac{B_h+8P_t}{A_{t-1}}$. Global optimality of the proximal update gives
\beq
\frac{\|\x_{t+1}-\x_t\|^2}{2\gamma_t}~~\leq~~ h(\x_t)-h(\x_{t+1})
-\langle\m_t,\x_{t+1}-\x_t\rangle~~\leq~~(K_h+\|\m_t\|)\|\x_{t+1}-\x_t\|.
\eeq
Thus $\|\r_t\|\leq2(K_h+\|\m_t\|)$, also when
$\x_{t+1}=\x_t$.
Since $\|\m_t\|\leq{{\widehat G}}+\|\s_t\|$ and
$\|\s_t\|^2=\sqrt{\alpha_t}P_t\leq P_t$, we have
$\|\r_t\|^2\leq B_h+8P_t$.
Using $\|\r_t\|^2=A_t-A_{t-1}$, we obtain
\[{\begin{array}{rcl}
\frac{\|\r_t\|^2}{A_{t-1}} ~~=~~\frac{\|\r_t\|^2}{A_t}
+\|\r_t\|^2 \left(\frac1{A_{t-1}}-\frac1{A_t}\right)~~\leq~~\ln A_t-\ln A_{t-1} +(B_h+8P_t) \left(\frac1{A_{t-1}}-\frac1{A_t}\right).
\end{array}}\]

On the other hand,
\beq
\ts \Psi_{t+1}-\Psi_t
=-(B_h+8P_t)
\left(\frac1{A_{t-1}}-\frac1{A_t}\right)
+\frac{8(P_{t+1}-P_t)}{A_t}.
\eeq
Adding these relations and using the tracking-error recursion yields
\beq
\ts \tfrac{\|\r_t\|^2}{A_{t-1}}+\Psi_{t+1}-\Psi_t
\leq\ln A_t-\ln A_{t-1}
+8\frac{\mathcal V_{t+1}+\mathcal E_{t+1}}{A_t}.
\eeq
Summing over $t=1,\ldots,T$, and using $\Psi_{T+1}\geq0$,
$\Psi_1=B_h+8P_1$, $P_1=\mathcal V_1$, and $\mathcal E_1=0$, gives
\beq \label{eq:MVR2:c:potential-bound}
\ts  H_T:=\theta^2\sum_{t=1}^T\tfrac{\|\r_t\|^2}{A_{t-1}}\leq\theta^2\left[
\ln A_T+B_h
+8\sum_{t=1}^{T+1}
\frac{\mathcal V_t+\mathcal E_t}{A_{t-1}}
\right].
\eeq

\noindent$\quad\bullet$\enspace
\emph{Bounding the cumulative innovations.} We first bound the nonnegative cumulative innovation
$\overline B_{T+1}$.
For $t\geq2$, the innovation estimate gives
\[
v_t\leq 11 L^2\theta^2
\frac{\|\r_{t-1}\|^2}{A_{t-2}A_{t-1}}
+\frac{2\dot{\mathcal E}_t}{A_{t-1}}.
\]
Since $\mathcal V_1=\dot{\mathcal E}_1$ and $A_{t-1}\geq1$,
\beq
\ts \overline B_{T+1}&\leq&\ts 11 L^2\theta^2\sum_{t=1}^T
\frac{\|\r_t\|^2}{A_{t-1}A_t}+2Q_{T+1}\nn\\
&=& \ts 11 L^2\theta^2\big(1-\tfrac1{A_T}\big)+2Q_{T+1}~~\leq~~ 11 L^2\theta^2+2Q_{T+1}.\nn
\eeq
The proof of part (b) established
$\mathbb E[Q_{T+1}^2]\leq C_q^2\ell_T^2$.
Therefore, Cauchy--Schwarz and $\ell_T\geq1$ yield
\begin{equation}
\mathbb E[\overline B_{T+1}^2]
\leq\left(11 L^2\theta^2
+2\sqrt{\mathbb E[Q_{T+1}^2]}\right)^2
\leq(11 L^2\theta^2+2C_q)^2\ell_T^2.
\label{eq:MVR2:c:innovation-second-moment}
\end{equation}

\noindent$\quad\bullet$\enspace
\emph{A direct second-moment estimate.}
We use partial sums as weights and summation by parts, as in the proof
of Lemma~\ref{lemma:MVR1:logTerm}(e).
For $t\geq2$, the tracking-error recursion, $\varpi_t^2\leq1$, and
$A_{t-1}\geq A_{t-2}$ imply
\[
\frac{P_t}{A_{t-1}}
\leq\frac{P_{t-1}}{A_{t-2}}+v_t+e_t.
\]
Since $P_1=\mathcal V_1$ and $e_1=0$, induction gives
\[
0\leq\frac{P_t}{A_{t-1}}
\leq\overline B_t+\overline Z_t.
\]
Moreover, $A_{t-1}$ is $\mathcal F_{t-1}$-measurable, so
$\mathbb E_t[e_t]=0$.
The pointwise estimate
$\mathcal E_t^2\leq4P_{t-1}\mathcal V_t$ gives, for $t\geq2$,
\[
e_t^2
\leq4\frac{P_{t-1}}{A_{t-1}}v_t
\leq4(\overline B_{t-1}+\overline Z_{t-1})v_t.
\]
The last inequality also holds at $t=1$ because $e_1=0$.
Expanding the square and using conditional mean zero, we obtain
\beq
\ts \mathbb E[\overline Z_{T+1}^2]
=\mathbb E\sum_{t=1}^{T+1}e_t^2
\leq4\mathbb E\sum_{t=1}^{T+1}
(\overline B_{t-1}+\overline Z_{t-1})v_t.\nn
\eeq
To control the right-hand side, summation by parts yields
\beq
\ts 2\sum_{t=1}^{T+1}\overline B_{t-1}v_t &=&\ts \overline B_{T+1}^2-\sum_{t=1}^{T+1}v_t^2, \nn \\
\ts \sum_{t=1}^{T+1}\overline Z_{t-1}v_t &=&\ts \overline B_{T+1}\overline Z_{T+1}
-\sum_{t=1}^{T+1}\overline B_{t-1}e_t
-\sum_{t=1}^{T+1}v_te_t.\nn
\eeq
Because $\mathbb E[\overline B_{t-1}e_t]=0$, it follows that
\beq
\ts \mathbb E[\overline Z_{T+1}^2]~~\leq~~2\mathbb E[\overline B_{T+1}^2]
+4\mathbb E[\overline B_{T+1}\overline Z_{T+1}]
-4\mathbb E\sum_{t=1}^{T+1}v_te_t.\nn
\eeq
Cauchy--Schwarz bounds the two cross terms as follows:
\beq
\ts \big|\mathbb E[\overline B_{T+1}\overline Z_{T+1}]\big|
~~\leq~~\sqrt{
\mathbb E[\overline B_{T+1}^2]
\mathbb E[\overline Z_{T+1}^2]},\nn
\eeq
and
\beq
\ts \left|\mathbb E\sum_{t=1}^{T+1}v_te_t\right|~~\leq~~\sqrt{
\big(\mathbb E\sum_{t=1}^{T+1}v_t^2\big)
\big(\mathbb E\sum_{t=1}^{T+1}e_t^2\big)}~~\leq~~\sqrt{
\mathbb E[\overline B_{T+1}^2]
\mathbb E[\overline Z_{T+1}^2]}.\nn
\eeq
Consequently, Young's inequality gives
\beq
\ts \mathbb E[\overline Z_{T+1}^2]~~\leq~~2\mathbb E[\overline B_{T+1}^2]
+8\sqrt{
\mathbb E[\overline B_{T+1}^2]
\mathbb E[\overline Z_{T+1}^2]}~~\leq~~34\mathbb E[\overline B_{T+1}^2]
+\frac12\mathbb E[\overline Z_{T+1}^2].\nn
\eeq
Thus $\mathbb E[\overline Z_{T+1}^2]
\leq68\mathbb E[\overline B_{T+1}^2]$.
Another application of Cauchy--Schwarz yields
\begin{equation} \label{eq:MVR2:c:weighted-increment-bound}
\mathbb E[(\overline B_{T+1}+\overline Z_{T+1})^2]
\leq(1+\sqrt{68})^2\mathbb E[\overline B_{T+1}^2]
\leq100\mathbb E[\overline B_{T+1}^2].
\end{equation}

\noindent$\quad\bullet$\enspace\emph{Establishing $\mathbb E[H_T^2]\leq C_h^2\ell_T^2$ and $\mathbb E[H_T]\leq C_h\ell_T$.} We combine the preceding estimates to bound $H_T$.
The preliminary estimate $\|\r_t\|\leq C_rt$ gives
\beq
\ln A_T\leq\ln(1+C_r^2T^3)\leq C_A\ell_T.\nn
\eeq
By~\eqref{eq:MVR2:c:potential-bound},
\beq
0~~\leq~~ H_T\leq\theta^2\left[C_A\ell_T+B_h+8(\overline B_{T+1}+\overline Z_{T+1})\right].\nn
\eeq
Expanding the square and applying Cauchy--Schwarz, together with
\eqref{eq:MVR2:c:innovation-second-moment} and
\eqref{eq:MVR2:c:weighted-increment-bound}, we obtain
\begin{align*}
\sqrt{\mathbb E[H_T^2]}
&\leq\theta^2\left[
C_A\ell_T+B_h
+8\sqrt{\mathbb E[
(\overline B_{T+1}+\overline Z_{T+1})^2]}
\right]\\
&\leq\theta^2\left[
C_A+B_h+80(11 L^2\theta^2+2C_q)
\right]\ell_T ~~=~~ \ddot C_h\ell_T ~~=~~ C_h\ell_T.
\end{align*}
Squaring this inequality and using
$\mathbb E[H_T]\leq\sqrt{\mathbb E[H_T^2]}$, we conclude that
\begin{equation}
\mathbb E[H_T^2]\leq C_h^2\ell_T^2,
\qquad
\mathbb E[H_T]\leq C_h\ell_T.
\label{eq:MVR2:c:nonconvex-conclusion}
\end{equation}

\paragraph{Part (d-i).} Define $\mathcal B_t:=\sum_{i=1}^t\mathcal V_i$ and $\mathcal Z_t:=\sum_{i=1}^t\mathcal E_i$, where $\mathcal B_0=\mathcal Z_0=0$.

The tracking recursion at iteration $t+1$ gives
\[
\alpha_{t+1}P_t
\leq P_t-P_{t+1}+\mathcal V_{t+1}+\mathcal E_{t+1}.
\]
Since $\alpha_t/\alpha_{t+1}\leq5/3<2$, summing over
$t=1,\ldots,T$ and using $P_1=\mathcal V_1$,
$\mathcal E_1=0$, and $P_{T+1}\geq0$, we obtain
\beq
\label{eq:MVR2:log:U-pathwise}
U_T &\leq& \ts 2\sum_{t=1}^T\alpha_{t+1}P_t~~\leq~~2\left[P_1-P_{T+1}
+\sum_{t=2}^{T+1}(\mathcal V_t+\mathcal E_t)\right] \nn\\
&=& \ts \ts 2(\mathcal B_{T+1}+\mathcal Z_{T+1}-P_{T+1})\leq2(\mathcal B_{T+1}+\mathcal Z_{T+1}).
\eeq
Furthermore,
\beq \label{eq:MVR2:log:B-bound}
\ts \mathcal B_{T+1} \leq\sum_{t=1}^{T+1}(\ddot{\mathcal E}_t+2\dot{\mathcal E}_t) = 11 L^2 H_T+2Q_{T+1}.
\eeq
Conditional mean zero gives $\mathbb E[\mathcal Z_{T+1}]=0$.
Taking expectations in~\eqref{eq:MVR2:log:U-pathwise}
and applying part (c) and the bound $\mathbb E[Q_{T+1}]\leq C_v\ell_T$ proved in part (b-i), we obtain the more precise first-moment estimate
\begin{equation} \label{eq:MVR2:log:U-first}
\ts \mathbb E[U_T] \leq2\mathbb E[\mathcal B_{T+1}] \leq2(11 L^2 C_h+2C_v)\ell_T.
\end{equation}

\paragraph{Part (d-ii).}
We again use partial sums and summation by parts, without
assuming a uniform bound on $P_t$.
The exact recursion gives the identity
\beq \label{eq:MVR2:log:nonnegative-sum}
\ts \mathcal B_t+\mathcal Z_t
=P_t+\sum_{i=1}^t(1-\varpi_i^2)P_{i-1}\geq P_t\geq0.
\eeq
By conditional mean zero,
\beq \label{eq:MVR2:log:Z-square}
\ts \mathbb E[\mathcal Z_{T+1}^2]~~=~~\mathbb E\sum_{t=1}^{T+1}\mathcal E_t^2
\leq4\mathbb E\sum_{t=1}^{T+1} P_{t-1}\mathcal V_t~~\leq~~4\mathbb E\sum_{t=1}^{T+1}
(\mathcal B_{t-1}+\mathcal Z_{t-1})\mathcal V_t.
\eeq
The partial-sum identities give
\beq
\ts 2\sum_{t=1}^{T+1}\mathcal B_{t-1}\mathcal V_t
&=& \ts \mathcal B_{T+1}^2-\sum_{t=1}^{T+1}\mathcal V_t^2, \nn \\
\ts \sum_{t=1}^{T+1}\mathcal Z_{t-1}\mathcal V_t
&=& \ts \mathcal B_{T+1}\mathcal Z_{T+1}
-\sum_{t=1}^{T+1}\mathcal B_{t-1}\mathcal E_t
-\sum_{t=1}^{T+1}\mathcal V_t\mathcal E_t. \nn
\eeq
Here $\mathbb E[\mathcal B_{t-1}\mathcal E_t]=0$
because $\mathcal B_{t-1}$ is predictable.
Combining these identities with~\eqref{eq:MVR2:log:Z-square}
and Cauchy--Schwarz gives
\begin{align*}
\ts \mathbb E[\mathcal Z_{T+1}^2]
& \ts \leq2\mathbb E[\mathcal B_{T+1}^2]
+4\mathbb E[\mathcal B_{T+1}\mathcal Z_{T+1}]
-4\mathbb E\sum_{t=1}^{T+1}\mathcal V_t\mathcal E_t\\
&\leq \ts 2\mathbb E[\mathcal B_{T+1}^2]
+8\sqrt{\mathbb E[\mathcal B_{T+1}^2]
                  \mathbb E[\mathcal Z_{T+1}^2]}\\
&\leq \ts 34\mathbb E[\mathcal B_{T+1}^2]
+\frac12\mathbb E[\mathcal Z_{T+1}^2].
\end{align*}
For the second inequality, the same-time product is bounded by
\[
\ts  \left|\mathbb E\sum_{t=1}^{T+1}\mathcal V_t\mathcal E_t\right|
\leq\sqrt{\left(\mathbb E\sum_{t=1}^{T+1}\mathcal V_t^2\right)
                   \left(\mathbb E\sum_{t=1}^{T+1}\mathcal E_t^2\right)}
\leq\sqrt{\mathbb E[\mathcal B_{T+1}^2]
                   \mathbb E[\mathcal Z_{T+1}^2]}.
\]
Consequently,
\begin{equation}
\label{eq:MVR2:log:sum-second}
\mathbb E[\mathcal Z_{T+1}^2]\leq68\mathbb E[\mathcal B_{T+1}^2],
\qquad
\sqrt{\mathbb E[(\mathcal B_{T+1}+\mathcal Z_{T+1})^2]}
\leq10\sqrt{\mathbb E[\mathcal B_{T+1}^2]}.
\end{equation}
Finally,~\eqref{eq:MVR2:log:B-bound}, part (c), and {the bound $\mathbb E[Q_{T+1}^2]\leq C_q^2\ell_T^2$ in~\eqref{eq:MVR2:noise-second-moment}} imply
\[
\sqrt{\mathbb E[\mathcal B_{T+1}^2]}
\leq 11 L^2\sqrt{\mathbb E[H_T^2]}
    +2\sqrt{\mathbb E[Q_{T+1}^2]}
\leq(11 L^2 C_h+2C_q)\ell_T.
\]
Since the right-hand side of~\eqref{eq:MVR2:log:U-pathwise}
is nonnegative, squaring that inequality and invoking
\eqref{eq:MVR2:log:sum-second} yields
\[
\sqrt{\mathbb E[U_T^2]}
\leq2\sqrt{\mathbb E[(\mathcal B_{T+1}+\mathcal Z_{T+1})^2]}
\leq20(11 L^2 C_h+2C_q)\ell_T
=C_u\ell_T.
\]
Together with~\eqref{eq:MVR2:log:U-first}, or simply
$\mathbb E[U_T]\leq\sqrt{\mathbb E[U_T^2]}$, this proves
both claims in part (d).

\end{proof}

With the first two moments of the tracking and normalized motion sums available, we can bound the descent energy. The following pathwise comparison then transfers this control to the residual
accumulator with the weight appropriate to MVR2.

\begin{lemma}
\label{lemma:MVR2:bridge}
(Descent energy and motion bridge for MVR2) Let $\Delta:=F(\x_1)-F_\star$ and $\ell_T:=\ln(e+T)$. Let $C_h$ and $C_u$ be defined in Lemma~\ref{lemma:MVR2:logTerm}. Then, for every integer $T\geq1$,
\begin{enumerate}[label=\textbf{(\alph*)},leftmargin=18pt,
itemsep=2pt,topsep=1pt]
\item Define $W_T:=\sum_{t=1}^T \gamma_t \|\r_t\|^2$. We have $\mathbb E[W_T]\leq C_w\ell_T$ and
$\mathbb E[W_T^2]\leq C_w^2\ell_T^2$,
where $C_w:=2\Delta+C_u+(1+L)C_h$.
\item $\mathbb E[\sqrt{\alpha_T}M_T]\leq C_m\ell_T^2$,
where $C_m:=C_w/\theta+C_w^2/(4\theta^2)$.
\end{enumerate}
\end{lemma}

\begin{proof}

Define $W_T:=\sum_{t=1}^T \gamma_t \|\r_t\|^2$. Let $L_+=L+1$.

\paragraph{Part (a-i).} Summing the first inequality of Lemma~\ref{lemma:MVR2:descent}
over $t=1,\ldots,T$ gives {$\frac12W_T\leq F(\x_1)-F(\x_{T+1})+\frac12U_T+\frac{1+L}{2}H_T$}. Since $F(\x_{T+1})\geq F_\star$, we obtain
\beq
0\leq W_T\leq2\Delta+U_T+L_+ H_T.
\label{eq:MVR2:bridge:pathwise-descent}
\eeq
Taking expectations and applying {Lemma~\ref{lemma:MVR2:logTerm}(c),(d)}, we have
\begin{align*}
\mathbb E[W_T]
~~\leq~~2\Delta+C_u\ell_T+L_+C_h\ell_T~~\leq~~\bigl(2\Delta+C_u+L_+ C_h\bigr)\ell_T~~
=~~C_w\ell_T.
\end{align*}

\paragraph{Part (a-ii).} Squaring~\eqref{eq:MVR2:bridge:pathwise-descent} gives
\beq
\mathbb E[W_T^2]
&\leq& 4\Delta^2+\mathbb E[U_T^2]+L_+^2\mathbb E[H_T^2] +4\Delta\mathbb E[U_T]
 +4\Delta L_+\mathbb E[H_T]
 +2L_+\mathbb E[U_TH_T]\nn\\
&\xxxx{1}{\leq}&4\Delta^2 + C_u^2\ell_T^2 + L_+^2 C_h^2 \ell_T^2 + 4\Delta C_u \ell_T
 +4\Delta L_+ C_h \ell_T + 2L_+\mathbb E[U_TH_T] \nn\\
&\xxxx{2}{\leq}&4\Delta^2 + C_u^2\ell_T^2 + L_+^2 C_h^2 \ell_T^2 + 4\Delta C_u \ell_T
 +4\Delta L_+ C_h \ell_T + 2L_+ C_u C_h \ell_T^2\nn\\
&=& (2\Delta + C_u \ell_T + L_+ C_h\ell_T )^2 \leq C_w^2 \ell_T^2, \nn
\eeq
\noi where step \step{1} uses Lemma \ref{lemma:MVR2:logTerm}; step \step{2} uses Cauchy--Schwarz that $\mathbb E[U_TH_T]
\leq\sqrt{\mathbb E[U_T^2]}\sqrt{\mathbb E[H_T^2]}
\leq C_uC_h\ell_T^2$, where the last step uses Lemma \ref{lemma:MVR2:logTerm}.

\paragraph{Part (b).} We obtain
\beq \label{eq:MVR2:bridge:root-motion}
\ts W_T &:=& \ts \sum_{t=1}^T \gamma_t \|\r_t\|^2 ~~\xxxx{1}{=}~~ \ts \theta\sum_{t=1}^T\alpha_t^{1/4}     \frac{A_t-A_{t-1}}{\sqrt{A_{t-1}}}~~\geq~~ \theta\alpha_T^{1/4}\sum_{t=1}^T\frac{A_t-A_{t-1}}{\sqrt{A_{t-1}}}\notag\\
&\xxxx{2}{\geq} &\ts 2\theta\alpha_T^{1/4}\sum_{t=1}^T
           (\sqrt{A_t}-\sqrt{A_{t-1}})~~=~~ 2\theta\alpha_T^{1/4}(\sqrt{1+M_T}-1),
\eeq
\noi where step \step{1} uses $\gamma_t=\theta\alpha_t^{1/4}/\sqrt{A_{t-1}}$ and $\|\r_t\|^2=A_t-A_{t-1}$; step \step{2} uses the fact that $\frac{A_t-A_{t-1}}{\sqrt{A_{t-1}}}
=(\sqrt{A_t}-\sqrt{A_{t-1}})
 \left(1+\sqrt{\frac{A_t}{A_{t-1}}}\right)
\geq2(\sqrt{A_t}-\sqrt{A_{t-1}})$.

We have from Inequality \eqref{eq:MVR2:bridge:root-motion}: $W_T + 2\theta \alpha_T^{1/4}\geq 2\theta \alpha_T^{1/4} \sqrt{1+M_T}$. Squaring both
nonnegative sides and dividing both sides by $4\theta^2$ gives $\tfrac{1}{4\theta^2}W^2_T + \alpha_T^{1/2} + \tfrac{1}{\theta} W_T\alpha_T^{1/4} \geq  \alpha_T^{1/2} (1+M_T)$, leading to:
\beq
\ts \sqrt{\alpha_T}M_T~~\leq~~\frac{\alpha_T^{1/4}}{\theta}W_T+\frac{W_T^2}{4\theta^2}
~~\leq~~\frac{W_T}{\theta}+\frac{W_T^2}{4\theta^2}.\nn
\eeq

Taking expectations and applying Part (a), we conclude that
\beq
\ts \mathbb E[\sqrt{\alpha_T}M_T]~~\leq~~\frac{C_w}{\theta}\ell_T+\frac{C_w^2}{4\theta^2}\ell_T^2~~\leq~~\left(\frac{C_w}{\theta}+\frac{C_w^2}{4\theta^2}\right)\ell_T^2~~=~~C_m\ell_T^2.\nn
\eeq

\end{proof}

We now assemble the MVR2 stationarity guarantee. Unlike MVR1, its discrepancy compares two samples at the same point, so the inverse-weight estimate depends directly on the noise variance;
combining this estimate with the motion and tracking bounds yields the stated rate.

\begin{theorem}\label{app:theorem:MVR2}
(Noise-adaptive composite stationarity of {\sf ALFPG-MVR2}) Suppose Assumptions \ref{ass:prox}, \ref{ass:oracle}, \ref{ass:bounded-oracle},
\ref{ass:smooth}, \ref{ass:regularizer}, and~\ref{ass:individual-smooth}
hold, and run {\sf ALFPG-MVR2}. Write $\ell_T:=\ln(e+T)$. Then there exists a constant $C_R>0$, independent of $T$ and $\sigma$, such that, for every $T\geq1$,
\beq\label{eq:main-rate}
\ts \frac1T
\sum_{t=1}^T\mathbb E[\mathcal{R}(\x_{t+1})] \leq C_R  \left(\frac{\ell_T}{T^{1/2}}+\frac{\sigma^{1/3}\ell_T}{T^{1/3}}\right).
\eeq
Consequently, if $R\sim{\rm Unif}\{1,\ldots,T\}$
is independent of the algorithmic randomness, then
$\mathbb E[\mathcal{R}(\x_{R+1})]
=\widetilde{\mathcal O}
(T^{-1/2}+\sigma^{1/3}T^{-1/3})$. Since {\sf ALFPG-MVR2} uses at most two stochastic-gradient
evaluations per iteration, its stochastic-gradient oracle complexity for obtaining
$\mathbb E[\mathcal{R}(\x_{R+1})]\leq\epsilon$ is $\widetilde{\mathcal O}(\epsilon^{-2}+\sigma\epsilon^{-3})$.

\end{theorem}

\begin{proof}

Define $\ell_T:=\ln(e+T)$, and $\mathcal{R}(\x_t):=\operatorname{dist} \bigl(0,\nabla f(\x_t)+\partial h(\x_t)\bigr)$.

Define $D_t=\sum_{i=1}^t\chi_i^2$ and $N_t=\max_{1\leq i\leq t}\chi_i^2$, where $\chi_1=0$.

Define the deterministic quantity $\mathcal B_T:=1+2(T-1)\sigma^2$.

Let $C_u$ and $C_m$ be the constants in
Lemmas~\ref{lemma:MVR2:logTerm} and~\ref{lemma:MVR2:bridge}, and set $C_R:=2(1+L\theta)\sqrt{C_m} +2\sqrt{C_u}$.

\paragraph{Bounding the term $\mathbb E[\alpha_T^{-1/2}]$.} For $T\geq2$, the following running-maximum expression is nonempty. Using the initial discrepancy $\chi_1=0$, we derive:
\beq
\ts \alpha_T^{-1/2} ~~=~~ \big(\max_{1\leq k\leq T-1}\frac{1+D_k}{1+N_k}\big)^{1/3} ~~\leq~~ \big(1+D_{T-1} \big)^{1/3} ~~\leq~~ \big(1 + \sum_{i=1}^{T-1} \chi_i^2 \big)^{1/3}.\nn
\eeq
Taking the expectations yields:
\beq \label{eq:MVR2:inverse-moment}
\ts \EEE[\alpha_T^{-1/2}] ~~\leq~~ \EEE[ \big(1 + \sum_{i=1}^{T-1} \chi_i^2 \big)^{1/3} ] ~~\xxxx{1}{\leq}~~ \big(1 + 2 (T-1)\sigma^2\big)^{1/3} ~~= ~~\mathcal B_T^{1/3},
\eeq
\noi where step \step{1} uses $\mathbb E[\chi_t^2]
=\mathbb E\|\epsilons_{t-1}\|^2
 +\mathbb E\|\epsilons_t^-\|^2
\leq2\sigma^2$. For $T=1$, Inequality (\ref{eq:MVR2:inverse-moment}) holds directly.

\paragraph{Bounding the composite residual.} We derive:
\beq \label{eq:MVR2:residual}
\ts \mathcal{R}(\x_{t+1}) &:=& \ts \operatorname{dist}
 \bigl(0,\nabla f(\x_{t+1})+\partial h(\x_{t+1})\bigr)\nn\\
&\xxxx{1}{\leq}& \ts \|\nabla f(\x_{t+1})+\r_t-\m_t\| ~~\xxxx{2}{\leq}~~ \|\r_t\|+\|\s_t\| +\|\nabla f(\x_{t+1})-\nabla f(\x_t)\| \nn\\
&\xxxx{3}{\leq}& \ts (1+L\gamma_t)\|\r_t\|+\|\s_t\|~~\xxxx{4}{\leq}~~(1+L\theta)\|\r_t\|+\|\s_t\|,
\eeq
\noi where step \step{1} uses the global optimality of the proximal update that $\r_t-\m_t\in\partial h(\x_{t+1})$; step \step{2} uses $\s_t=\m_t-\nabla f(\x_t)$; step \step{3} uses $L$-smoothness of $f(\cdot)$; step \step{4} uses $\gamma_t=\theta\alpha_t^{1/4}\beta_t\leq\theta$.

\paragraph{Removing the adaptive weights.} We obtain
\beq \label{eq:MVR2:sqrtMT}
\mathbb E[\sqrt{M_T}]~~=~~\mathbb E\big[\sqrt{\sqrt{\alpha_T}M_T}\,\alpha_T^{-1/4}\big]~~\xxxx{1}{\leq}~~\sqrt{\mathbb E[\sqrt{\alpha_T}M_T]}\sqrt{\mathbb E[\alpha_T^{-1/2}]}~~\leq~~\sqrt{C_m}\ell_T\mathcal B_T^{1/6},
\eeq
\noi where step \step{1} uses Cauchy--Schwarz inequality; step \step{2} uses Lemma~\ref{lemma:MVR2:bridge}(b), and Inequality (\ref{eq:MVR2:inverse-moment}).

Another application of Cauchy--Schwarz and Lemma~\ref{lemma:MVR2:logTerm}(d) gives:
\beq \label{eq:MVR2:sqrtST}
\ts \mathbb E[\sqrt{S_T}]~~\xxxx{1}{\leq}~~\mathbb E[\sqrt{U_T}\,\alpha_T^{-1/4}]~~\xxxx{2}{\leq}~~\sqrt{\mathbb E[U_T]}\sqrt{\mathbb E[\alpha_T^{-1/2}]}~~\xxxx{3}{\leq}~~\sqrt{C_u \ell_T} \mathcal B_T^{1/6},
\eeq
\noi where step \step{1} uses $\sqrt{\alpha_T}S_T\leq\sum_{t=1}^T\sqrt{\alpha_t}\|\s_t\|^2:=U_T$ since $\alpha_t\geq\alpha_T$ for $t\leq T$; step \step{2} uses Cauchy--Schwarz inequality; step \step{3} uses Inequality (\ref{eq:MVR2:inverse-moment}).

\paragraph{Establishing the noise-adaptive rate.} We derive:
\beq
\ts && \ts\frac1T\sum_{t=1}^T\mathbb E[\mathcal{R}(\x_{t+1})]~~\xxxx{1}{\leq}~~\frac{1+L\theta}{T}\mathbb E\sum_{t=1}^T\|\r_t\|+\frac1T\mathbb E\sum_{t=1}^T\|\s_t\| \nn \\
&\leq& \ts \frac{1+L\theta}{\sqrt T}\mathbb E[\sqrt{M_T}]+\frac1{\sqrt T}\mathbb E[\sqrt{S_T}]~~\xxxx{2}{\leq}~~\frac{(1+L\theta)\sqrt{C_m}\ell_T+\sqrt{C_u}\ell_T^{1/2}}{\sqrt T} \mathcal B_T^{1/6}\nn\\
&\xxxx{3}{\leq}& \ts \underbrace{\big(2(1+L\theta)\sqrt{C_m} +2\sqrt{C_u}\big)}_{:=C_R}\cdot \frac{ \ell_T}{2\sqrt T} (1+2 T^{1/6}\sigma^{1/3})~~\leq~~\ts C_R \big( \frac{\ell_T}{\sqrt T} +  \frac{\ell_T}{T^{1/3}}\sigma^{1/3}\big),\nn
\eeq
\noi where step \step{1} uses Inequality (\ref{eq:MVR2:residual}); step \step{2} uses Inequalities (\ref{eq:MVR2:sqrtMT}) and (\ref{eq:MVR2:sqrtST}); step \step{3} uses $\mathcal B_T^{1/6} \leq1+\bigl(2(T-1)\sigma^2\bigr)^{1/6}
\leq1+2^{1/6}T^{1/6}\sigma^{1/3} \leq 1+2 T^{1/6}\sigma^{1/3}$.
\end{proof}

\section{Experiments (Further Details)}
\label{app:sect:exp}

This section evaluates the proposed {\sf ALFPG} methods on three applications: sparse phase retrieval (SPR), linear eigenvalue problems (LEV), and neural network regression (NNR).

$\blacktriangleright$ \textbf{General Experimental Settings.} All methods are implemented in MATLAB and evaluated on a machine with an Intel 2.6 GHz CPU and 64 GB of RAM. Within each comparison, all methods use
the same initialization and nominal time budget. We report the shifted
best-so-far objective $\min_{0\leq j\leq t}F(\x_j)-F_{\rm best}+10^{-4}$, where $F_{\rm best}$ is the smallest objective value recorded across all compared methods on the same instance. The positive offset enables logarithmic-scale visualization. The figures plot the recorded objective traces against elapsed time.

$\blacktriangleright$ \textbf{Datasets}. We construct experimental instances from seven public benchmark
datasets (TDT2, 20Newsgroups, Sector, MNIST, CIFAR, Gisette, Caltech, and E2006) and from synthetic Gaussian data. For each public dataset, we select $\dot m$ samples and $\dot n$ features uniformly without replacement
from the stored feature matrix. For a synthetic instance, we generate
a $\dot m\times\dot n$ matrix with independent standard normal entries.
In both cases, denoting the resulting matrix by $\A_{\rm raw}$, we rescale it as $\A=\frac{\A_{\rm raw}}{\|\A_{\rm raw}\|}$, $\A\in\mathbb{R}^{\dot m\times\dot n}$. Subfigure labels use the format
\texttt{dataset-$\dot m$-$\dot n$}. For example,
\texttt{TDT2-$\dot m$-$\dot n$} identifies a submatrix of TDT2 with
$\dot m$ samples and $\dot n$ features, while
\texttt{Gaussian-$\dot m$-$\dot n$} identifies a synthetic instance of
the same size. In the problem formulations, we write $N=\dot m$ and
$n=\dot n$, and use $\a_i^\top$ for the $i$th row of $\A$.
The three tasks use these matrices as inputs and construct their
observations or targets as described below.

\subsection{Problem Formulations}

$\blacktriangleright$ \textbf{Sparse Phase Retrieval (SPR)}.
Sparse phase retrieval aims to recover a sparse signal from squared linear
measurements. For each sensing matrix $\A\in\mathbb{R}^{N\times n}$, let
$\a_i^\top$ denote its $i$th row. We generate a ground-truth signal
$\x_\star\in\mathbb{R}^n$ by choosing its support $\mathcal{S}$ uniformly
from all $k$-element subsets of $\{1,\ldots,n\}$, where
$k=\operatorname{round}(0.1n)$. We draw the entries on $\mathcal{S}$
independently from $\mathcal{N}(0,10^2)$ and set the remaining entries to
zero. We form the noiseless intensities $\u$ and the noisy observations
$\b$ as
\beq
\ts \u = (\A\x_\star)^{\odot 2},
\qquad
\b = \u + 10^{-3}\|\u\|_2\boldsymbol{\varepsilon},
\qquad
\boldsymbol{\varepsilon}\sim\mathcal{N}(\zero,\I_N),
\nn
\eeq
where $(\cdot)^{\odot 2}$ denotes entrywise squaring. Conditional on $\u$,
each observation has noise standard deviation $10^{-3}\|\u\|_2$. We minimize $F(\x)=f(\x)+h(\x)$, where
\beq
\ts f(\x)=\frac{1}{N}\sum_{i=1}^{N}
      \bigl((\a_i^\top\x)^2-{[\b]_i}\bigr)^2,
\qquad
h(\x)=\dot\lambda\sum_{j=1}^{n}\min\{\dot\delta,|{[\x]_j}|\}
      +\iota_{\{\|\x\|_\infty\leq\dot r\}}(\x).
\nn
\eeq
The first term in $h$ is the capped-$\ell_1$ penalty, with weight
$\dot\lambda>0$ and cap $\dot\delta>0$. The indicator imposes the box
constraint with radius $\dot r>0$. The Hessian of the data-fitting loss is $\nabla^2 f(\x)
=
\frac{4}{N}\A^\top
\operatorname{diag}\!\bigl(3(\A\x)^{\odot 2}-\b\bigr)\A$. Its spectral norm is bounded on the feasible box, so $\nabla f$ is
Lipschitz continuous there.

$\blacktriangleright$ \textbf{Linear Eigenvalue Problems (LEV)}. We use the preprocessed feature matrix $\A$ to construct
$\C=-\A^\top\A\preceq 0$. Each feature vector $\a_i$ serves as both
the input and its reconstruction target. For a matrix
$\V\in\mathbb{R}^{n\times\dot r}$ with orthonormal columns,
$\V\V^\top\a_i$ is the orthogonal projection of $\a_i$ onto the
column space of $\V$. The deterministic experiment solves
\beq
\min_{\V\in\mathbb{R}^{n\times\dot r},~~\V^\top\V=\I_{\dot r}}
    \operatorname{tr}(\V^\top\C\V)+|\operatorname{tr}(\C)|. \nn
\eeq
The columns of a minimizer span a subspace associated with the $\dot r$ smallest eigenvalues of $\C$, counted with multiplicity. For every feasible $\V$, the objective is the total squared reconstruction
error $\operatorname{tr}(\V^\top\C\V)+|\operatorname{tr}(\C)| = \|\A\V\V^\top-\A\|_{\fro}^2 =\sum_{i=1}^{N}
      \|\V\V^\top\a_i-\a_i\|_2^2$. In the expectation setting, we sample rows uniformly with replacement
and use the corresponding per-sample reconstruction losses.
The reported objective is their mean over the full data matrix,
$N^{-1}\|\A\V\V^\top-\A\|_{\fro}^2$.

$\blacktriangleright$ \textbf{Three-Layer Neural Network Regression (NNR)}. We use the benchmark feature matrices to construct supervised regression instances. A fixed teacher network then assigns one scalar regression target to each normalized input. The teacher and student have the
same architecture, with two $\tanh$ hidden layers of widths $16$ and
$8$ and a scalar linear output.

Let $\W_1\in\mathbb{R}^{n\times16}$,
$\W_2\in\mathbb{R}^{16\times8}$, and
$\W_3\in\mathbb{R}^{8\times1}$ be the weight matrices, and let
$\b_1\in\mathbb{R}^{16}$, $\b_2\in\mathbb{R}^{8}$, and
$b_3\in\mathbb{R}$ be the biases. With $\w$ denoting the vector of all
weights and biases, the forward map is
\beq
\h_{i,1} = \tanh\!\bigl(\W_1^\top\widetilde{\a}_i+\b_1\bigr),\quad  \h_{i,2} =\tanh\!\bigl(\W_2^\top\h_{i,1}+\b_2\bigr),\quad  \phi_{\w}(\widetilde{\a}_i) =\W_3^\top\h_{i,2}+b_3,\nn
\eeq
where $\tanh$ is applied entrywise. The network has $16n+161$
trainable parameters, including biases.

To generate the teacher parameters $\w_\star$, we first draw the entries of its first-layer weight matrix independently from $\mathcal{N}(0,3^2)$. We then apply an independent Bernoulli mask with
retention probability $0.1$, setting the masked entries to zero.
The second- and third-layer weights are drawn independently from
$\mathcal{N}(0,1/16)$ and $\mathcal{N}(0,1/8)$, respectively. We evaluate this fixed network on the
normalized inputs to obtain the noise-free scalar labels $y_i=\phi_{\w_\star}(\widetilde{\a}_i),~i=1,\ldots,N$. The resulting pairs $(\widetilde{\a}_i,y_i)$ define the
regression instance. The finite-sum and expectation experiments use
the same empirical objective
\beq
\ts \min_{\|\w\|_\infty\leq \dot{r}}
    \frac{1}{2N}\sum_{i=1}^{N}
    \bigl(\phi_{\w}(\widetilde{\a}_i)-y_i\bigr)^2
    +\dot\lambda\|\w\|_1,\nn
\eeq
where $\dot\lambda\in\{0.01,0.001\}$ and $\dot{r} = 10$.

For the student network, we draw independent zero-mean Gaussian weights
with standard deviations $1/\sqrt n$, $1/\sqrt{16}$, and $1/\sqrt8$
in the first, second, and third layers, respectively. 

\subsection{Proximal Operators}

$\blacktriangleright$ \textbf{Projection onto Orthogonality Constraints}. When $h(\x)=\iota_{\mathcal{M}}(\mat(\x))$ with
$\mathcal{M}\triangleq
\{\V\in\mathbb{R}^{d'\times r'}\,|\,\V\trans\V=\I_{r'}\}$
and $d'\geq r'$, the computation of the generalized proximal
operator with a scalar metric $\mu\I$, $\mu>0$, reduces to
solving the following optimization problem:
\beq
\textstyle
\bar{\x}\in\arg\min_{\x}
\frac{\mu}{2}\|\x-\x'\|_2^2,
\quad\st\quad \mat(\x)\in\mathcal{M}.\nn
\eeq
\noi This corresponds to finding the nearest matrix with
orthonormal columns. A global optimal solution is given by
$\bar{\x}=\vecc(\hat{\mathbf{U}}\hat{\mathbf{V}}\trans)$,
where
$\mat(\x')=\hat{\mathbf{U}}\Diag(\mathbf{s})
\hat{\mathbf{V}}\trans$
is a thin singular value decomposition (SVD) of $\mat(\x')$.
Here, $\vecc(\mathbf{V})\in\mathbb{R}^{d'r'}$ denotes the
vector formed by stacking the columns of
$\mathbf{V}\in\mathbb{R}^{d'\times r'}$, and $\mat(\x)$
reshapes $\x\in\mathbb{R}^{d'r'}$ into a matrix in
$\mathbb{R}^{d'\times r'}$, satisfying
$\mat(\vecc(\mathbf{V}))=\mathbf{V}$.

$\blacktriangleright$ \textbf{Generalized Proximal Operator for the
Capped-$\ell_1$ Penalty}. When
$h(\x)=\dot{\lambda}\|\min(|\x|,\dot{\tau})\|_1
+\iota_{\Omega}(\x)$,
where
$\Omega\triangleq\{\x\mid\|\x\|_\infty\leq\dot r\}$
and $\dot{\lambda},\dot{\tau},\dot r>0$,
the generalized proximal operator reduces to solving
the following nonconvex optimization problem:
\beq
\textstyle
\bar{\x}\in\arg\min_{\x\in\mathbb{R}^n}
\dot{\lambda}\|\min(|\x|,\dot{\tau})\|_1
+\frac12\|\x-\a\|_{\mathbf c}^2,
\quad\st\quad -\dot r\leq\x\leq\dot r.\nn
\eeq
Here,
$\|\z\|_{\mathbf c}^2=\sum_{i=1}^n{[\mathbf c]_i} [\z]_i^2$
with ${[\mathbf c]_i}>0$.
This problem decomposes into $n$ independent subproblems:
\beq
\textstyle
{[\bar{\x}]_i}\in
\arg\min_{-\dot r\leq x\leq\dot r}q_i(x),
\qquad
q_i(x)\triangleq
\frac{{[\mathbf c]_i}}{2}(x-{[\a]_i})^2
+\dot{\lambda}\min(|x|,\dot{\tau}).
\label{eq:one:dim:q_i}
\eeq

\noi To simplify, we write
$q_i(x)=\min\{\dot{q}_i(x),\ddot{q}_i(x)\}$, where
\beq
\textstyle
\dot{q}_i(x)\triangleq
\frac{{[\mathbf c]_i}}{2}(x-{[\a]_i})^2+\dot{\lambda}|x|,
\qquad
\ddot{q}_i(x)\triangleq
\frac{{[\mathbf c]_i}}{2}(x-{[\a]_i})^2
+\dot{\lambda}\dot{\tau}.\nn
\eeq
{The two coordinate-wise minimizers define candidate vectors
$\dot{\x}$ and $\ddot{\x}$ through}
\beq
{[\dot{\x}]_i}~=~\ts \operatorname{sign}({[\a]_i})
  \min\big(\dot r,\,
  \max(0,|{[\a]_i}|-\frac{\dot{\lambda}}{{[\mathbf c]_i}})
  \big), \quad {[\ddot{\x}]_i} ~=~ \ts \max(-\dot r,\min(\dot r,{[\a]_i})). \nn
\eeq
Indeed, the first candidate is obtained by soft-thresholding
followed by projection onto the box, while the second is
the projection of ${[\a]_i}$ onto the box. Since $\textstyle
\min_{|x|\leq\dot r}q_i(x)
=
\min\big\{
\min_{|x|\leq\dot r}\dot{q}_i(x),\,
\min_{|x|\leq\dot r}\ddot{q}_i(x)
\big\}$, a global optimal solution of
Problem~(\ref{eq:one:dim:q_i}) is given by
\beq
{[\bar{\x}]_i}=
\begin{cases}
{[\dot{\x}]_i}, &
\dot{q}_i({[\dot{\x}]_i})
\leq \ddot{q}_i({[\ddot{\x}]_i}),\\
{[\ddot{\x}]_i}, &
\rm{else}.
\end{cases}\nn
\eeq
Thus, the generalized proximal operator is computed exactly
by comparing two candidates for each coordinate. 
 
\subsection{Results of {\sf ALFPG-FS} on SPR and NNR}
\label{app:exp:phase}

$\blacktriangleright$ \textbf{Experimental Settings}.
For both SPR and NNR, we set the box radius to $\dot r=10$ and
test $\dot\lambda\in\{0.01,0.001\}$. For SPR, we set the capped-$\ell_1$
threshold to $\dot\delta=0.1$. All methods start from the same Gaussian
initialization on each instance. We set $\widehat L=10\|\A\|$ and use
the stepsize $1/\widehat L$ for {\sf ProxSARAH}, {\sf SpiderBoost},
{\sf SpiderBoost-M}, and {\sf SGP-SPIDER}. For the quartic SPR loss,
$\widehat L$ serves as a numerical scale for stepsize selection.
For {\sf ALFPG-FS}, we initialize the diagonal metric as
$\v_0=(\widehat L/2)\mathbf{1}$ and use
$(\alpha,\beta)=(0.01,1)$ with $\theta\in\{0,0.1,0.5,0.9\}$.
The mini-batch size and the refresh period are both set to
$\operatorname{round}(\sqrt N)$.

\begin{figure*}[!t]
\centering
\begin{subfigure}{.24\textwidth}\centering\includegraphics[width=\linewidth]{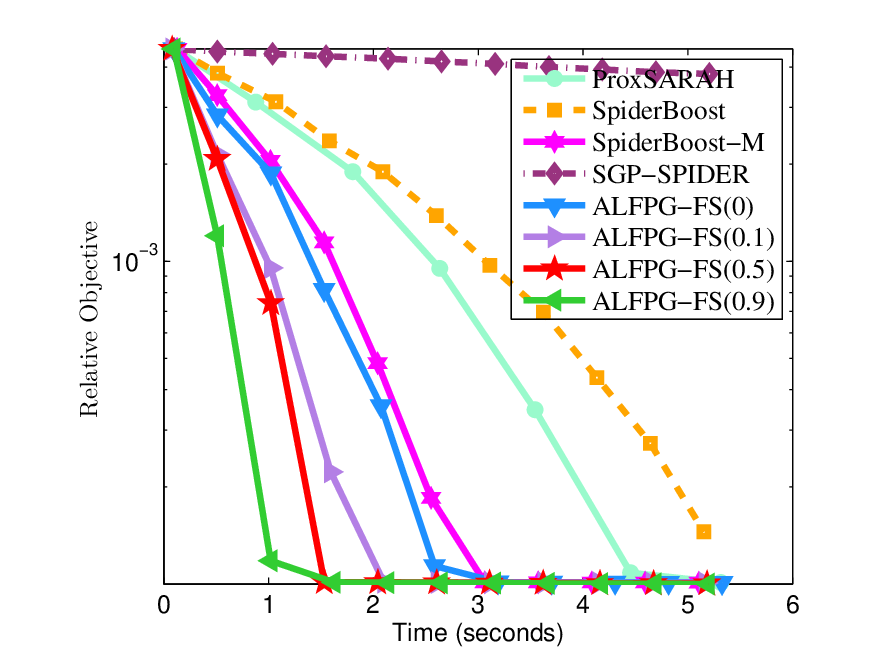}\caption{{\scriptsize TDT2{-9000-1000}}}\label{fig:fs:phase:1:1}\end{subfigure}
\begin{subfigure}{.24\textwidth}\centering\includegraphics[width=\linewidth]{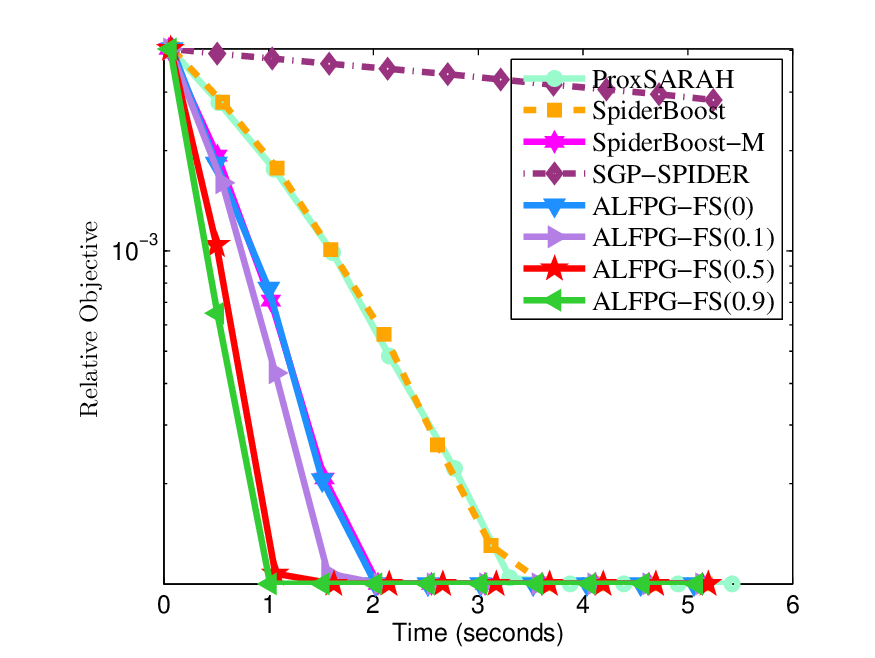}\caption{{\scriptsize 20Newsgroups{-10000-1000}}}\label{fig:fs:phase:1:2}\end{subfigure}
\begin{subfigure}{.24\textwidth}\centering\includegraphics[width=\linewidth]{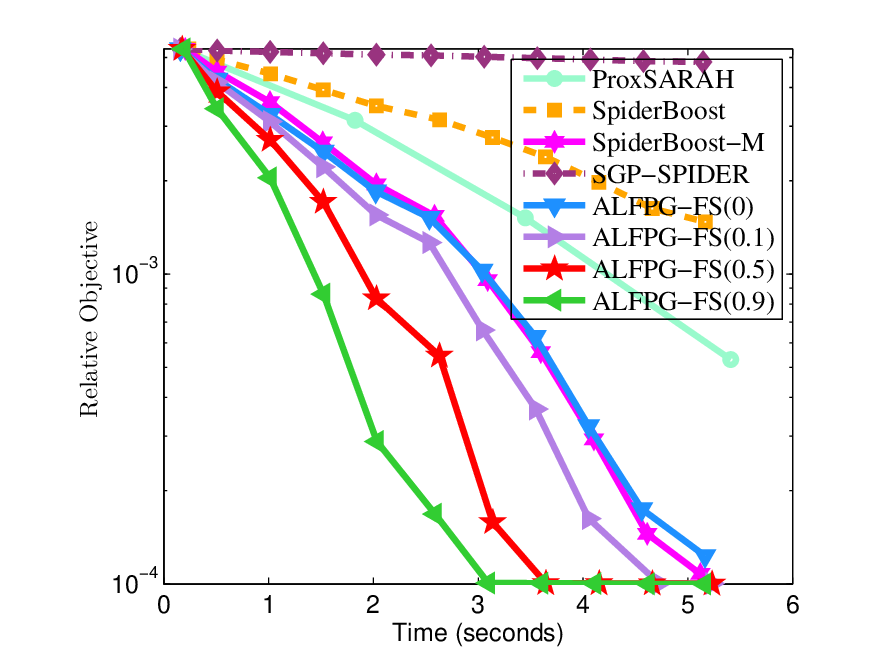}\caption{{\scriptsize Sector{-6412-1000}}}\label{fig:fs:phase:1:3}\end{subfigure}
\begin{subfigure}{.24\textwidth}\centering\includegraphics[width=\linewidth]{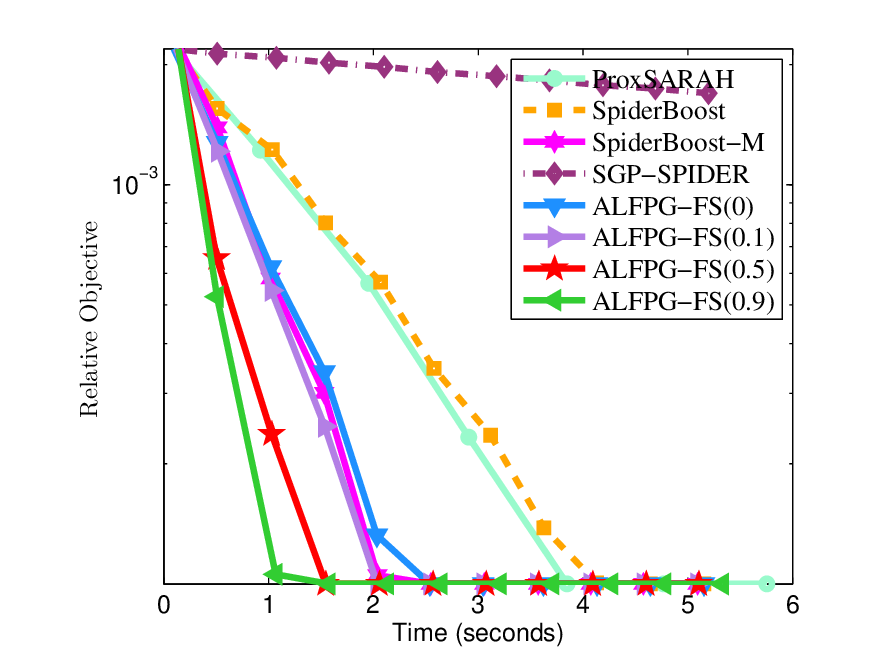}\caption{{\scriptsize MNIST{-2000-784}}}\label{fig:fs:phase:1:4}\end{subfigure}

\begin{subfigure}{.24\textwidth}\centering\includegraphics[width=\linewidth]{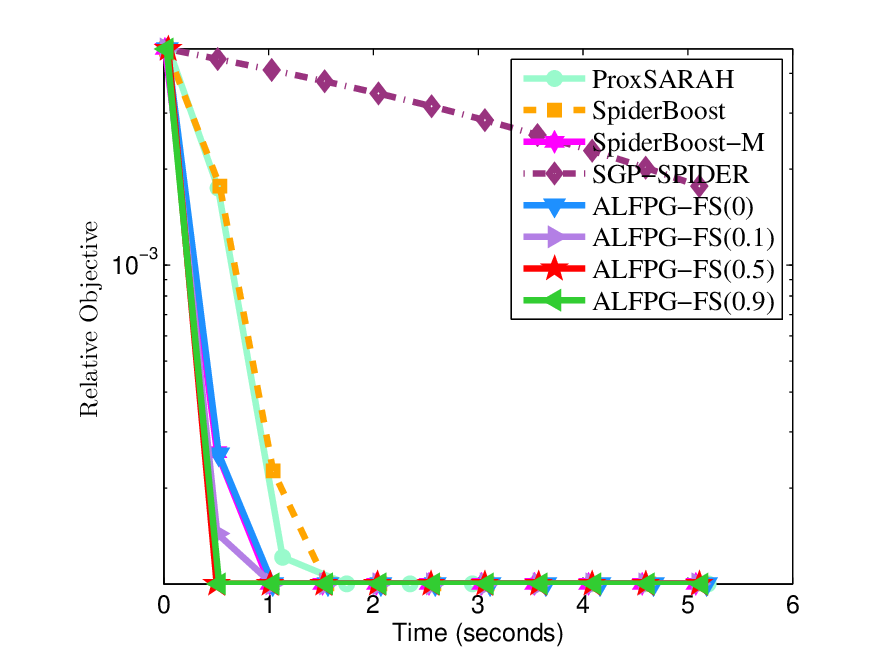}\caption{{\scriptsize CIFAR{-10000-1000}}}\label{fig:fs:phase:3:5}\end{subfigure}
\begin{subfigure}{.24\textwidth}\centering\includegraphics[width=\linewidth]{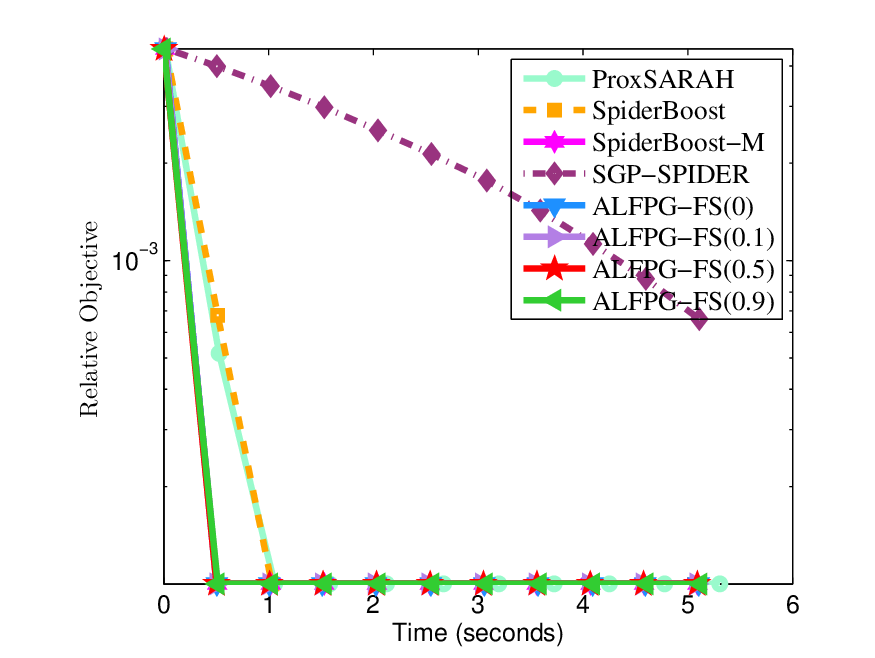}\caption{{\scriptsize Gisette{-5000-1000}}}\label{fig:fs:phase:3:6}\end{subfigure}
\begin{subfigure}{.24\textwidth}\centering\includegraphics[width=\linewidth]{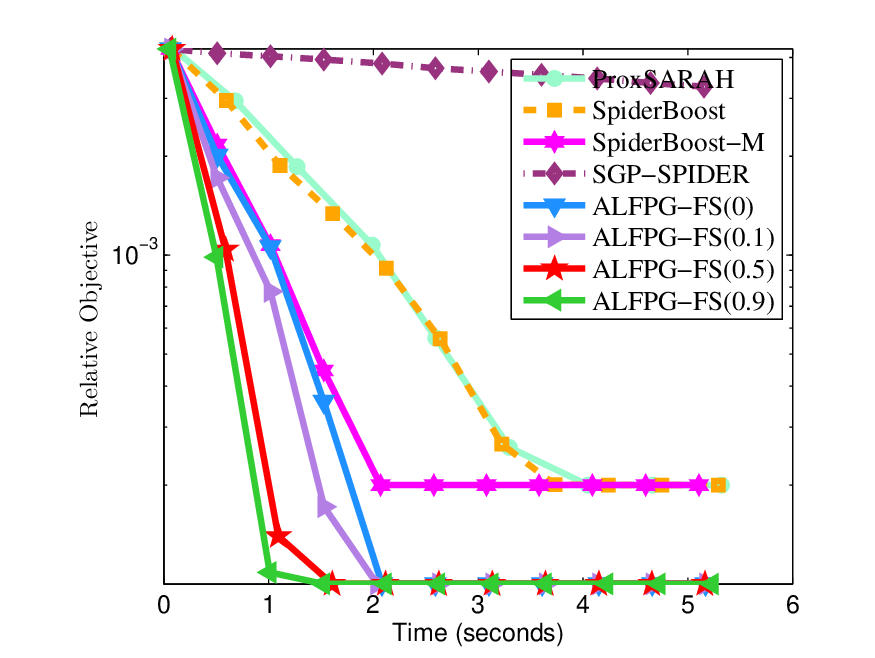}\caption{{\scriptsize Caltech{-3000-1000}}}\label{fig:fs:phase:3:7}\end{subfigure}
\begin{subfigure}{.24\textwidth}\centering\includegraphics[width=\linewidth]{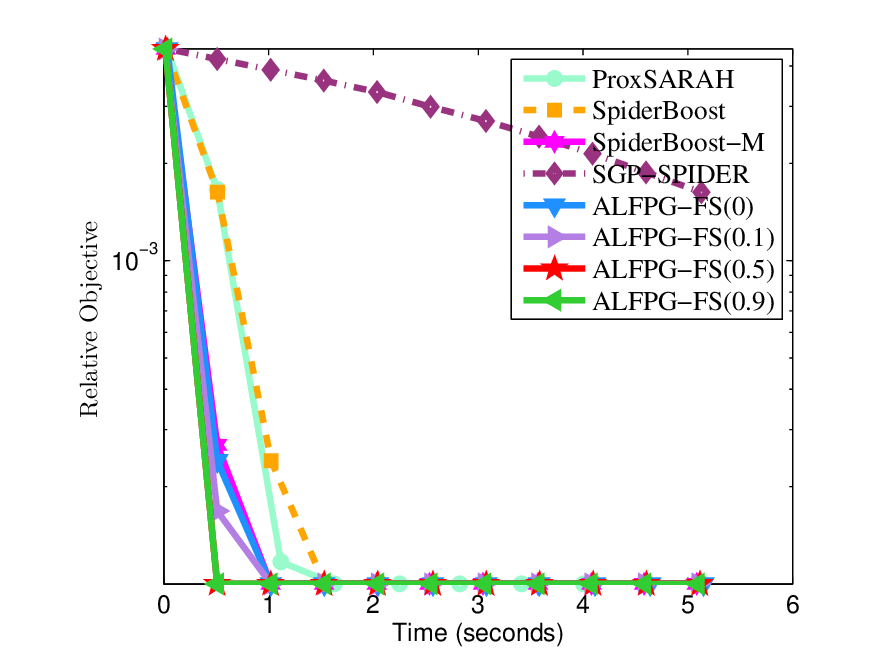}\caption{{\scriptsize Gaussian{-10000-1000}}}\label{fig:fs:phase:3:8}\end{subfigure}

\caption{Sparse phase retrieval with $\dot\lambda=0.01$ in the finite-sum setting.}
\label{fig:fs:phase:1}
\end{figure*}

\begin{figure*}[!h]
\centering
\begin{subfigure}{.24\textwidth}\centering\includegraphics[width=\linewidth]{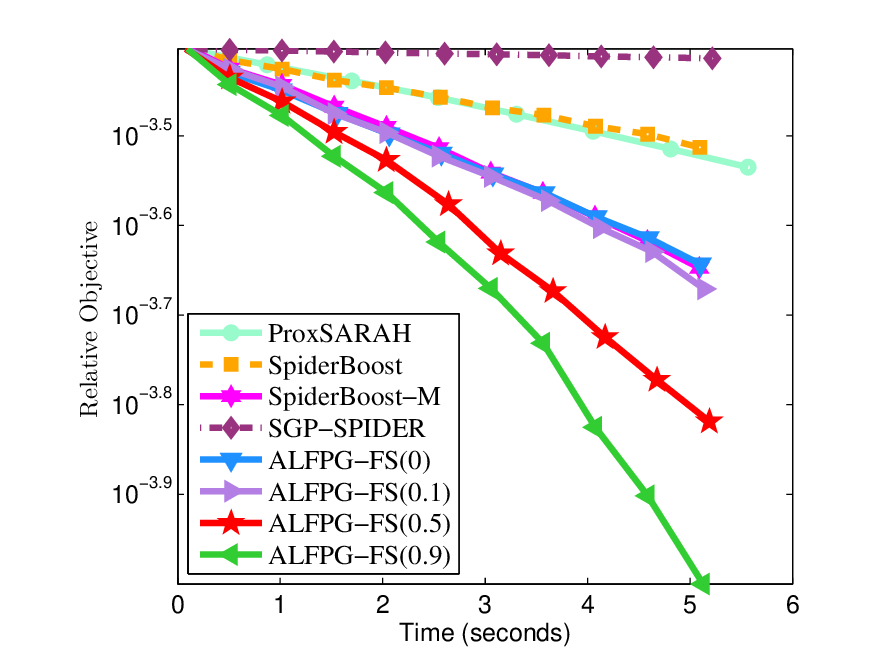}\caption{{\scriptsize TDT2{-9000-1000}}}\label{fig:fs:phase:2:1}\end{subfigure}
\begin{subfigure}{.24\textwidth}\centering\includegraphics[width=\linewidth]{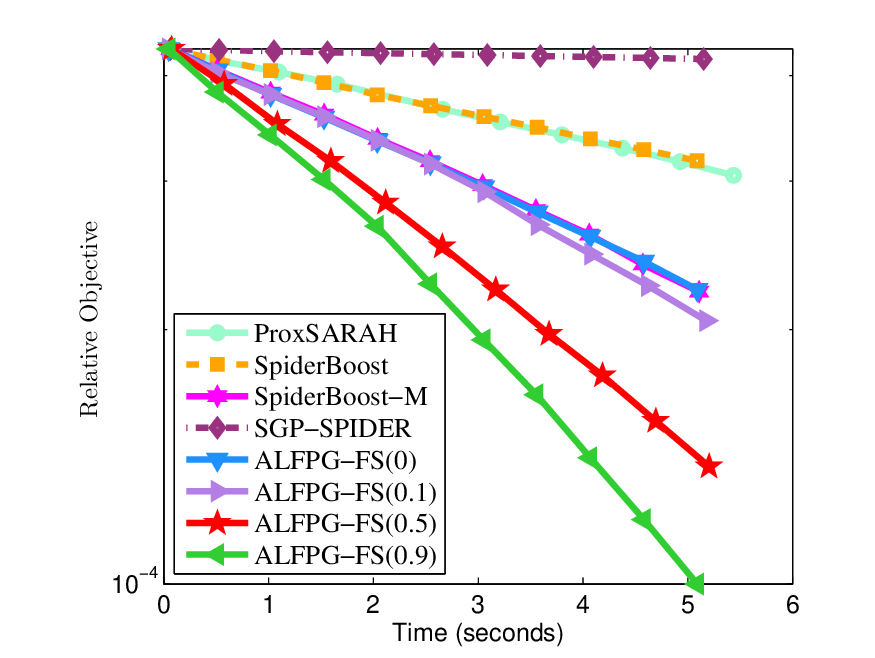}\caption{{\scriptsize 20Newsgroups{-10000-1000}}}\label{fig:fs:phase:2:2}\end{subfigure}
\begin{subfigure}{.24\textwidth}\centering\includegraphics[width=\linewidth]{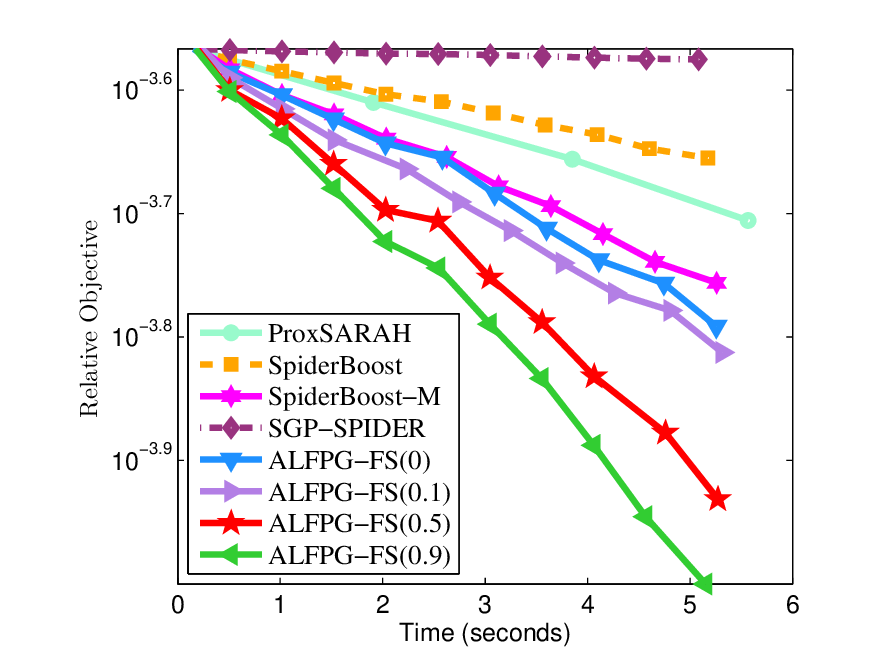}\caption{{\scriptsize Sector{-6412-1000}}}\label{fig:fs:phase:2:3}\end{subfigure}
\begin{subfigure}{.24\textwidth}\centering\includegraphics[width=\linewidth]{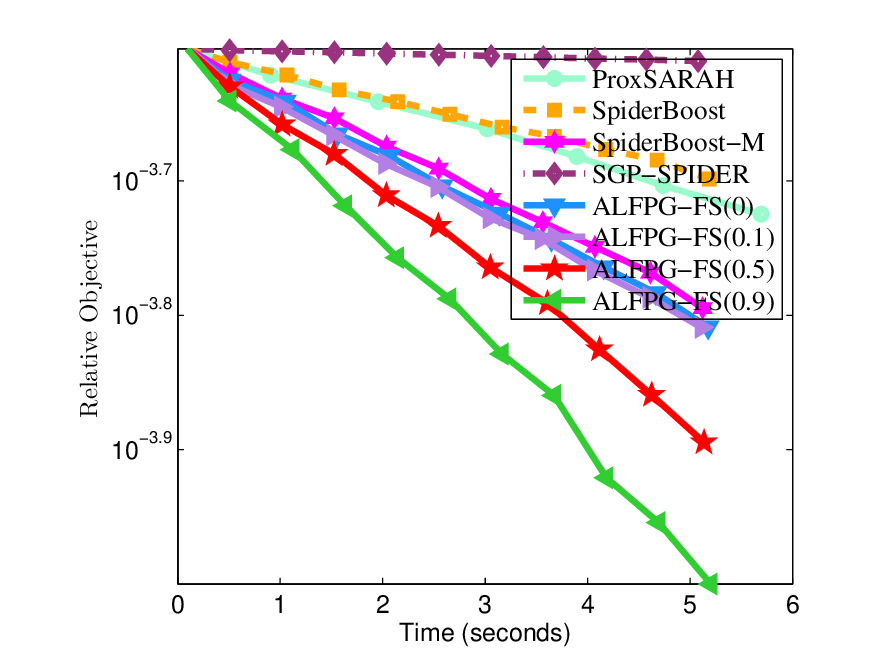}\caption{{\scriptsize MNIST{-2000-784}}}\label{fig:fs:phase:2:4}\end{subfigure}

\centering
\begin{subfigure}{.24\textwidth}\centering\includegraphics[width=\linewidth]{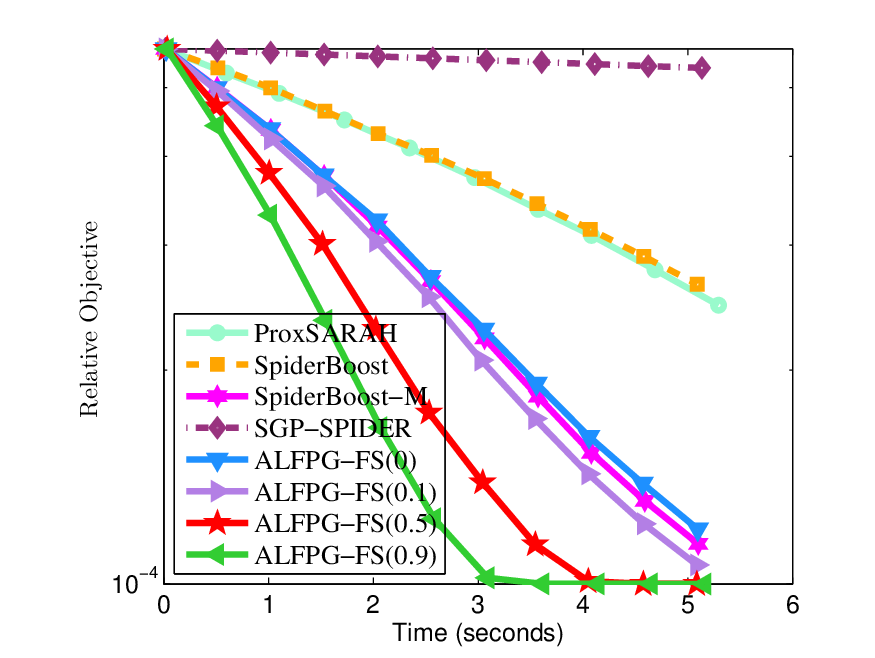}\caption{{\scriptsize CIFAR{-10000-1000}}}\label{fig:fs:phase:2:5}\end{subfigure}
\begin{subfigure}{.24\textwidth}\centering\includegraphics[width=\linewidth]{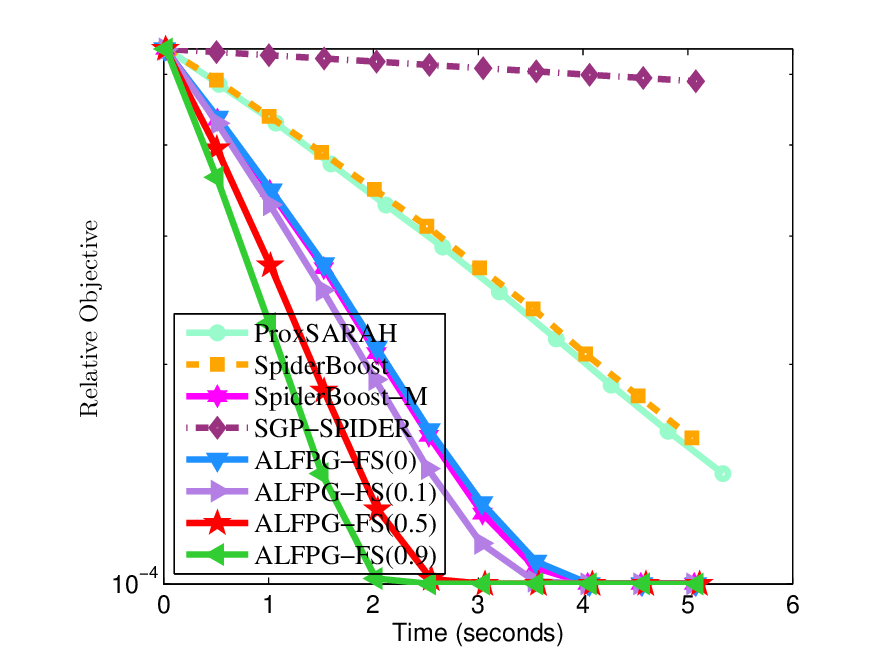}\caption{{\scriptsize Gisette{-5000-1000}}}\label{fig:fs:phase:2:6}\end{subfigure}
\begin{subfigure}{.24\textwidth}\centering\includegraphics[width=\linewidth]{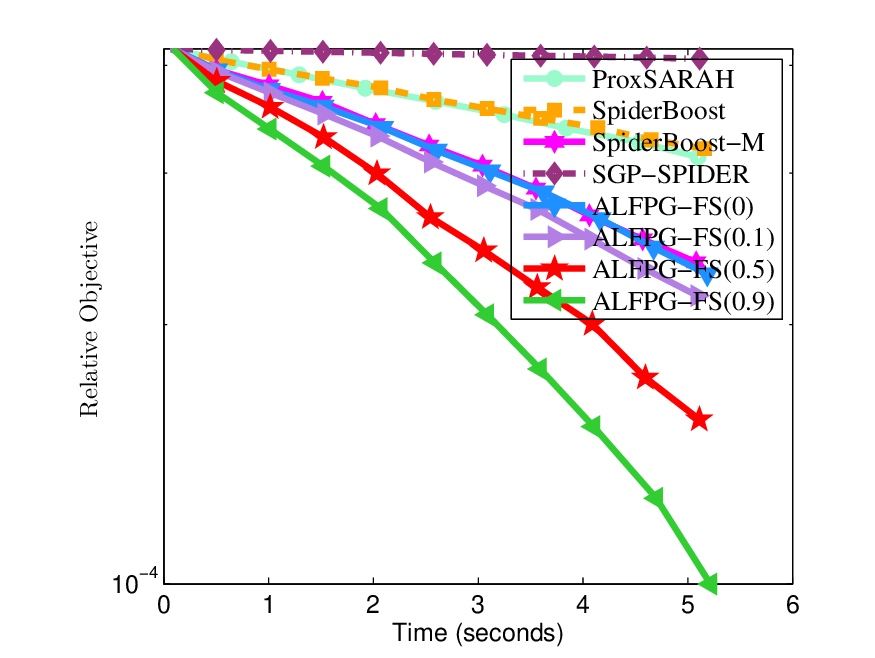}\caption{{\scriptsize Caltech{-3000-1000}}}\label{fig:fs:phase:2:7}\end{subfigure}
\begin{subfigure}{.24\textwidth}\centering\includegraphics[width=\linewidth]{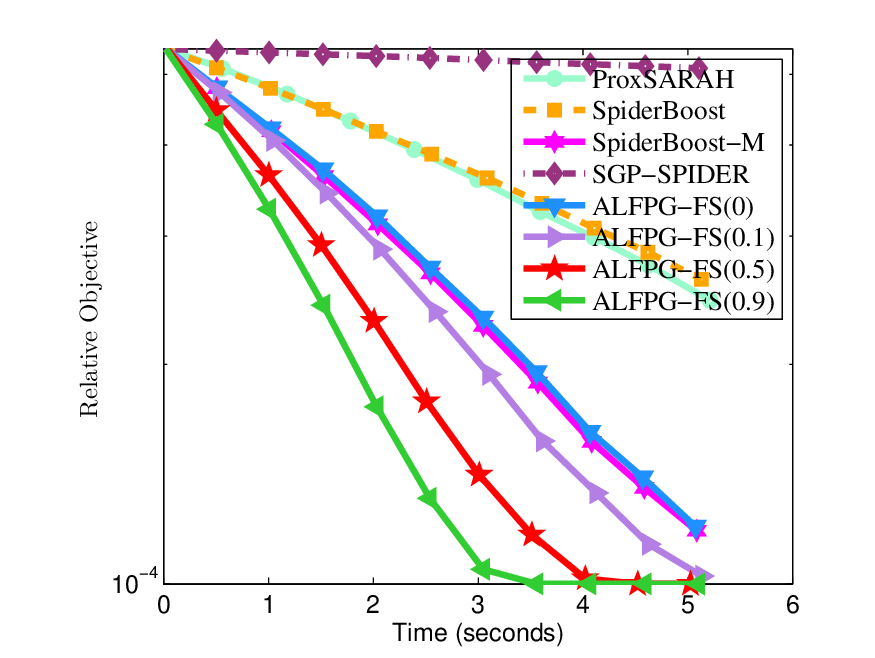}\caption{{\scriptsize Gaussian{-10000-1000}}}\label{fig:fs:phase:2:8}\end{subfigure}
\caption{Sparse phase retrieval with $\dot\lambda=0.001$ in the finite-sum setting.}
\label{fig:fs:phase:2}
\end{figure*}

\begin{figure*}[!t]
\centering
\begin{subfigure}{.24\textwidth}\centering\includegraphics[width=\linewidth]{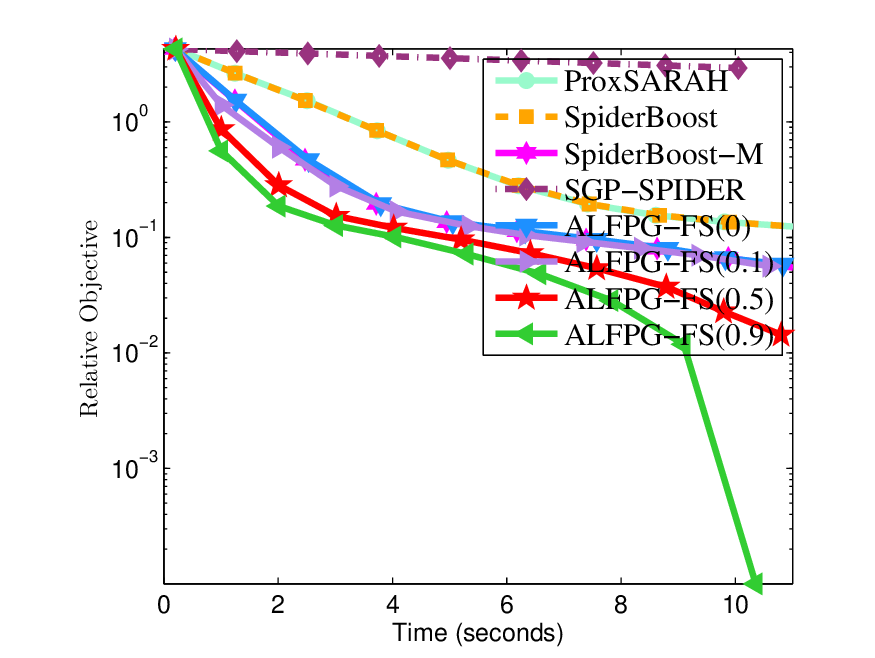}\caption{{\scriptsize TDT2{-3000-3000}}}\label{fig:fs:nn:1:1}\end{subfigure}
\begin{subfigure}{.24\textwidth}\centering\includegraphics[width=\linewidth]{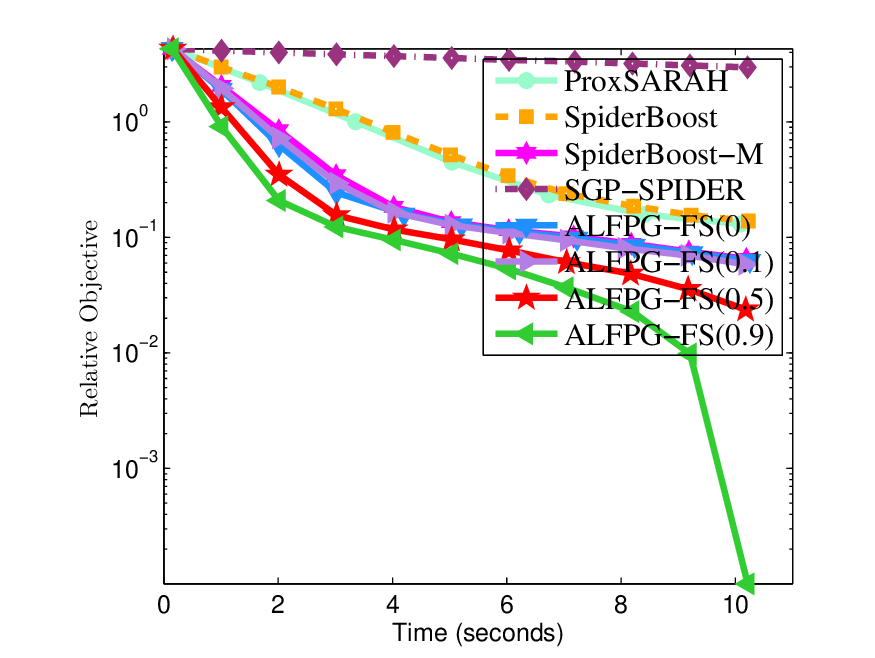}\caption{{\scriptsize 20Newsgroups{-5000-1000}}}\label{fig:fs:nn:1:2}\end{subfigure}
\begin{subfigure}{.24\textwidth}\centering\includegraphics[width=\linewidth]{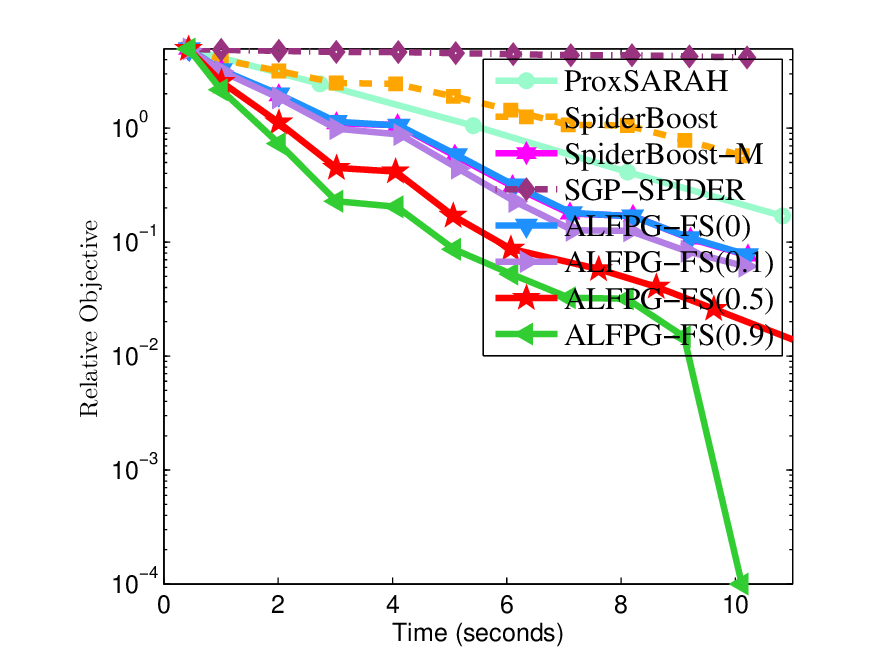}\caption{{\scriptsize Sector{-6000-1000}}}\label{fig:fs:nn:1:3}\end{subfigure}
\begin{subfigure}{.24\textwidth}\centering\includegraphics[width=\linewidth]{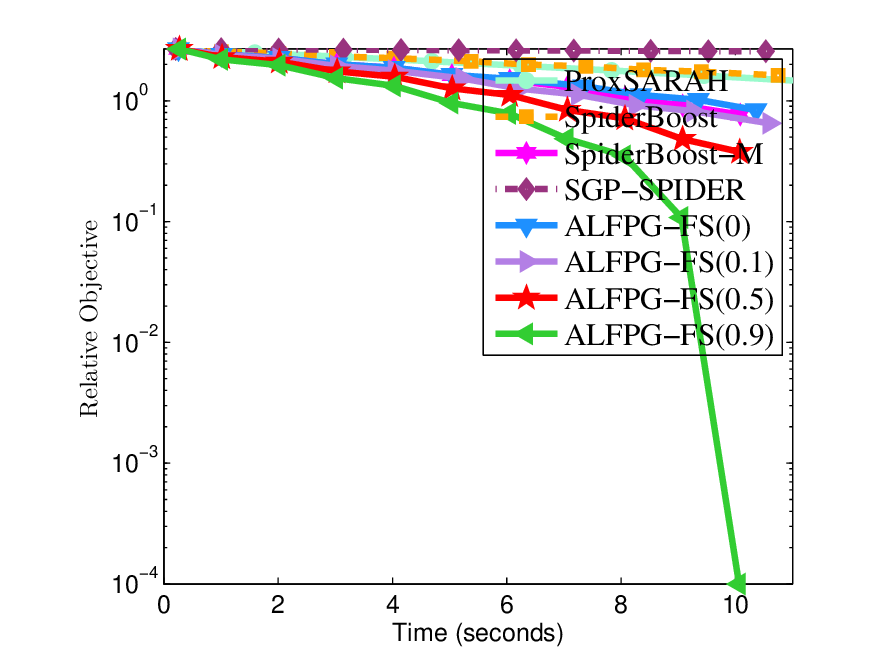}\caption{{\scriptsize MNIST{-5000-784}}}\label{fig:fs:nn:1:4}\end{subfigure}

\begin{subfigure}{.24\textwidth}\centering\includegraphics[width=\linewidth]{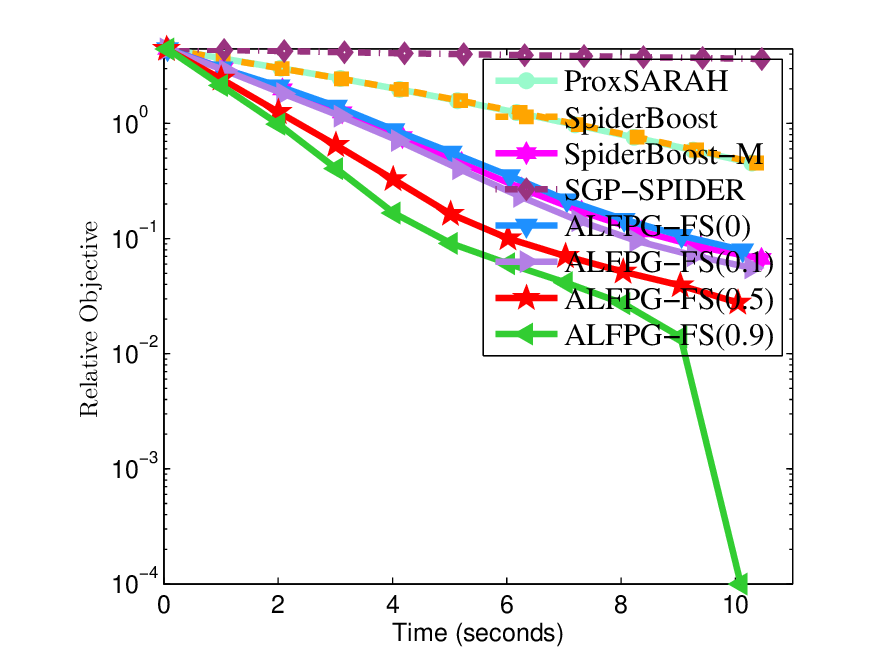}\caption{{\scriptsize CIFAR{-5000-1000}}}\label{fig:fs:nn:1:5}\end{subfigure}
\begin{subfigure}{.24\textwidth}\centering\includegraphics[width=\linewidth]{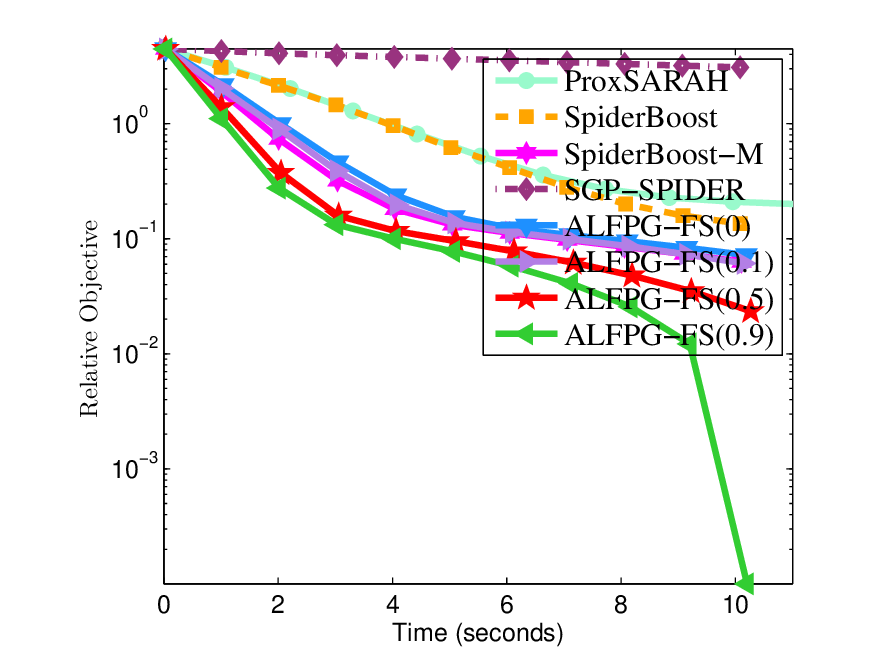}\caption{{\scriptsize Gisette{-6000-3000}}}\label{fig:fs:nn:1:6}\end{subfigure}
\begin{subfigure}{.24\textwidth}\centering\includegraphics[width=\linewidth]{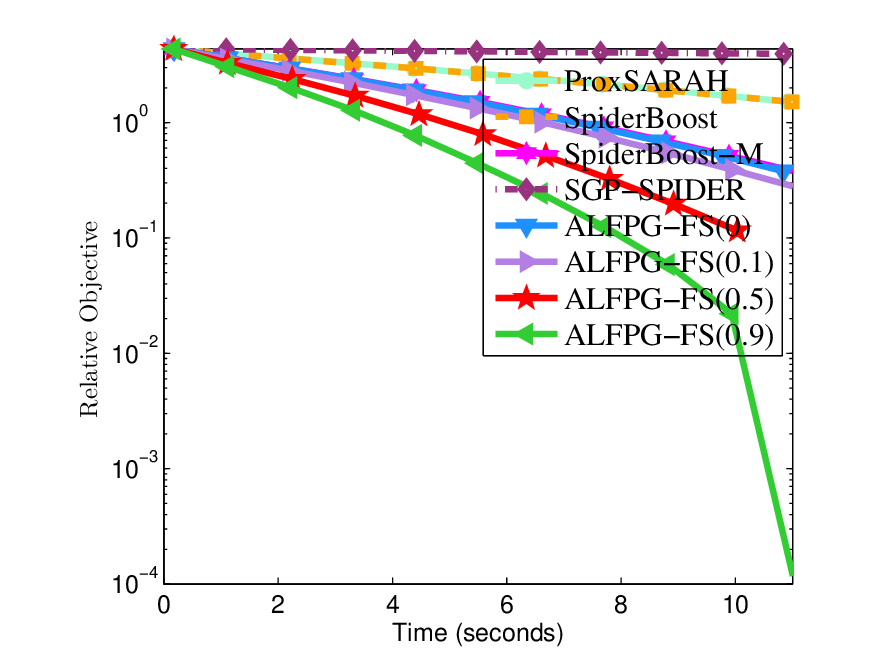}\caption{{\scriptsize Caltech{-2000-1000}}}\label{fig:fs:nn:1:7}\end{subfigure}
\begin{subfigure}{.24\textwidth}\centering\includegraphics[width=\linewidth]{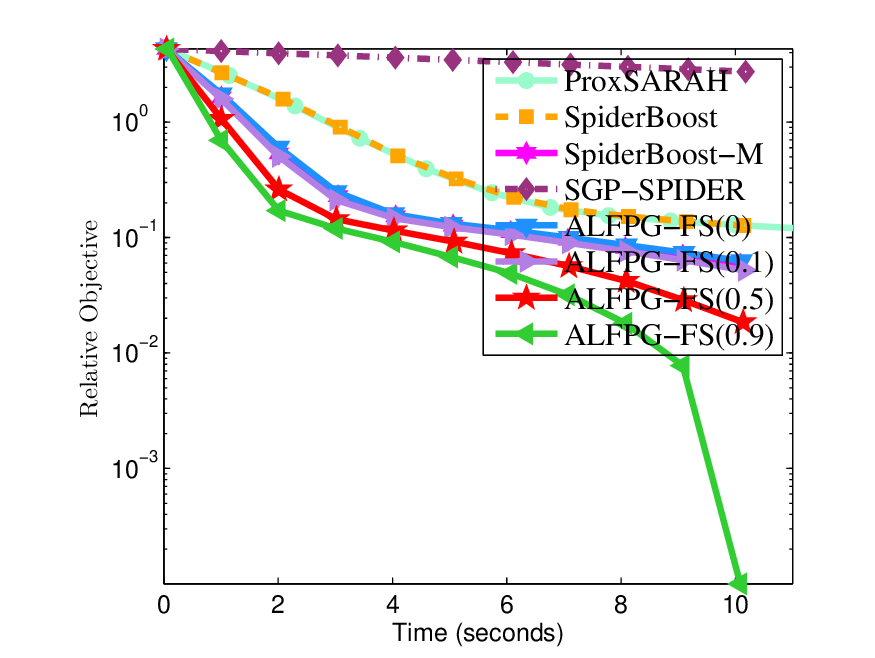}\caption{{\scriptsize Gaussian{-3000-1000}}}\label{fig:fs:nn:1:8}\end{subfigure}

\caption{Three-layer neural network regression with $\dot\lambda=0.01$ in the finite-sum setting.}
\label{fig:fs:nn:1}
\end{figure*}

\begin{figure*}[!h]
\centering
\begin{subfigure}{.24\textwidth}\centering\includegraphics[width=\linewidth]{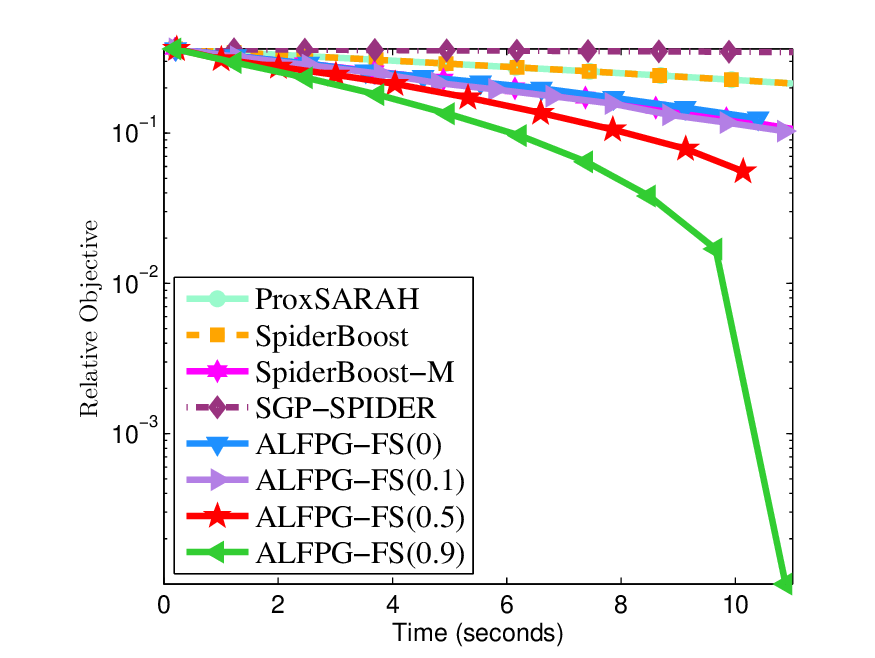}\caption{{\scriptsize TDT2{-3000-3000}}}\label{fig:fs:nn:2:1}\end{subfigure}
\begin{subfigure}{.24\textwidth}\centering\includegraphics[width=\linewidth]{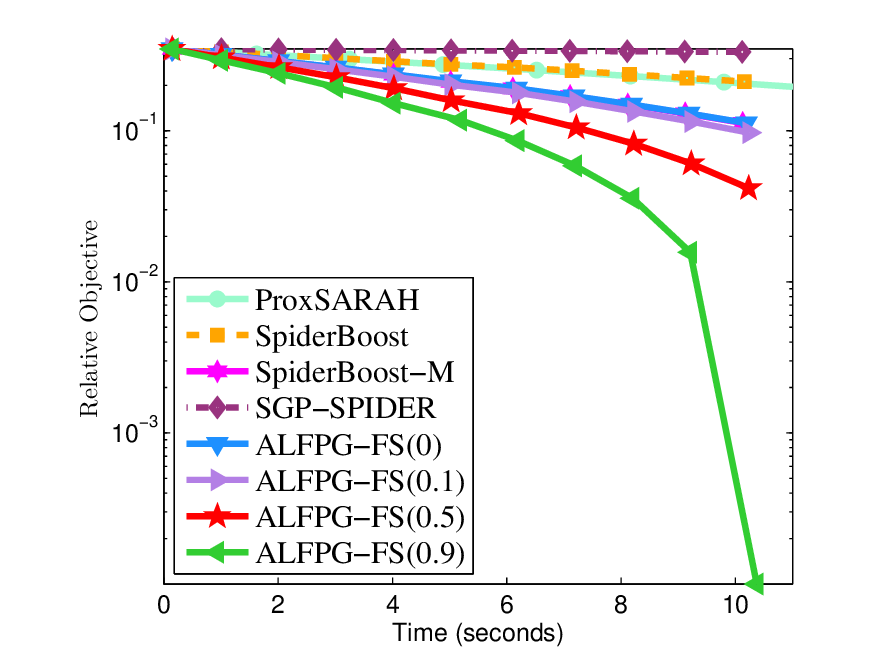}\caption{{\scriptsize 20Newsgroups{-5000-1000}}}\label{fig:fs:nn:2:2}\end{subfigure}
\begin{subfigure}{.24\textwidth}\centering\includegraphics[width=\linewidth]{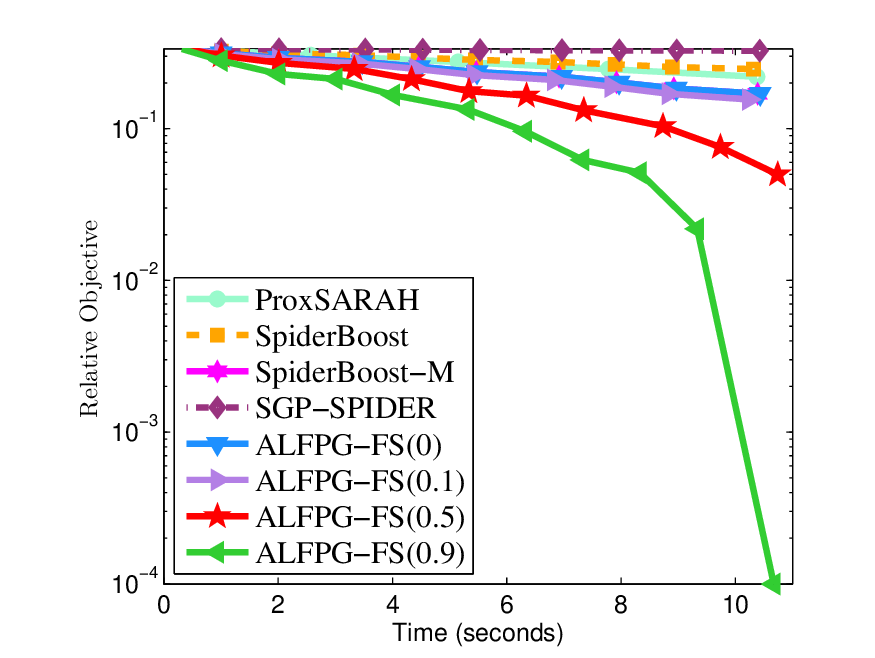}\caption{{\scriptsize Sector{-6000-1000}}}\label{fig:fs:nn:2:3}\end{subfigure}
\begin{subfigure}{.24\textwidth}\centering\includegraphics[width=\linewidth]{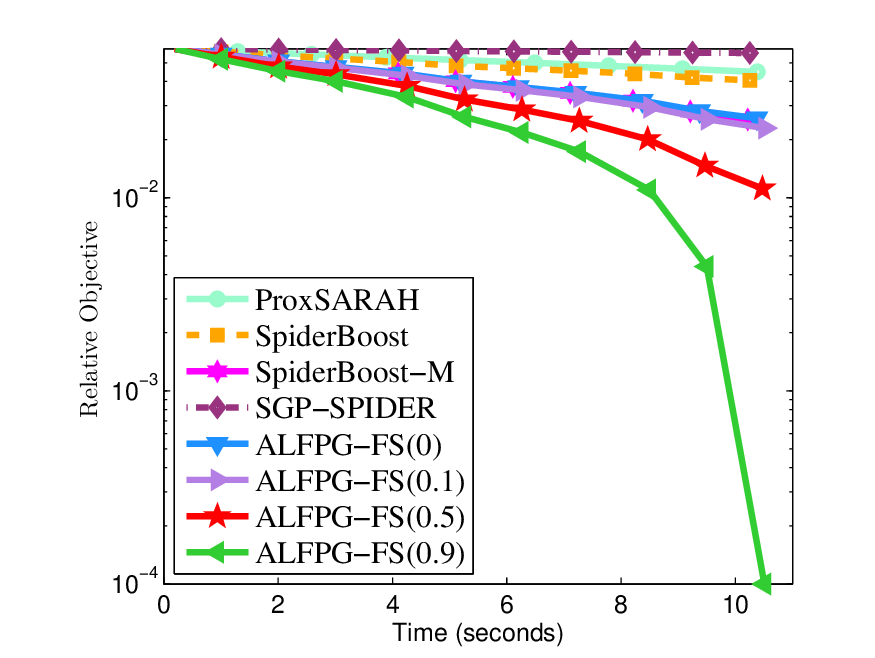}\caption{{\scriptsize MNIST{-5000-784}}}\label{fig:fs:nn:2:4}\end{subfigure}

\begin{subfigure}{.24\textwidth}\centering\includegraphics[width=\linewidth]{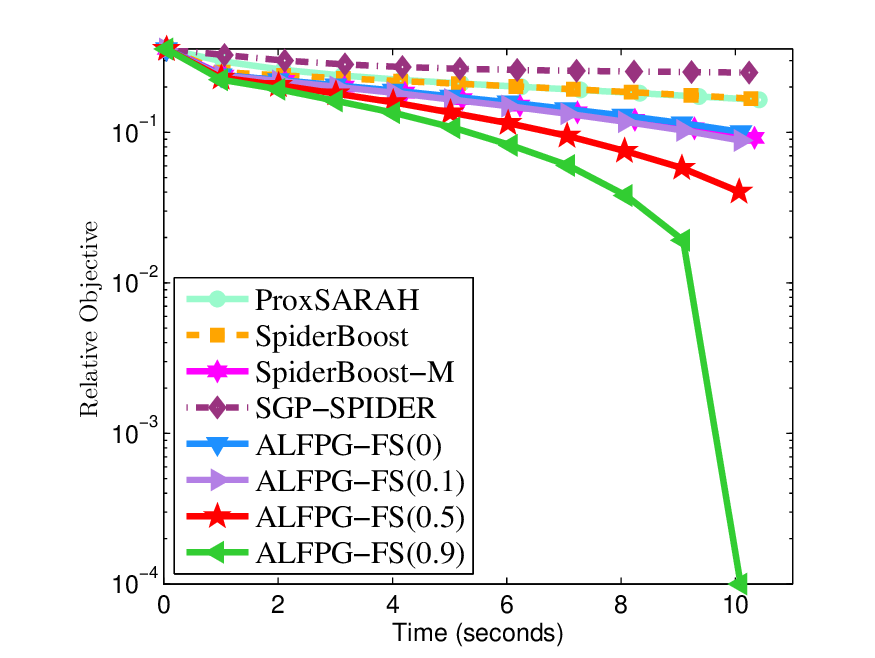}\caption{{\scriptsize CIFAR{-5000-1000}}}\label{fig:fs:nn:2:5}\end{subfigure}
\begin{subfigure}{.24\textwidth}\centering\includegraphics[width=\linewidth]{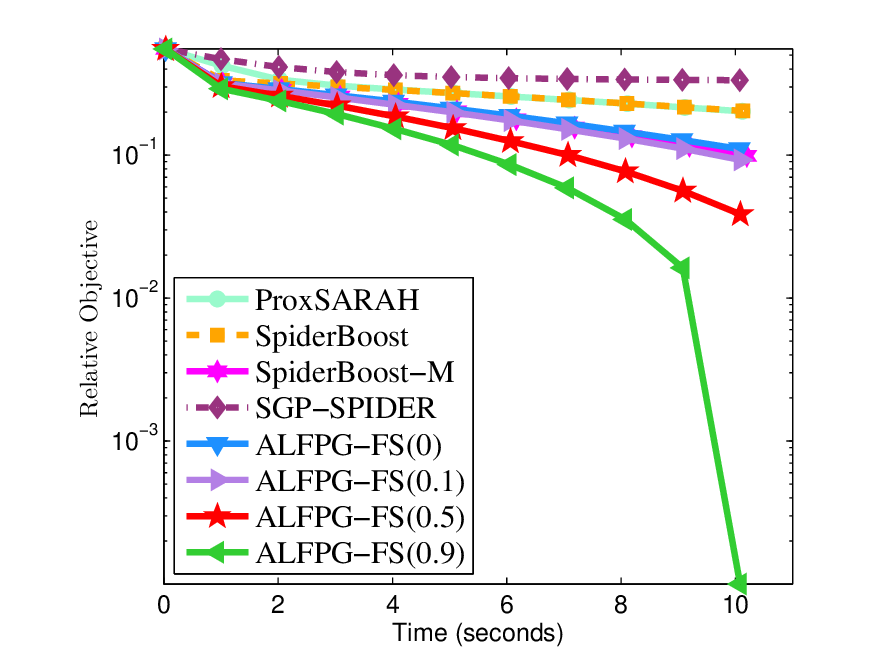}\caption{{\scriptsize Gisette{-6000-3000}}}\label{fig:fs:nn:2:6}\end{subfigure}
\begin{subfigure}{.24\textwidth}\centering\includegraphics[width=\linewidth]{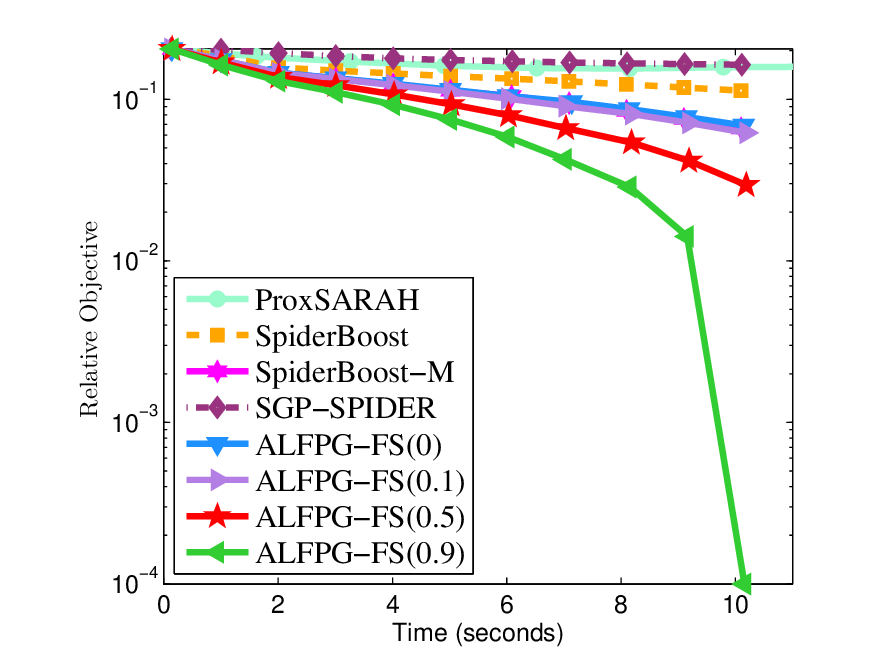}\caption{{\scriptsize Caltech{-2000-1000}}}\label{fig:fs:nn:2:7}\end{subfigure}
\begin{subfigure}{.24\textwidth}\centering\includegraphics[width=\linewidth]{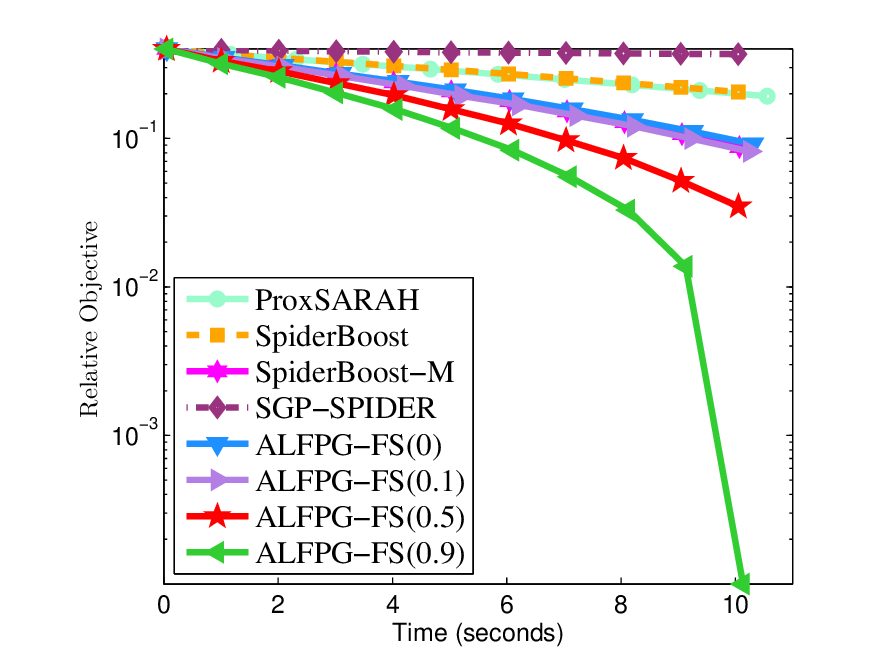}\caption{{\scriptsize Gaussian{-3000-1000}}}\label{fig:fs:nn:2:8}\end{subfigure}
\caption{Three-layer neural network regression with $\dot\lambda=0.001$ in the finite-sum setting.}
\label{fig:fs:nn:2}
\end{figure*}

$\blacktriangleright$ \textbf{Experimental Results}.
\textbf{(i)} For SPR, Figures~\ref{fig:fs:phase:1}--\ref{fig:fs:phase:2} show that {\sf ALFPG-FS} with $\theta=0.9$ often reduces the objective faster than
with $\theta=0$. {\sf ALFPG-FS} attains some of the lowest recorded
objective values, and {\sf SpiderBoost-M} is also competitive on
several instances. \textbf{(ii)} For NNR, Figures~\ref{fig:fs:nn:1}--\ref{fig:fs:nn:2} show that
$\theta=0.9$ reaches a lower objective than $\theta=0$ in many panels.
{\sf ALFPG-FS} often attains the lowest recorded objective among the
compared methods on this task.

\subsection{Results of {\sf ALFPG-D} on LEV}
\label{app:exp:orth}

$\blacktriangleright$ \textbf{Experimental Settings}.
We use the authors' implementations of {\sf OptM} and {\sf FOForth}
with their default solver settings. For {\sf ALFPG-D}, we set
$(\underline{\rm v},\alpha,\beta)=(0.001,0.001,0)$ and test
$\theta\in\{0,0.1,0.5,0.9\}$. We evaluate all methods at ranks
$\dot r\in\{20,50\}$. On each instance, all methods start from the
same point, obtained by polar projection of a Gaussian matrix.

\begin{figure*}[!h]
\centering
\begin{subfigure}{.24\textwidth}\centering\includegraphics[width=\linewidth]{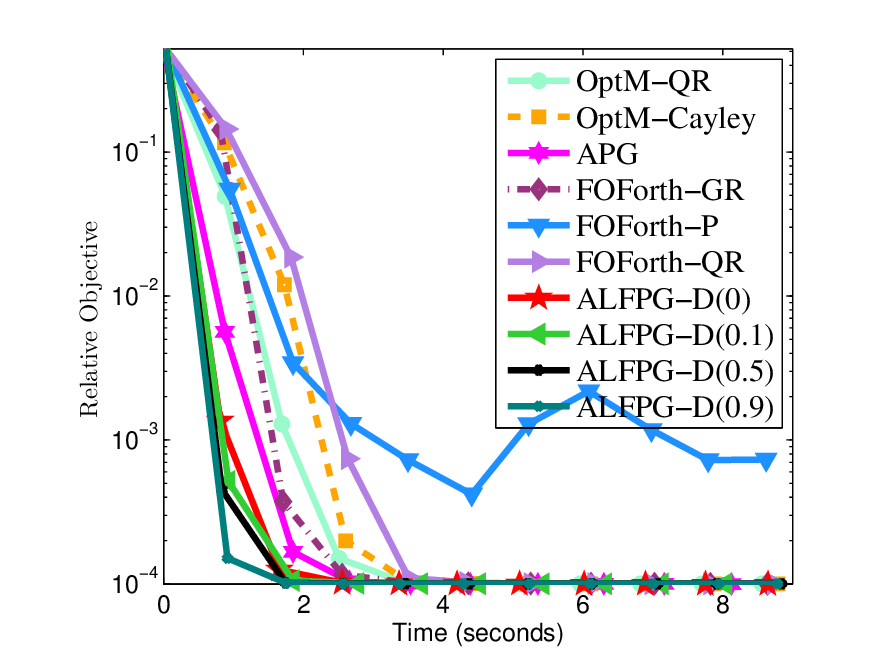}\caption{{\scriptsize TDT2{-3000-3000}}}\label{fig:d:orth:1:1}\end{subfigure}
\begin{subfigure}{.24\textwidth}\centering\includegraphics[width=\linewidth]{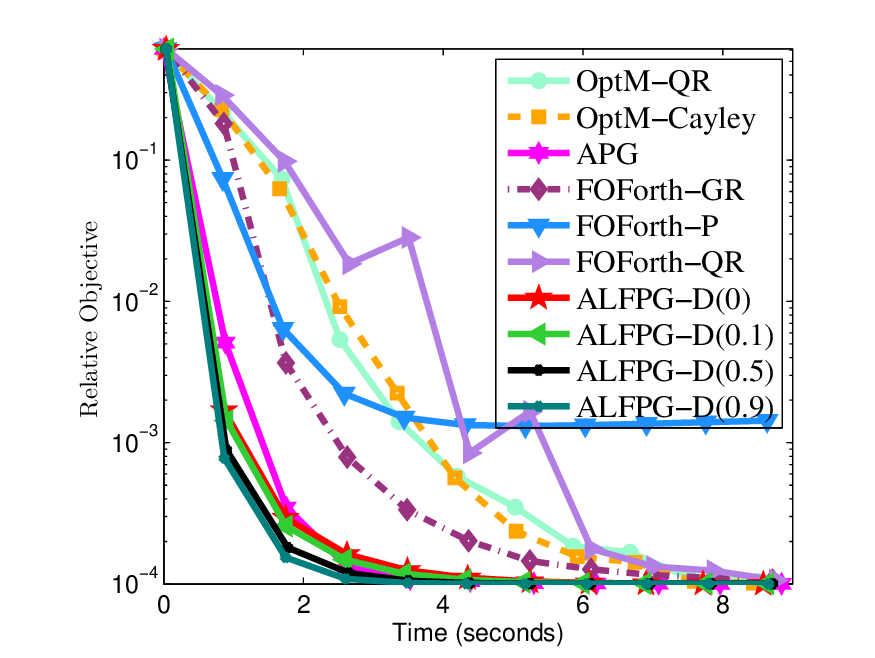}\caption{{\scriptsize 20Newsgroups{-5000-1000}}}\label{fig:d:orth:1:2}\end{subfigure}
\begin{subfigure}{.24\textwidth}\centering\includegraphics[width=\linewidth]{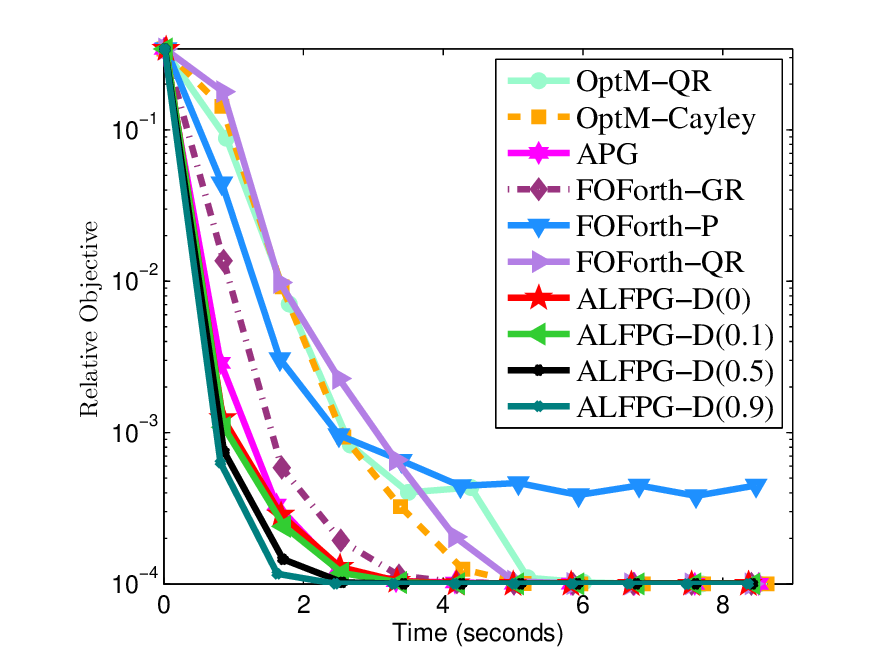}\caption{{\scriptsize Sector{-6000-1000}}}\label{fig:d:orth:1:3}\end{subfigure}
\begin{subfigure}{.24\textwidth}\centering\includegraphics[width=\linewidth]{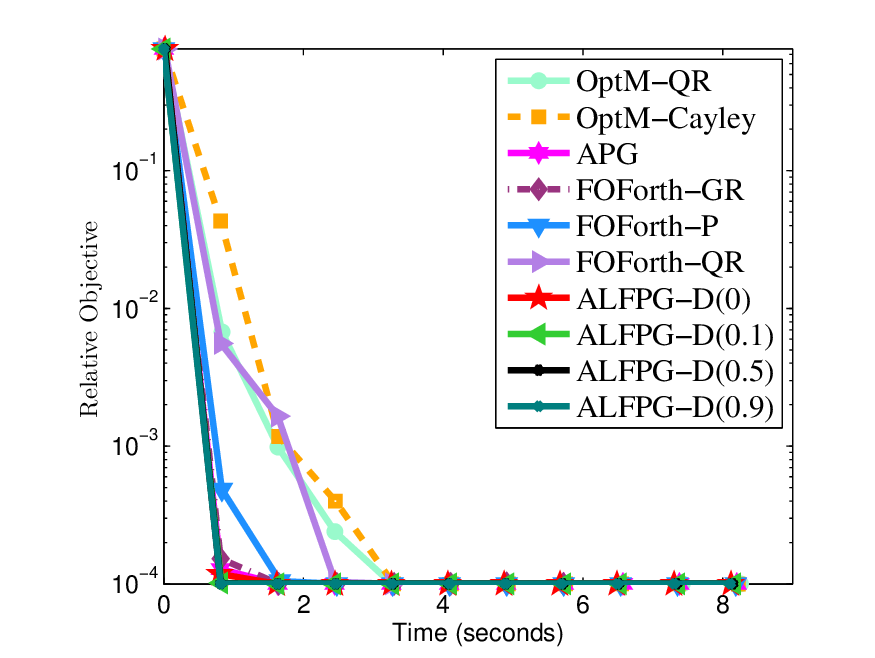}\caption{{\scriptsize MNIST{-5000-784}}}\label{fig:d:orth:1:4}\end{subfigure}

\begin{subfigure}{.24\textwidth}\centering\includegraphics[width=\linewidth]{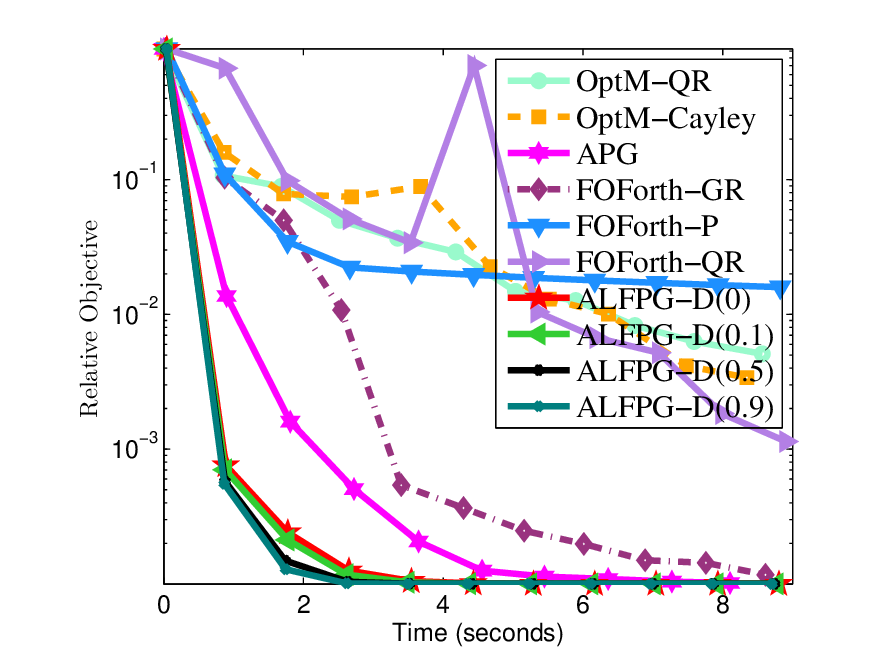}\caption{{\scriptsize CIFAR{-5000-1000}}}\label{fig:fs:orth:1:5}\end{subfigure}
\begin{subfigure}{.24\textwidth}\centering\includegraphics[width=\linewidth]{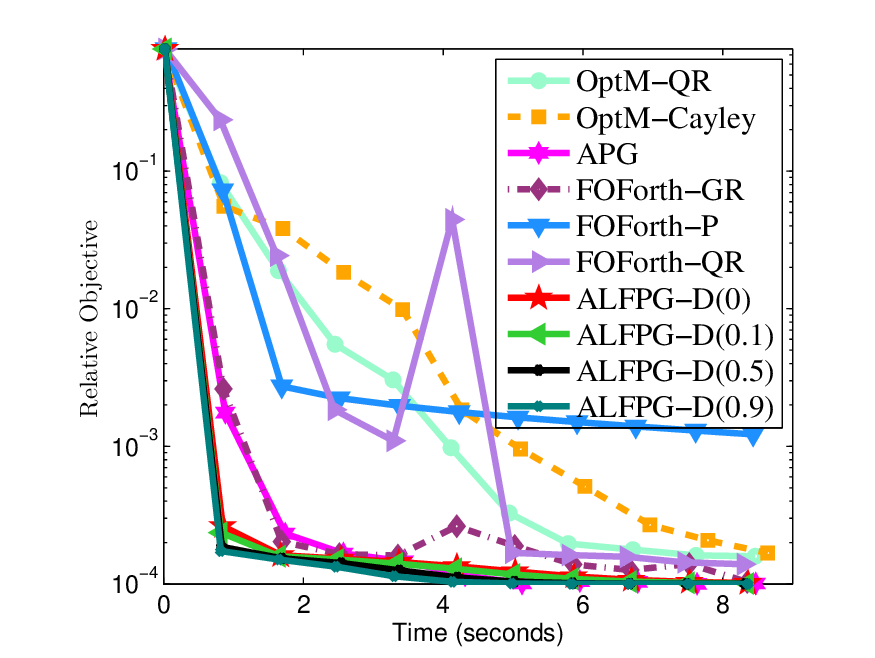}\caption{{\scriptsize Gisette{-6000-3000}}}\label{fig:fs:orth:1:6}\end{subfigure}
\begin{subfigure}{.24\textwidth}\centering\includegraphics[width=\linewidth]{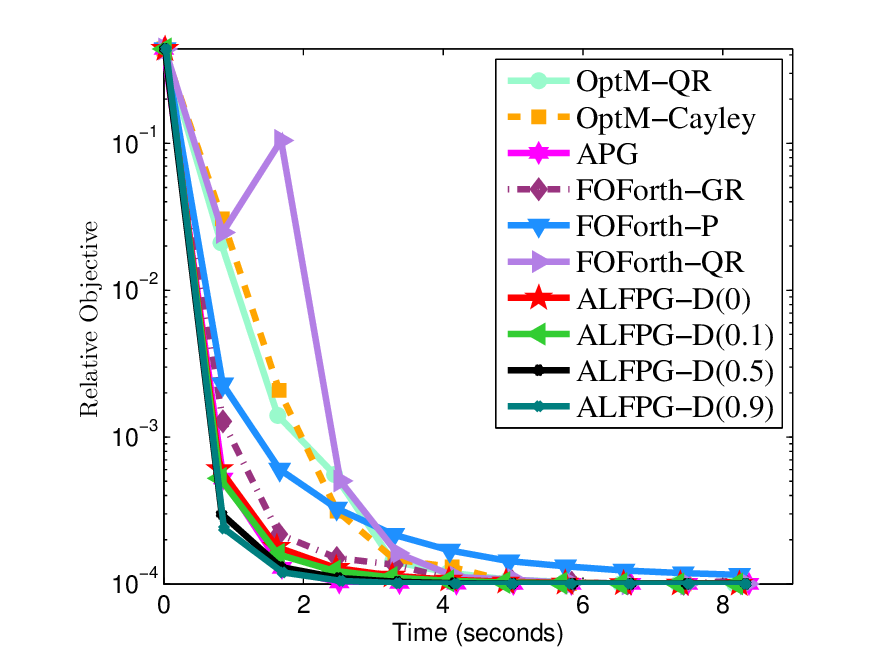}\caption{{\scriptsize Caltech{-2000-1000}}}\label{fig:fs:orth:1:7}\end{subfigure}
\begin{subfigure}{.24\textwidth}\centering\includegraphics[width=\linewidth]{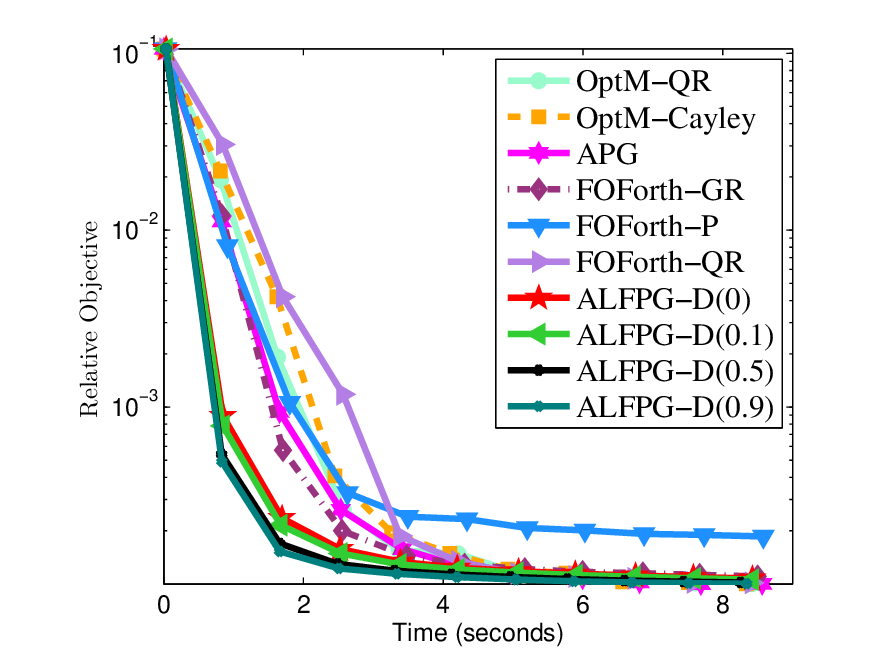}\caption{{\scriptsize Gaussian{-3000-1000}}}\label{fig:fs:orth:1:8}\end{subfigure}
\caption{Linear eigenvalue problems with $\dot r=20$ in the deterministic setting.}
\label{fig:d:orth:1}
\end{figure*}

\begin{figure*}[!h]
\centering

\begin{subfigure}{.24\textwidth}\centering\includegraphics[width=\linewidth]{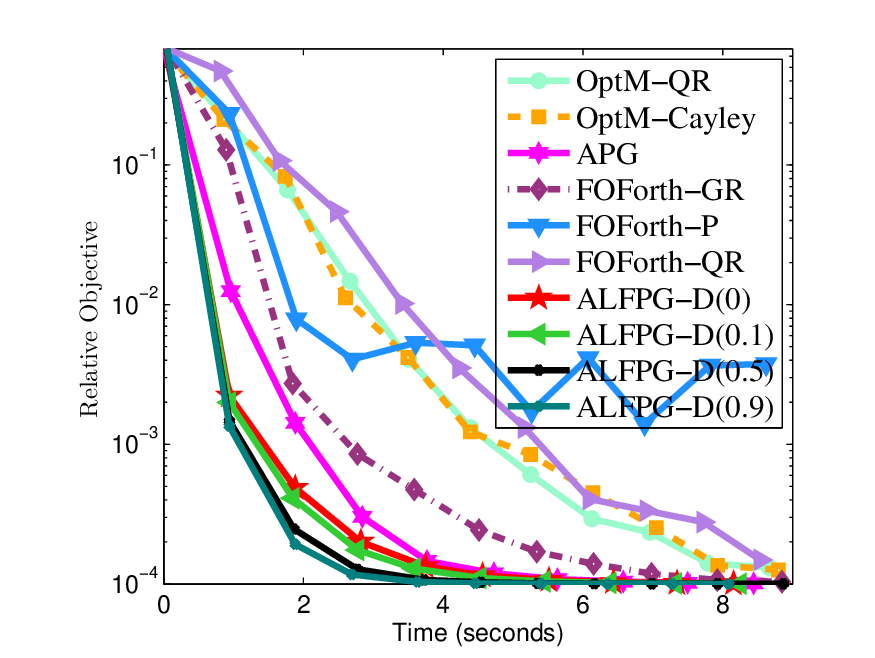}\caption{{\scriptsize TDT2{-3000-3000}}}\label{fig:d:orth:2:1}\end{subfigure}
\begin{subfigure}{.24\textwidth}\centering\includegraphics[width=\linewidth]{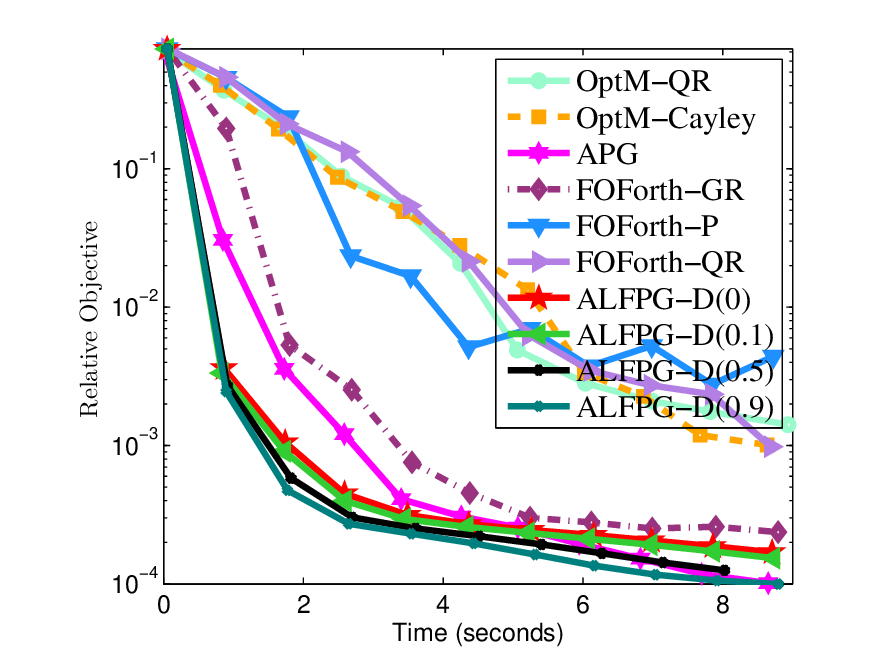}\caption{{\scriptsize 20Newsgroups{-5000-1000}}}\label{fig:d:orth:2:2}\end{subfigure}
\begin{subfigure}{.24\textwidth}\centering\includegraphics[width=\linewidth]{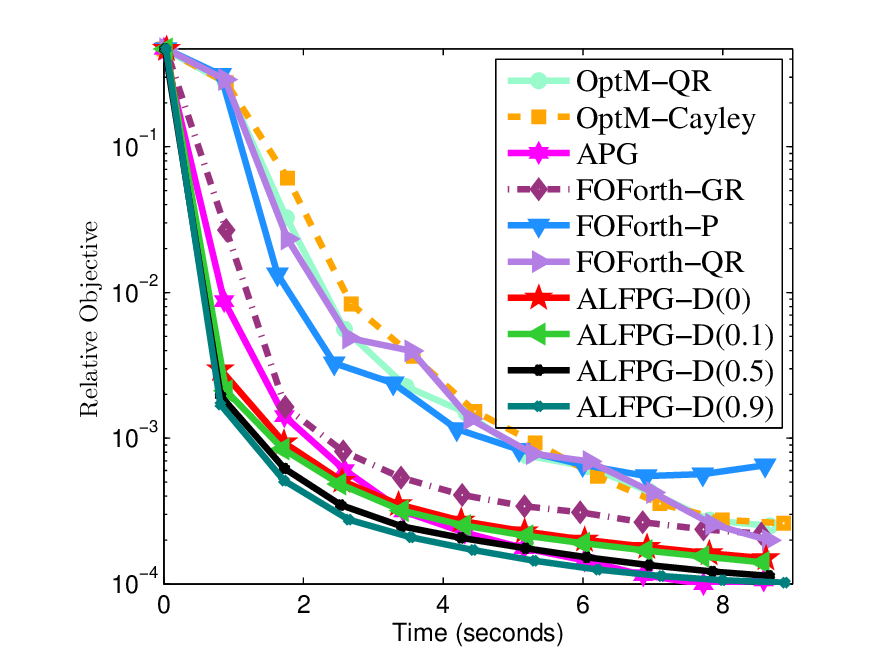}\caption{{\scriptsize Sector{-6000-1000}}}\label{fig:d:orth:2:3}\end{subfigure}
\begin{subfigure}{.24\textwidth}\centering\includegraphics[width=\linewidth]{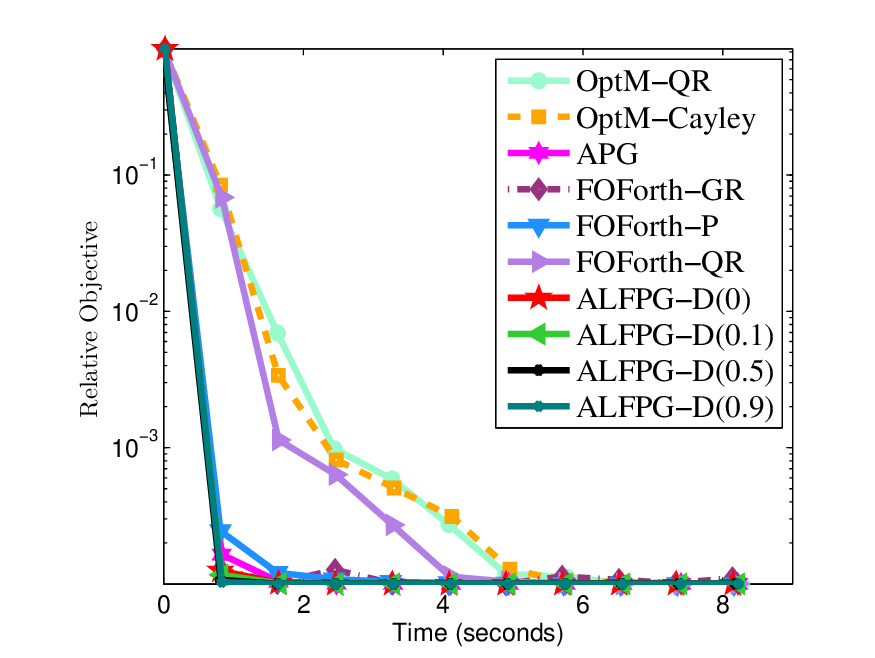}\caption{{\scriptsize MNIST{-5000-784}}}\label{fig:d:orth:2:4}\end{subfigure}

\begin{subfigure}{.24\textwidth}\centering\includegraphics[width=\linewidth]{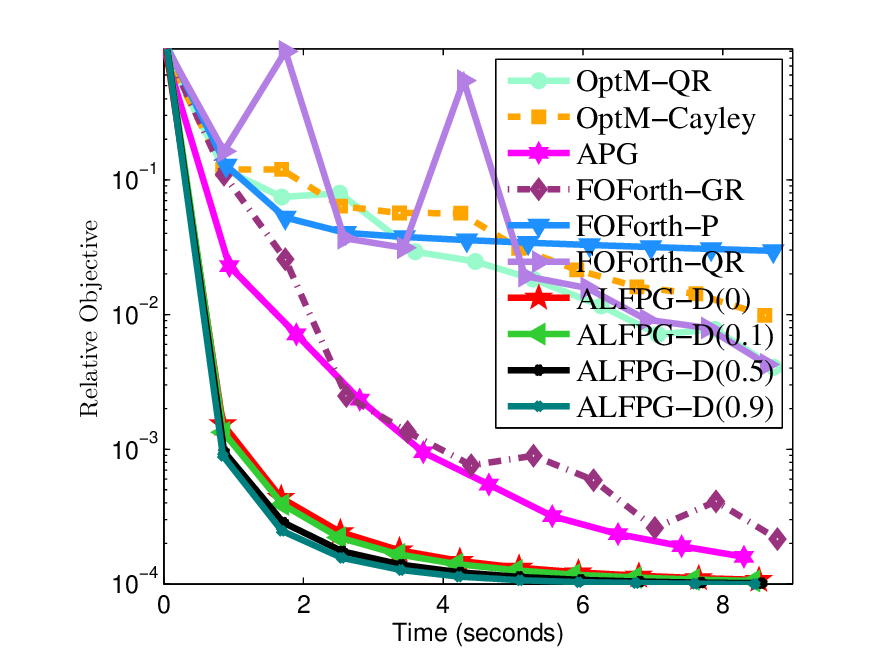}\caption{{\scriptsize CIFAR{-5000-1000}}}\label{fig:fs:orth:2:5}\end{subfigure}
\begin{subfigure}{.24\textwidth}\centering\includegraphics[width=\linewidth]{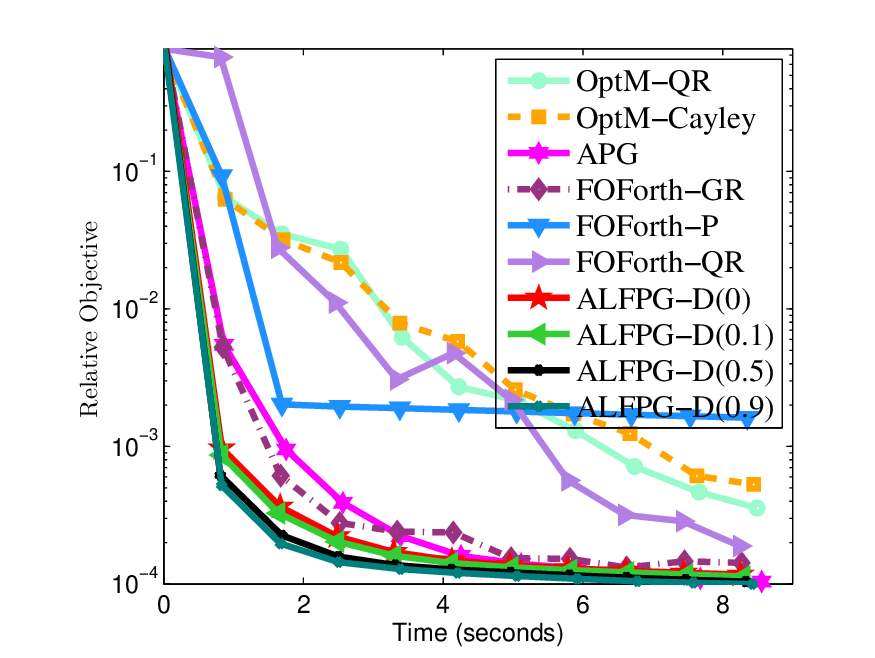}\caption{{\scriptsize Gisette{-6000-3000}}}\label{fig:fs:orth:2:6}\end{subfigure}
\begin{subfigure}{.24\textwidth}\centering\includegraphics[width=\linewidth]{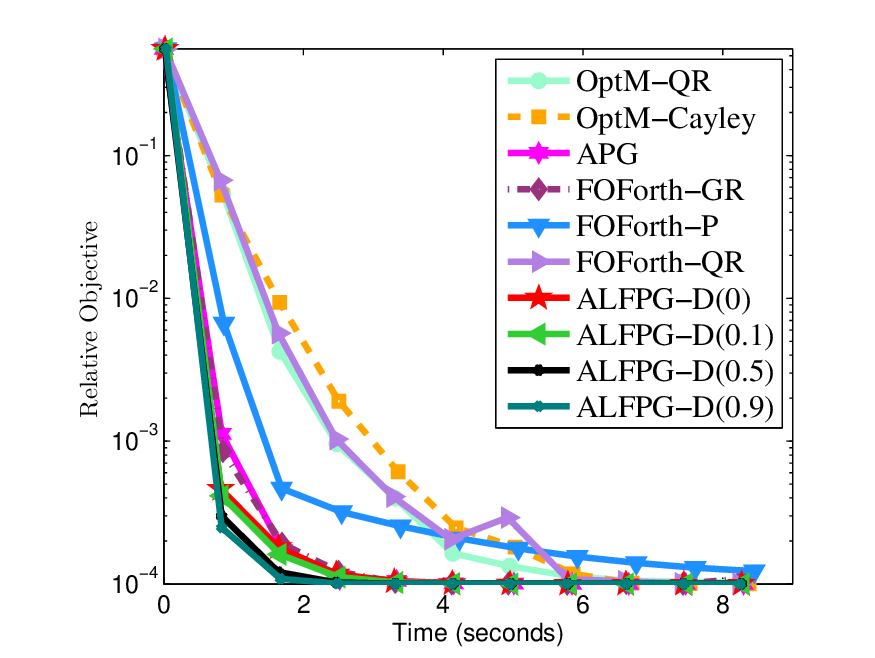}\caption{{\scriptsize Caltech{-2000-1000}}}\label{fig:fs:orth:2:7}\end{subfigure}
\begin{subfigure}{.24\textwidth}\centering\includegraphics[width=\linewidth]{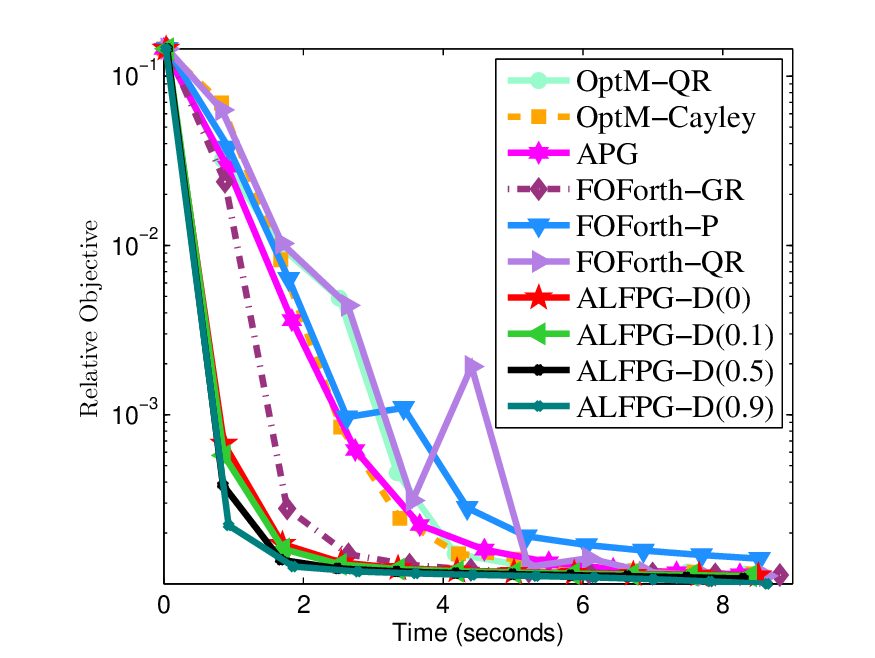}\caption{{\scriptsize Gaussian{-3000-1000}}}\label{fig:fs:orth:2:8}\end{subfigure}

\caption{Linear eigenvalue problems with $\dot r=50$ in the deterministic setting.}
\label{fig:d:orth:2}
\end{figure*}

$\blacktriangleright$ \textbf{Experimental Results}.
Figures~\ref{fig:d:orth:1}--\ref{fig:d:orth:2} show that the
extrapolated variants of {\sf ALFPG-D} often reduce the objective
more rapidly than the variant with $\theta=0$. The setting
$\theta=0.9$ typically gives the lowest recorded objective among
the {\sf ALFPG-D} variants. {\sf ALFPG-D} attains the lowest recorded
objective on many instances. Several baselines also reach comparable
objective values, and the ranking varies with the dataset and rank.

\subsection{Results of {\sf ALFPG-MVR1} and {\sf ALFPG-MVR2}
on SPR, LEV, and NNR}
\label{app:exp:expectation}

$\blacktriangleright$ \textbf{Experimental Settings}. On each instance, all methods start from the same feasible point and use mini-batches of $256$ samples drawn uniformly with replacement. For SPR, we use box-constrained $\ell_1$ regularization with box radius
$\dot r=1000$ and $\dot\lambda\in\{10^{-3},10^{-5}\}$. For LEV, we test ranks $\dot r\in\{20,50\}$ and initialize with the polar projection of a Gaussian matrix. For NNR, we use two $\tanh$
hidden layers of widths $16$ and $8$, a scalar linear output, and
noise-free targets generated by a teacher network. We set the box
radius to $\dot r=10$ and use $\dot\lambda\in\{10^{-2},10^{-3}\}$.

For {\sf SpiderBoost}, {\sf ProxSARAH}, {\sf ProxHSGD}, {\sf MB-SPG}, and {\sf Prox-AdaGrad-Norm},
we use a constant proximal stepsize of $0.1$ for SPR and $0.01$ for
LEV and NNR. {\sf PStorm} uses the decaying stepsize
$\eta(t+1)^{-1/3}$ for iterations $t\geq 1$, with $\eta=0.1$ for SPR
and $\eta=0.01$ for LEV and NNR. {\sf ProxGen} uses the adaptive
stepsize $\eta/(\sqrt{s_t}+10^{-8})$, where $s_t$ is a scalar
exponential moving average of the mean squared gradient entries,
initialized at zero with decay $0.999$. For {\sf ProxGen},
{\sf ALFPG-MVR1}, and {\sf ALFPG-MVR2}, we set the stepsize coefficient
$\eta$ to $100$ for SPR and $1$ for LEV and NNR.

\begin{figure*}[!t]
\centering
\begin{subfigure}{.24\textwidth}\centering\includegraphics[width=\linewidth]{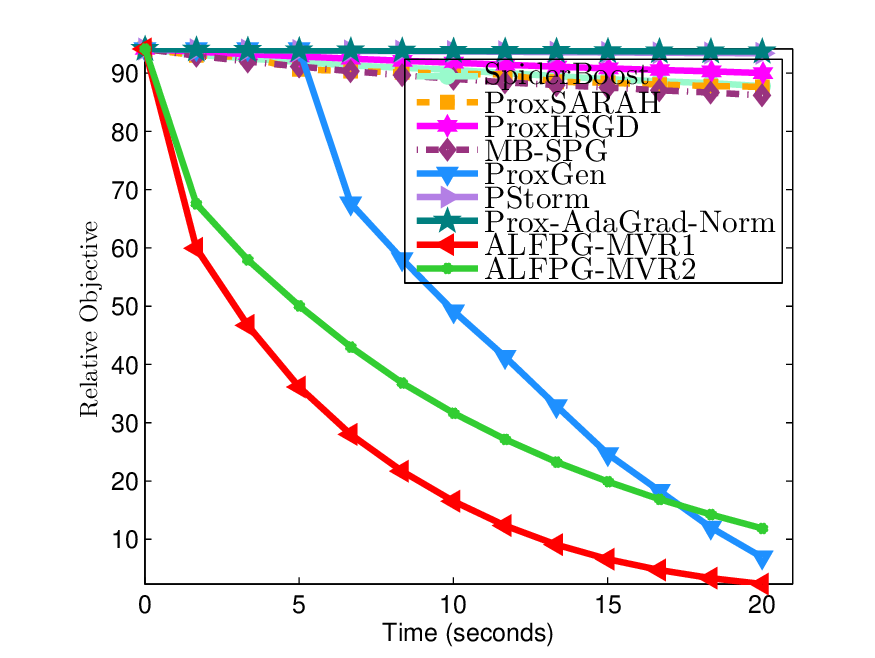}\caption{{\scriptsize TDT2{-9000-1000}}}\label{fig:ex:phase:1:1}\end{subfigure}
\begin{subfigure}{.24\textwidth}\centering\includegraphics[width=\linewidth]{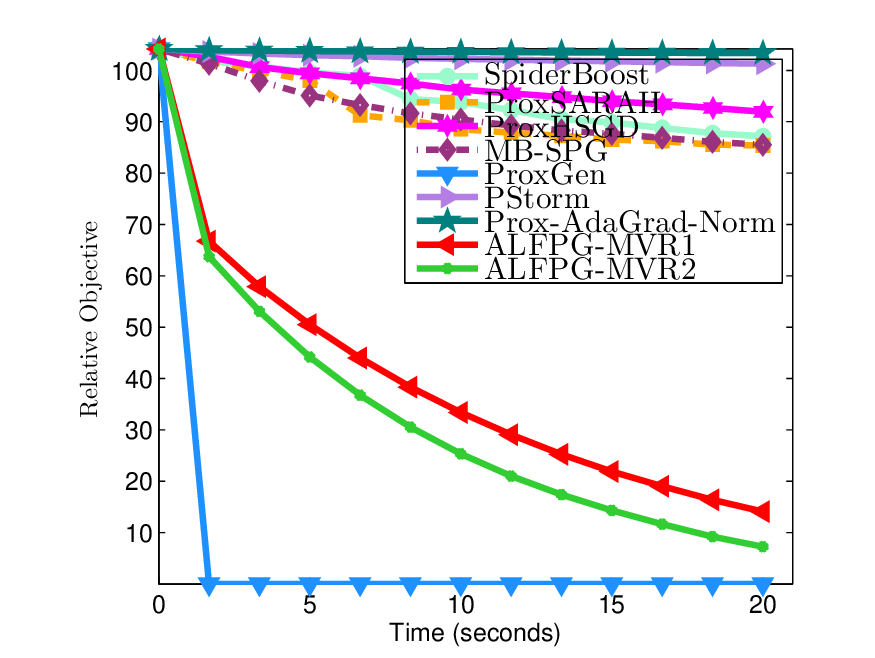}\caption{{\scriptsize 20Newsgroups{-10000-1000}}}\label{fig:ex:phase:1:2}\end{subfigure}
\begin{subfigure}{.24\textwidth}\centering\includegraphics[width=\linewidth]{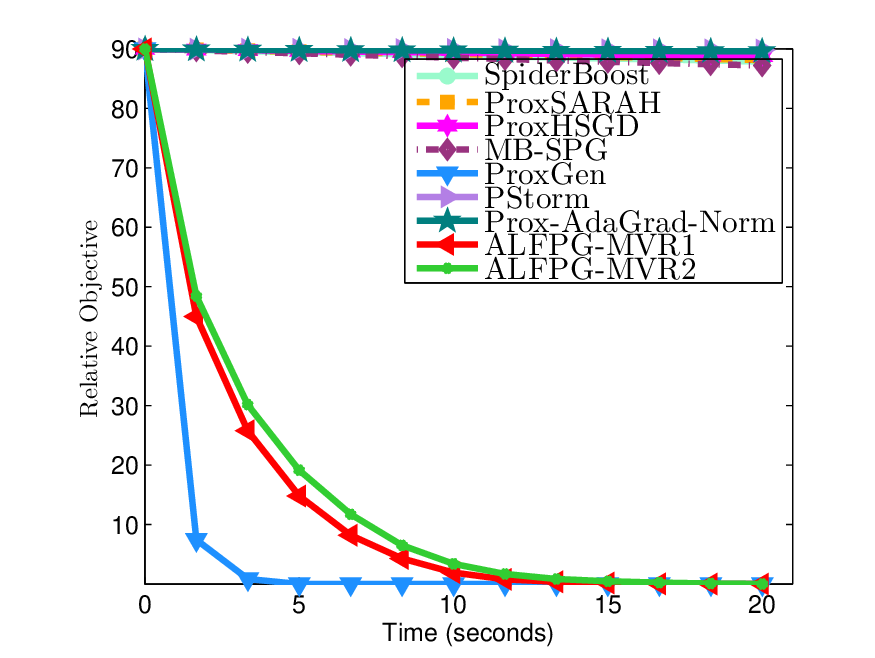}\caption{{\scriptsize Sector{-6412-1000}}}\label{fig:ex:phase:1:3}\end{subfigure}
\begin{subfigure}{.24\textwidth}\centering\includegraphics[width=\linewidth]{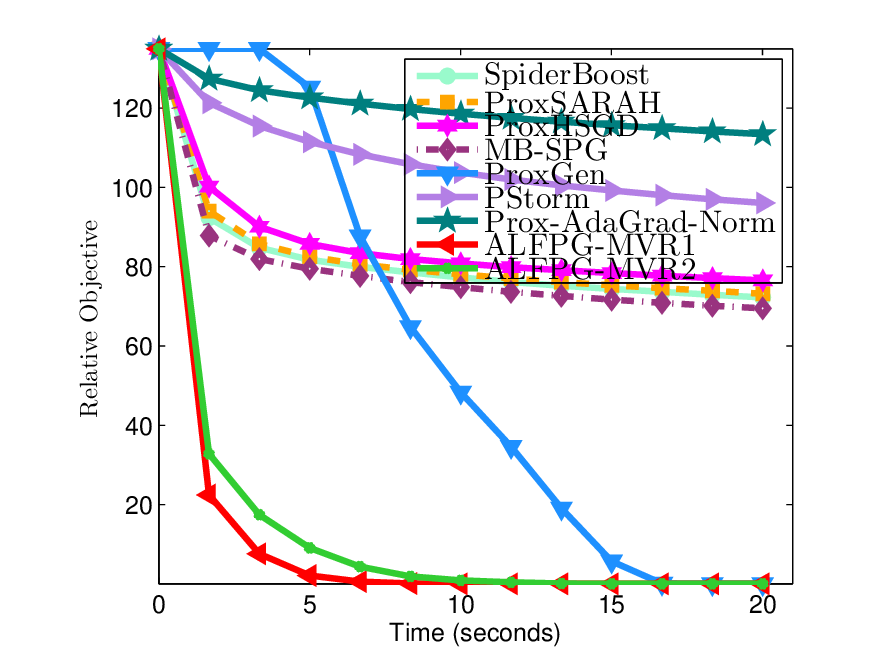}\caption{{\scriptsize MNIST{-2000-784}}}\label{fig:ex:phase:1:4}\end{subfigure}

\begin{subfigure}{.24\textwidth}\centering\includegraphics[width=\linewidth]{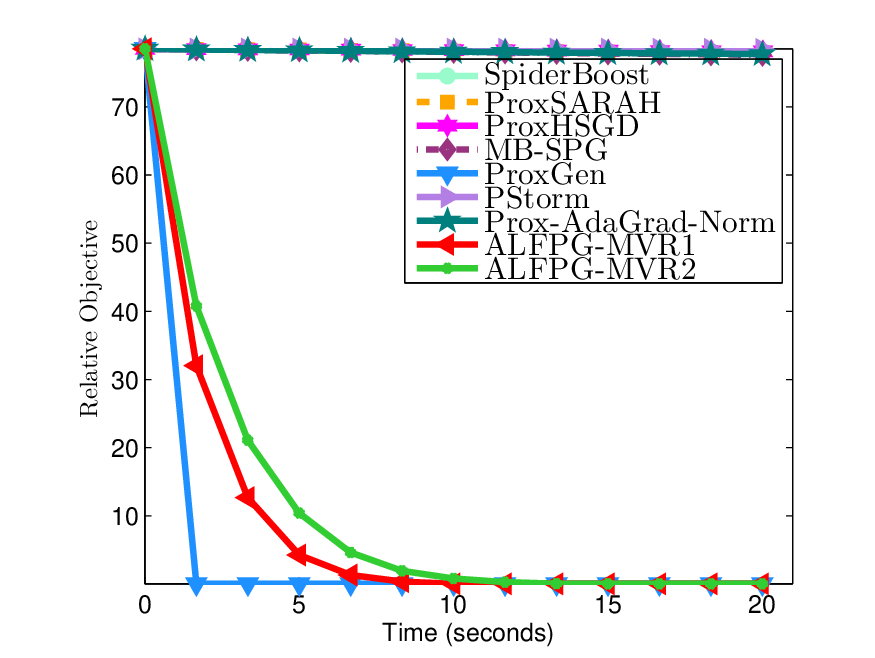}\caption{{\scriptsize CIFAR{-10000-1000}}}\label{fig:ex:phase:1:5}\end{subfigure}
\begin{subfigure}{.24\textwidth}\centering\includegraphics[width=\linewidth]{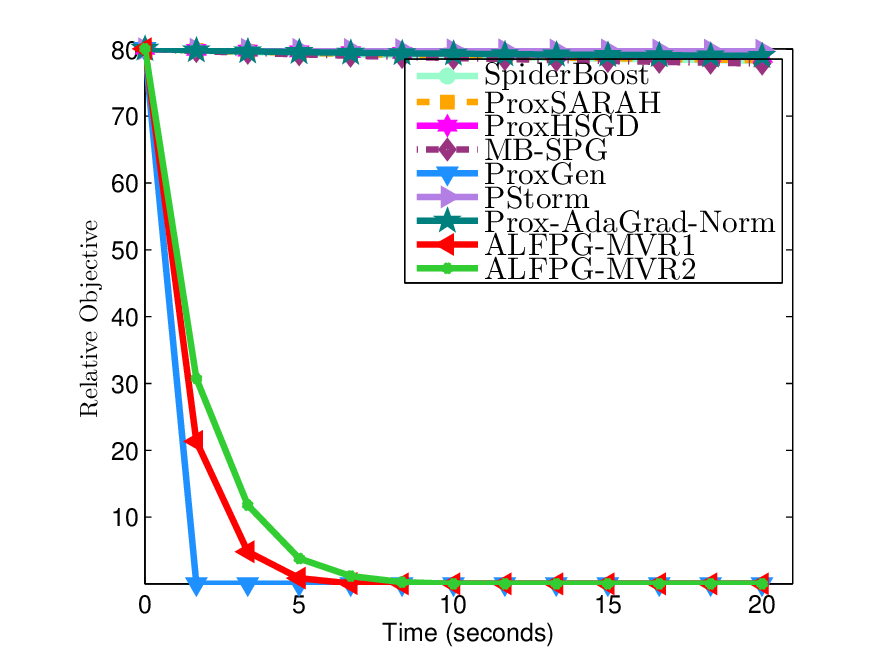}\caption{{\scriptsize Gisette{-5000-1000}}}\label{fig:ex:phase:1:6}\end{subfigure}
\begin{subfigure}{.24\textwidth}\centering\includegraphics[width=\linewidth]{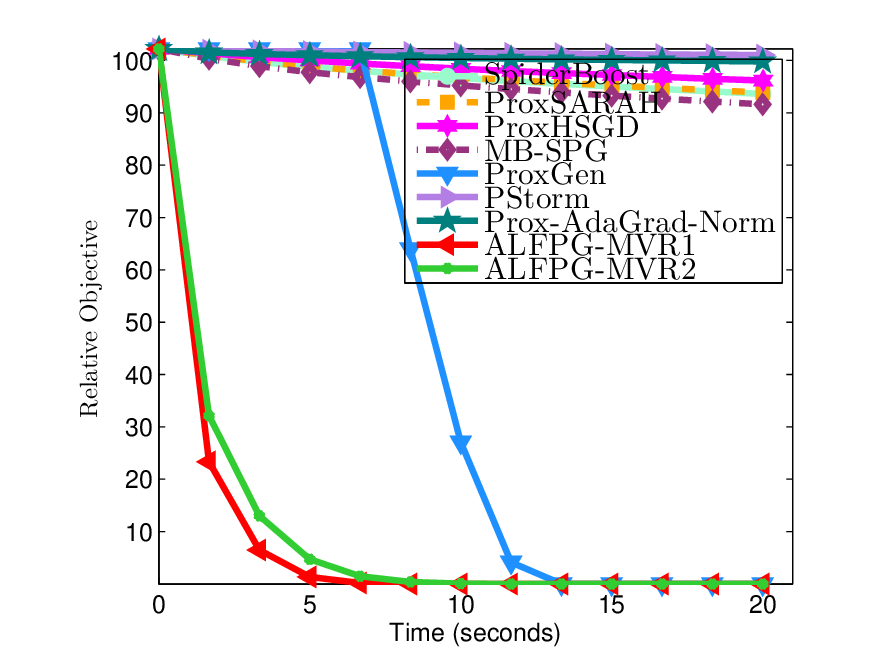}\caption{{\scriptsize Caltech{-3000-1000}}}\label{fig:ex:phase:1:7}\end{subfigure}
\begin{subfigure}{.24\textwidth}\centering\includegraphics[width=\linewidth]{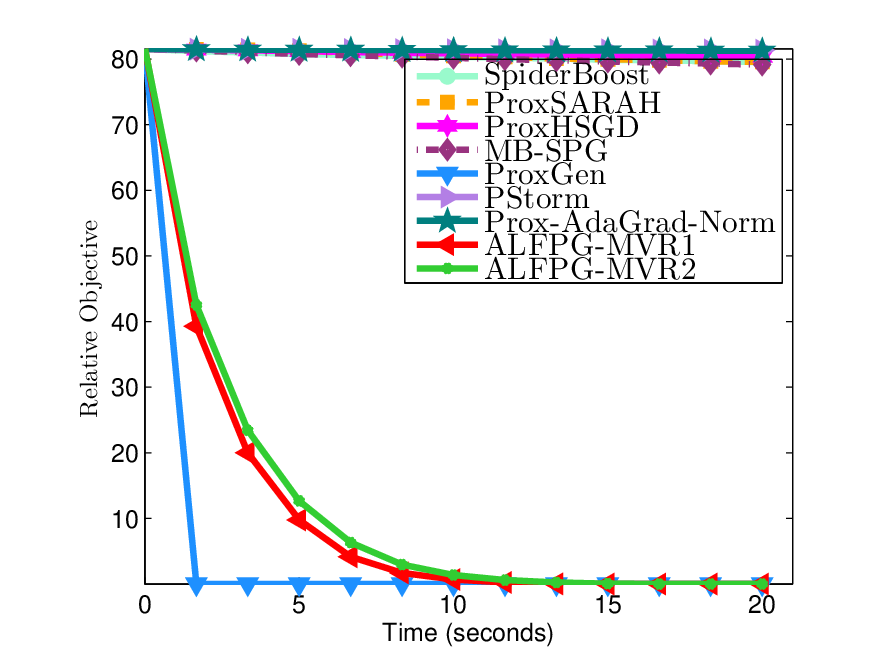}\caption{{\scriptsize Gaussian{-10000-1000}}}\label{fig:ex:phase:1:8}\end{subfigure}

\caption{Sparse phase retrieval with $\dot\lambda=10^{-3}$ in the expectation setting.}
\label{fig:ex:phase:1}
\end{figure*}

\begin{figure*}[!h]
\centering
\begin{subfigure}{.24\textwidth}\centering\includegraphics[width=\linewidth]{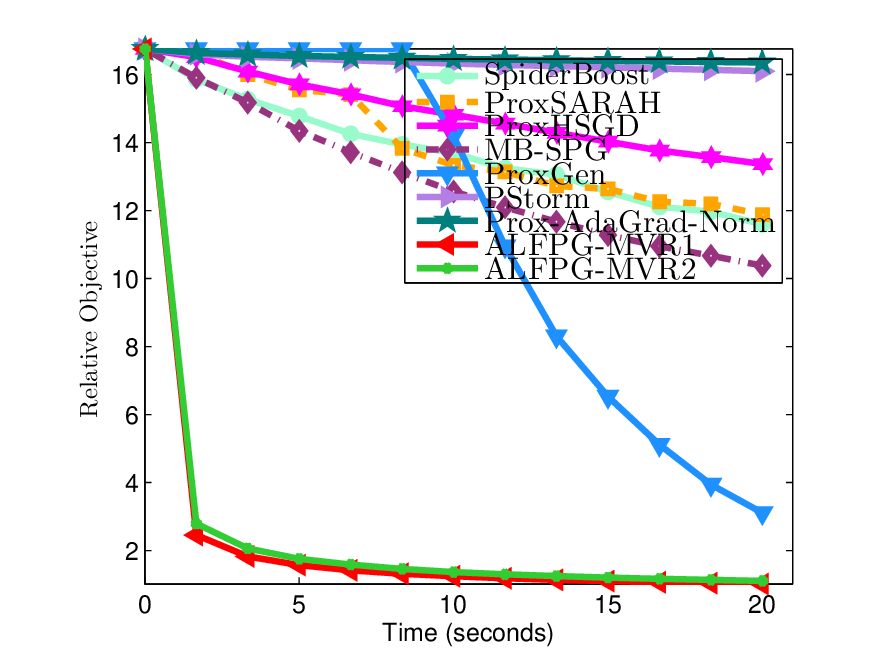}\caption{{\scriptsize TDT2{-9000-1000}}}\label{fig:ex:phase:2:1}\end{subfigure}
\begin{subfigure}{.24\textwidth}\centering\includegraphics[width=\linewidth]{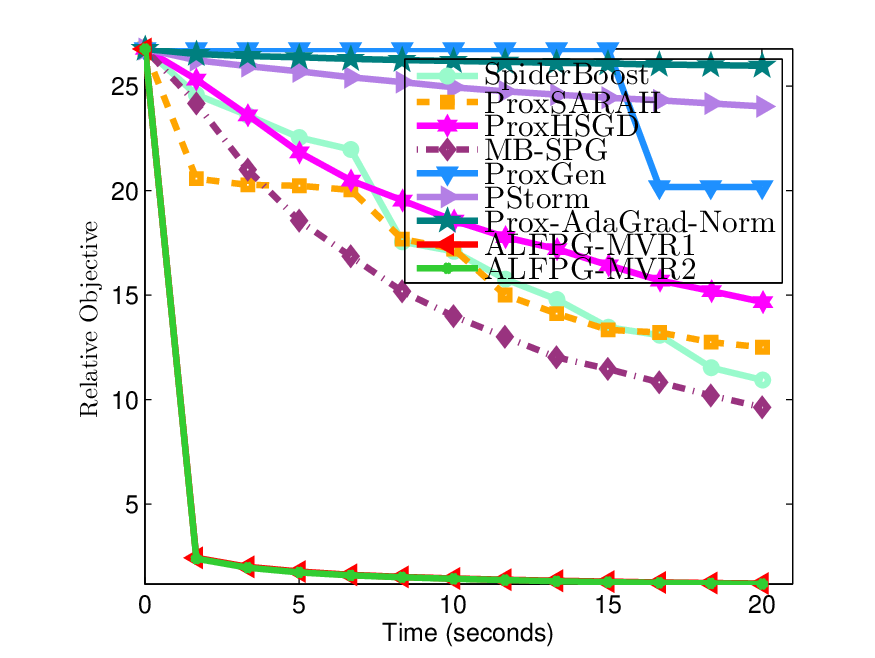}\caption{{\scriptsize 20Newsgroups{-10000-1000}}}\label{fig:ex:phase:2:2}\end{subfigure}
\begin{subfigure}{.24\textwidth}\centering\includegraphics[width=\linewidth]{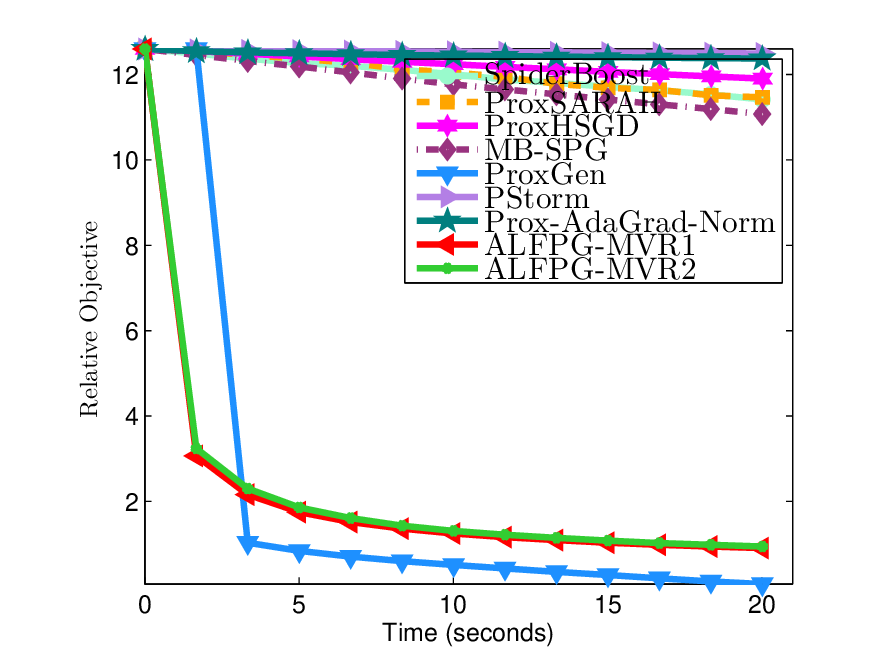}\caption{{\scriptsize Sector{-6412-1000}}}\label{fig:ex:phase:2:3}\end{subfigure}
\begin{subfigure}{.24\textwidth}\centering\includegraphics[width=\linewidth]{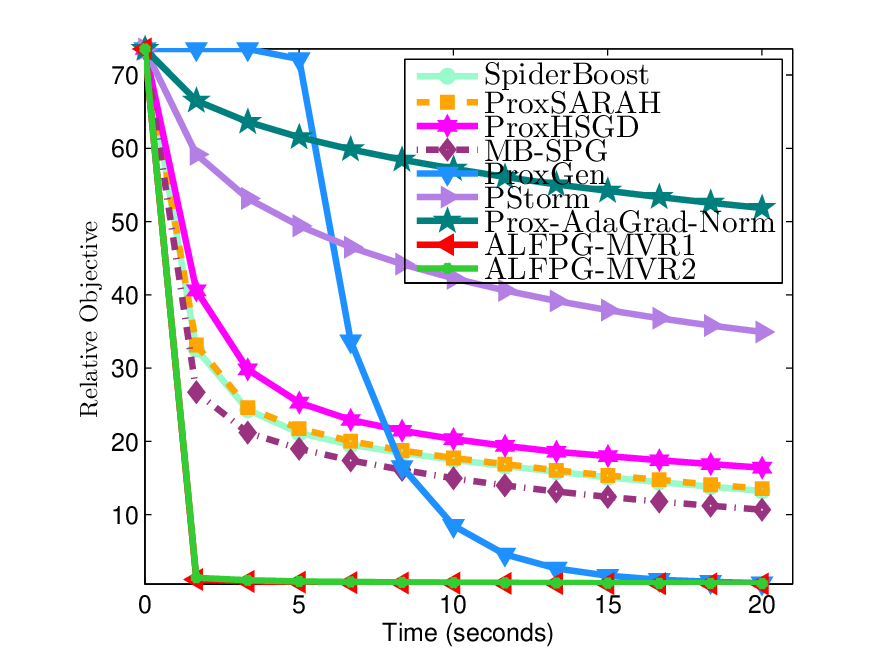}\caption{{\scriptsize MNIST{-2000-784}}}\label{fig:ex:phase:2:4}\end{subfigure}

\begin{subfigure}{.24\textwidth}\centering\includegraphics[width=\linewidth]{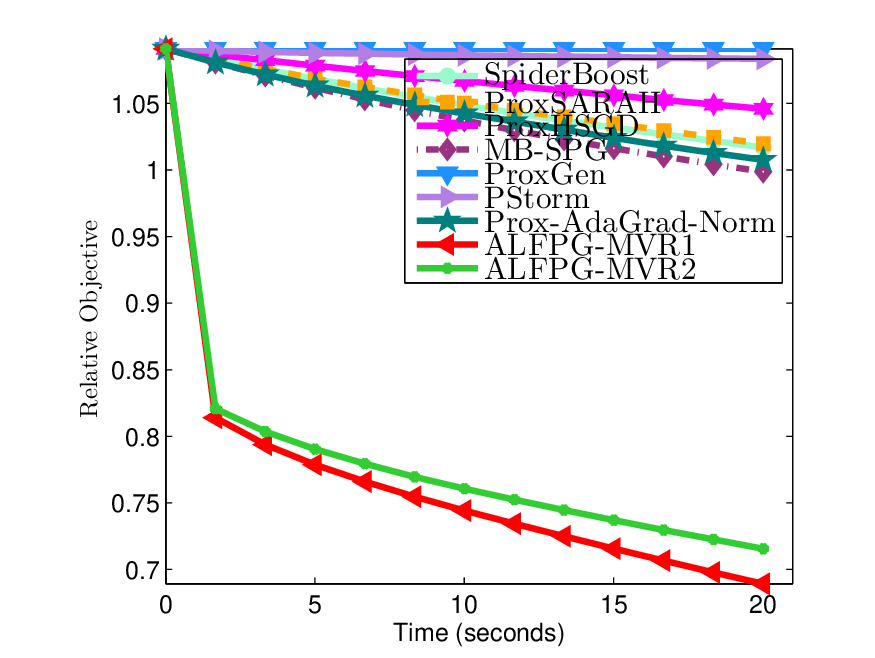}\caption{{\scriptsize CIFAR{-10000-1000}}}\label{fig:ex:phase:2:5}\end{subfigure}
\begin{subfigure}{.24\textwidth}\centering\includegraphics[width=\linewidth]{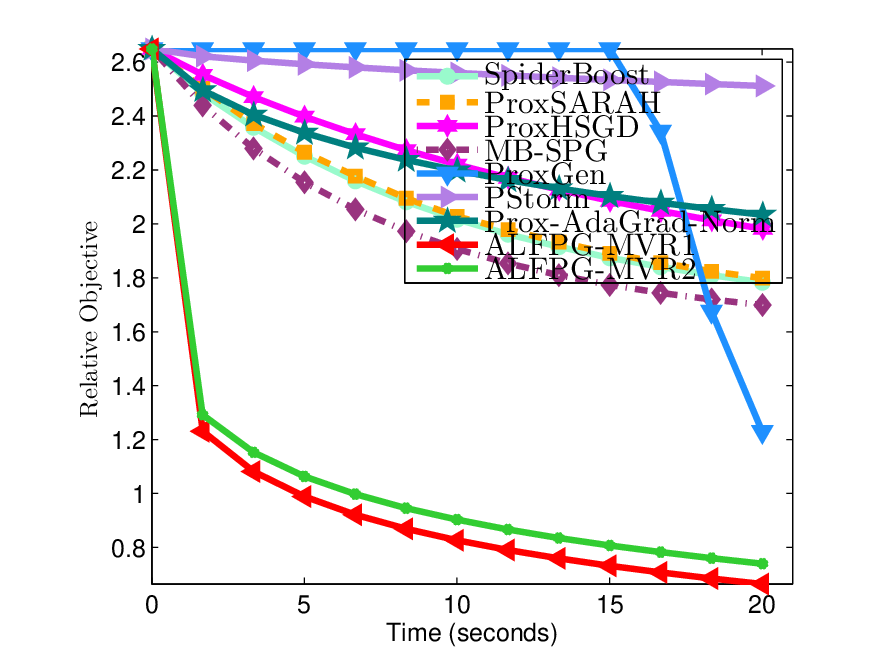}\caption{{\scriptsize Gisette{-5000-1000}}}\label{fig:ex:phase:2:6}\end{subfigure}
\begin{subfigure}{.24\textwidth}\centering\includegraphics[width=\linewidth]{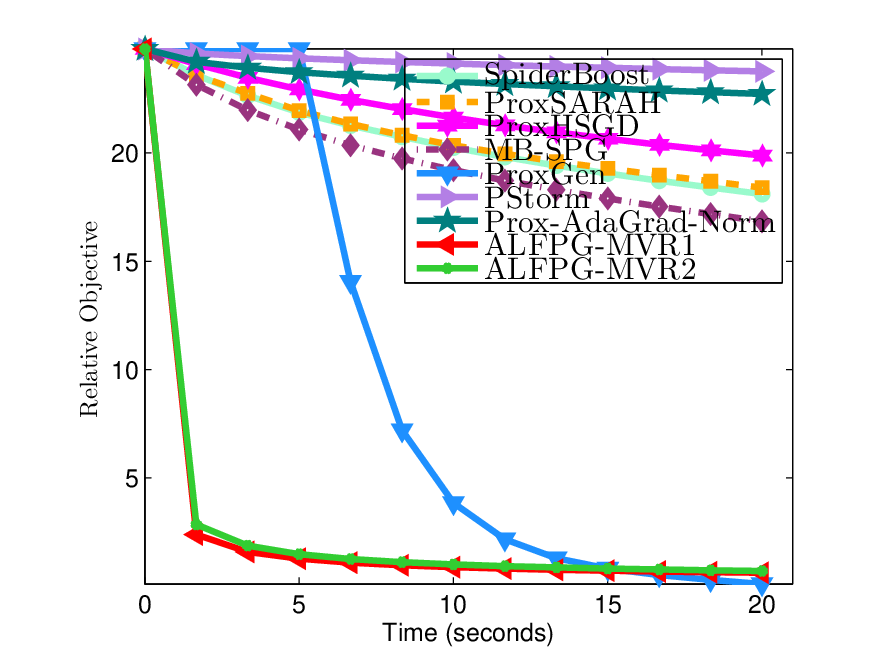}\caption{{\scriptsize Caltech{-3000-1000}}}\label{fig:ex:phase:2:7}\end{subfigure}
\begin{subfigure}{.24\textwidth}\centering\includegraphics[width=\linewidth]{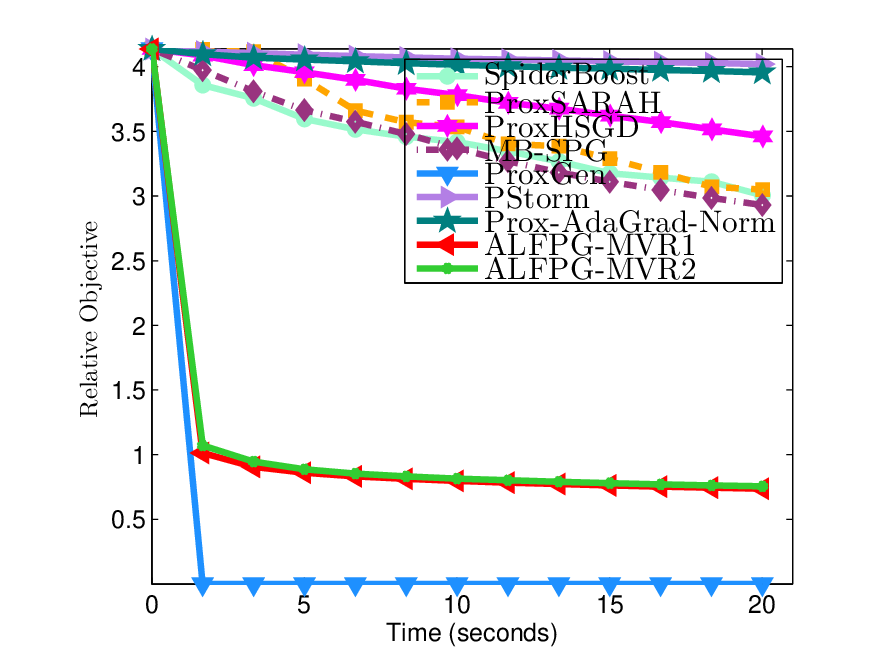}\caption{{\scriptsize Gaussian{-10000-1000}}}\label{fig:ex:phase:2:8}\end{subfigure}
\caption{Sparse phase retrieval with $\dot\lambda=10^{-5}$ in the expectation setting.}
\label{fig:ex:phase:2}
\end{figure*}

\begin{figure*}[!t]
\centering
\begin{subfigure}{.24\textwidth}\centering\includegraphics[width=\linewidth]{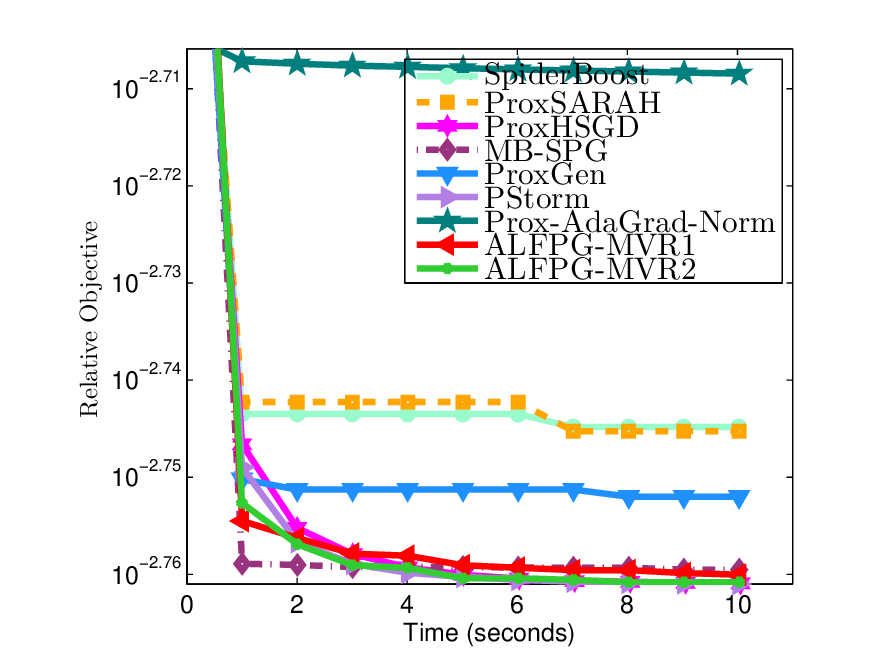}\caption{{\scriptsize TDT2{-36771-1000}}}\label{fig:ex:orth:1:1}\end{subfigure}
\begin{subfigure}{.24\textwidth}\centering\includegraphics[width=\linewidth]{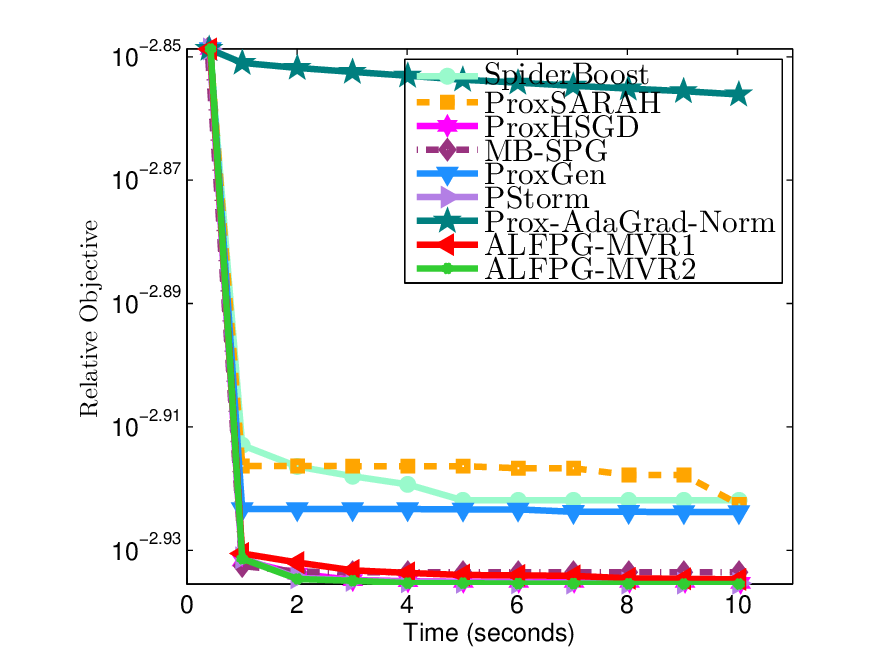}\caption{{\scriptsize 20Newsgroups{-26214-1000}}}\label{fig:ex:orth:1:2}\end{subfigure}
\begin{subfigure}{.24\textwidth}\centering\includegraphics[width=\linewidth]{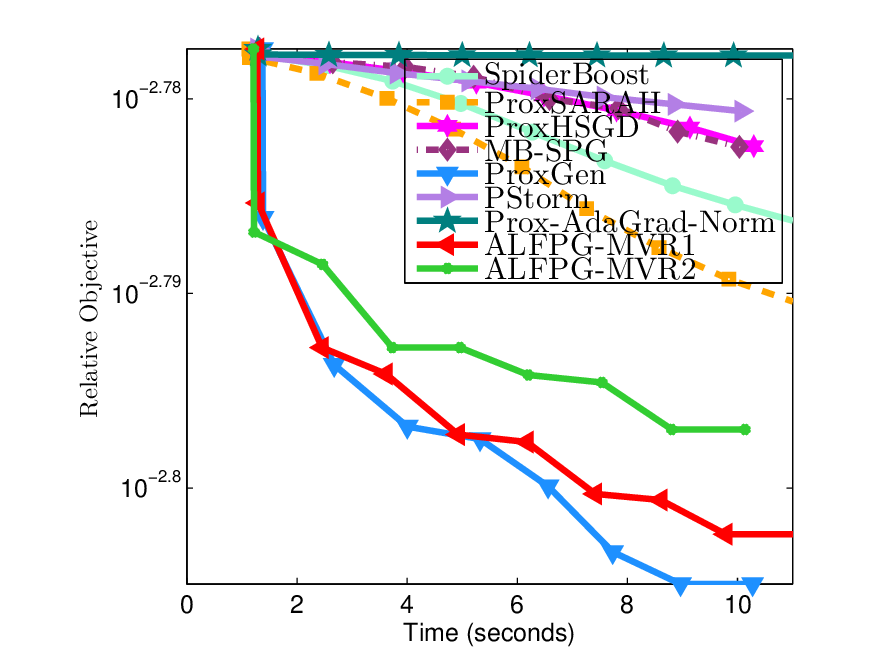}\caption{{\scriptsize Sector{-55197-1412}}}\label{fig:ex:orth:1:3}\end{subfigure}
\begin{subfigure}{.24\textwidth}\centering\includegraphics[width=\linewidth]{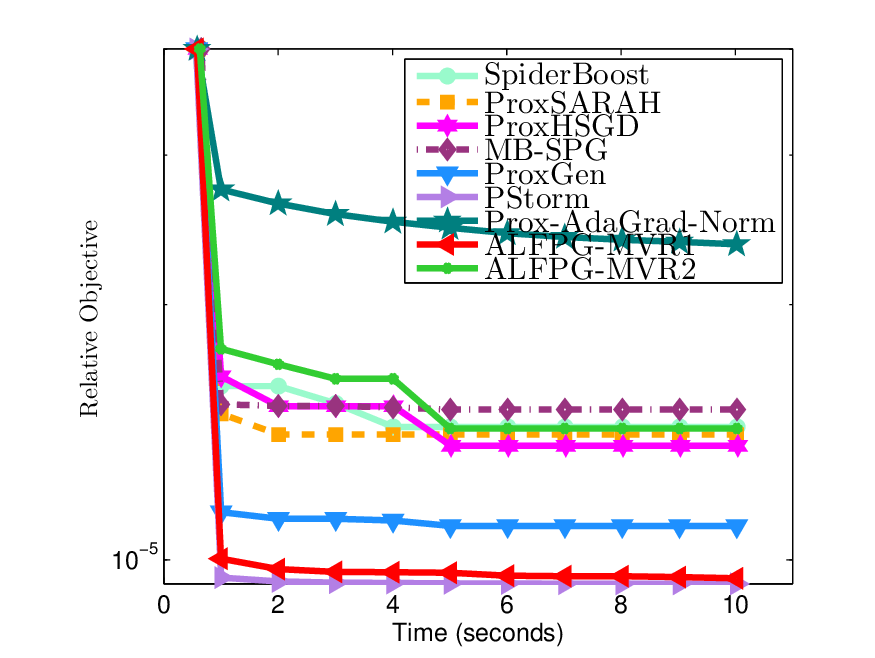}\caption{{\scriptsize MNIST{-60000-784}}}\label{fig:ex:orth:1:4}\end{subfigure}

\begin{subfigure}{.24\textwidth}\centering\includegraphics[width=\linewidth]{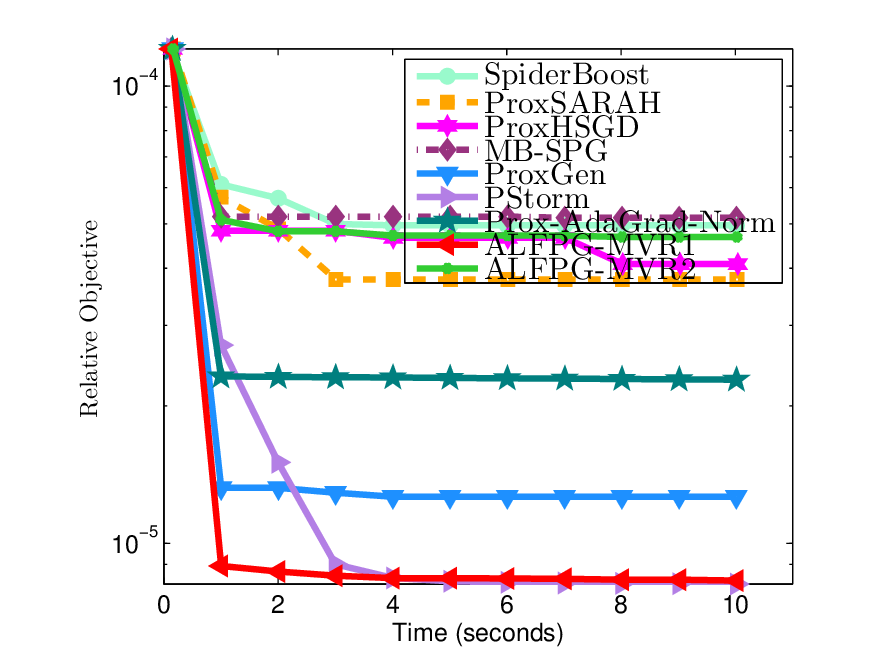}\caption{{\scriptsize CIFAR{-10000-1072}}}\label{fig:ex:orth:1:5}\end{subfigure}
\begin{subfigure}{.24\textwidth}\centering\includegraphics[width=\linewidth]{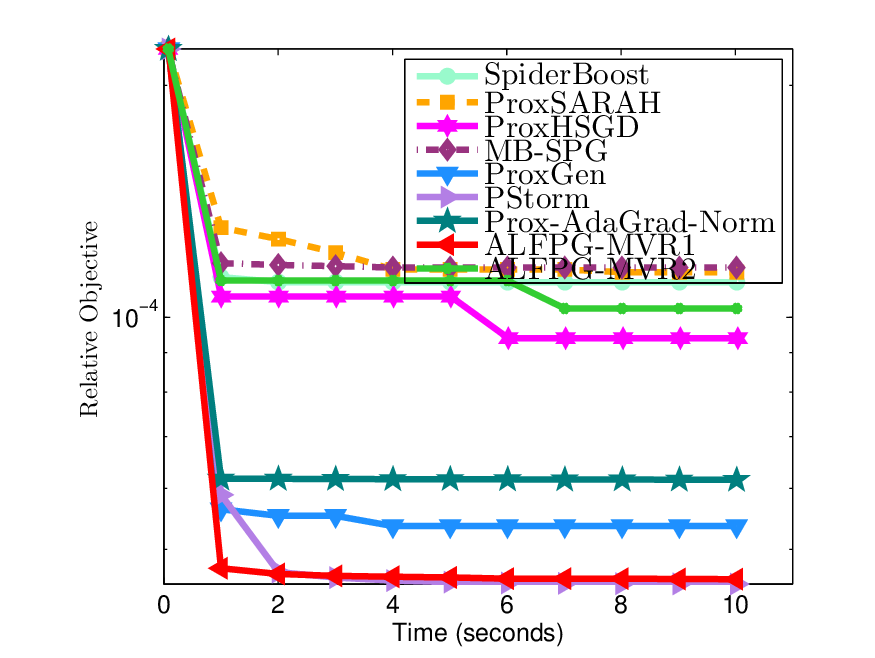}\caption{{\scriptsize Gisette{-6000-1000}}}\label{fig:ex:orth:1:6}\end{subfigure}
\begin{subfigure}{.24\textwidth}\centering\includegraphics[width=\linewidth]{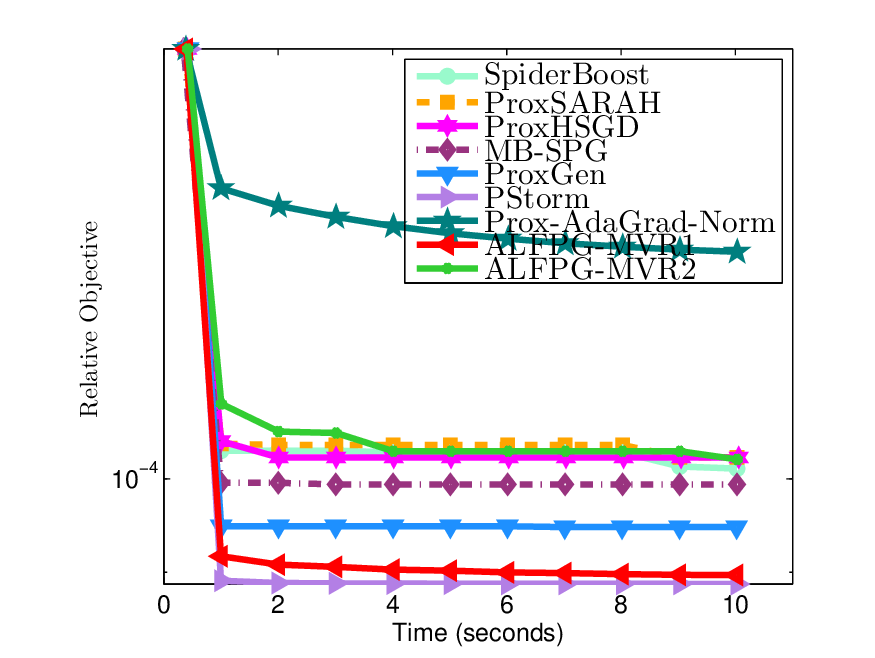}\caption{{\scriptsize Caltech{-29358-1096}}}\label{fig:ex:orth:1:7}\end{subfigure}
\begin{subfigure}{.24\textwidth}\centering\includegraphics[width=\linewidth]{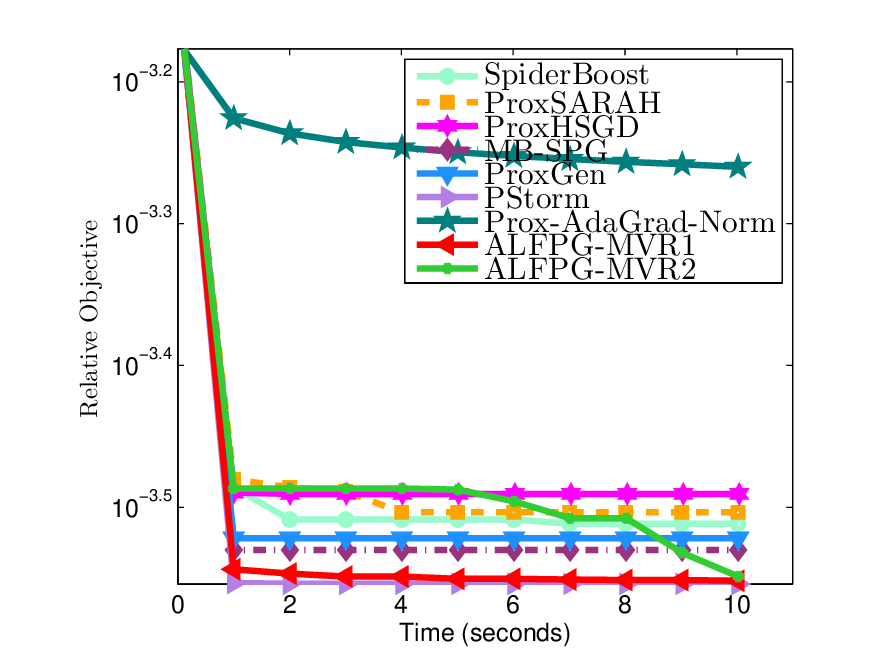}\caption{{\scriptsize E2006{-10000-1000}}}\label{fig:ex:orth:1:8}\end{subfigure}

\caption{Linear eigenvalue problems with $\dot r=20$ in the expectation setting.}
\label{fig:ex:orth:1}
\end{figure*}

\begin{figure*}[!h]
\centering
\begin{subfigure}{.24\textwidth}\centering\includegraphics[width=\linewidth]{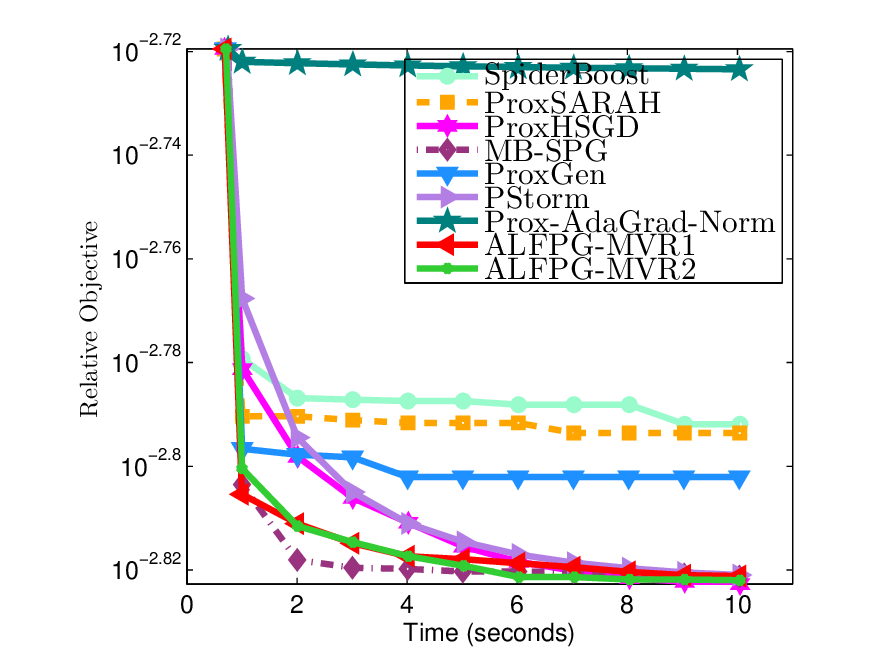}\caption{{\scriptsize TDT2{-36771-1000}}}\label{fig:ex:orth:2:1}\end{subfigure}
\begin{subfigure}{.24\textwidth}\centering\includegraphics[width=\linewidth]{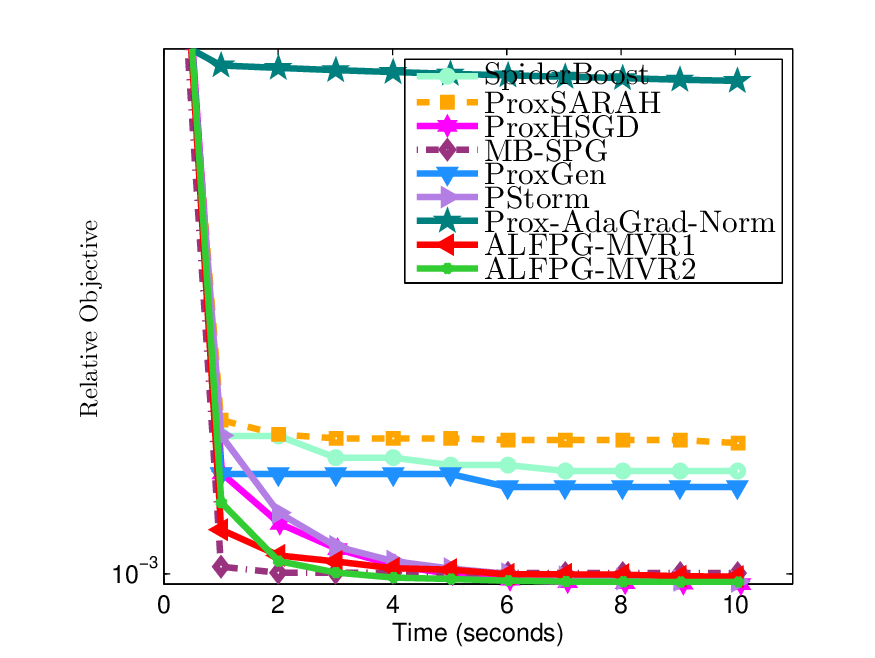}\caption{{\scriptsize 20Newsgroups{-26214-1000}}}\label{fig:ex:orth:2:2}\end{subfigure}
\begin{subfigure}{.24\textwidth}\centering\includegraphics[width=\linewidth]{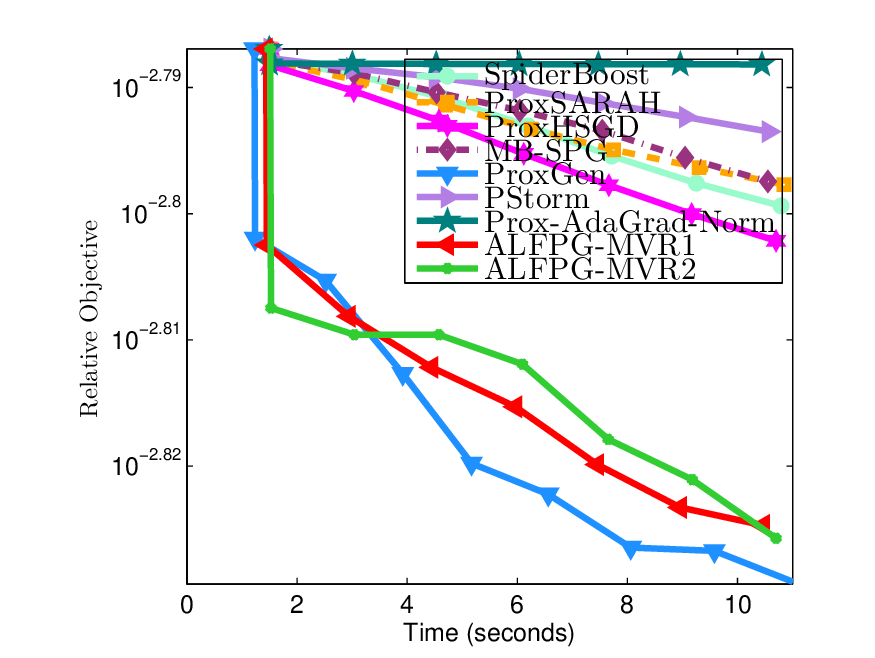}\caption{{\scriptsize Sector{-55197-1412}}}\label{fig:ex:orth:2:3}\end{subfigure}
\begin{subfigure}{.24\textwidth}\centering\includegraphics[width=\linewidth]{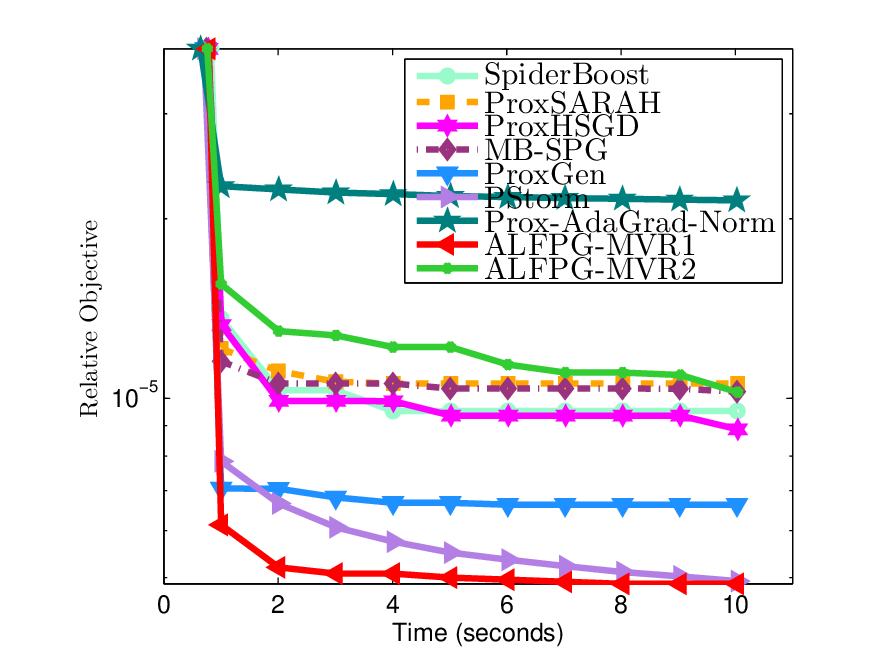}\caption{{\scriptsize MNIST{-60000-784}}}\label{fig:ex:orth:2:4}\end{subfigure}

\begin{subfigure}{.24\textwidth}\centering\includegraphics[width=\linewidth]{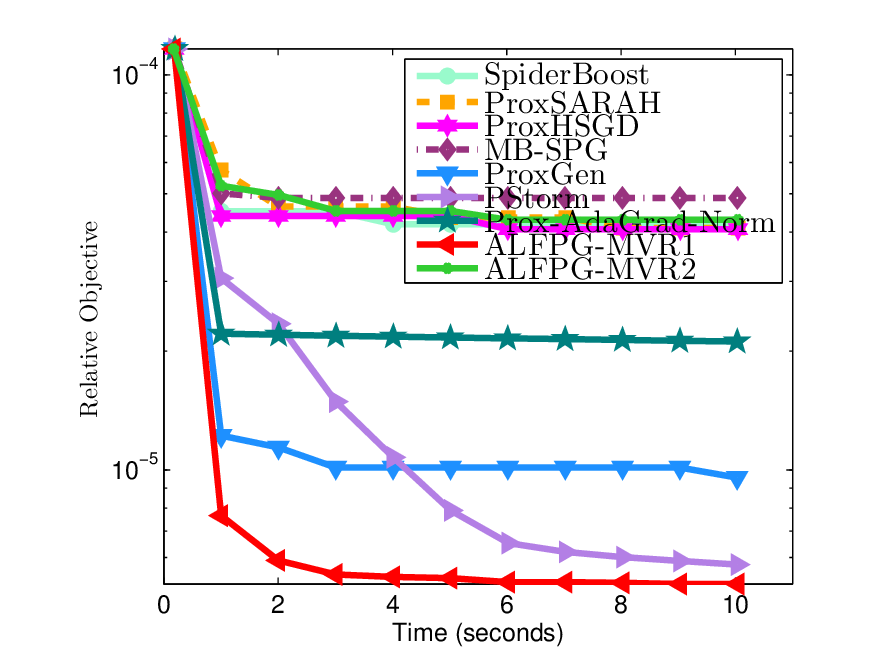}\caption{{\scriptsize CIFAR{-10000-1072}}}\label{fig:ex:orth:2:5}\end{subfigure}
\begin{subfigure}{.24\textwidth}\centering\includegraphics[width=\linewidth]{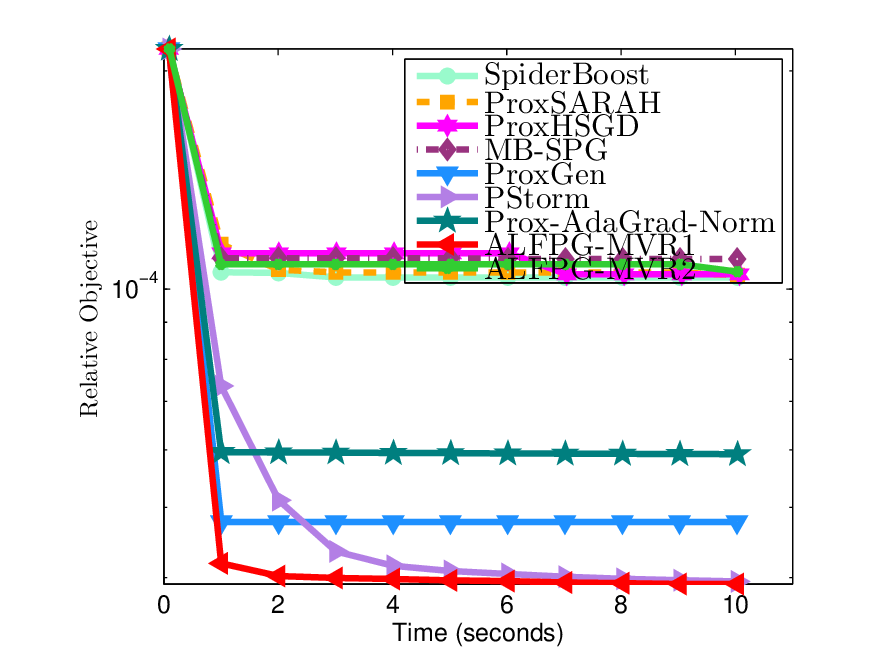}\caption{{\scriptsize Gisette{-6000-1000}}}\label{fig:ex:orth:2:6}\end{subfigure}
\begin{subfigure}{.24\textwidth}\centering\includegraphics[width=\linewidth]{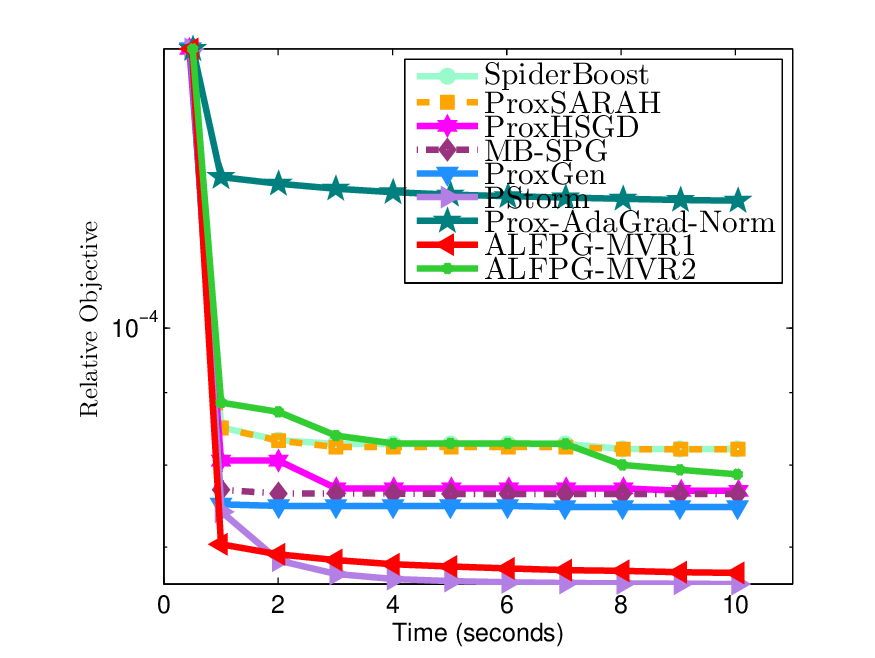}\caption{{\scriptsize Caltech{-29358-1096}}}\label{fig:ex:orth:2:7}\end{subfigure}
\begin{subfigure}{.24\textwidth}\centering\includegraphics[width=\linewidth]{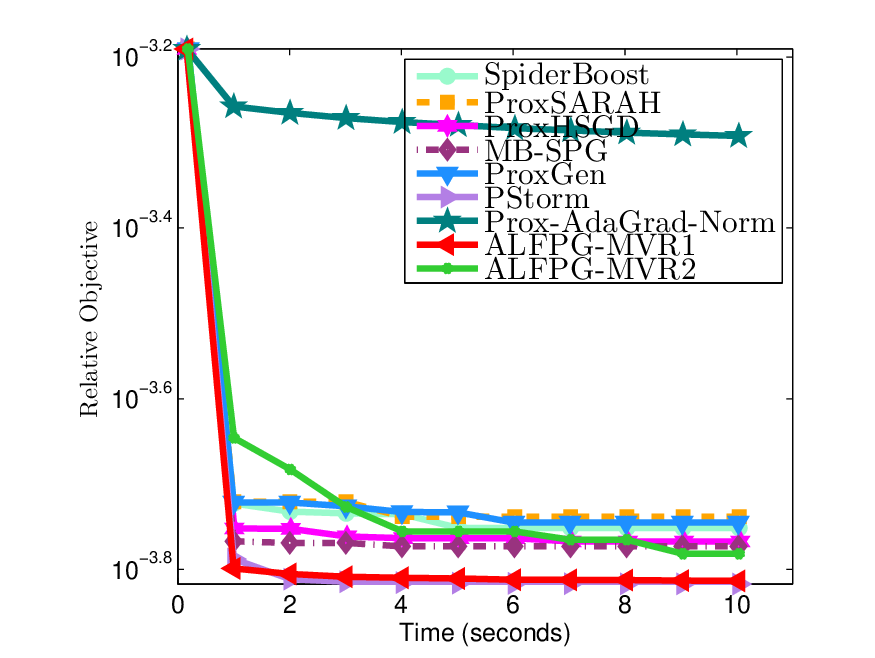}\caption{{\scriptsize E2006{-10000-1000}}}\label{fig:ex:orth:2:8}\end{subfigure}
\caption{Linear eigenvalue problems with $\dot r=50$ in the expectation setting.}
\label{fig:ex:orth:2}

\end{figure*}

\begin{figure*}[!t]
\centering
\begin{subfigure}{.24\textwidth}\centering\includegraphics[width=\linewidth]{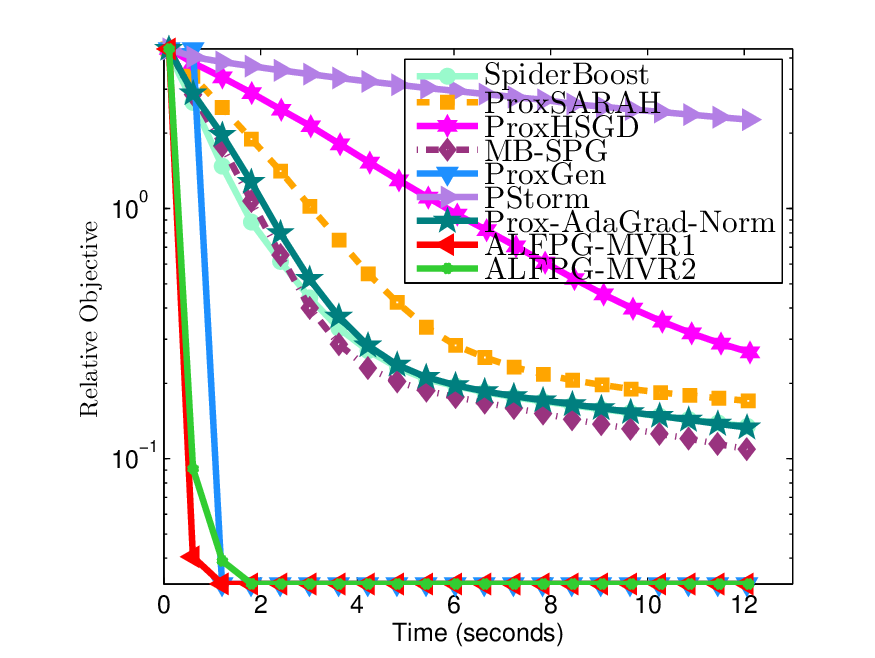}\caption{{\scriptsize TDT2{-36771-1000}}}\label{fig:ex:nn:1:1}\end{subfigure}
\begin{subfigure}{.24\textwidth}\centering\includegraphics[width=\linewidth]{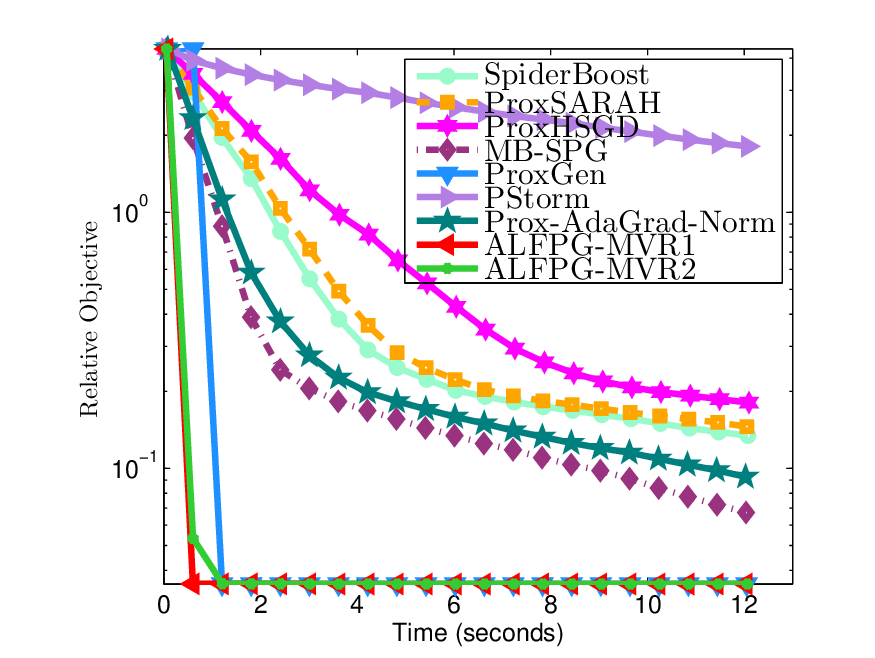}\caption{{\scriptsize 20Newsgroups{-26214-1000}}}\label{fig:ex:nn:1:2}\end{subfigure}
\begin{subfigure}{.24\textwidth}\centering\includegraphics[width=\linewidth]{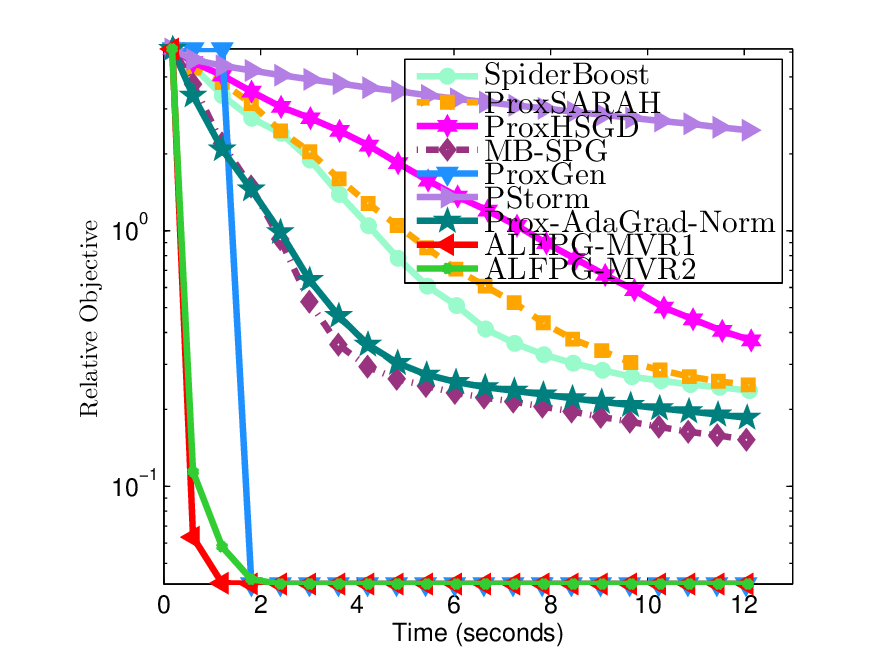}\caption{{\scriptsize Sector{-55197-1412}}}\label{fig:ex:nn:1:3}\end{subfigure}
\begin{subfigure}{.24\textwidth}\centering\includegraphics[width=\linewidth]{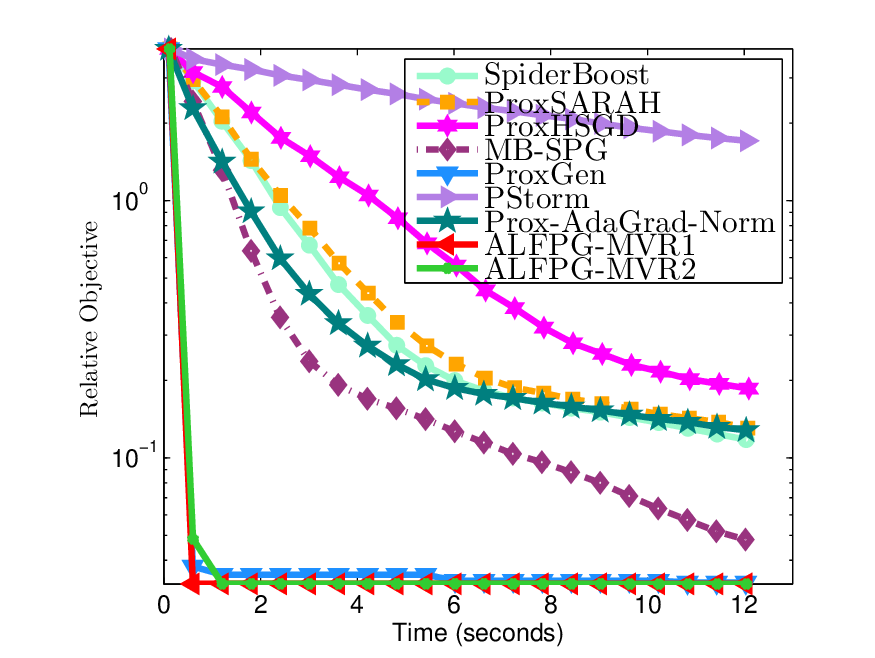}\caption{{\scriptsize MNIST{-60000-784}}}\label{fig:ex:nn:1:4}\end{subfigure}

\begin{subfigure}{.24\textwidth}\centering\includegraphics[width=\linewidth]{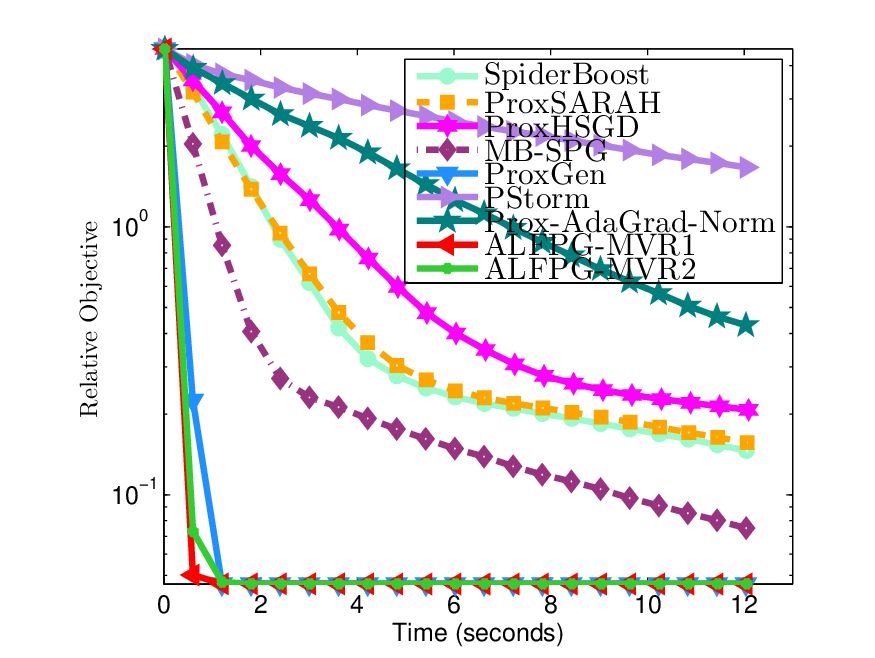}\caption{{\scriptsize CIFAR{-10000-1072}}}\label{fig:ex:nn:1:5}\end{subfigure}
\begin{subfigure}{.24\textwidth}\centering\includegraphics[width=\linewidth]{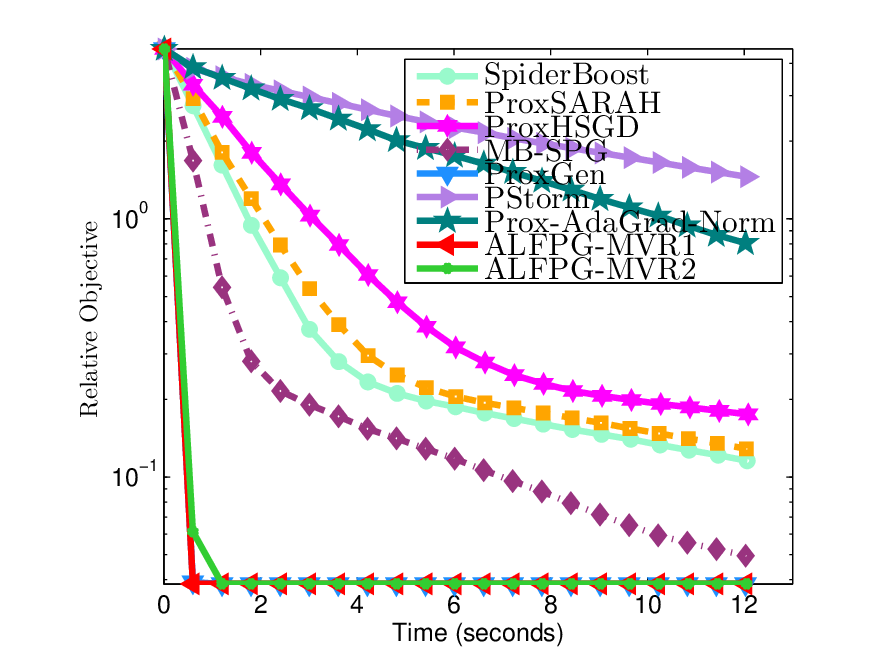}\caption{{\scriptsize Gisette{-6000-1000}}}\label{fig:ex:nn:1:6}\end{subfigure}
\begin{subfigure}{.24\textwidth}\centering\includegraphics[width=\linewidth]{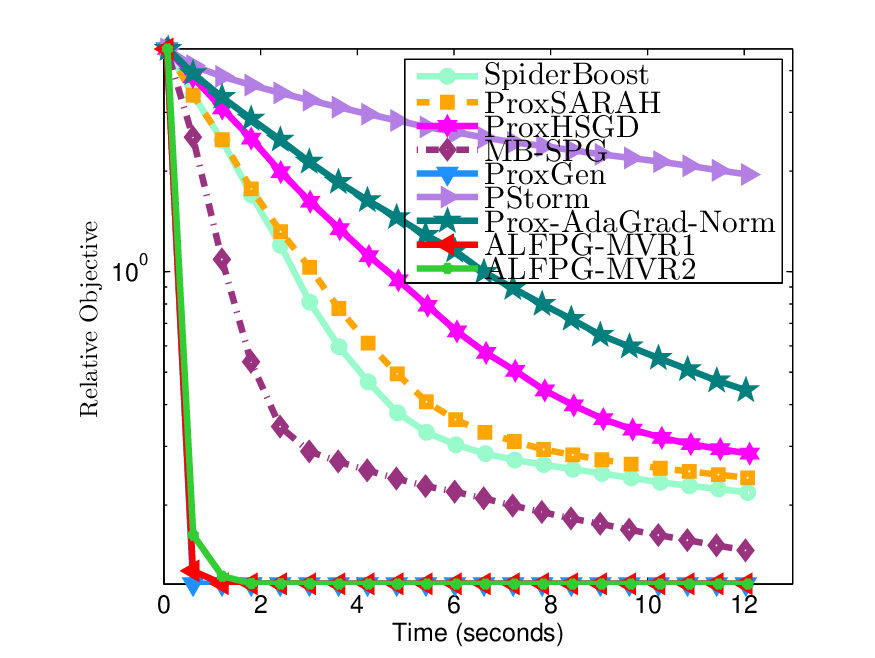}\caption{{\scriptsize Caltech{-29358-1096}}}\label{fig:ex:nn:1:7}\end{subfigure}
\begin{subfigure}{.24\textwidth}\centering\includegraphics[width=\linewidth]{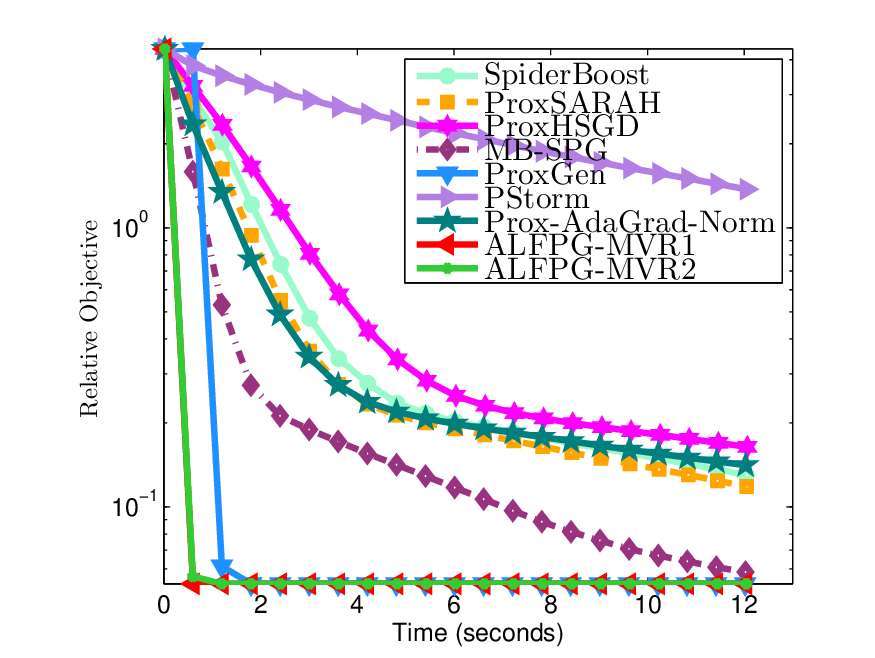}\caption{{\scriptsize E2006{-10000-1000}}}\label{fig:ex:nn:1:8}\end{subfigure}

\caption{Three-layer neural network regression with $\dot\lambda=0.01$ in the expectation setting.}
\label{fig:ex:nn:1}
\end{figure*}

\begin{figure*}[!h]
\centering
\begin{subfigure}{.24\textwidth}\centering\includegraphics[width=\linewidth]{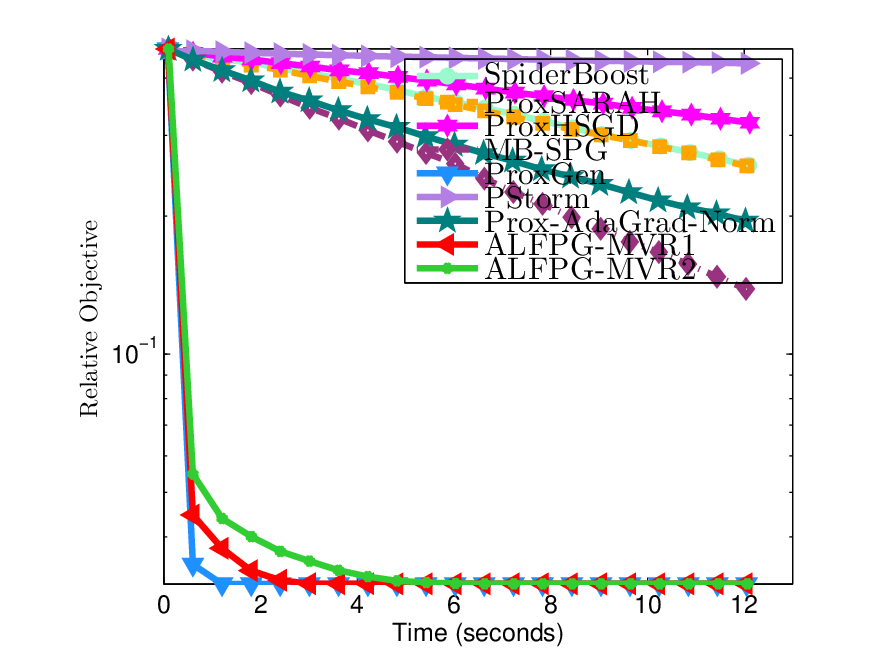}\caption{{\scriptsize TDT2{-36771-1000}}}\label{fig:ex:nn:2:1}\end{subfigure}
\begin{subfigure}{.24\textwidth}\centering\includegraphics[width=\linewidth]{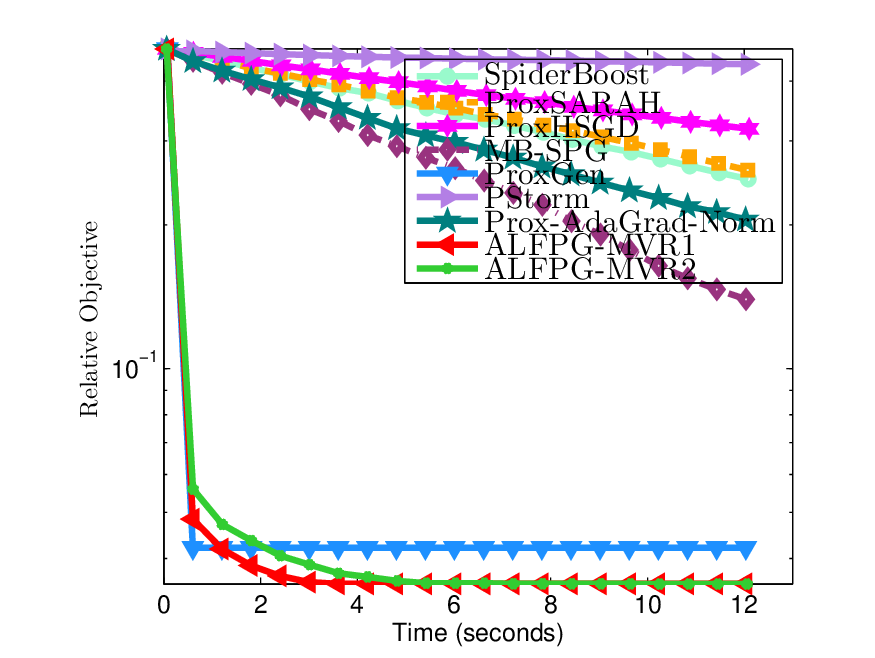}\caption{{\scriptsize 20Newsgroups{-26214-1000}}}\label{fig:ex:nn:2:2}\end{subfigure}
\begin{subfigure}{.24\textwidth}\centering\includegraphics[width=\linewidth]{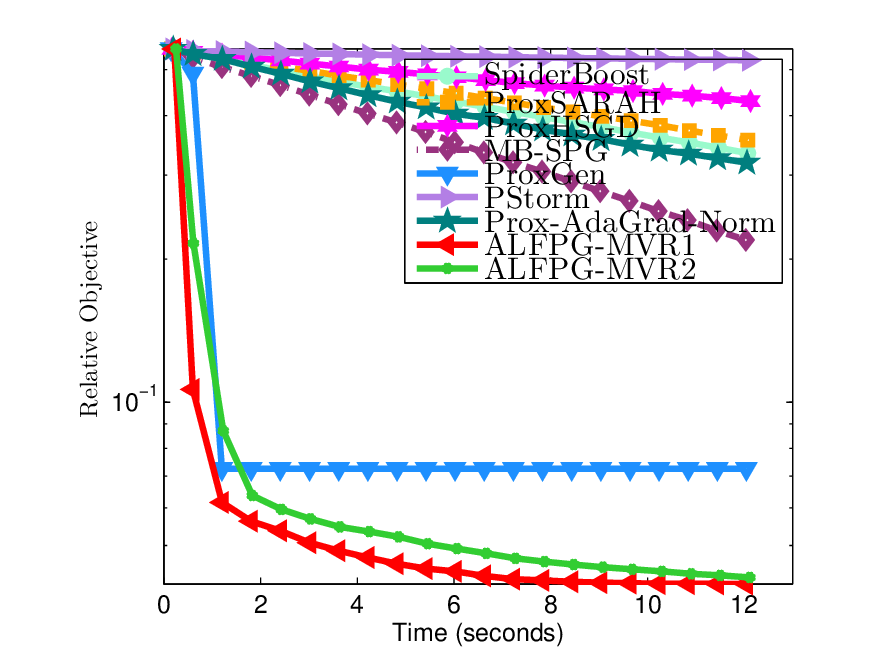}\caption{{\scriptsize Sector{-55197-1412}}}\label{fig:ex:nn:2:3}\end{subfigure}
\begin{subfigure}{.24\textwidth}\centering\includegraphics[width=\linewidth]{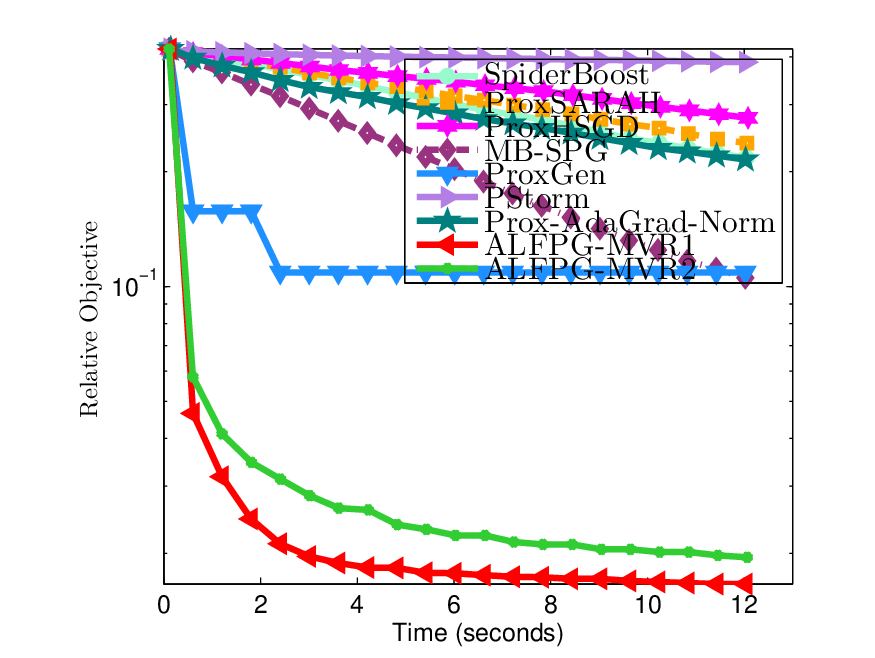}\caption{{\scriptsize MNIST{-60000-784}}}\label{fig:ex:nn:2:4}\end{subfigure}

\begin{subfigure}{.24\textwidth}\centering\includegraphics[width=\linewidth]{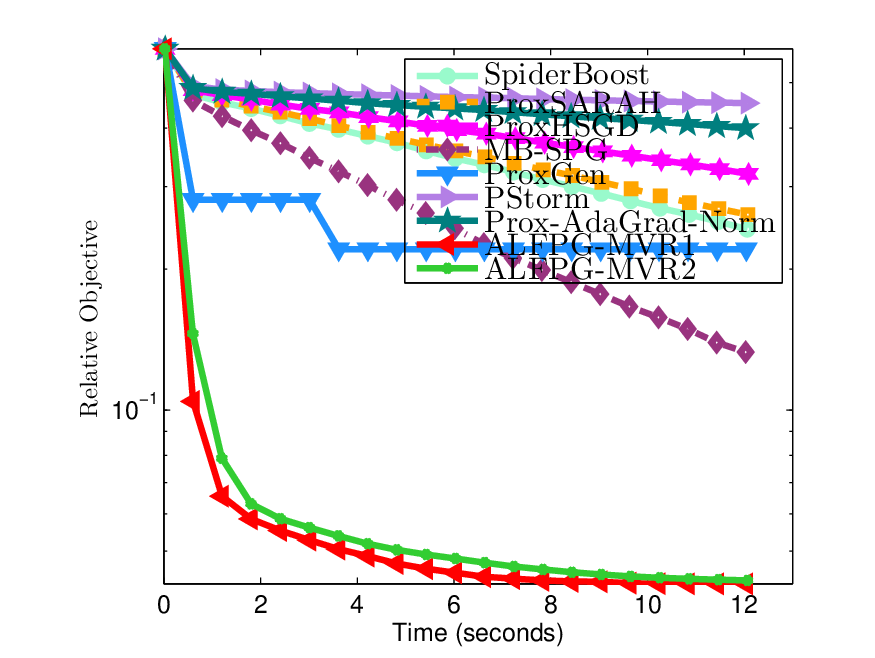}\caption{{\scriptsize CIFAR{-10000-1072}}}\label{fig:ex:nn:2:5}\end{subfigure}
\begin{subfigure}{.24\textwidth}\centering\includegraphics[width=\linewidth]{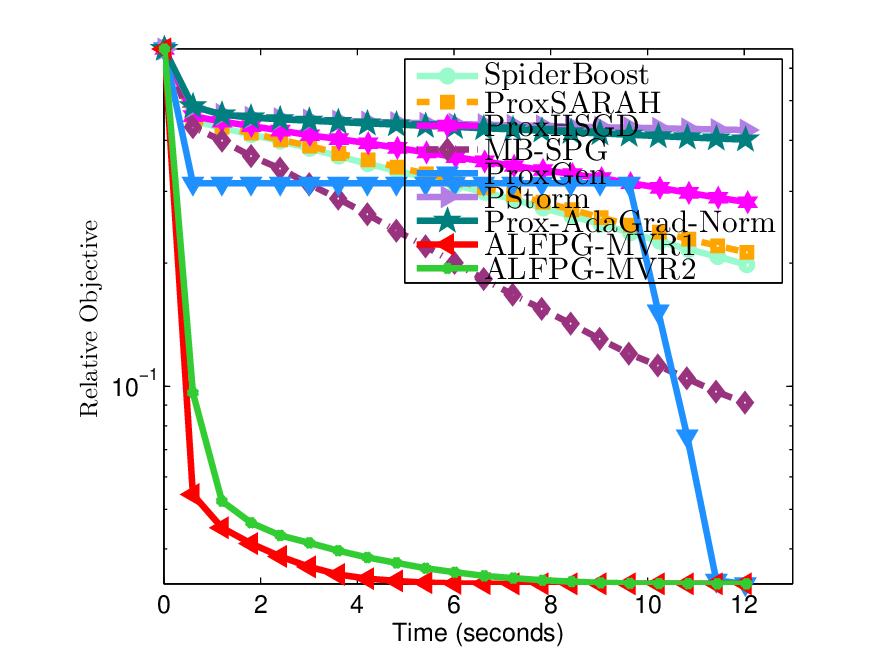}\caption{{\scriptsize Gisette{-6000-1000}}}\label{fig:ex:nn:2:6}\end{subfigure}
\begin{subfigure}{.24\textwidth}\centering\includegraphics[width=\linewidth]{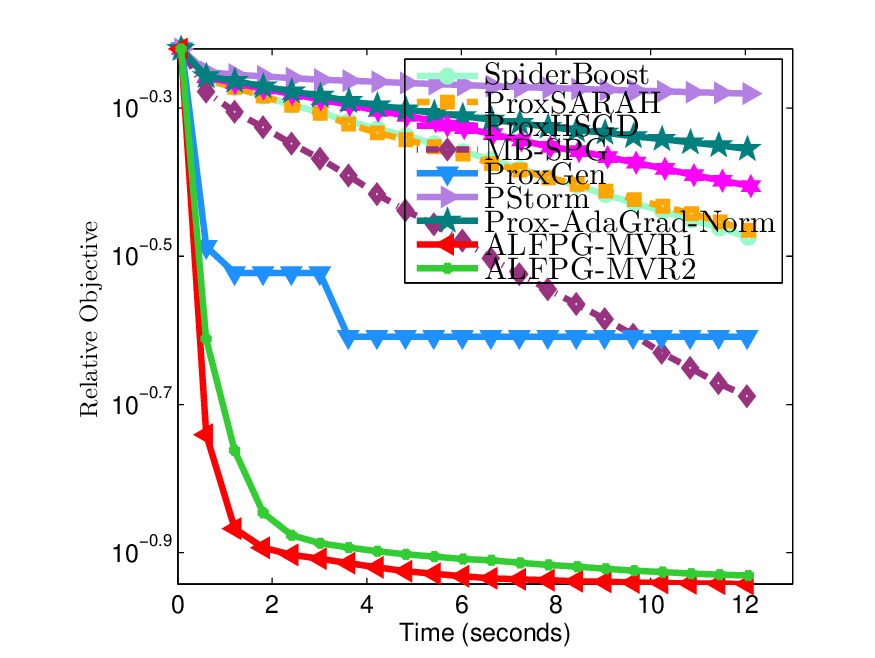}\caption{{\scriptsize Caltech{-29358-1096}}}\label{fig:ex:nn:2:7}\end{subfigure}
\begin{subfigure}{.24\textwidth}\centering\includegraphics[width=\linewidth]{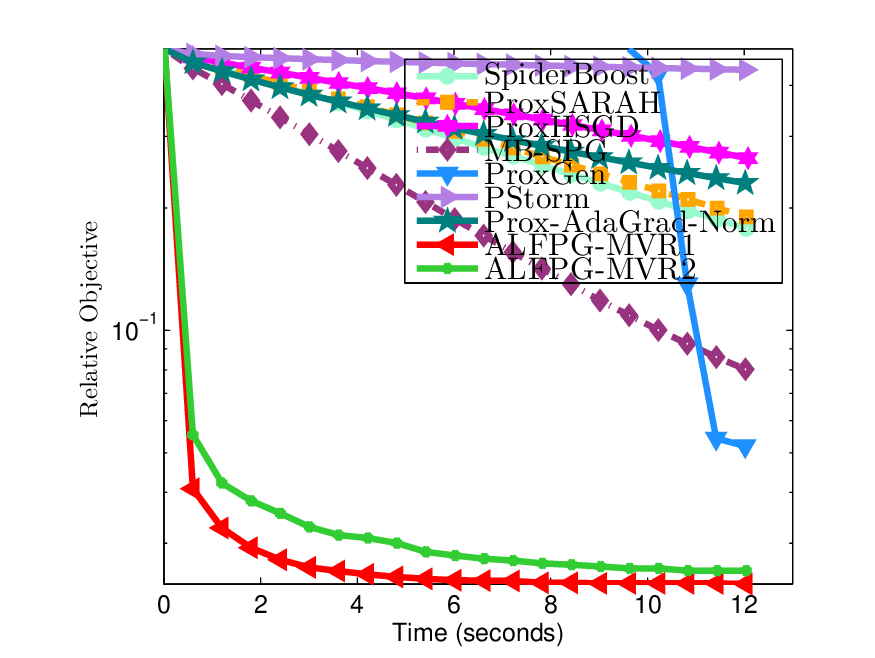}\caption{{\scriptsize E2006{-10000-1000}}}\label{fig:ex:nn:2:8}\end{subfigure}

\caption{Three-layer neural network regression with $\dot\lambda=0.001$ in the expectation setting.}
\label{fig:ex:nn:2}
\end{figure*}

$\blacktriangleright$ \textbf{Experimental Results}. \textbf{(i)} Figures~\ref{fig:ex:phase:1}--\ref{fig:ex:phase:2} show the experiment results on SPR. {\sf ALFPG-MVR1} and {\sf ALFPG-MVR2} reduce the objective rapidly on most datasets at both regularization levels. At \(\dot\lambda=10^{-3}\), MVR2 attains a substantially lower objective than MVR1 on 20Newsgroups. ProxGen reaches low objective values sooner on Sector, CIFAR, Gisette, and Gaussian.  \textbf{(ii)} Figures~\ref{fig:ex:orth:1}--\ref{fig:ex:orth:2} show the experiment results on LEV. At both ranks, {\sf ALFPG-MVR1} quickly reaches the lowest or nearly the lowest objective on most datasets, with lower final values than MVR2 on MNIST, Gisette, and Caltech. ProxGen is particularly competitive on Sector at rank 50. \textbf{(iii)} Figures~\ref{fig:ex:nn:1}--\ref{fig:ex:nn:2} show the experiment results on NNR. At \(\dot\lambda=0.01\), both {\sf ALFPG} variants quickly reach similarly low objective values. At \(\dot\lambda=0.001\), {\sf ALFPG-MVR1} attains a comparable or lower final objective than {\sf ALFPG-MVR2} on most datasets, with a clear improvement on MNIST.

\clearpage

\clearpage

\clearpage

\clearpage


\clearpage


\clearpage

\clearpage

\end{document}